\documentclass[11pt,a4paper,usenames,dvipsnames]{article}

\usepackage{url,amsmath,amssymb,latexsym,pstricks,mathrsfs,comment,amsthm,graphicx,tikz,tikz-cd,enumerate,accents,pgffor,cite,wrapfig,multicol,float,calc,geometry,tikz-3dplot}
\usepackage[colorlinks]{hyperref}
\usepackage[subnum]{cases}

\usepackage[margin=10pt,font=small,labelfont=bf, labelsep=period]{caption}

\newcommand{\nc}{\newcommand}
\nc{\rnc}{\renewcommand}

\nc\editonlyPre[1]{#1}
\nc\editonlyIE[1]{#1}
\nc\editonlyC[1]{#1}
\nc\editonlyT[1]{#1}
\nc\editonlyP[1]{#1}
\nc\editonlyB[1]{#1}
\nc\editonlyTL[1]{#1}
\nc\editonlyJ[1]{#1}
\nc\editonlyL[1]{#1}
\nc\editonlyIIE[1]{#1}
\nc\editonlyIB[1]{#1}

\usepackage{xcolor}

\usepackage{enumitem}

\usepackage[T1]{fontenc}
\usepackage{ wasysym }
\usepackage{stmaryrd}

\usepackage{tikz}
\usetikzlibrary{decorations.markings}
\usetikzlibrary{arrows,matrix}
\usepgflibrary{arrows}
  \usepgflibrary{fadings}
\usetikzlibrary{matrix,arrows,positioning,automata}
\pgfdeclarelayer{background layer}
\pgfsetlayers{background layer,main}

\tikzset{->-/.style={decoration={
  markings,
  mark=at position #1 with {\arrow{>}}},postaction={decorate}}}

\usetikzlibrary{intersections,through,hobby}

\usetikzlibrary{calc}

\renewcommand{\arraystretch}{1.2}

\allowdisplaybreaks

\DeclareMathSymbol{\widetildesym}{\mathord}{largesymbols}{"65}
\newcommand\lowerwidetildesym{%
  \text{\smash{\raisebox{-1.3ex}{%
    $\widetildesym$}}}}
\newcommand\fixwidetilde[1]{%
  \mathchoice
    {\accentset{\displaystyle\lowerwidetildesym}{#1}}
    {\accentset{\textstyle\lowerwidetildesym}{#1}}
    {\accentset{\scriptstyle\lowerwidetildesym}{#1}}
    {\accentset{\scriptscriptstyle\lowerwidetildesym}{#1}}
}

\nc\LiftCong[2]{\Cong^{#2}\hspace{-0.6truemm}(#1)}
\nc\colvert[3]{\fill (#1,#2)circle(.13); \fill[#3] (#1,#2)circle(.1);}
\nc\colvertx[4]{\fill (#1,#2)circle(#4+.03); \fill[#3] (#1,#2)circle(#4);}
\nc\smallcolvert[3]{\fill[#3] (#1,#2)circle(.1);}
\nc\redcolvert[3]{\fill[red] (#1,#2)circle(.13); \fill[#3] (#1,#2)circle(.1);}

\nc\OrdP[1]{\mathscr P^+(#1)}
\nc\OrdPP[1]{\mathscr P_{{\circ}}^+(#1)} 
\nc\OrdR[1]{\mathscr P^-(#1)}
\nc\OrdPR[1]{\mathscr P^\pm(#1)}
\nc\OriP[1]{\mathscr A^+(#1)}
\nc\OriR[1]{\mathscr A^-(#1)}
\nc\OriPR[1]{\mathscr A^\pm(#1)}

\nc\bd{{\bf d}}
\nc\br{{\bf r}}
\nc{\Sone}{S^1}
\nc\comp{\sim}
\nc\wt\fixwidetilde
\nc\nab{\wt\nabla}
\nc\itemit[1]{\item[\emph{(#1)}]}
\nc\ve\varepsilon
\nc\Lam\Lambda
\nc\IB{\mathcal{IB}}
\nc\ab{\overline{\al}}
\nc\bb{\overline{\be}}
\nc\Th\Theta
\nc\rev{\operatorname{rev}}
\nc\RR{\mathcal R}
\nc\E{\mathcal E}
\nc\G{\mathcal G}
\nc\Q{\mathcal Q}
\nc\PT{\mathcal{PT}}
\nc\A{\mathcal A}

\newcommand{\dottedintervalx}[5]{\draw(#1,#2+#5)--(#1,#2+#3) (#1,#2+#3+#4)--(#1,#2+#3+#4+#3-#5);\draw[dotted](#1,#2+#3)--(#1,#2+#3+#4);}
\nc\LLL{\mathbb L}
\nc\RRR{\mathbb R}
\nc\bn{{\bf n}}
\nc\btwo{{\bf 2}}
\nc\bthree{{\bf 3}}
\nc\Mod[1]{\ (\operatorname{mod}\ #1)}
\nc\even{{\operatorname{even}}}
\nc\odd{{\operatorname{odd}}}
\nc\TL{\mathcal T\!\mathcal L}
\nc\TLpm{\TL^\pm}
\nc\TLm{\TL^-}
\nc\Jpm{\JJ^\pm}
\newcommand{\darcxx}[4]{\draw[#4](#1,0)arc(180:90:#3) (#1+#3,#3)--(#2-#3,#3) (#2-#3,#3) arc(90:0:#3);}
\newcommand{\uarcxx}[4]{\draw[#4](#1,2)arc(180:270:#3) (#1+#3,2-#3)--(#2-#3,2-#3) (#2-#3,2-#3) arc(270:360:#3);}
\newcommand{\stlinex}[3]{\draw[#3](#1,2)--(#2,0);}
\newcommand{\alt}{\widetilde{\al}}
\newcommand{\stlines}[1]{{\foreach \x/\y in {#1} { \stline{\x}{\y} }}}
\nc\LL{\mathcal L}
\nc\Fx{\F^\times}
\nc\PL{P\LL}
\nc\Ker{\operatorname{Ker}}
\nc\vu{{\bf u}}
\nc\vv{{\bf v}}
\nc\vw{{\bf w}}
\nc\vzero{{\bf 0}}
\nc\GL{\operatorname{GL}}
\nc\SL{\operatorname{SL}}
\nc\PGL{\operatorname{PGL}}
\nc\tr{{\operatorname{T}}}
\nc\BB{\mathscr B}
\nc{\mat}[4]{\left[\begin{matrix}#1&#2\\#3&#4\end{matrix}\right]}
\nc{\tmat}[4]{\left[\begin{smallmatrix}#1&#2\\#3&#4\end{smallmatrix}\right]}
\nc\sicongnode[6] {\node[rounded corners,rectangle,draw,fill=blue!20] (#1#2) at (#3,#4) {$\si_{#5,#6}$};}
\nc\redsicongnode[6] {\node[rounded corners,rectangle,draw,fill=red!20] (#1#2) at (#3,#4) {$\si_{#5,#6}$};}
\nc\taucongnode[4] {\node[rounded corners,rectangle,draw,fill=red!20] (#1) at (#2,#3) {$\tau_{#4}$};}
\nc\congnode[4] {\node[rounded corners,rectangle,draw,fill=blue!20] (#3) at (#1,#2) {$#4$};}
\nc\redcongnode[4] {\node[rounded corners,rectangle,draw,fill=red!20] (#3) at (#1,#2) {$#4$};}
\nc\thickcongnode[4] {\node[rounded corners,rectangle,draw,fill=blue!20,ultra thick] (#3) at (#1,#2) {$#4$};}
\nc\thickredcongnode[4] {\node[rounded corners,rectangle,draw,fill=red!20,ultra thick] (#3) at (#1,#2) {$#4$};}
\nc\thicksicongnode[6] {\node[rounded corners,rectangle,draw,fill=blue!20,ultra thick] (#1#2) at (#3,#4) {$\si_{#5,#6}$};}
\nc\thickredsicongnode[6] {\node[rounded corners,rectangle,draw,fill=red!20,ultra thick] (#1#2) at (#3,#4) {$\si_{#5,#6}$};}
\nc\thicktaucongnode[4] {\node[rounded corners,rectangle,draw,fill=red!20,ultra thick] (#1) at (#2,#3) {$\tau_{#4}$};}
\nc\Znode[4] {\node[rounded corners,rectangle,draw,fill=blue!20] (Z#1#2) at (#3,#4) {$Z_{#1,#2}$};}
\nc\Snode[6] {\node[rounded corners,rectangle,draw,fill=green!20] (S#1#2) at (#3+#5,#4+#6) {$S_{#1,#2}$};}
\nc\squarelattice[2]{
\draw (#1,#2)--(#1+1,#2+1)--(#1,#2+2)--(#1-1,#2+1)--(#1,#2);
\foreach \x/\y in {#1/#2,#1+1/#2+1,#1/#2+2,#1-1/#2+1} {\colvert{\x}{\y}{white}}
}
\nc\understring[2]{\draw[thick] (#1,4) .. controls (#1,2) and (#2,2) .. (#2,0);}
\nc\overstring[2]{\draw[white,line width=2mm] (#1,4) .. controls (#1,2) and (#2,2) .. (#2,0); \draw[thick] (#1,4) .. controls (#1,2) and (#2,2) .. (#2,0);}
\nc\redoverstring[2]{\draw[white,line width=2mm] (#1,4) .. controls (#1,2) and (#2,2) .. (#2,0); \draw[red,thick] (#1,4) .. controls (#1,2) and (#2,2) .. (#2,0);}
\nc\halfunderstring[2]{\draw[thick] (#1,2) .. controls (#1,1) and (#2,1) .. (#2,0);}
\nc\halfoverstring[2]{\draw[white,line width=2mm] (#1,2) .. controls (#1,1) and (#2,1) .. (#2,0); \draw[thick] (#1,2) .. controls (#1,1) and (#2,1) .. (#2,0);}
\nc\redhalfoverstring[2]{\draw[white,line width=2mm] (#1,2) .. controls (#1,1) and (#2,1) .. (#2,0); \draw[red,thick] (#1,2) .. controls (#1,1) and (#2,1) .. (#2,0);}
\nc\understringx[4]{\draw[thick] (#1,#2) .. controls (#1,#2/2+#4/2) and (#3,#2/2+#4/2) .. (#3,#4);}
\nc\overstringx[4]{\draw[white,line width=2mm] (#1,#2) .. controls (#1,#2/2+#4/2) and (#3,#2/2+#4/2) .. (#3,#4); \draw[thick] (#1,#2) .. controls (#1,#2/2+#4/2) and (#3,#2/2+#4/2) .. (#3,#4);}
\nc\redoverstringx[4]{\draw[white,line width=2mm] (#1,#2) .. controls (#1,#2/2+#4/2) and (#3,#2/2+#4/2) .. (#3,#4); \draw[red,thick] (#1,#2) .. controls (#1,#2/2+#4/2) and (#3,#2/2+#4/2) .. (#3,#4);}
\nc\dotts[3]{\draw[dotted](#1+.35,#3)--(#2-.35,#3);}
\nc\braidrest[1]{^{\downarrow(#1)}}
\nc\V{\mathcal V}
\nc\DD{\mathcal D}
\nc\CC{\mathcal C}
\nc\bN{{\bf N}}
\nc\tran\ol
\nc\CongH{\Cong_{\H}}

\nc\overuarcx[3]{
\draw[white,line width=2mm](#1,4)arc(180:270:#3) (#1+#3,4-#3)--(#2-#3,4-#3) (#2-#3,4-#3) arc(270:360:#3);
\draw[thick](#1,4)arc(180:270:#3) (#1+#3,4-#3)--(#2-#3,4-#3) (#2-#3,4-#3) arc(270:360:#3);
}

\nc\overdarcx[3]{
\draw[white,line width=2mm](#1,0)arc(180:90:#3) (#1+#3,#3)--(#2-#3,#3) (#2-#3,#3) arc(90:0:#3);
\draw[thick](#1,0)arc(180:90:#3) (#1+#3,#3)--(#2-#3,#3) (#2-#3,#3) arc(90:0:#3);
}

\nc{\set}[2]{\{#1:#2\}}
\nc{\bigset}[2]{\big\{#1:#2\big\}}
\nc\sub\subseteq
\nc\mt\mapsto
\nc\al\alpha
\nc\be\beta
\nc\ga\gamma
\nc\de\delta
\nc\si\sigma
\nc\ep\epsilon
\nc\lam\lambda
\nc\ze\zeta
\nc\vk\varkappa
\nc\vt{\widetilde{\nu}}
\nc\vs\varsigma
\nc\ka\kappa
\rnc\th\theta
\nc\om\omega
\nc\ups\upsilon
\nc\Ga\Gamma
\nc\De\Delta
\nc\Si\Sigma
\nc\Om\Omega
\rnc\O{\mathcal O}
\nc\bit{\begin{itemize}[leftmargin=6mm,itemsep=1mm]}
\nc\eit{\end{itemize}}
\nc\ben{\begin{enumerate}[label=\textup{(\roman*)},leftmargin=7mm]}
\nc\bena{\begin{enumerate}[label=\textup{(\alph*)},leftmargin=7mm]}
\nc\een{\end{enumerate}}
\nc{\leqR}{\leq_{\R}}
\nc{\leqL}{\leq_{\L}}
\nc{\leqJ}{\leq_{\J}}
\nc{\leqK}{\leq_{\K}}
\nc{\geqR}{\geq_{\R}}
\nc{\geqL}{\geq_{\L}}
\nc{\geqJ}{\geq_{\J}}

\nc\bp{{\bf p}}
\nc\bq{{\bf q}}
\rnc\iff{\ \Leftrightarrow\ }
\rnc\implies{\ \Rightarrow\ }

\nc\pf{\begin{proof}}
\nc\epf{\end{proof}}
\nc\epfres{\hfill\qed}
\nc\epfreseq{\tag*{\qed}}

\let\oldproofname=\proofname
\renewcommand{\proofname}{\rm\bf{\oldproofname}}

\nc\AND{\qquad\text{and}\qquad}
\nc\WHERE{\qquad\text{where}\qquad}
\nc\ANd{\quad\text{and}\quad}
\nc\anD{\ \ \ \text{and}\ \ \ }
\nc\ANDSIM{\qquad\text{and similarly}\qquad}
\nc{\COMMA}{,\qquad}
\nc{\COMMa}{,\quad}
\nc\ul\underline
\nc\ol\overline
\nc\permdec[1]{#1^{\natural}}
\nc\ext[1]{#1^\textup{E}}
\nc\wh\widehat
\nc\normal\unlhd
\rnc\emptyset\varnothing
\nc{\pfitem}[1]{\medskip\noindent #1}
\nc{\firstpfitem}[1]{#1}
\nc{\pfcase}[1]{\medskip\noindent {\bf Case #1.}}
\nc{\pfstep}[1]{\medskip\noindent {\bf Step #1.}}
\nc\aftercases{\medskip\noindent}
\nc{\pfclaim}[1]{\medskip\noindent{\bf Claim #1.}}
\nc{\pfclaimnn}{\medskip\noindent{\bf Claim.} } 
\nc\afterclaim{\medskip}
\nc{\pfsubcase}[1]{\medskip\noindent {\bf Subcase #1.}}
\nc\im{\operatorname{im}}
\nc\LSUB{\operatorname{LSUB}}

\nc\pre\preceq
\nc\Y{\mathcal Y}
\nc\B{\mathcal B}
\nc\Z{\mathcal Z}
\nc\ZZ{\mathbb Z}
\nc\F{\mathbb F}
\nc\T{\mathcal T}
\nc\TT{\mathscr T}
\nc\PP{\mathscr P\mathcal P}
\nc\C{\mathscr C}
\nc\I{\mathcal I}
\nc\Eq{\mathfrak{Eq}}
\nc\Part{\mathbb{P}}
\nc\cg[2]{(#1,#2)^\sharp}
\nc\Rev{\operatorname{Rev}}
\nc\cR{\mathcal R}
\nc\Ptop{P^\top}
\nc\Qtop{Q^\top}
\nc\la\langle
\nc\ra\rangle
\nc\tb[1]{\operatorname{Seq}(#1)}
\nc\RevX{{\Rev}\big([0,|X|],[0,|X|^+]\big)}

\makeatletter
\DeclareRobustCommand\widecheck[1]{{\mathpalette\@widecheck{#1}}}
\def\@widecheck#1#2{%
    \setbox\z@\hbox{\m@th$#1#2$}%
    \setbox\tw@\hbox{\m@th$#1%
       \widehat{%
          \vrule\@width\z@\@height\ht\z@
          \vrule\@height\z@\@width\wd\z@}$}%
    \dp\tw@-\ht\z@
    \@tempdima\ht\z@ \advance\@tempdima2\ht\tw@ \divide\@tempdima\thr@@
    \setbox\tw@\hbox{%
       \raise\@tempdima\hbox{\scalebox{1}[-1]{\lower\@tempdima\box
\tw@}}}%
    {\ooalign{\box\tw@ \cr \box\z@}}}
\makeatother

\nc\ch\widecheck

\newcommand{\uv}[1]{\fill (#1,2)circle(.17);}
\newcommand{\uuv}[1]{\fill (#1,4)circle(.17);}
\newcommand{\uuuv}[1]{\fill (#1,3)circle(.17);}
\newcommand{\lv}[1]{\fill (#1,0)circle(.17);}
\newcommand{\uvs}[1]{{\foreach \x in {#1} { \uv{\x}}}}
\newcommand{\uuvs}[1]{{\foreach \x in {#1} { \uuv{\x}}}}
\newcommand{\uuuvs}[1]{{\foreach \x in {#1} { \uuuv{\x}}}}
\newcommand{\lvs}[1]{{\foreach \x in {#1} { \lv{\x}}}}
\newcommand{\darcx}[3]{\draw(#1,0)arc(180:90:#3) (#1+#3,#3)--(#2-#3,#3) (#2-#3,#3) arc(90:0:#3);}
\newcommand{\darc}[2]{\darcx{#1}{#2}{.4}}
\newcommand{\uarcx}[3]{\draw(#1,2)arc(180:270:#3) (#1+#3,2-#3)--(#2-#3,2-#3) (#2-#3,2-#3) arc(270:360:#3);}
\newcommand{\uarc}[2]{\uarcx{#1}{#2}{.4}}
\newcommand{\stline}[2]{\draw(#1,2)--(#2,0);}

\newcommand{\uustline}[2]{\draw(#1,3)--(#2,0);}

\nc\uur[3]{\fill[blue!20] (#1,4)--(#2,4)--(#2,4-#3)--(#1,4-#3); \draw[dotted] (#1,4)--(#2,4);}
\nc\uuur[3]{\fill[blue!20] (#1,3)--(#2,3)--(#2,3-#3)--(#1,3-#3); \draw[dotted] (#1,3)--(#2,3);}
\nc\ur[3]{\fill[blue!20] (#1,2)--(#2,2)--(#2,2-#3)--(#1,2-#3); \draw[dotted] (#1,2)--(#2,2);}
\nc\lr[3]{\fill[blue!20] (#1,0)--(#2,0)--(#2,#3)--(#1,#3); \draw[dotted] (#1,0)--(#2,0);}
\nc\lrr[3]{\fill[red!20] (#1,0)--(#2,0)--(#2,#3)--(#1,#3); \draw[dotted] (#1,0)--(#2,0);}

\nc{\buv}[1]{\fill (#1,2)circle(.18);}
\nc{\buvs}[1]{{
\foreach \x in {#1}
{ \buv{\x}}
}}
\nc{\blv}[1]{\fill (#1,0)circle(.18);}
\nc{\blvs}[1]{{
\foreach \x in {#1}
{ \blv{\x}}
}}

\nc{\uarcs}[1]{
{\foreach \x/\y in {#1}
{ \uarc{\x}{\y} }
}
}

\nc{\darcs}[1]{
{\foreach \x/\y in {#1}
{ \darc{\x}{\y} }
}
}
\nc{\darcxhalf}[3]{\draw(#1,0)arc(180:90:#3) (#1+#3,#3)--(#2,#3) ;}
\nc{\darchalf}[2]{\darcxhalf{#1}{#2}{.4}}
\nc{\uarcxhalf}[3]{\draw(#1,2)arc(180:270:#3) (#1+#3,1.5-#3)--(#2,1.5-#3) ;}
\nc{\uarchalf}[2]{\uarcxhalf{#1}{#2}{.4}}
\nc{\colv}[3]{\fill[#3] (#1,#2)circle(.17);}

\nc{\uvert}[1]{\fill (#1,2)circle(.2);}
\rnc{\lvert}[1]{\fill (#1,0)circle(.2);}

\nc{\custpartn}[3]{{\lower1.4 ex\hbox{
\begin{tikzpicture}[scale=.3]
\foreach \x in {#1}
{ \uvert{\x}  }
\foreach \x in {#2}
{ \lvert{\x}  }
#3 \end{tikzpicture}
}}}

\renewcommand{\P}{\mathcal P} 
\newcommand{\PB}{\mathcal{PB}} 
\newcommand{\JJ}{\mathcal{J}} 
\renewcommand{\S}{\mathcal{S}}
\newcommand{\M}{\mathcal{M}}

\nc\MYZ{\mathcal M_{Y\cup Z}}

\renewcommand{\H}{\mathscr H}
\renewcommand{\L}{\mathscr L}
\newcommand{\R}{\mathscr R}
\newcommand{\D}{\mathscr D}
\newcommand{\J}{\mathscr J}
\newcommand{\K}{\mathscr K}

\newcommand{\N}{\mathbb{N}}

\newcommand{\NN}{\mathcal{N}}
\nc\HH{\mathcal H}

\newcommand{\Cong}{\operatorname{Cong}}
\newcommand{\coker}{\operatorname{coker}}
\newcommand{\dom}{\operatorname{dom}} 
\newcommand{\codom}{\operatorname{codom}}
\newcommand{\rank}{\operatorname{rank}}

\newcommand{\id}{\operatorname{id}}

\newcommand{\Dmax}[1]{\mathcal{D}\textup{\textsf{max}}(#1)}
\newcommand{\ngen}[1]{\nabla\textup{\textsf{gen}}(#1)}

\newcommand{\Sep}[1]{\textup{\textsf{Sep}}(#1)}
\newcommand{\Sepb}[1]{\textup{\textsf{Sep}}^\flat(#1)}
\newcommand{\Sepz}[1]{\textup{\textsf{Sep}}^\ze(#1)}
\newcommand{\Mult}[1]{\textup{\textsf{Mult}}(#1)}
\newcommand{\Multb}[1]{\textup{\textsf{Mult}}^\flat(#1)}
\newcommand{\Multz}[1]{\textup{\textsf{Mult}}^\ze(#1)}

\renewcommand{\c}{@{}c@{}}
\newcommand{\cend}{@{}c@{\hspace{1.5truemm}}}
\newcommand{\cstart}{@{\hspace{1.5truemm}}c@{}}
\newcommand{\cstartend}{@{\hspace{1.5truemm}}c@{\hspace{1.5truemm}}}

\newcommand{\partn}[4]{
\Big(   \hspace{-1.5 truemm}
{ \scriptsize \renewcommand*{\arraystretch}{1}
\begin{array} {\cstart|\cend}
 #1 \:&\: #2 \\ \cline{2-2}
 #3 \:&\: #4 
\rule[0mm]{0mm}{2.7truemm}
\end{array} 
}
\hspace{-1.5 truemm} \Big)
}

\newcommand{\partI}[8]{
\Big( 
{ \scriptsize \renewcommand*{\arraystretch}{1}
\begin{array} {\c|\c|\c|\c|\c|\cend}
 #1 \:&\: \cdots \:&\: #2 \:&\: #3 \:&\: \cdots \:&\: #4 \\ \cline{4-6}
 #5 \:&\: \cdots \:&\: #6 \:&\: #7 \:&\: \cdots \:&\: #8 
\rule[0mm]{0mm}{2.7mm}
\end{array} 
}
\hspace{-1.5 truemm} \Big)
}

\newcommand{\sqpartI}[8]{
\Big[ 
{ \scriptsize \renewcommand*{\arraystretch}{1}
\begin{array} {\c|\c|\c|\c|\c|\cend}
 #1 \:&\: \cdots \:&\: #2 \:&\: #3 \:&\: \cdots \:&\: #4 \\ \cline{4-6}
 #5 \:&\: \cdots \:&\: #6 \:&\: #7 \:&\: \cdots \:&\: #8 
\rule[0mm]{0mm}{2.7mm}
\end{array} 
}
\hspace{-0.8 truemm} \Big]
}

\newcommand{\sqpartII}[8]{
\Big[ 
{ \scriptsize \renewcommand*{\arraystretch}{1}
\begin{array} {\c|\c|\c|\c|\cend}
 #1 \:&\: \cdots \:&\: #2 \:&\: #3 \:&\: #4 \\ \cline{5-5}
 #5 \:&\: \cdots \:&\: #6 \:&\: #7 \:&\: #8 
\rule[0mm]{0mm}{2.7mm}
\end{array} 
}
\hspace{-0.8 truemm} \Big]
}

\newcommand{\partII}[6]{
\Big(  \hspace{-1.5 truemm}
{\scriptsize \renewcommand*{\arraystretch}{1} \begin{array} {\cstart|\c|\cend}
#1 \:&\: #2 \:&\: #3  \\ \cline{1-3}
#4 \:&\: #5 \:&\: #6 
\rule[0mm]{0mm}{2.7mm}
\end{array}  }
\hspace{-1.5 truemm} \Big) 
}

\newcommand{\partIII}[8]{
\Big(  \hspace{-1.5 truemm}
{\scriptsize \renewcommand*{\arraystretch}{1} \begin{array} {\cstart|\c|\c|\cend}
#1 \:&\: #2 \:&\: #3 \:&\: #4 \\ \cline{1-4}
#5 \:&\: #6 \:&\: #7 \:&\: #8
\rule[0mm]{0mm}{2.7mm}
\end{array}  }
\hspace{-1.5 truemm} \Big) 
}

\newcommand{\partIV}[8]{
\Big(  
{\scriptsize \renewcommand*{\arraystretch}{1} \begin{array} {\c|\c|\c|\cend}
#1 \:&\: #2 \:&\: #3 \:&\: #4 \\ \cline{2-4}
#5 \:&\: #6 \:&\: #7 \:&\: #8
\rule[0mm]{0mm}{2.7mm}
\end{array}  }
\hspace{-1.5 truemm} \Big) 
}

\newcommand{\partV}[8]{
\Big(  
{\scriptsize \renewcommand*{\arraystretch}{1} \begin{array} {\c|\c|\c|\c|\cend}
#1 \:&\: #2 \:&\: #3 \:&\: \cdots \:&\: #4 \\ \cline{3-5}
#5 \:&\: #6 \:&\: #7 \:&\: \cdots \:&\: #8
\rule[0mm]{0mm}{2.7mm}
\end{array}  }
\hspace{-1.5 truemm} \Big) 
}

\newcommand{\partVI}[8]{
\Big(  
{\scriptsize \renewcommand*{\arraystretch}{1} \begin{array} {\c|\c|\c|\c|\cend}
#1 \:&\: #2 \:&\: #3 \:&\: \cdots \:&\: #4 \\ \cline{2-5}
#5 \:&\: #6 \:&\: #7 \:&\: \cdots \:&\: #8
\rule[0mm]{0mm}{2.7mm}
\end{array}  }
\hspace{-1.5 truemm} \Big) 
}

\newcommand{\partVII}[8]{
\Big(  \hspace{-1.5 truemm}
{\scriptsize \renewcommand*{\arraystretch}{1} \begin{array} {\cstart|\c|\c|\c|\cend}
#1 \:&\: #2 \:&\: #3 \:&\: \cdots \:&\: #4 \\ \cline{1-5}
#5 \:&\: #6 \:&\: #7 \:&\: \cdots \:&\: #8
\rule[0mm]{0mm}{2.7mm}
\end{array}  }
\hspace{-1.5 truemm} \Big) 
}

\newcommand{\sqpartXX}[2]{
\Big[  \hspace{-0.8 truemm}
{\scriptsize \renewcommand*{\arraystretch}{1} \begin{array} {\cstartend}
#1   \\ \cline{1-1}
#2
\rule[0mm]{0mm}{2.7mm}
\end{array}  }
\hspace{-0.8 truemm} \Big]
}

\newcommand{\sqpartXXIII}[3]{
\Big[  \hspace{-0.9 truemm}
{\scriptsize \renewcommand*{\arraystretch}{1} \begin{array} {\cstart|\c|\cend}
#1 \:&\: #2 \:&\: #3  \\ \cline{1-3}
\multicolumn{3}{c}{}
\rule[0mm]{0mm}{2.7mm}
\end{array}  }
\hspace{-0.7 truemm} \Big] 
}

\newcommand{\partXXIV}[6]{
\Big( 
{ \scriptsize \renewcommand*{\arraystretch}{1}
\begin{array} {\c|\c|\cend}
 #1 \:&\: #2 \:&\: #3 \\ \cline{3-3}
 #4 \:&\: #5 \:&\: #6
\rule[0mm]{0mm}{2.7mm}
\end{array} 
}
\hspace{-1.5 truemm} \Big)
}

\newcommand{\partpermIII}[6]{
\Big(
{ \scriptsize \renewcommand*{\arraystretch}{1}
\begin{array} {\c|\c|\c}
 #1 \:&\: #2 \:&\: #3 \\ 
 #4 \:&\: #5 \:&\: #6 
\rule[0mm]{0mm}{2.7mm}
\end{array} 
}
\Big)
}

\newcommand{\sqpartpermIII}[6]{
\Big[
{ \scriptsize \renewcommand*{\arraystretch}{1}
\begin{array} {\c|\c|\c}
 #1 \:&\: #2 \:&\: #3 \\ 
 #4 \:&\: #5 \:&\: #6 
\rule[0mm]{0mm}{2.7mm}
\end{array} 
}
\Big]
}

\newcommand{\partpermIV}[6]{
\Big(
{ \scriptsize \renewcommand*{\arraystretch}{1}
\begin{array} {\c|\c|\c|\c}
 #1 \:&\: \cdots \:&\: #2 \:&\: #3 \\ 
 #4 \:&\: \cdots \:&\: #5 \:&\: #6 
\rule[0mm]{0mm}{2.7mm}
\end{array} 
}
\Big)
}

\newcommand{\partpermV}[6]{
\Big(
{ \scriptsize \renewcommand*{\arraystretch}{1}
\begin{array} {\c|\c|\c|\c}
 #1 \:&\: #2 \:&\: \cdots \:&\: #3 \\ 
 #4 \:&\: #5 \:&\: \cdots \:&\: #6 
\rule[0mm]{0mm}{2.7mm}
\end{array} 
}
\Big)
}

\nc\congsquare[4]
{\fill[black!15] (#1*6,#2*6,#3*8)--({#1*6+(#4)*6},#2*6,#3*8)--({#1*6+(#4)*6},{#2*6+(#4)*6},#3*8)--(#1*6,{#2*6+(#4)*6},#3*8)--(#1*6,#2*6,#3*8);
\foreach \x in {0,1,...,#4} {\draw (#1*6,#2*6+\x*6,#3*8)--({#1*6+(#4)*6},#2*6+\x*6,#3*8); \draw (#1*6+\x*6,#2*6,#3*8)--(#1*6+\x*6,{#2*6+(#4)*6},#3*8);}}

\nc\congsquareconnections[7]
{\foreach \x in {0,...,#7} \foreach \y in {0,...,#7} {\draw(#1*6+6*\x,#2*6+6*\y,#3*8)--(#4*6+6*\x,#5*6+6*\y,#6*8);}}

\numberwithin{equation}{section}

\newtheorem{thm}[equation]{Theorem}
\newtheorem{lemma}[equation]{Lemma}
\newtheorem{cor}[equation]{Corollary}
\newtheorem{prop}[equation]{Proposition}

\theoremstyle{definition}

\newtheorem{rem}[equation]{Remark}
\newtheorem{defn}[equation]{Definition}

\newcommand{\sm}{\setminus}
\newcommand{\restr}{{\restriction}}

\begin{document}

\title{\vspace{-1.4cm}Congruence lattices of ideals in categories and (partial) semigroups}
\author{James East\footnote{Centre for Research in Mathematics, School of Computing, Engineering and Mathematics, Western Sydney University, Locked Bag 1797, Penrith NSW 2751, Australia. {\it Email:} {\tt j.east\,@\,westernsydney.edu.au}} 
\ and
Nik Ru\v{s}kuc\footnote{Mathematical Institute, School of Mathematics and Statistics, University of St Andrews, St Andrews, Fife KY16 9SS, UK. {\it Email:} {\tt nik.ruskuc\,@\,st-andrews.ac.uk}}}
\date{}

\maketitle

\vspace{-1cm}
\begin{abstract}
\noindent
This paper presents a unified framework for determining the congruences on a number of monoids and categories of transformations, diagrams, matrices and braids, and on all their ideals.
The key theoretical advances present an iterative process of stacking certain normal subgroup lattices on top of each other to successively build congruence lattices of a chain of ideals.  This is applied to several specific categories of: transformations; order/orientation preserving/reversing transformations; partitions; planar/annular partitions;
Brauer, Temperley--Lieb and Jones partitions; linear and projective linear transformations; and partial braids.
Special considerations are needed for certain small ideals, and technically more intricate theoretical underpinnings for the linear and partial braid categories.
\medskip

\noindent
\emph{Keywords}: 
Categories; semigroups; congruences; $\H$-congruences; lattices; ideals; chains of ideals; diagram categories; partition categories; Brauer categories; Temperley-Lieb categories; Jones categories; transformation categories; linear categories; projective linear categories; partial braid categories.
\medskip

\noindent
MSC: 20M20 (primary), 08A30, 15A04, 18B40, 20M10, 20M50.

\end{abstract}

\setcounter{tocdepth}{1}

\tableofcontents

\section{Introduction}\label{sect:intro}

A congruence on an algebraic structure is an equivalence relation  compatible with its operations.  Congruences are therefore in one-one correspondence with quotient structures, and are the kernels of homomorphisms/representations.  The set $\Cong(S)$ of all congruences on a structure $S$ is a (complete, algebraic) lattice, called the congruence lattice of $S$.  For some structures, congruences can be described by means of special substructures, such as normal subgroups of groups, or ideals of rings, but this is not the case in general.

The structures we are mostly concerned with here are (small) categories, viewed as sets of morphisms equipped with the (partial) operation of composition.  In fact, it is both convenient and necessary to broaden our scope to the wider class of partial semigroups, which may be thought of as categories in which identities need not exist at every object.  As special cases, monoids and semigroups occur when there is only a single object.

The last 70 years have seen many important studies of congruence lattices of monoids and semigroups.  Arguably the most influential is Mal'cev's seminal paper \cite{Malcev1952}, which completely classifies the congruences on the full transformation monoid $\T_X$: i.e., the monoid of all self-maps of the set $X$, under composition.  In the case of finite $X$, congruences are elegantly described in terms of the ideals of $\T_X$, and the normal subgroups of the symmetric groups $\S_n$, $n\leq|X|$.  Mal'cev's follow-up paper \cite{Malcev1953} carries out the same program for the full linear monoid $\M_n(\F)$, consisting of all $n\times n$ matrices over the field $\F$ (equivalently, linear transformations of an $n$-dimensional vector space over $\F$).  The description here is more elaborate than for $\T_X$, but still involves ideals and normal subgroups, this time of appropriate general linear groups.  A crucial structural feature shared by the monoids $\T_X$ and $\M_n(\F)$ is that their ideals form chains under inclusion.

Mal'cev's two papers \cite{Malcev1952,Malcev1953} set the direction for a number of subsequent studies, and his ideas and techniques have been developed and extended to many other important families of monoids; see for example \cite{Aizenstat1962,FGJ2009,FGJ2005,ABG2016,Fernandes2001, Liber1953,MSS2000,Scheiblich1973,Sutov1961,Sutov1961_2,EMRT2018,ER2018,Klimov1977}.  To the authors' knowledge, much less is known about congruences on categories.  For example, it seems that the congruences of categories of (finite) sets and mappings have not previously been classified.  
Accordingly, one of the main motivations for the current work is to initiate the study of congruence lattices of categories and partial semigroups,
including categories of (partial) mappings, partitions/diagrams, linear transformations, and (partial) braids.  Categories of mappings and (linear) transformations are of course among the most fundamental and well-studied categories.  
Diagram and braid categories play a crucial role across an entire spectrum of mathematical and scientific disciplines, including classical groups, representation theory, topology, theoretical physics, knot theory, logic, invariant theory, combinatorics, theoretical computer science
and many more \cite{BH2019,HR2005,Brauer1937,TL1971,Jones1994,Jones2001,Jones1994_a,Jones1987,Jones1983_2,GL1998,Abramsky2008,DP2013,DW2000, BS2011, KX2001,Xi1999,Wenzl1988,Wilcox2007,Martin2015, Kauffman1990,Kauffman1987,BD1995,Stroppel2005,LZ2017,LZ2012,LZ2015,Martin1991,Martin1994,Martin2008,
JS1993,Baez1994,VZ1998,FY1989,BN1996}.

The second source of motivation for the present work is to do with ideals.
We have already noted that the ideals of a semigroup are crucially involved in its congruences.  On the other hand, any such ideal is a semigroup in its own right, and many interesting studies have investigated the ideals of important families of monoids
 \cite{EG2017,DEG2017,Gray2008,Howie1990,Gray2007,YY2004,Howie1978,Klimov1977,Malcev1953}.  
In fact, Mal'cev's paper \cite{Malcev1953} classified the congruences on all ideals of the linear monoid $\M_n(\F)$, including $\M_n(\F)$ itself.  Congruences on ideals of the full transformation semigroups $\T_X$ were described some years later by Klimov \cite{Klimov1977}.
Formally, Mal'cev's classification of the congruences on~$\T_X$ \cite{Malcev1952} is a special case of Klimov's result.
In actual fact, however,
Klimov's proof relies crucially on Mal'cev's classification: to understand the congruences 
on the ideals, it was necessary to first understand the congruences on the semigroup itself.

To the best of our knowledge, no further results exist classifying the congruences on ideals of natural semigroups.
Here we give such classifications for (many) other families of semigroups.  In fact, we do so in far greater generality, classifying the congruences on all ideals of many \emph{categories} of mappings, partitions/diagrams, linear transformations, and braids.
In several cases, the congruences on the whole category were not previously known, even in the special case of monoids.

The approach we develop works in the opposite direction to Klimov's, as described above.  
In all of our applications, the category $U$ in question is the union of an ascending chain of ideals, $I_0\subset I_1\subset\cdots$, of order-type at most~$\om$.
The goal is to classify the congruences on each ideal~$I_r$, and to describe the lattice $\Cong(I_r)$. 

However, rather than starting with $U$ itself, and then (somehow) deducing results about the ideals $I_r$, we start at the other end with the smallest ideal(s).  In each case, we find that the congruences on the first few ideals are characterised by great freedom; in some extreme cases, every equivalence is a congruence.  There is then a certain ``cut-off'' ideal $I_k$ (generally with $k\leq2$), whose congruences are far more restricted and well behaved, and from which point it is possible to iteratively compute the congruences of the ideals $I_{k+1},I_{k+2},\ldots$, building up eventually to the entire category $U$.  
In this way, classifying congruences of ideals can be seen as 
integral to computing the congruences of $U$ itself, and to providing an explanation for the structure of $\Cong(U)$. 

In Sections \ref{sect:T}--\ref{sect:J}, \ref{sect:L} and \ref{sect:IB}, we carry out the above program for 
some
thirty important families of categories and monoids.  In each case we classify the congruences on every ideal, including the category/monoid itself, and completely describe the corresponding lattices.  Each such classification is obtained by applying general theoretical results, proved in Sections \ref{sect:Pre}--\ref{sect:chains} and~\ref{sect:H}.

Let us now give a more detailed summary of the contents of the paper.

As mentioned above,
Sections \ref{sect:Pre}--\ref{sect:chains} are theoretical in nature, and underpin 
the computations undertaken in the remainder of the paper.
Section \ref{sect:Pre} contains background on categories and partial semigroups, and their ideals and congruences.  Basic definitions are given in Section \ref{subsect:S}, and then a number of general congruence constructions are outlined in Section \ref{subsect:EMRT}, adapting and extending the theory developed in \cite{EMRT2018} to the categorical context.  Several useful lemmas are proved in Section \ref{subsect:tricks}.

Section \ref{sect:IE} concerns ideal extensions, the main result being Theorem~\ref{thm:E6}.  This theorem shows how the congruence lattice of an extension $T$ of $S$ may be described in terms of the so-called liftable congruences on $S$ and the normal subgroups of a certain group in $T\sm S$, under a series of natural, 
albeit technical, assumptions.  We introduce and analyse these assumptions in Section~\ref{subsect:properties}, and discuss liftable congruences in Section \ref{subsect:lift}.  Theorem \ref{thm:E6} is then proved in Section~\ref{subsect:IE}, where we also describe how to visualise the lattice $\Cong(T)$ in terms of the above-mentioned data; see in particular Remarks~\ref{rem:E1} and~\ref{rem:E3}, and Figures~\ref{fig:E1}--\ref{fig:E32}.

The results of Section \ref{sect:IE} are then applied in Section \ref{sect:chains} to prove a number of results concerning congruences on certain general classes of partial semigroups, and on their ideals.  Such a partial semigroup $U$ is the union of an ascending chain of ideals $I_0\subset I_1\subset\cdots$, of order-type at most~$\om$.
The main results of Section \ref{sect:chains} describe the congruences on each ideal $I_r$, including $U$ itself.  Sections \ref{subsect:PRB}--\ref{subsect:small} do this for ``small''~$r$; see Propositions~\ref{prop:PRB},~\ref{prop:small_ideal},~\ref{prop:small01} and~\ref{prop:small012}.
Section \ref{subsect:chains} deals with ``large'' $r$, the main result being Theorem~\ref{thm:E2}.  This theorem assumes that $\Cong(I_k)$ is known, and well-behaved, for a certain ``cut-off'' value $k$ (which is always very small in our applications), and shows how to construct $\Cong(I_r)$, for $r\geq k$, out of $\Cong(I_k)$ and a sequence of groups $G_{k+1},G_{k+2},\ldots$ outside of $I_k$.  A crucial consequence of the theorem is that the lattice $\Cong(I_r)$, for $r>k$, is obtained by stacking on top of the lattice $\Cong(I_k)$ copies of the lattices of normal subgroups of $G_{k+1},G_{k+2},\ldots$ in a way made precise in Remark \ref{rem:E2}; see especially Figure~\ref{fig:E2}.  This ``stacking mechanism'' is one of the most important advances of the paper, and it can be seen in many of the lattice diagrams, such as Figures \ref{fig:T}, \ref{fig:QT}, \ref{fig:P}, \ref{fig:QP}, \ref{fig:TL}, \ref{fig:J} and more.  The meaning of the colours in such lattice diagrams, and other visualisation conventions, are discussed in Section~\ref{subsect:viscon}.

Sections \ref{sect:T}--\ref{sect:J} present a sequence of applications of the general theory developed in Sections~\mbox{\ref{sect:Pre}--\ref{sect:chains}}.  We begin in Section \ref{sect:T} with the transformation category $\T=\T(\C)$.  Here $\C$ is a set of non-empty finite sets, and $\T$ consists of all mappings between members of $\C$.  We give the definitions in Section \ref{subsect:defnT}, and then in Section \ref{subsect:prelimT} establish certain multiplicative properties that will be used to verify the relevant assumptions listed in Sections~\ref{sect:IE} and \ref{sect:chains}; most importantly,~$\T$ is a chain of ideals $I_1\subset I_2\subset\cdots$ as in Definition \ref{defn:chain}.
All of this will enable an application of Theorem \ref{thm:E2} (and Propositions \ref{prop:PRB} and \ref{prop:small01}) in the category $\T$, culminating in Theorem \ref{thm:T}, which completely describes the congruences on $\T$, and on all of its ideals.  The congruences of the minimal ideal $I_1$ have rather different behaviour to those of the larger ideals $I_r$ ($r\geq2$); this can be seen by comparing parts \ref{T1} and \ref{T2} of Theorem \ref{thm:T}, and also by examining Figure \ref{fig:T}, which pictures the congruence lattices of relatively small ideals of $\T$.  
In the above language, $I_2$ is the cut-off ideal from which the stacking mechanism kicks in.
In the case that $\C$ consists of a single set of size $n$, the category $\T$ is simply the full transformation semigroup of degree $n$, and so we recover the results of Mal'cev \cite{Malcev1952} for $r=n$, and of Klimov \cite{Klimov1977} for arbitrary $r$.  In Section~\ref{subsect:QT} we consider several subcategories of $\T$, which we call ``planar and annular reducts''; these consist of transformations preserving and/or reversing the order or orientation on the underlying sets; such mappings may also be characterised geometrically in terms of certain graphs with no crossings in a given planar or cylindrical region.  The main result here is Theorem \ref{thm:QT}, which describes the congruences on these reduct categories and on their ideals; see also Figure \ref{fig:QT}.  In Section \ref{subsect:OtherT} we discuss how to obtain similar results for several other categories of (partial) mappings.

The applications considered in Sections \ref{sect:P}--\ref{sect:J} all concern certain diagram categories.  The largest of these, the partition category $\P=\P(\C)$, is treated in Section \ref{sect:P}.  The structure of Section \ref{sect:P} closely parallels that of Section \ref{sect:T}.  We give the definitions in Section \ref{subsect:defnP}, prove the requisite multiplicative properties in Section \ref{subsect:prelimP}, showing in particular that $\P$ is again a chain of ideals $I_0\subset I_1\subset\cdots$.  In Section~\ref{subsect:P} we prove the main result, Theorem \ref{thm:P}, which classifies the congruences on $\P$ and on all of its ideals; Figure~\ref{fig:P} depicts the congruence lattice of a sufficiently large ideal.  This time 
the congruences for
the smallest \emph{two} ideals, $I_0$ and $I_1$, must be described separately (using Propositions \ref{prop:PRB} and \ref{prop:small_ideal}), and the stacking mechanism kicks in at $I_2$.  In Section~\ref{subsect:QP} we classify the congruences of planar and annular reducts of~$\P$, and of the ideals of these reducts; see Theorem \ref{thm:QP} and Figure \ref{fig:QP}.  In Section \ref{subsect:PB} we indicate how the methods of Section~\ref{sect:P} also apply to other subcategories of $\P$, namely the partial Brauer category and its planar and annular reducts, including the Motzkin category.

Section \ref{sect:B} concerns the Brauer category $\B=\B(\C)$, another natural subcategory of the partition category.  The ideals of $\B$ do not form a chain in general, so the general machinery developed in Sections \ref{sect:IE} and \ref{sect:chains} does not immediately apply.  However, it is possible to split the category $\B$ into even and odd parts, in a way made precise in Section \ref{subsect:defnB}, and we may then restrict our attention to these even and odd subcategories, in which ideals do form chains, $I_0\subset I_2\subset\cdots$ and $I_1\subset I_3\subset\cdots$.  After proving the required multiplicative properties in Section~\ref{subsect:prelimB}, we treat the odd and even cases in Sections \ref{subsect:Bodd} and \ref{subsect:Beven}, respectively, the main results being Theorems~\ref{thm:Bodd} and~\ref{thm:Beven}; see also~Figures \ref{fig:Bodd}--\ref{fig:B4}.  
There are substantial differences between the even and odd cases, including the height at which the stacking mechanism kicks in.

The planar and annular reducts of $\B$ include well-known categories such as the Temperley--Lieb and Jones categories, $\TL$ and $\JJ$.  Unlike the situation with the transformation and partition categories, however, it is not possible to simply modify the arguments for $\B$ to work for its reducts.  Accordingly, we treat the planar and annular reducts of $\B$ in Sections \ref{sect:TL} and \ref{sect:J}, respectively.  The main results of Section \ref{sect:TL} are Theorems~\ref{thm:TLodd},~\ref{thm:TLeven} and~\ref{thm:TLpmeven}, concerning the Temperley--Lieb and anti-Temperley--Lieb categories; see also Figure \ref{fig:TL}.  
The main results of Section \ref{sect:J} are Theorems~\ref{thm:Jodd} and~\ref{thm:Jeven}, concerning the Jones and anti-Jones categories; see also Figures \ref{fig:J}--\ref{fig:Jpmeven}.  Again, each of these categories splits into even and odd parts, which exhibit rather different behaviour.  Many of the arguments of Sections~\ref{sect:TL} and~\ref{sect:J} require a fairly intricate geometrical analysis, due to the planarity constraints.

By taking $\C$ to be a single set of size $n$, and taking $r=n$, we obtain new proofs of the main results of \cite{EMRT2018} as special cases of the theorems presented in Sections \ref{sect:P}--\ref{sect:TL}.  The (annular) Jones monoids were not considered in \cite{EMRT2018}.

While the remaining categories considered in the paper are still chains of ideals, some of them fail to satisfy certain assumptions built into the general machinery developed in Section~\ref{sect:chains}.  Thus, in Section \ref{sect:H} we undertake another theoretical development, proving more general results that will allow us to treat our remaining examples.  The main idea is that of a so-called $\H$-congruence, by which we mean a congruence contained in Green's $\H$-relation.
The $\H$-congruences of a certain class of partial semigroups are classified in Section \ref{subsect:H}; see Proposition~\ref{prop:ThN}.  Returning to chains of ideals in Section \ref{subsect:chains2}, we prove the main result, Theorem~\ref{thm:E4}, which is a strengthening of Theorem \ref{thm:E2} to a wider class of 
chains of ideals.
We also prove Theorem~\ref{thm:E5}, which is a simplification of Theorem \ref{thm:E4} particularly suited to our remaining examples.

Section \ref{sect:L} concerns certain categories of linear transformations, namely the linear and projective linear categories, $\LL=\LL(\C)$ and $\PL=\PL(\C)$, defined with respect to a set $\C$ of finite-dimensional vector spaces over an arbitrary 
fixed field $\F$.  Following the structure of previous sections, we give definitions in Section \ref{subsect:defnL}, establish multiplicative properties in Section~\ref{subsect:prelimL}, and then prove the main results in Sections \ref{subsect:L} and \ref{subsect:PL}.  Theorem \ref{thm:L} describes the congruences on the linear category $\LL$, and on all of its ideals; Theorem \ref{thm:PL} does the same for the projective linear category $\PL$.  While the general results of Section \ref{sect:chains} are sufficient to treat~$\PL$, we need the stronger results of Section \ref{sect:H} for $\LL$.  Crucially, the $\H$-congruences on $\LL$ and its ideals involve certain sequences of normal subgroups of (centres of) general linear groups of degrees bounded by the supremum of the dimensions of the spaces from $\C$.
Section~\ref{subsect:visual_L} explains how to visualise the congruence lattices of ideals of~$\LL$.  These lattices are far more complicated than for any of the categories considered in Sections~\ref{sect:T}--\ref{sect:J}; see especially Figures~\ref{fig:GLM27} and~\ref{fig:GLM37}, and compare these to Figures \ref{fig:T}, \ref{fig:QT}, \ref{fig:P}, \ref{fig:QP}--\ref{fig:B4}, and so on.  As a special case of our results on~$\LL$, we obtain Mal'cev's classification \cite{Malcev1953} of the congruences on the ideals of the full linear monoid $\M_n(\F)$.

Finally, Section \ref{sect:IB} applies the general results of Section \ref{sect:H} to categories $\IB=\IB(\C)$ of partial braids, the main result being Theorem \ref{thm:IB}.  This time, $\H$-congruences correspond to sequences of normal subgroups of (pure) braid groups of appropriate dimensions.  As a special case, we obtain a classification of the congruences of the inverse braid monoid $\IB_n$ \cite{EL2004}; to the best of our knowledge, such a classification has never before been established.

Before we begin, we briefly fix notation for some standard mathematical concepts.
We work in standard ZFC set theory; see for example \cite{Jech2003}.  We write $\N=\{0,1,2,\ldots\}$ for the set of natural numbers, and $\om$ for the first infinite ordinal/cardinal (which may of course be identified with~$\N$).  For $n\in\N$ we write $[n]=\{1,\ldots,n\}$, interpreting $[0]=\emptyset$.  

For a group $G$, we write $H\leq G$ and $N\normal G$ to indicate that $H$ is a subgroup and $N$ a normal subgroup.  The set $\NN(G)$ of all normal subgroups of $G$ is a complete lattice under inclusion.  The meet of a family $N_q\normal G$, $q\in Q$, is their intersection, $\bigwedge_qN_q=\bigcap_qN_q$; the join $\bigvee_qN_q=\big\la\bigcup_qN_q\big\ra$ is the (normal) subgroup generated by their union.  

If $\si$ is an equivalence relation on a set $X$, we will write $x\mathrel\si y$ or $(x,y)\in\si$ to indicate~$\si$-equivalence of $x,y\in X$.  If $Y\sub X$ is a union of $\si$-classes, we write $Y/\si$ for the set of all $\si$-classes contained in~$Y$.  For any $Y\sub X$, we write $\si\restr_Y=\si\cap(Y\times Y)$ for the restriction of $\si$ to~$Y$.  The set~$\Eq(X)$ of all equivalences on $X$ is a lattice under inclusion.  The meet of a family $\si_q\in\Eq(X)$, $q\in Q$, is their intersection; the join of the family is the transitive closure of their union.  When $X=[n]$ for some $n\in\N$, we will sometimes write $\Eq_n=\Eq(X)$.  We adopt similar conventions for other families of lattices, groups, etc.  Thus, for example, we will write $\S_X$ for the symmetric group on the set $X$, and in particular write $\S_n=\S_X$ when $X=[n]$; we write $\A_n$ for the alternating group.

\section{Categories and partial semigroups}\label{sect:Pre}

This section contains preliminary material on categories and partial semigroups, and their ideals and congruences.  The basic definitions are given in Section \ref{subsect:S}, where we also prove some preliminary results.  In Section \ref{subsect:EMRT} we adapt and extend a number of constructions and results of \cite[Section~3]{EMRT2018} to the context of partial semigroups.  In Section \ref{subsect:tricks}, we prove several lemmas that will be used frequently in the sequel.

\subsection{Preliminaries}\label{subsect:S}

We follow \cite{DE2018}, slightly adapting notation to suit our purposes.  
A \emph{partial semigroup} is a 5-tuple $(S,\C,\bd,\br,\cdot)$, where $S$ and $\C$ are sets, $\bd$ and $\br$ are mappings $S\to \C$, and $\cdot$ is a partial binary operation $(x,y)\mt x\cdot y$ such that, for all $x,y,z\in S$,~
\bit
\item $x\cdot y$ is defined if and only if $\br(x)=\bd(y)$, in which case $\bd(x\cdot y)=\bd(x)$ and $\br(x\cdot y)=\br(y)$,
\item if $x\cdot y$ and $y\cdot z$ are defined, then  $(x\cdot y)\cdot z=x\cdot (y\cdot z)$.
\eit
We usually abbreviate $x\cdot y$ to $xy$.  
The elements of $\C$ are called \emph{objects}, and the elements of~$S$ \emph{morphisms}.
For $A,B\in \C$ we write $S_{A,B}=\set{x\in S}{\bd(x)=A,\ \br(x)=B}$, and call these \emph{hom-sets}.  For $A\in \C$, we write $S_A=S_{A,A}$, noting that these are (possibly empty) semigroups under the induced operation; we call these \emph{endomorphism semigroups}.  Any semigroup is trivially a partial semigroup with $|\C|=1$.  
If the set~$\C$, the functions $\bd,\br$, and the operation~$\cdot$ are clear from context, we will often refer to ``the partial semigroup $S$'' instead of ``the partial semigroup $(S,\C,\bd,\br,\cdot)$''.

We say that a partial semigroup $S$ is a \emph{(small) category} if there exist elements $e_A\in S_A$~(${A\in\C}$) such that, for all $A,B\in\C$ and $x\in S_{A,B}$, $e_A\cdot x=x=x\cdot e_B$.  
In a category, the endomorphism semigroups~$S_A$ are monoids (with identities~$e_A$), and the elements of $S_A$ that are invertible with respect to the identity~$e_A$ are called \emph{units} (or \emph{automorphisms}).  As explained in \cite{DE2018}, any partial semigroup $S$ embeds in a category $\Sone$ in an obvious way; the category $\Sone$ is simply $S$ with as many identities adjoined as required.

An element $x$ of a partial semigroup $S$ is an \emph{idempotent} if $x=x^2$; note that this forces $\bd(x)=\br(x)$.  An element $x\in S$ is \emph{regular} if $x=xax$ for some $a\in S$; note that this forces $\bd(a)=\br(x)$ and $\br(a)=\bd(x)$, and note that with $b=axa$ we have $x=xbx$ and $b=bxb$.  We say $S$ itself is regular if all of its elements are regular; all of the partial semigroups we consider in Sections \ref{sect:T}--\ref{sect:IB} will be regular categories.

For subsets $U,V$ of a partial semigroup $S$, we write
\[
UV = \set{uv}{u\in U,\ v\in V,\ \br(u)=\bd(v)}
\]
for the (possibly empty) set of all well-defined products $uv$, where $u\in U$ and $v\in V$.  We extend this in the obvious way to products of more than two sets, and write $uV$ for $\{u\}V$, and so on.

A non-empty subset $I$ of a partial semigroup~$S$ is an \emph{ideal} if $IS$ and $SI$ are both contained in~$I$, in which case we also say that $S$ is an \emph{ideal extension} of $I$.  
For $x\in S$, the set $S^1xS^1$ is an ideal of $S$, called the \emph{principal ideal generated by $x$}.  
An ideal is certainly a partial semigroup in its own right, with the same set of objects as $S$.
Ideals of categories need not be categories themselves due to the absence of some identities at some objects.

Green's $\R$, $\L$, $\J$, $\H$ and $\D$ relations on a partial semigroup $S$ will play an important role in all that follows.  For $x,y\in S$, we have
\[
(x,y)\in\R \iff xS^1=yS^1 \COMMa
(x,y)\in\L \iff S^1x=S^1y \COMMa
(x,y)\in\J \iff S^1xS^1=S^1yS^1 .
\]
Note that the first of these says that $x$ and $y$ are $\R$-related if and only if $x=ya$ and $y=xb$ for some $a,b\in\Sone$.  
We also define $\H=\R\cap\L$ and $\D=\R\vee\L$ to be the meet and join of $\R$ and~$\L$ in the lattice of equivalence relations on $S$; as in \cite[Lemma 2.6]{DE2018}, we have $\D=\R\circ\L=\L\circ\R$.  If $x\in S$, and if $\K$ is any of $\R$, $\L$, $\J$, $\H$ or $\D$, we write $K_x = \set{y\in S}{(x,y)\in\K}$ for the $\K$-class of $x$; note that $K_x$ was denoted $[x]_\K$ in \cite{DE2018}.  We will sometimes write $x\mathrel\K y$ in place of $(x,y)\in\K$.  We will also frequently use the pre-order $\leqJ$ defined for $x,y\in S$ by
\[
x\leqJ y \iff x\in S^1yS^1.
\]
Since $\J = {\leqJ}\cap{\geqJ}$, this induces a partial order $\leq$ on the set $S/\J$ of $\J$-classes of $S$:
\[
J_x \leq J_y \iff x\leqJ y.
\]

It will often be important to distinguish Green's relations, classes, etc., on different partial semigroups.  Thus, we will at times write $\R^S$, $R_x^S$, 
and so on, to emphasise the partial semigroup~$S$.  
Of special significance will be the case in which $S$ is contained in a partial semigroup~$T$.  Here, if $\K$ is any of Green's relations, we have $\K^S\sub\K^T\restr_S$, but the reverse inclusion need not hold in general.  If the reverse inclusion does hold, so that in fact $\K^S=\K^T\restr_S$, we say that Green's $\K$ relation on $S$ is \emph{simply inherited} from $T$.  If $S$ is regular, then $\R$, $\L$ and~$\H$ are simply inherited from any $T$ containing $S$ (cf.~\cite[Lemma~1.8]{Sandwich1}), but $\D$ and $\J$ might not be.  We will see in Lemma~\ref{lem:SD}, however, that all of Green's relations are simply inherited if~$S$ is a regular \emph{ideal} of~$T$.

As in \cite[Section 2.3]{CPbook}, if $x\in S$ is regular, then so too is every element of $D_x$, so we refer to the latter as a \emph{regular $\D$-class}.  If $x=xax$, then $x$ is $\R$-related to $xa$, and $\L$-related to $ax$, with $xa$ and $ax$ both being idempotents.  The $\H$-class of an idempotent is a group; group $\H$-classes of $\D$-related idempotents are all isomorphic to each other.

The next simple lemma will be useful frequently.  For the statement and proof, we say an element $x$ of a partial semigroup is \emph{self-divisible} if $x\in xS\cap Sx$: i.e., if $x=ux=xv$ for some $u,v\in S$.  For example, any regular element of a partial semigroup is self-divisible, and every element of a category is self-divisible.

\begin{lemma}\label{lem:SD}
If $S$ is a partial semigroup in which every element is self-divisible (in particular, if $S$ is regular or is a category), and if~$T$ is an ideal extension of $S$, then the following hold:
\ben
\item \label{SD1} If $\K$ denotes any of Green's relations, then $\K^S=\K^T\restr_S$.
\item \label{SD2} If $\K$ denotes any of Green's relations, then every $\K^S$-class of $S$ is a $\K^T$-class of $T$.
\item \label{SD3} Every ideal of $S$ is an ideal of $T$.
\een
\end{lemma}

\pf
We begin with \ref{SD2}. Let $x\in S$.  Since $K_x^S\sub K_x^T$ for all $\K$, we need only show the reverse inclusion.

Beginning with $\K=\R$, let $y\in R_x^T$.  So $x=ya$ and $y=xb$ for some $a,b\in T^1$.  Since $S$ is an ideal, $y=xb$ gives $y\in S$.  Let $u,v\in S$ be such that $x=xu$ and $y=yv$.  Then $x=y(va)$ and $y=x(ub)$.  Since $va,ub\in S$ (as $S$ is an ideal), this shows that $(x,y)\in\R^S$: i.e., that $y\in R_x^S$.

The proof for $\K=\L$ is dual, and the proof for $\K=\J$ is similar.  For $\H$, we have $H_x^T=R_x^T\cap L_x^T=R_x^S\cap L_x^S=H_x^S$, and for $\D=\R\circ\L$, 
\[
D_x^T = \bigcup_{y\in R_x^T}L_y^T = \bigcup_{y\in R_x^S}L_y^T = \bigcup_{y\in R_x^S}L_y^S = D_x^S.
\]

\pfitem{\ref{SD1}}  This follows from \ref{SD2}.

\pfitem{\ref{SD3}}  Let $I$ be an ideal of $S$.  Since every element of $S$ is self-divisible, we have $I=IS$ and $S=ST$, from which it follows that $IT=(IS)T=I(ST)=IS=I$.  Similarly, $TI=I$.
\epf

\subsection{Congruences}\label{subsect:EMRT}

Throughout this section, $(S,\C,\bd,\br,\cdot)$ will denote a fixed partial semigroup.  
Following \cite[p.~52]{MacLane1998} (see also \cite{Kelly1964}), a \emph{congruence} on $S$ is an equivalence relation $\si$ on $S$ such that for all $x,y,a\in S$
\begin{enumerate}[label=\textup{(C\arabic*)},leftmargin=12mm]
\item \label{C1} $\phantom{[}(x,y)\in\si \implies [\bd(x)=\bd(y)$ and $\br(x)=\br(y)]$,
\item \label{C2} $[(x,y)\in\si$ and $\br(a)=\bd(x)] \implies (ax,ay)\in\si$,
\item \label{C3} $[(x,y)\in\si$ and $\bd(a)=\br(x)] \implies (xa,ya)\in\si$.
\end{enumerate}
Condition \ref{C1} says that $\si$ only identifies elements within fixed hom-sets $S_{A,B}$; \ref{C2} says that~$\si$ is \emph{left-compatible} with the partial product in $S$; \ref{C3} is dual to \ref{C2}, and concerns \emph{right-compatibility}.  For each $A\in\C$, the restriction $\si\restr_{S_A}$ of $\si$ to an endomorphism semigroup $S_A$ is a semigroup congruence on $S_A$.  The set $\Cong(S)$ of all congruences on $S$ is a complete lattice under containment, called the \emph{congruence lattice of $S$}, and is a sublattice of $\Eq(S)$, the lattice of all equivalence relations on $S$.  In particular, the meet of any family of congruences is their intersection, and the join is the transitive closure of their union.

Congruences govern the formation of quotients, are the kernels of morphisms/functors, and form the basis of presentations by generators and relations; again see \cite{MacLane1998,Kelly1964}.  If $\si$ is a congruence on the partial semigroup~$S$, then the set $S/\si$ of all $\si$-classes is itself a partial semigroup under the induced operation and $\bd,\br$ maps.  The kernel of a morphism $\phi:S\to T$ of partial semigroups is the congruence $\ker(\phi)=\bigset{(x,y)\in S\times S}{\bd(x)=\bd(y),\ \br(x)=\br(y),\ x\phi=y\phi}$.  If $\phi$ is injective on objects in the sense that $x{\phi=y\phi\implies \bd(x)=\bd(y)}$ and $\br(x)=\br(y)$, for all $x,y\in S$, then $S/\ker(\phi)\cong\im(\phi)$.  Presentations arise when $S$ is a suitable free category.

The trivial relation
\[
\De_S=\bigset{(x,x)}{x\in S}
\]
is clearly a congruence.  
The universal relation $\nabla_S=S\times S$ is not a congruence in general, as condition \ref{C1} will be violated if $S$ has more than one non-empty hom-set.  Instead, the largest congruence on $S$ is
\[
\nab_S = \bigset{(x,y)\in S\times S}{\bd(x)=\bd(y),\ \br(x)=\br(y)}.
\]
It is called the \emph{universal} or \emph{full} congruence, and its equivalence classes are the hom-sets $S_{A,B}$ ($A,B\in\C$).  
For $x,y\in S$, we will also write $x\comp y$ to indicate that $(x,y)\in\nab_S$: i.e., that $x$ and~$y$ belong to a common hom-set.
For any relation $\si$ on $S$, we define
\[
\wt\si = \si\cap\nab_S = \bigset{(x,y)\in\si}{x\sim y}.
\]
For example, two elements of $S$ are $\wt\R$-related if they are $\R$-related and belong to a common hom-set.

For a set of pairs $\Om\sub \nab_S$, we write $\Om_S^\sharp$ (or $\Om^\sharp$ if the context is clear) for the congruence on~$S$ generated by $\Om$: i.e., the smallest congruence on $S$ containing $\Om$.  If $\Om=\big\{(x,y)\big\}$ contains just a single pair, we write $\cg xy_S$ for $\Om_S^\sharp$; again, we will often abbreviate this to $\cg xy$.

Another simple family of congruences comes from ideals.  If $I$ is an ideal of $S$, then the \emph{Rees congruence} $R_I=R_I^S$ is defined by
\[
R_I = \nab_I\cup\De_S = \bigset{(x,y)\in S\times S}{x=y\text{ or } [x,y\in I \text{ and } x\comp y]}.
\]
The $R_I^S$ notation will be used when there may be potential for confusion as to the semigroup in which $I$ is an ideal, and $R_I$ when no such confusion is likely.   

Further congruences arise from a number of constructions described in \cite[Section 3]{EMRT2018}.  These were originally formulated in the context of (ordinary) semigroups, but the constructions and proofs work without difficulty for partial semigroups.

We say a $\J$-class $J$ of $S$ is \emph{stable} if for all $x\in J$ and $a\in S$,
\[
(x,xa)\in\J \implies (x,xa)\in\R \AND (x,ax)\in\J \implies (x,ax)\in\L.
\]
As in \cite[Lemma 3.10]{EMRT2018}, any stable $\J$-class is in fact a $\D$-class.  

Suppose $D$ is a stable, regular $\J$-class of $S$.  Let $G$ be a group $\H$-class contained in $D$, and let $N\normal G$ be a normal subgroup.  We define the relation
\begin{align*}
\nu_N = S^1(N\times N)S^1\cap(D\times D) 
&= \bigset{(axb,ayb)}{x,y\in N,\ a,b\in S^1,\ axb,ayb\in D}.
\end{align*}
As in \cite[Lemma 3.17]{EMRT2018}, $\nu_N$ is an equivalence on $D$, and is contained in $\H$.
We claim that in fact,
\begin{equation}\label{eq:nu}
\nu_N = D(N\times N)D\cap(D\times D).
\end{equation}
Indeed, with ($\supseteq$) being obvious, suppose $x,y\in N$ and $a,b\in S^1$ are such that $axb,ayb\in D$.  Then, writing $e$ for the identity element of $G$, we have
\[
J_x = J_{axb} = J_{aexb} \leq J_{ae} \leq J_e = J_x,
\]
so that $J_{ae}=J_x$; this gives $ae\in J_x=D$.  A similar argument shows that $eb\in D$.  Since $(axb,ayb)=\big((ae)x(eb),(ae)y(eb)\big)$, this completes the proof of \eqref{eq:nu}.  

The significance of \eqref{eq:nu} is that while $D$ could be a $\J$-class of multiple partial semigroups, the relation $\nu_N$ depends only on $N$ and $D$, and not on the containing partial semigroup.  

If $D$ (the stable, regular $\D$-class of $S$, as above) happens to be a $\J$-class of a partial semigroup $T$ containing $S$, then $D$ is still stable in $T$.  Indeed, suppose $x\in D$ and $a\in T$ are such that $(x,xa)\in\J^T$ (i.e., $xa\in D$).  
Fix any right identity $e\in D$ for $x$ (which exists because~$D$ is regular).  Since
\[
J_x^T = J_{xa}^T = J_{xea}^T \leq J_{ea}^T \leq J_e^T = J_x^T \implies J_{ea}^T=J_x^T=D,
\]
we have $ea\in D\sub S$.  So, since $(x,x(ea))=(x,xa)\in\J^S$ (as $xa\in D$), stability of $D$ in $S$ gives $(x,xa)=(x,x(ea))\in\R^S\sub\R^T$.  A dual argument shows that $(x,ax)\in\J^T\implies(x,ax)\in\L^T$.  

As in \cite[Lemma~3.16]{EMRT2018}, the relations $\nu_N$ do not depend on the choice of the group $\H$-class~$G$ in $D$; i.e., if $G_1$ and $G_2$ are group $\H$-classes in $D$, and if $N_1\normal G_1$, then there exists $N_2\normal G_2$ such that~$\nu_{N_1}=\nu_{N_2}$.  

As in \cite[Definition 3.18]{EMRT2018}, an \emph{IN-pair} in $S$ is a pair $(I,N)$, where $I$ is an ideal of $S$, and $N$ is a normal subgroup of a group $\H$-class contained in some stable, regular $\D$-class that is minimal in the set $(S\setminus I)/\J$.  To such an IN-pair, we may associate the relation
\[
R_{I,N} = R_I \cup \nu_N.
\]
As in \cite[Propostion 3.23]{EMRT2018}, $R_{I,N}$ is a congruence on $S$.  As with Rees congruences, we will write~$R_{I,N}^S$ for the congruence $R_{I,N}$ when it is important to emphasise the semigroup $S$.  

If the intersection of all the ideals of $S$ is non-empty, then this intersection is the unique minimal (indeed, \emph{minimum}) ideal of $S$.  Not every partial semigroup (or even ordinary semigroup) has a minimal ideal, but if it exists, the minimal ideal is always a single $\J$-class.  The following constructions from \cite{EMRT2018} apply only to partial semigroups that have a minimal ideal, so for the rest of Section \ref{subsect:EMRT} we assume $S$ has a minimal ideal $M$.  We will also assume that $M$ is stable and regular, so in fact,~$M$ is a single $\D$-class.   (Note that stability of $M$ does not automatically imply regularity of~$M$, as is the case for (ordinary) semigroups; cf.~\cite[Lemma 3.11]{EMRT2018}.  For example, if $S$ has a single element $x$ with $\bd(x)\not=\br(x)$, then $M=S$ is not regular, but it is vacuously stable.)

A congruence $\tau$ on $M$ (considered as a partial semigroup in its own right) is called \emph{liftable} (to~$S$) if $\tau\cup\De_S$ is a congruence on $S$.   As in \cite[Lemma 3.12]{EMRT2018}, the following relations are all liftable congruences on $M$:
\begin{align*}
\wt\L^M &= \bigset{(x,y)\in M\times M}{x\sim y,\ (x,y)\in\L}, \\
\wt\R^M &= \bigset{(x,y)\in M\times M}{x\sim y,\ (x,y)\in\R}, \\
\wt\H^M &= \bigset{(x,y)\in M\times M}{x\sim y,\ (x,y)\in\H}.
\end{align*}
(Since $M$ is a regular ideal, Lemma \ref{lem:SD}\ref{SD1} gives $\K^M=\K\restr_M$ for any $\K$.)
Clearly $\De_M$ and~${\nab_M=\wt\D^M=\wt\J^M}$ are also liftable, with
$\De_M\cup \De_S=\De_S$ and $\nab_M\cup\De_S=R_M$.  If $S$ is contained in a partial semigroup $T$ also with $M$ as its minimal ideal, then a congruence on~$M$ liftable to $T$ must be liftable to $S$, but the converse need not hold in general.  Liftability will be treated in greater generality in Section \ref{subsect:lift}, and the preceding facts can be viewed as special cases of Lemma \ref{lem:E2}.

Let $I$ be an ideal of $S$.  A map $f:I\to M$ is a \emph{retraction} if
\bit
\item $x\sim xf$ for all $x\in I$,
\item $(xy)f=(xf)(yf)$ for all $x,y\in I$ with $\br(x)=\bd(y)$, and
\item $xf=x$ for all $x\in M$.
\eit
The ideal $I$ is \emph{retractable} if such a retraction exists, in which case there is a unique such retraction (since $M$ is regular); cf.~\cite[Corollary 3.4]{EMRT2018}.  Moreover, if $x\in I$ and $a,b\in S^1$ are such that the product $axb$ is defined, then as in~\cite[Lemma 3.3]{EMRT2018}, we have
\begin{equation}\label{eq:axb}
(axb)f = a(xf)b.
\end{equation}  
It follows quickly from \eqref{eq:axb} that if $\K$ is any of Green's relations, $(x,y)\in\K\implies(xf,yf)\in\K$ for any $x,y\in I$.

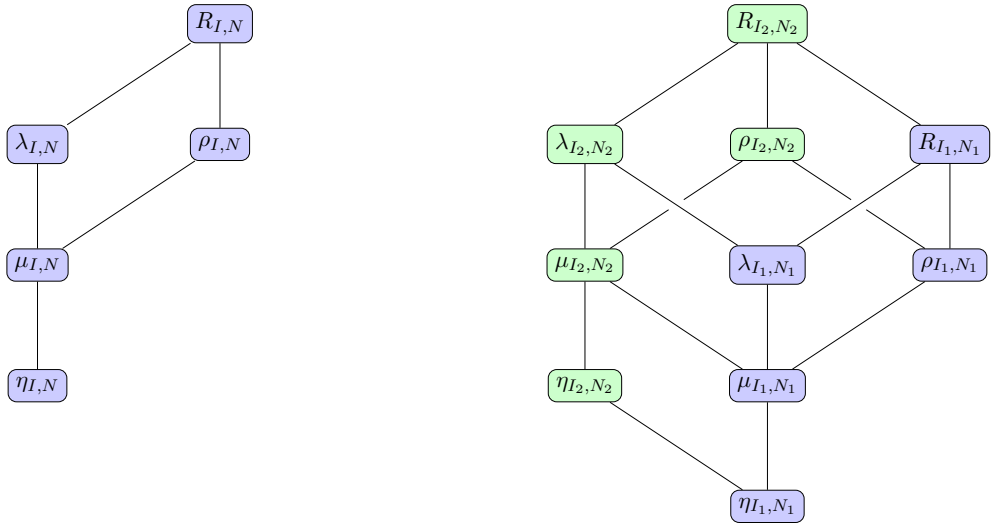
\begin{figure}[ht]
\begin{center}
\scalebox{0.8}{
\begin{tikzpicture}[scale=1]
\node[rounded corners,rectangle,draw,fill=blue!20] (rl) at (6,4) {$R_{I,N}$};
\node[rounded corners,rectangle,draw,fill=blue!20] (r) at (3,2) {$\lambda_{I,N}$};
\node[rounded corners,rectangle,draw,fill=blue!20] (l) at (6,2) {$\rho_{I,N}$};
\node[rounded corners,rectangle,draw,fill=blue!20] (d) at (3,0) {$\mu_{I,N}$};
\node[rounded corners,rectangle,draw,fill=blue!20] (dd) at (3,-2) {$\eta_{I,N}$};
\draw
(dd)--(d)--(l)
(d)--(r)
(l)--(rl)
(r)--(rl)
;
\begin{scope}[shift={(12,-2)}]
\node[rounded corners,rectangle,draw,fill=green!20] (mrl) at (3,6) {$R_{I_2,N_2}$};
\node[rounded corners,rectangle,draw,fill=green!20] (mr) at (0,4) {$\lambda_{I_2,N_2}$};
\node[rounded corners,rectangle,draw,fill=green!20] (ml) at (3,4) {$\rho_{I_2,N_2}$};
\node[rounded corners,rectangle,draw,fill=blue!20] (rl) at (6,4) {$R_{I_1,N_1}$};
\node[rounded corners,rectangle,draw,fill=green!20] (m) at (0,2) {$\mu_{I_2,N_2}$};
\node[rounded corners,rectangle,draw,fill=green!20] (mm) at (0,0) {$\eta_{I_2,N_2}$};
\node[rounded corners,rectangle,draw,fill=blue!20] (r) at (3,2) {$\lambda_{I_1,N_1}$};
\node[rounded corners,rectangle,draw,fill=blue!20] (l) at (6,2) {$\rho_{I_1,N_1}$};
\node[rounded corners,rectangle,draw,fill=blue!20] (d) at (3,0) {$\mu_{I_1,N_1}$};
\node[rounded corners,rectangle,draw,fill=blue!20] (dd) at (3,-2) {$\eta_{I_1,N_1}$};
\draw
(dd)--(d)--(l) (m)--(ml) 
(d)--(m) (l)--(ml) (rl)--(mrl)
(d)--(r) (dd)--(mm)--(m)--(mr)
(l)--(rl) (ml)--(mrl)
;
\fill[white] (1.5,5)circle(.15);
\fill[white] (1.5,3)circle(.15);
\fill[white] (4.5,3)circle(.15);
\draw
(r)--(rl) (mr)--(mrl)  (r)--(mr)
;
\end{scope}
\end{tikzpicture}
}
\caption{Left: the sublattice of $\Cong(S)$ formed by the five congruences $R_{I,N}$,~$\lam_{I,N}$,~$\rho_{I,N}$,~$\mu_{I,N}$ and~$\eta_{I,N}$, associated to a retractable ideal IN-pair $(I,N)$,
or, ignoring the subscript $N$, the congruences
$R_I$, $\lam_I$, $\rho_I$, $\mu_I$ and $\eta_I$, associated to a retractable ideal $I$.  
Right: the sublattice formed by the ten such congruences associated to retractable IN-pairs $(I_1,N_1)$ and $(I_2,N_2)$, in the cases that $I_1\cup N_1\sub I_2$, or $I_1=I_2$ and $N_1\sub N_2$, or, again ignoring subscripts $N_1,N_2$,
the congruences associated to retractable ideals $I_1,I_2$ with $I_1\sub I_2$.}
\label{fig:diamonds2}
\end{center}
\end{figure}

If $I$ is retractable, and if $\tau$ is a liftable congruence on~$M$, then as in \cite[Proposition 3.8]{EMRT2018}, and writing $f:I\to M$ for the unique retraction, the relation
\begin{align*}
\th_{I,\tau} &=\bigset{(x,y)\in I\times I}{(xf,yf)\in\tau}\cup\De_S
\intertext{is a congruence on $S$, which we will again denote by $\th_{I,\tau}^S$ when necessary.  In the case that $I=M$ (which is obviously retractable), we have $\th_{M,\tau}=\tau\cup\De_S$.  With respect to the liftable congruences $\wt\L^M$, $\wt\R^M$, $\wt\H^M$ and~$\De_M$ on $M$ mentioned above, there is a family of congruences associated to a retractable ideal $I$, defined by}
\lam_I = \th_{I,\wt\L^M} &=  \bigset{(x,y)\in I\times I}{x\sim y,\ (xf,yf)\in\L}\cup\De_S, \\
\rho_I = \th_{I,\wt\R^M} &=  \bigset{(x,y)\in I\times I}{x\sim y,\ (xf,yf)\in\R}\cup\De_S, \\
\mu_I = \th_{I,\wt\H^M} &=  \bigset{(x,y)\in I\times I}{x\sim y,\ (xf,yf)\in\H}\cup\De_S, \\
\eta_I = \th_{I,\De_M} &=  \bigset{(x,y)\in I\times I}{xf=yf}\cup\De_S.
\end{align*}
Note that $\th_{I,\nab_M}$ is simply the Rees congruence $R_I$.  If $M$ is $\H$-trivial, then $\mu_I=\eta_I$; this will be the case in all the examples we consider in Sections \ref{sect:T}--\ref{sect:IB}.  As usual, we will denote the above congruences by $\lam_I^S$, $\rho_I^S$, etc., when it is important to emphasise the partial semigroup~$S$.  Figure~\ref{fig:diamonds2}, ignoring the subscripts $N,N_1,N_2$, shows the relative inclusions of these congruences for a single retractable ideal $I$ (left), and for two retractable ideals $I_1,I_2$ with $I_1\sub I_2$ (right); the joins and meets of these relations may be inferred from the diagrams as well.
  Note that Figure \ref{fig:diamonds2} pictures the generic situation in which all of these congruences are distinct; this is not always the case, however; it depends on whether~$M$ is $\R$- and/or $\L$-trivial.

We say an IN-pair $(I,N)$ is \emph{retractable} if $I$ is retractable, and if $Nx$ and $xN$ both have size at most $1$ for all $x\in M$.  Since $N$ is contained in a single group $\H$-class, we have 
$N\times N\subseteq \nab_S$, so that for any $x$ the condition $|Nx|\leq1$ means that either no product $gx$ ($g\in N$) is defined, or else all these products are defined and equal to each other.
If $(I,N)$ is retractable, and if $\tau$ is a liftable congruence on $M$, then as in \cite[Proposition 3.23]{EMRT2018}, and writing $f:I\to M$ for the unique retraction, the relation
\begin{align*}
\th_{I,N,\tau} = \th_{I,\tau}\cup\nu_N &=  \bigset{(x,y)\in I\times I}{(xf,yf)\in\tau}\cup\nu_N \cup \De_S
\intertext{is a congruence on $S$.  When $N$ is trivial, we have $\th_{I,N,\tau}=\th_{I,\tau}$.  When $I=M$, we have $\th_{M,N,\tau}=\tau\cup\nu_N \cup \De_S$.  With respect to the liftable congruences $\wt\L^M$, $\wt\R^M$, $\wt\H^M$ and~$\De_M$ on~$M$ mentioned above, there is a family of congruences associated to a retractable IN-pair $(I,N)$, defined by}
\lam_{I,N} = \th_{I,N,\wt\L^M} &=  \bigset{(x,y)\in I\times I}{x\sim y,\ (xf,yf)\in\L}\cup\nu_N \cup \De_S , \\
\rho_{I,N} = \th_{I,N,\wt\R^M} &=  \bigset{(x,y)\in I\times I}{x\sim y,\ (xf,yf)\in\R}\cup\nu_N \cup \De_S , \\
\mu_{I,N} = \th_{I,N,\wt\H^M} &=  \bigset{(x,y)\in I\times I}{x\sim y,\ (xf,yf)\in\H}\cup\nu_N \cup \De_S , \\
\eta_{I,N} = \th_{I,N,\De_M} &=  \bigset{(x,y)\in I\times I}{xf=yf}\cup\nu_N \cup \De_S .
\end{align*}
Note that $\th_{I,N,\nab_M}$ is simply the relation $R_{I,N}$.  
If $M$ is $\H$-trivial, then $\mu_{I,N}=\eta_{I,N}$.  
Note that $\eta_{M,N}=\nu_N \cup \De_S$ is contained in $\H$.
Again, we will also denote these congruences by $\lam_{I,N}^S$, $\rho_{I,N}^S$, and so on, when required.
The various relationships (in terms of inclusions, meets and joins) between the congruences defined above were described in \cite[Propositions 3.25--3.29]{EMRT2018}, and these are summarised in Figure~\ref{fig:diamonds2} for IN-pairs $(I_1,N_1)$ and $(I_2,N_2)$ in the cases that $I_1\cup N_1\sub I_2$, or $I_1=I_2$ and $N_1\sub N_2$.

\subsection{Useful lemmas}\label{subsect:tricks}

We conclude Section \ref{sect:Pre} with several basic results concerning the kinds of relations and congruences discussed in Section \ref{subsect:EMRT}.  They will be used frequently.  

\begin{lemma}\label{lem:D_trick}
Let $\si$ be a congruence on a partial semigroup $S$.  Suppose $R_I\sub\si$ for some ideal~$I$ of $S$, and that there exists a pair $(x,y)\in\si$ with $y\in I$.  Then, writing $K=S^1xS^1$ for the principal ideal of $S$ generated by $x$, we have $R_K\sub\si$.
\end{lemma}

\pf
Since $R_I\sub\si$, it suffices to show that every element of $K$ is $\si$-related to an element of~$I$.  To do so, let $u\in K$ be arbitrary.  Then $u=axb$ for some $a,b\in S^1$.  We then have $(u,ayb)\in\si$, with $ayb\in I$.
\epf

The next three lemmas concern retractable ideals, IN-pairs, and liftable congruences.

\begin{lemma}\label{lem:Sretract}
Suppose $S$ is a partial semigroup with a stable, regular minimal ideal $M$.  If $S$ is itself retractable, then every congruence on $M$ is liftable to $S$.
\end{lemma}

\pf
Let $\si\in\Cong(M)$, and let $f:S\to M$ be the retraction.  Let $(x,y)\in\si$ and $a\in S$ be arbitrary.  We must show that
\[
\br(a)=\bd(x)\implies(ax,ay)\in\si \AND \bd(a)=\br(x)\implies(xa,ya)\in\si.
\]
By symmetry, we may just do the first, so suppose $\br(a)=\bd(x)$.  Since $x,ax\in M$, we have ${ax=(ax)f=(af)(xf)=(af)x}$; similarly $ay=(af)y$.  But then ${(ax,ay)=\big((af)x,(af)y\big)\in\si}$ since $af\in M$ and $\si$ is a congruence on $M$.
\epf

\begin{lemma}\label{lem:INretract}
Suppose $S$ is a partial semigroup with a stable, regular, $\H$-trivial minimal ideal~$M$.  If an IN-pair $(I,N)$ is such that $I\cup N$ is contained in a retractable ideal of $S$, then $(I,N)$ is retractable.
\end{lemma}

\pf
Let $D$ be the (stable, regular) $\D$-class of $S$ containing $N$.  Then $I\cup D$ is an ideal of~$S$ (by minimality of $D$ in $(S\sm I)/\J$), and is the smallest ideal of $S$ containing $I\cup N$.  Thus, by assumption, $I\cup D$ is contained in a retractable ideal, and is hence itself retractable.  Let $f:I\cup D\to M$ be the retraction.  Let $x\in M$.  We must show that $Nx$ and $xN$ both have size at most $1$.  By symmetry, we just do this for $Nx$, so suppose this set is non-empty.  Let $a,b\in N$.  We must show that $ax=bx$.  Since $(a,b)\in\H$, it follows from \eqref{eq:axb} that $(af,bf)\in\H$, so that $af=bf$ since $M$ is $\H$-trivial.  Using this, and \eqref{eq:axb} again, we have
\[
ax = (ax)f = (af)x= (bf)x = (bx)f = bx. \qedhere
\]
\epf

\begin{lemma}\label{lem:mu_trick}
Suppose $S$ is a partial semigroup with a stable, regular minimal ideal $M$, and that~$I$ is a retractable ideal of $S$.  Suppose $\si$ is a congruence on~$S$, and that there exists $(x,y)\in\si$ with $x\in I$ and $y\in M$.  Suppose also that the principal ideal $K=S^1xS^1$ is regular.  Then $\eta_K\sub\si$.  In particular, if $M$ is $\H$-trivial, then $\mu_K\sub\si$.
\end{lemma}

\pf
Let $u\in K$ with $u=axb$.  Since $\mu_K=\eta_K$ if $M$ is $\H$-trivial, we may prove the lemma by showing that $(u,uf)\in\si$, where $f:I\to M$ is the retraction.  Let $z=ayb$, noting that $z\in M$ (as $y\in M$), and that $(u,z)=(axb,ayb)\in\si$.  Since $u\in K$ is regular, we have $u=uhu$ for some $h\in S$.  Let $v=uhzhu$, again noting that $v\in M$.  Since $(u,z)\in\si$, we have $(u,v)=(uhuhu,uhzhu)\in\si$, so we may complete the proof by showing that $(uf,v)\in\sigma$.  

Since $uf$ and $v$ both belong to $M$, and since $u$ and $hzhu$ both belong to $I$, we have
\[
uf\mathrel\J v = vf = (uhzhu)f = (uf)\cdot(hzhu)f.
\]
Stability then gives $uf\mathrel\R (uf)\cdot(hzhu)f = v$.  
Now let $e\in M$ be a left identity for $uf$ (which exists because $M$ is regular).  Since $(v,uf)\in\R$, as just shown, $e$ is also a left identity for $v$; i.e., $v=ev$.  But then $uf = e(uf) = (ef)(uf) = (eu)f = eu$, and so $(uf,v)=(eu,ev)\in\si$, as $(u,v)\in\si$.  As noted above, this completes the proof.
\epf

The remaining lemmas concern the $\nu_N$ relations.

\begin{lemma}\label{lem:nu_trick}
Suppose $\si$ is a congruence on a partial semigroup $S$, and that $\si\restr_D\sub\H$ for some stable, regular $\J$-class $D$ of $S$.  Let $G$ be a group $\H$-class contained in $D$, and let $e$ be the identity of $G$.  Then $\si\restr_D=\nu_N$ for the normal subgroup $N=\set{g\in G}{(e,g)\in\si}$ of $G$.
\end{lemma}

\pf
It is a routine matter to check that $N\normal G$, so the proof will be complete if we can show that $\si\restr_D=\nu_N$.  In what follows, if $g\in G$, we will write $g^{-1}$ for the inverse of $g$ in $G$.  

Suppose first that $(u,v)\in\nu_N$.  So $u,v\in D$ and $(u,v)=(cpd,cqd)$ for some $c,d\in S^1$ and $p,q\in N$.  
By definition of $N$, we have $(e,qp^{-1})\in\si$, so also $(p,q)=(e\cdot p,qp^{-1}\cdot p)\in\si$.  But then $(u,v)=(cpd,cqd)\in\si$ as well, and so $(u,v)\in\si\restr_{D}$.

Conversely, suppose $(u,v)\in\si\restr_{D}$.  
By assumption, we have $(u,v)\in\H$.  
Let~$H$ be the~$\H$-class containing $u$ and $v$.  By Green's Lemma (see \cite[Lemma 2.7]{DE2018}), there exist $c,d,c',d'\in S$ such that the maps $H\to G:h\mt chd$ and $G\to H:g\mt c'gd'$
are mutually inverse bijections, and we may assume (since $G=H_e$ and $(u,e)\in\D$) that $cud=e$.  From $(u,v)\in\si$ it follows that $(e,cvd)=(cud,cvd)\in\si$, and since $cvd\in G$, it follows from the definition of $N$ that $cvd\in N$.  But then $(u,v) = (c'\cdot e\cdot d',c'\cdot cvd\cdot d')  \in\nu_N$.
\epf

We will generally use the next lemma without explicit reference.

\begin{lemma}\label{lem:nuN}
Let $G$ be a group $\H$-class in some stable, regular $\J$-class of a partial semigroup~$S$, and let $e$ be the identity of $G$.
\ben
\item \label{nuN1} For any $N\normal G$, we have $N=\set{g\in G}{(e,g)\in\nu_N}$.
\item \label{nuN2} For any $N_1,N_2\normal G$, we have $N_1\leq N_2 \iff \nu_{N_1}\sub\nu_{N_2}$.
\een
\end{lemma}

\pf
\firstpfitem{\ref{nuN1}}
Let $D$ be the (stable, regular) $\J$-class of $S$ containing $G$, and let ${I=S^1DS^1\sm D}$ be the ideal consisting of all elements of $S$ strictly $\leqJ$-below $D$.  Then $(I,N)$ is an IN-pair, so we may consider the congruence $R_{I,N}$.  Since $R_{I,N}\restr_D=\nu_N\sub\H$, Lemma~\ref{lem:nu_trick} gives ${N=\set{g\in G}{(e,g)\in R_{I,G}}=\set{g\in G}{(e,g)\in R_{I,G}\restr_D}=\set{g\in G}{(e,g)\in\nu_N}}$.  

In the previous argument, it is possible that $I=\emptyset$, which occurs when $D$ is in fact the minimal ideal of $S$.  In this case, we replace $R_{I,N}=\nab_I\cup\nu_N\cup\De_S$ by $\nu_N\cup\De_S$, which is itself a congruence; cf.~\cite[Remark 3.22]{EMRT2018}.

\pfitem{\ref{nuN2}}
The forward implication follows immediately from the definitions.  For the reverse, using part~\ref{nuN1}, $g\in N_1 \implies (e,g)\in\nu_{N_1}\sub\nu_{N_2} \implies g\in N_2$.
\epf

Recall that $\NN(G)$ denotes the (complete) lattice of normal subgroups of a group $G$.  If we have some collection $N_q$ ($q\in Q$) of normal subgroups of $G$, then their join $\bigvee_{q\in Q}N_q$ in $\NN(G)$ is the least common upper bound of the $N_q$, and is precisely the subgroup of $G$ generated by the union $\bigcup_{q\in Q}N_q$.  When $Q=\{q_1,\ldots,q_k\}$ is finite, this join is just the product $N_{q_1}\cdots N_{q_k}$.  The meet of the $N_q$ (the greatest common lower bound) is simply their intersection.

\begin{lemma}\label{lem:nu_joins}
Let $N_q$ ($q\in Q$) be a collection of normal subgroups of a group $\H$-class $G$ in some stable, regular $\J$-class $D$ of a partial semigroup $S$.  Then
\[
\bigvee_{q\in Q} \nu_{N_q}=\nu_N \AND \bigcap_{q\in Q} \nu_{N_q}=\nu_{N'},
\]
where $N=\bigvee_{q\in Q}N_q$ and $N'=\bigcap_{q\in Q}N_q$.
\end{lemma}

\pf
We just prove the statement concerning joins, as the other is similar.  Again let $I$ be the ideal $S^1DS^1\sm D$, and note that $(I,N_q)$ is an IN-pair for each $q$.  Since each $R_{I,N_q}$ is contained in $R_{I,G}=\nab_I\cup\H\restr_D\cup\De_{S\sm(I\cup D)}$, so too is the join
\[
\si=\bigvee_{q\in Q}R_{I,N_q}=\nab_I\cup\bigvee_{q\in Q}\nu_{N_q}\cup\De_{S\sm(I\cup D)}.
\]
It follows from Lemma \ref{lem:nu_trick} that $\bigvee_{q\in Q}\nu_{N_q}=\si\restr_D=\nu_N$ for some $N\normal G$.  We must show that $N=\bigvee_{q\in Q}N_q$: i.e., that $N$ is the least upper bound of the $N_q$ in $\NN(G)$.  The following argument uses Lemma \ref{lem:nuN}\ref{nuN2} several times.

Since $\nu_N=\bigvee_{q\in Q}\nu_{N_q}$, we have $\nu_{N_q}\sub\nu_N$, and hence $N_q\leq N$, for all $q\in Q$; i.e., $N$ is an upper bound for the $N_q$.  
Now suppose $N''\normal G$ is an arbitrary upper bound for the $N_q$; the proof will be complete if we can show that $N\leq N''$.  But for any $q\in Q$ we have $N_q\sub N''$, and so $\nu_{N_q}\sub\nu_{N''}$.  It follows that $\nu_N=\bigvee_{q\in Q}\nu_{N_q}\sub\nu_{N''}$, which gives $N\leq N''$, as required.
\epf

\begin{rem}\label{rem:nu_lattice}
It follows from Lemmas \ref{lem:nuN}\ref{nuN2} and \ref{lem:nu_joins} that the set $\set{\nu_N}{N\normal G}$ is a (complete) sublattice of~$\Eq(D)$, and is isomorphic to $\NN(G)$.  Moreover, if $I$ is any ideal of $S$ for which $D$ is minimal in $(S\sm I,\leqJ)$, then the sublattice $\set{R_{I,N}}{N\normal G}$ of $\Cong(S)$ is also isomorphic to~$\NN(G)$.  
\end{rem}

\section{Ideal extensions}\label{sect:IE}

This section and the next are central to the paper.
In them we develop a general iterative mechanism that will enable us to compute the congruence lattices in a swift and uniform way for a number of concrete classes of partial semigroups and all their ideals.
The key to our approach is viewing our partial semigroups as constructed by a succession of ideal extensions.
In the current section we focus on the general situation of a partial semigroup $S$ and an ideal extension~$T$ of~$S$
by a single $\J$-class, and we ask how to obtain $\Cong(T)$ from $\Cong(S)$.
In Section \ref{subsect:properties} we abstract a number of technical properties shared by most of our subsequent examples (but see also Section \ref{sect:H} for more).
The key idea with these properties is that, while technical in nature, they ought to be concerning objects that are simpler than congruences; typically, and ideally, they will be phrased just in terms of multiplication.
In Section \ref{subsect:lift} we introduce the concept of a liftable congruence, and prove that the set $\LiftCong ST$ of all congruences of $S$ liftable to $T$ is a complete sublattice of $\Cong(S)$.
The culmination is then Theorem~\ref{thm:E6} in Section~\ref{subsect:IE}, which describes $\Cong(T)$ as an
extension of $\LiftCong ST$.
In Section \ref{sect:chains} we turn to the iterative aspect of the construction, and show what happens when
we repeatedly apply Theorem~\ref{thm:E6} to a partial semigroup that is a chain of ideals.

\subsection{Generation and separation properties}\label{subsect:properties}

We begin with our programme outline above, and identify a handful of abstract technical properties that the majority of our subsequent applications will turn out to possess; further such properties will be discussed in Section \ref{sect:H}.
The first requires the existence of a maximum $\J$-class:

\begin{defn}
\label{defn:Dmax}
We say that a partial semigroup $T$ satisfies the property $\Dmax{T}$ if it has a maximum (greatest) $\J$-class, and if this $\J$-class is regular and stable.  This $\J$-class is then in fact a $\D$-class, and will be denoted by $D_T$.
\end{defn}

Our second property concerns the generation of the full congruence:

\begin{defn}
\label{defn:ngen}
A partial semigroup $T$ satisfies the property $\ngen{T}$ if $T$ satisfies $\Dmax{T}$, and if for any $x\in D_T$ and any $y\in T\sm H_x^T$ with $x\sim y$, we have $\cg xy_T=\nab_T$.
\end{defn}

The property $\ngen{T}$ can be interpreted as saying that every pair that could theoretically generate $\nab_T$ in fact does so.
Indeed, up to symmetry, and writing $S$ for the ideal $T\sm D_T$, the pairs $(x,y)$ not covered by the definition are those where $x,y\in S$, or where $(x,y)\in\H\restr_T$.
Either of those generate a congruence contained in $R_{S,G}^T\subsetneq\nab_T$, where $G$ is some group $\H$-class contained in $D_T$.
More formally, we have:

\begin{lemma}\label{lem:dmax_ngen}
Suppose the partial semigroup $T$ satisfies $\Dmax{T}$, let $S$ be the ideal $T\sm D_T$, and let~$G$ be a group $\H$-class contained in $D_T$.  Then $\ngen T$ holds if and only if every non-universal congruence is contained in $R_{S,G}^T$; i.e., $R_{S,G}^T$ is the top element of $\Cong(T)\sm\{\nab_T\}$.
\epfres
\end{lemma}

 Next we have the first of three \emph{separation properties}; the second will be given shortly in Definition \ref{defn:Sepb}, and the third much later in Definition \ref{defn:Sepz}.  Roughly speaking, we consider a pair ${(x,y)\in S\times S}$ with $x\sim y$ to be \emph{separated} if $x$ and $y$ are not related under one or more of Green's relations.  Each separation property stipulates that certain separated pairs generate congruences containing other pairs separated in specified ways.

\begin{defn}
\label{defn:Sep}
A partial semigroup $T$ satisfies the property $\Sep{T}$ if it satisfies $\Dmax{T}$, the ideal $S=T\setminus D_T$ satisfies $\Dmax{S}$, and the following hold:
\begin{enumerate}[label=\textup{\textsf{(S\arabic*)}},leftmargin=10mm]
\item \label{S1} for any $x\in D_T$ and $y\in S$ with $x\sim y$, there exists $(x',y')\in\cg xy_T$ such that $x'\in D_S$ and $y'\in S\sm H_{x'}^S$, 
\item \label{S2} for any $x\in D_T$ and $y\in D_T\setminus H_x^T$ with $x\sim y$, there exists $(x',y')\in\cg xy_T$ such that $x'\in D_T$ and $y'\in S$,
\item \label{S3} for any $x\in D_T$ and $y\in H_x^T\sm\{x\}$, there exists $(x',y')\in\cg xy_T$ such that $x'\in D_S $ and $y'\in S\sm H_{x'}^S$.
\end{enumerate}
\end{defn}

As it stands, the property $\Sep{T}$ does not fulfil our requirement that it concerns objects of simpler nature than congruences. Fortunately,
in all of the applications studied in Sections~\mbox{\ref{sect:T}--\ref{sect:J}}, we have a property stronger than $\Sep{T}$, where the pair $(x',y')\in\cg xy_T$ can be assumed to have one of the forms $(ax,ay)$, $(ay,ax)$, $(xa,ya)$ or $(ya,xa)$.

\begin{defn}
\label{defn:Mult}
A partial semigroup $T$ satisfies the property $\Mult{T}$ if it satisfies $\Dmax{T}$, the ideal $S=T\setminus D_T$ satisfies $\Dmax{S}$, and the following hold:
\begin{enumerate}[label=\textup{\textsf{(M\arabic*)}},leftmargin=10mm]
\item \label{M1} for any $x\in D_T$ and $y\in S$ with $x\sim y$, there exist $a,b\in T^1$ such that (swapping $x,y$ if necessary) $axb\in D_S$ and $ayb\in S\sm H_{axb}^S$, 
\item \label{M2} for any $x\in D_T$ and $y\in D_T\setminus H_x^T$ with $x\sim y$, there exists $a,b\in T^1$ such that (swapping $x,y$ if necessary) $axb\in D_T$ and $ayb\in S$,
\item \label{M3} for any $x\in D_T$ and $y\in H_x^T\sm\{x\}$, there exists $a,b\in T^1$ such that (swapping $x,y$ if necessary) $axb\in D_S $ and $ayb\in S\sm H_{axb}^S$.
\end{enumerate}
\end{defn}

The categories studied in Sections \ref{sect:L} and \ref{sect:IB} do not have the full $\Sep T$ property.  Nevertheless, they do satisfy the following weaker version.

\begin{defn}
\label{defn:Sepb}
A partial semigroup $T$ satisfies the property $\Sepb{T}$ if it satisfies $\Dmax{T}$, the ideal $S=T\setminus D_T$ satisfies $\Dmax{S}$, and conditions \ref{S1} and \ref{S2} from Definition \ref{defn:Sep} hold.
\end{defn}

We also have the corresponding weaker version of $\Mult T$:

\begin{defn}
\label{defn:Multb}
A partial semigroup $T$ satisfies the property $\Multb{T}$ if it satisfies $\Dmax{T}$, the ideal $S=T\setminus D_T$ satisfies $\Dmax{S}$, and conditions \ref{M1} and \ref{M2} from Definition \ref{defn:Mult} hold.
\end{defn}

In Section \ref{sect:H}, we will introduce intermediate conditions $\Sepz T$ and $\Multz T$; cf.~Definitions~\ref{defn:Sepz} and~\ref{defn:Multz}.  The relative strengths of these separation and multiplication properties are summarised as follows:
\[
\begin{tikzpicture}[xscale=1.2]
\node[above] () at (0,1.4) {$\Mult T$};
\node[above] () at (2,1.4) {$\Multz T$};
\node[above] () at (4,1.4) {$\Multb T$};
\node[above] () at (0,0) {$\Sep T$};
\node[above] () at (2,0) {$\Sepz T$};
\node[above] () at (4,0) {$\Sepb T$};
\foreach \x/\y in {1/0,3/0,1/1.4,3/1.4} {\node[above] () at (\x,\y+.1) {$\Rightarrow$};}
\foreach \x/\y in {0/.65,2/.65,4/.65} {\node[above] () at (\x,\y+.1) {$\Downarrow$};}
\end{tikzpicture}
\]
The next lemma gives a sufficient condition for $\ngen T$ to hold.  The condition is not necessary, however; cf.~Propositions \ref{prop:small01} and \ref{prop:small012}.

\begin{lemma}\label{lem:sepb_ngen}
Let $T$ be a partial semigroup.  If $\Sepb T$ and $\ngen{T\sm D_T}$ both hold, then so too does $\ngen T$.
\end{lemma}

\pf
Write $S=T\sm D_T$.  Let $x\in D_T$ and $y\in T\sm H_x^T$ with $x\sim y$, and let $\si=\cg xy_T$.  We must show that $\si=\nab_T$.

\pfcase1  Suppose first that $y\in S$.  By \ref{S1}, there exists $(x',y')\in\si$ such that $x'\in D_S$ and $y'\in S\sm H_{x'}^S$.  By $\ngen S$, we have $\nab_S=\cg{x'}{y'}_S\sub\cg{x'}{y'}_T\sub\si$, and so $R_S^T=\nab_S\cup\De_{D_T}\sub\si$.  Since $(x,y)\in\si$ with $x\in D_T$ and $y\in S$, Lemma \ref{lem:D_trick} gives $\nab_T=R_T^T\sub\si$, so that $\si=\nab_T$.

\pfcase2  If $y\not\in S$, then in fact $y\in D_T\sm H_x^T$.  By \ref{S2} there exists $(x'',y'')\in\si$ such that $x''\in D_T$ and $y''\in S$.  The pair $(x'',y'')\in\si$ then satisfies the assumptions of Case 1, so by that case we have $\nab_T=\cg{x''}{y''}_T\sub\si$, whence $\si=\nab_T$.
\epf

\subsection{Liftable congruences}\label{subsect:lift}

The next key idea is that of a liftable congruence:

\begin{defn}
\label{defn:lift}
If $S$ is a partial semigroup, and if $T$ is an ideal extension of $S$, we say a congruence~$\si$ on $S$ is \emph{liftable} to $T$ if $\si\cup\De_T$ is a congruence on $T$, or equivalently if $\si=\tau\restr_S$ for some congruence $\tau$ on $T$.  We write $\LiftCong ST$ for the set of all congruences on $S$ liftable to $T$.
\end{defn}

Clearly $\De_S$ and $\nab_S$ are liftable to any ideal extension, and we have $\LiftCong SS=S$.
First we show that liftable congruences form a sublattice of all congruences:

\begin{lemma}\label{lem:lift}
If $S$ is a partial semigroup, and if $T$ is an ideal extension of $S$, then $\LiftCong ST$ is a complete sublattice of $\Cong(S)$.
\end{lemma}

\pf
Let $\si_i$ ($i\in I$) be an arbitrary family of congruences on $S$, each liftable to $T$.  Then each $\si_i\cup\De_T=\si_i\cup\De_{T\sm S}$ is a congruence on $T$, and 
\[
\Big( \bigcap_{i\in I}\si_i \Big) \cup \De_{T\sm S} = \bigcap_{i\in I}(\si_i \cup \De_{T\sm S}) \AND  
\Big( \bigvee_{i\in I}\si_i \Big) \cup \De_{T\sm S} = \bigvee_{i\in I}(\si_i \cup \De_{T\sm S})
\]
are both congruences on $T$.
\epf

We conclude by showing that certain congruences described in Section \ref{subsect:EMRT} are always liftable.

\begin{lemma}\label{lem:E2}
Let $S$ be a partial semigroup with a stable, regular minimal ideal, and suppose $T$ is an ideal extension of $S$ such that every ideal of $S$ is an ideal of $T$.  Then each of the following congruences of~$S$ is liftable to $T$:
\ben
\item \label{c1} $R_I^S$, for some ideal $I$ of $S$,
\item \label{c2} $R_{I,N}^S$, for some IN-pair $(I,N)$ of $S$,
\item \label{c3} $\lam_I^S$, $\rho_I^S$, $\mu_I^S$ or $\eta_I^S$, for some retractable ideal $I$ of $S$,
\item \label{c4} $\lam_{I,N}^S$, $\rho_{I,N}^S$, $\mu_{I,N}^S$ or $\eta_{I,N}^S$, for some retractable IN-pair $(I,N)$ of $S$.
\een
In particular, if every congruence on $S$ has one of the above forms, then $\LiftCong ST=\Cong(S)$.
\end{lemma}

\begin{rem}
By Lemma \ref{lem:SD}\ref{SD3}, the 
assumption that every ideal of $S$ is an ideal of $T$ holds if every element of $S$ is self-divisible, and so in particular if $S$ is regular or a category.
\end{rem}

\begin{proof}[\bf Proof of Lemma \ref{lem:E2}]
\firstpfitem{\ref{c1}}  By assumption, $I$ is an ideal of $T$, and we have $R_I^S\cup\De_T=R_I^T$.

\pfitem{\ref{c2}}  Let $D$ be the (stable, regular) $\J^S$-class of $S$ containing $N$, and let $G$ be the group $\H^S$-class with $N\normal G$.  Now, $I$ and $I\cup D$ are both ideals of $S$ (since $D$ is minimal in $(S\sm I)/\J^S$), and hence these are also ideals of $T$.  It follows that $I$ and $I\cup D$ are both unions of $\J^T$-classes; so too therefore is $D$.  But $D$ is a single $\J^S$-class; thus, since $\J^S\sub\J^T$, it follows that $D$ is a single $\J^T$-class.  As we observed in Section \ref{subsect:EMRT}, $D$ is then a stable, regular $\J^T$-class, so it follows that $(I,N)$ is an IN-pair in $T$ as well, so that $\nu_N^S=\nu_N^T$; cf.~\eqref{eq:nu}.  It follows that $R_{I,N}^S\cup\De_T=R_{I,N}^T$ is a congruence on $T$.

\pfitem{\ref{c3} and \ref{c4}}  These proofs follow the same pattern as \ref{c2}, making use of the following facts: the minimal ideal of $S$ is also the minimal ideal of $T$; any retractable ideal (or IN-pair) of $S$ is also retractable in $T$.  We have $\lam_I^S\cup\De_T=\lam_I^T$ and $\lam_{I,N}^S\cup\De_T=\lam_{I,N}^T$, as appropriate, with similar expressions for the other kinds of congruences.
\end{proof}

\subsection{Congruences on ideal extensions}\label{subsect:IE}

We are now ready to prove the main result of this section, Theorem \ref{thm:E6} below, which concerns an ideal extension $T$ of a partial semigroup $S$ by a single stable, regular $\D$-class $D_T$.
Part \ref{E6i} of the theorem describes all of the congruences contained in $R_{S,G}^T$, where $G$ is a group $\H$-class in~$D_T$.  Parts \ref{E6ii} and \ref{E6iii} use this to describe the entire lattice $\Cong(T)$ under increasingly stronger assumptions.

The statement of Theorem \ref{thm:E6} refers to congruences $\tau_N=(\nu_N)_T^\sharp\restr_S$,
where $N$ is a normal subgroup of $G$. 
Thus, $\tau_N$ is the restriction to $S$ of the congruence on~$T$ generated by $\nu_N$.
From
\[
\nu_N\sub(\nu_N)_T^\sharp\restr_{D_T}\sub (R_{S,N}^T)_T^\sharp\restr_{D_T}= R_{S,N}^T\restr_{D_T}=\nu_N,
\]
we have $(\nu_N)_T^\sharp\restr_{D_T}=\nu_N$.
Since also $(\tau_N)_T^\sharp \sub R_{S,N}^T \sub \nab_S\cup\H\restr_{D_T}$, it follows that
\begin{equation}\label{eq:tau_N}
(\nu_N)_T^\sharp = (\nu_N)_T^\sharp\restr_S\cup(\nu_N)_T^\sharp\restr_{D_T} = \tau_N\cup\nu_N.
\end{equation}

\begin{thm}\label{thm:E6}
Suppose $T$ is a partial semigroup such that $\Dmax{T}$ holds.  Let $S$ be the ideal $T\setminus D_T$, let $G$ be a group $\H^T$-class contained in $D_T$, and for $N\normal G$ define $\tau_N=(\nu_N)_T^\sharp\restr_S$.  
\ben
\item \label{E6i} The interval $[\De_T,R_{S,G}^T]$ in the lattice $\Cong(T)$ is given by
\bit
\item[] \hspace{-1.8mm}\hspace{0.7mm}$\phantom{\Cong(T)}[\De_T,R_{S,G}^T] = \set{\si\cup\nu_N}{\si\in\LiftCong ST,\ N\normal G,\ \tau_N\sub\si}$.
\eit
\item \label{E6ii} If $\ngen{T}$ holds, then
\bit
\item[] \hspace{-1.8mm}$\phantom{[\De_T,R_{S,G}^T]}\Cong(T) = \set{\si\cup\nu_N}{\si\in\LiftCong ST,\ N\normal G,\ \tau_N\sub\si} \cup \{\nab_T\}$.
\eit
\item \label{E6iii} If $\Sep{T}$ and $\ngen{S}$ hold, then
\bit
\item[] \hspace{-1.8mm}$\phantom{[\De_T,R_{S,G}^T]}\Cong(T) = \bigset{\si\cup\De_{D_T}}{\si\in\LiftCong ST} \cup \bigset{R_{S,N}^T}{N\normal G} \cup \{\nab_T\}$.
\eit
\een
\end{thm}

\pf
\firstpfitem{\ref{E6i}}
Beginning with the backwards inclusion, suppose $\si\in\LiftCong ST$ and $N\normal G$ are such that $\tau_N\sub\si$, and put $\th=\si\cup\nu_N$.  Since $\th$ is the union of equivalence relations on disjoint sets, it is an equivalence.  Let $(x,y)\in\th$, and suppose $a\in T$.  To complete the proof that $\th$ is a congruence, we must show that
\[
\br(a)=\bd(x) \implies (ax,ay)\in\th \AND \bd(a)=\br(x) \implies (xa,ya)\in\th.
\]
By symmetry it is enough to do only the first of these, so suppose $\br(a)=\bd(x)$.  If $(x,y)\in\si$, then ${(ax,ay)\in\si\cup\De_T\sub\th}$ since $\si\cup\De_T$ is a congruence on $T$.  If $(x,y)\in\nu_N$, then by \eqref{eq:tau_N} we have $(ax,ay)\in(\nu_N)_T^\sharp=\tau_N\cup\nu_N\sub\si\cup\nu_N=\th$.

Conversely, suppose $\th$ is a congruence on $T$, and that $\th\sub R_{S,G}^T$.  
Let $\si=\th\restr_S$.  Since $\si\cup\De_{D_T}=\th\cap R_S^T$ is a congruence on~$T$ (as both $\th$ and $R_S^T$ are), it follows that $\si\in\LiftCong ST$.    Since $\th\sub R_{S,G}^T=\nab_S\cup\H\restr_{D_T}$, we have $\th=\th\restr_S\cup\th\restr_{D_T}=\si\cup\nu_N$ for some $N\normal G$, using Lemma~\ref{lem:nu_trick} in the last step.  We also have $\tau_N=(\nu_N)_T^\sharp\restr_S\sub\th_T^\sharp\restr_S=\th\restr_S=\si$.

\pfitem{\ref{E6ii}}  
This follows immediately from \ref{E6i} and Lemma \ref{lem:dmax_ngen}.

\pfitem{\ref{E6iii}}  
Since $\Sep T$ implies $\Sepb T$, and since $\Sepb T$ and $\ngen S$ together imply $\ngen T$ (cf.~Lemma \ref{lem:sepb_ngen}), part \ref{E6ii} applies.  To complete the proof, it is enough to show that for $N\normal G$,
\[
\tau_N = \begin{cases}
\De_S &\text{if $N=\{\id_G\}$}\\
\nab_S &\text{if $N\not=\{\id_G\}$.}
\end{cases}
\]
The first is clear, since $\nu_{\{\id_G\}}=\De_{D_T}$.  For the second, suppose $N\not=\{\id_G\}$, and let $x,y\in N$ be distinct, noting that $(x,y)\in\nu_N$.  Then $x,y\in D_T$, $x\sim y$ and $y\in G\sm\{x\}=H_x^T\sm\{x\}$.  By~\ref{S3}, there exists $(x',y')\in\cg xy_T$ such that $x'\in D_S$ and $y'\in S\sm H_{x'}^S$.  By $\ngen S$, we have $\nab_S=\cg{x'}{y'}_S\sub\cg{x'}{y'}_T\sub\cg xy_T\sub(\nu_N)_T^\sharp$, so that $\nab_S\sub(\nu_N)_T^\sharp\restr_S=\tau_N$, whence $\tau_N=\nab_S$.
\epf

\begin{rem}\label{rem:tau_N}
With the notation of Theorem~\ref{thm:E6}, note that for $N\normal G$, the congruence $(\nu_N)_T^\sharp$ is the transitive closure of the relation
\[
\bigset{(axb,ayb)}{(x,y)\in\nu_N,\ a,b\in T,\ \br(a)=\bd(x),\ \bd(b)=\br(x)}.
\]
In fact, since every pair from $\nu_N$ is of the form $(cud,cvd)$ for some $u,v\in N$ and $c,d\in T$, it quickly follows that $(\nu_N)_T^\sharp$ is the transitive closure of the set
\[
T(N\times N)T = \bigset{(axb,ayb)}{x,y\in N,\ a,b\in T,\ \br(a)=\bd(x),\ \bd(b)=\br(x)}.
\]
\end{rem}

\begin{figure}[ht]
\begin{center}
\begin{tikzpicture}[scale=0.9]
\draw[fill=black!25] (0,0) to [bend left=80] (0,4) to [bend left=80] (0,0);
\draw[fill=green!25] (0,4) to [bend left=80] (0,7) to [bend left=80] (0,4);
\draw (0,7) -- (0,8);
\foreach \x in {0,4,7,8} {\draw[-{latex}] (-2,\x)--(-.5,\x);}
\foreach \x in {0,4,8} {\colvert{0}{\x}{black}}
\foreach \x in {7} {\colvert{0}{\x}{white}}
\foreach \x/\y in {0/\De_T,4/R_S^T,7/R_{S,G}^T,8/\nab_T} {\node[left] () at (-2,\x) {$\y$};}
\draw[|-|] (2,0)--(2,4);
\draw[|-] (2,7)--(2,4);
\node[right] () at (2,2) {$\Si_1\cong\LiftCong ST$};
\node[right] () at (2,5.5) {$\Si_2\cong\NN(G)$};
\end{tikzpicture}
\caption{Hasse diagram of $\Cong(T)$, for $T$ as in Theorem \ref{thm:E6}\ref{E6iii}; cf.~Remark \ref{rem:E1}.}
\label{fig:E1}
\end{center}
\end{figure}
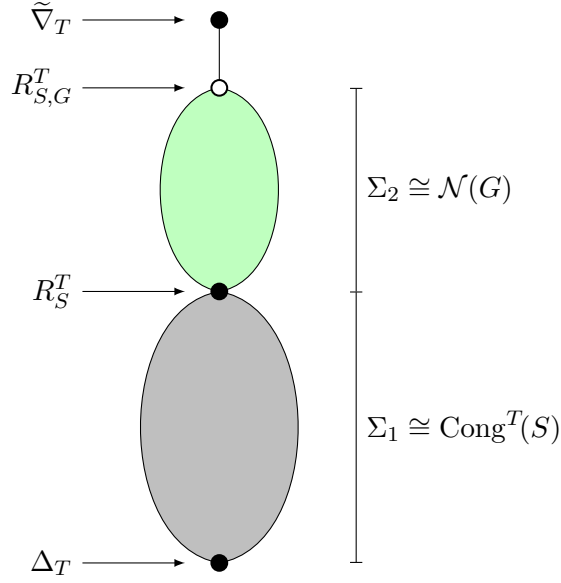

\begin{rem}\label{rem:E1}
Consider $T$ satisfying the assumptions of Theorem \ref{thm:E6}\ref{E6iii}, and let
\[
\Si_1 = \bigset{\si\cup\De_{D_T}}{\si\in\LiftCong ST}
\AND
\Si_2 = \bigset{R_{S,N}^T}{N\normal G}.
\]
Then the following are all easily checked:
\bit
\item $\Si_1$ is a sublattice of $\Cong(T)$ and is isomorphic to $\LiftCong ST$.
\item $\Si_2$ is a sublattice of $\Cong(T)$ and is isomorphic to $\NN(G)$, the lattice of normal subgroups of~$G$ (cf.~Remark \ref{rem:nu_lattice}).  
\item $\Si_1$ has bottom element $\De_T=\De_S\cup\De_{D_T}$ and top element $R_S^T=\nab_S\cup\De_{D_T}$.
\item $\Si_2$ has bottom element $R_S^T=R_{S,\{\id_G\}}^T$ and top element $R_{S,G}^T$.
\item The only remaining congruence on $T$ is $\nab_T$, which sits directly above $R_{S,G}^T$ in the lattice (cf.~Lemma \ref{lem:dmax_ngen}).
\eit
Informally, $\Cong(T)$ is obtained by taking $\LiftCong ST$, stacking $\NN(G)$ on top of it and identifying their respective top and bottom elements, and finally adding
$\nab_T$ on top of that.  This is illustrated in Figure \ref{fig:E1}.  In this figure and coming ones, the trivial congruence and all Rees congruences are depicted as solid vertices, while all other congruences have a white background.  The choice of colours will be explained in Section \ref{subsect:viscon}.
\end{rem}

\begin{figure}[ht]
\begin{center}
\begin{tikzpicture}[line width=0.2mm]
\nc\aaa{9}
\nc\bbb0
\nc\ccc{-.66}
\nc\ddd2
\nc\lele{-3}
\nc\lelele{1.5}
\nc\riri3
\nc\ririri{1.5}
\small
\draw[fill=gray!40, opacity=0.7] (2,0)--(4,2)--(4,4)--(2.5,2.5)--(2,3)--(1,2)--(0,3)--(0,2)--(2,0);
\draw[fill=gray!60, opacity=0.7] (1,2)--(2,3)--(.5,4.5)--(0,4)--(0,3)--(1,2);
\draw[fill=gray!60, opacity=0.7] (2,3)--(3.5,4.5)--(4,4)--(2.5,2.5)--(2,3);
\draw[fill=gray!80, opacity=0.7] (2,3)--(3.5,4.5)--(3,5)--(2,4)--(1,5)--(.5,4.5)--(2,3);
\draw[fill=gray!100, opacity=0.7] (2,4)--(3,5)--(2,6)--(1,5)--(2,4);
\draw (2,0)--(4,2)--(4,4)--(2,6)--(0,4)--(0,2)--(2,0);
\draw (1,2)--(3.5,4.5)--(2,6)--(0,4)--(0,3)--(1,2);
\draw (2.5,2.5)--(4,4)--(2,6)--(.5,4.5)--(2.5,2.5);
\draw (2,4)--(3,5)--(2,6)--(1,5)--(2,4);
\node () at (2.2,-.3) {$\De_S=\tau_{\{\id_G\}}$};
\node () at (1.3,1.7) {$\tau_{N_1}$};
\node () at (2.9,2.3) {$\tau_{N_2}$};
\node () at (1.7,3.8) {$\tau_{G}$};
\node () at (2.3,6.2) {$\nab_S$};
\foreach \x/\y in {1/2,2/4,2.5/2.5} {\colvertx{\x}{\y}{white}{0.08}}
\foreach \x/\y in {2/0,2/6} {\colvertx{\x}{\y}{black}{0.08}}
\begin{scope}[shift={(\aaa,\bbb)}]
\begin{scope}[shift={(\lele+\riri,\lelele+\ririri)}]
\draw[fill=gray!100, opacity=0.7] (2,4)--(3,5)--(2,6)--(1,5)--(2,4);
\draw(2,6)--(2,7);
\draw (2,4)--(3,5)--(2,6)--(1,5)--(2,4);
\foreach \x/\y in {2/4,1/5,3/5,2/6} {\draw(\x-\riri,\y-\ririri)--(\x,\y);}
\fill[gray!70] (2,4.5) circle(.1);
\foreach \x/\y in {2/4,1/5,3/5,2/6} {\draw(\x,\y)--(\x-\lele,\y-\lelele);}
\fill[white] (1,3.5) circle(.1);
\fill[white] (.5,3.75) circle(.1);
\fill[white] (.5,3.25) circle(.1);
\fill[white] (-.5,4.25) circle(.1);
\fill[white] (4-1,3.5) circle(.1);
\fill[white] (4-.5,3.75) circle(.1);
\fill[white] (4-.5,3.25) circle(.1);
\fill[white] (4+.5,4.25) circle(.1);
\end{scope}
\begin{scope}[shift={(\lele,\lelele)}]
\draw[fill=gray!60, opacity=0.7] (1,2)--(3.5,4.5)--(3,5)--(2,4)--(1,5)--(0,4)--(0,3)--(1,2);
\draw[fill=gray!100, opacity=0.7] (2,4)--(3,5)--(2,6)--(1,5)--(2,4);
\draw (2,4)--(3,5)--(2,6)--(1,5)--(2,4);
\draw (1,2)--(3.5,4.5)--(2,6)--(0,4)--(0,3)--(1,2);
\end{scope}
\begin{scope}[shift={(\riri,\ririri)}]
\draw[fill=gray!60, opacity=0.7] (2.5,2.5)--(4,4)--(3,5)--(2,4)--(1,5)--(.5,4.5)--(2.5,2.5);
\draw[fill=gray!100, opacity=0.7] (2,4)--(3,5)--(2,6)--(1,5)--(2,4);
\draw (2,4)--(3,5)--(2,6)--(1,5)--(2,4);
\draw (2.5,2.5)--(4,4)--(2,6)--(.5,4.5)--(2.5,2.5);
\end{scope}
\foreach \x/\y in {1/2,0/3,0/4,2/6,3.5/4.5} {\draw(\x,\y)--(\x+\lele,\y+\lelele);}
\foreach \x/\y in {2.5/2.5,4/4,2/6,.5/4.5} {\draw(\x,\y)--(\x+\riri,\y+\ririri);}
\foreach \x/\y in {2/4,1/5,3/5} {\draw(\x+\riri,\y+\ririri)--(\x,\y)--(\x+\lele,\y+\lelele);}
\draw[fill=gray!40, opacity=0.7] (2,0)--(4,2)--(4,4)--(2.5,2.5)--(2,3)--(1,2)--(0,3)--(0,2)--(2,0);
\draw[fill=gray!60, opacity=0.7] (1,2)--(2,3)--(.5,4.5)--(0,4)--(0,3)--(1,2);
\draw[fill=gray!60, opacity=0.7] (2,3)--(3.5,4.5)--(4,4)--(2.5,2.5)--(2,3);
\draw[fill=gray!80, opacity=0.7] (2,3)--(3.5,4.5)--(3,5)--(2,4)--(1,5)--(.5,4.5)--(2,3);
\draw[fill=gray!100, opacity=0.7] (2,4)--(3,5)--(2,6)--(1,5)--(2,4);
\draw (1,2)--(3.5,4.5)--(2,6)--(0,4)--(0,3)--(1,2);
\draw (2.5,2.5)--(4,4)--(2,6)--(.5,4.5)--(2.5,2.5);
\draw (2,4)--(3,5)--(2,6)--(1,5)--(2,4);
\draw (2,0)--(4,2)--(4,4)--(2,6)--(0,4)--(0,2)--(2,0);
\foreach \x/\y in {1/2,2/4,2.5/2.5} {\colvertx{\x}{\y}{white}{0.08}}
\foreach \x/\y in {2/0,2/6} {\colvertx{\x}{\y}{black}{0.08}}
\begin{scope}[shift={(\lele,\lelele)}]
\foreach \x/\y in {1/2,2/4,2/6} {\colvertx{\x}{\y}{white}{0.08}}
\end{scope}
\begin{scope}[shift={(\riri,\ririri)}]
\foreach \x/\y in {2/4,2.5/2.5,2/6} {\colvertx{\x}{\y}{white}{0.08}}
\end{scope}
\begin{scope}[shift={(\lele+\riri,\lelele+\ririri)}]
\foreach \x/\y in {2/4,2/6} {\colvertx{\x}{\y}{white}{0.08}}
\foreach \x/\y in {2/7} {\colvertx{\x}{\y}{black}{0.08}}
\end{scope}
\end{scope}
\end{tikzpicture}
\caption{
Left: the relative positions of $\tau_{\{\id_G\}}$, $\tau_{N_1}$, $\tau_{N_2}$ and $\tau_G$ in the lattice $\LiftCong ST$.  
Right: the lattice $\Cong(T)$.  See Remark \ref{rem:E3}; cf.~Theorem \ref{thm:E6}\ref{E6ii}.}
\label{fig:E32}
\end{center}
\end{figure}
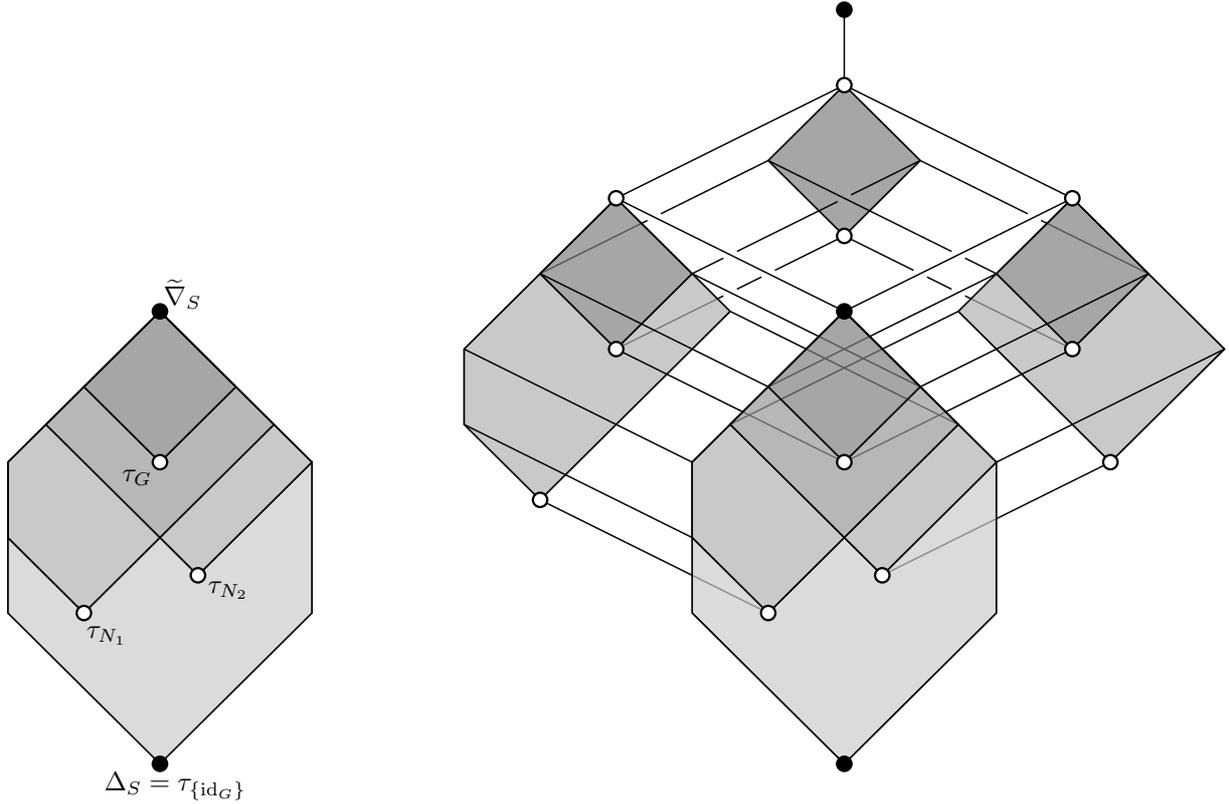

\begin{rem}\label{rem:E3}
When $T$ satisfies the weaker assumptions of Theorem \ref{thm:E6}\ref{E6ii}, as opposed to those of part \ref{E6iii}, the structure of the lattice $\Cong(T)$ tends to be more complex.  Nevertheless, much can still be said.

First, the order in the lattice $\Cong(T)$ is easily described in terms of the orders in $\LiftCong ST$ and $\NN(G)$, as we clearly have
\[
\si_1\cup\nu_{N_1} \sub \si_2\cup\nu_{N_2} \iff \si_1\sub\si_2 \text{ and } N_1\leq N_2.
\]
In other words, $\Cong(T)\sm\{\nab_T\}$ is a sublattice of the direct product $\LiftCong ST\times\NN(G)$.
This sublattice is proper in general, as not all pairs $(\sigma, N)$ give rise to a congruence $\si\cup\nu_N$ on $T$.
It is therefore possible to view $\Cong(T)\sm\{\nab_T\}$ as consisting of pairwise-disjoint \emph{layers}
\[
\Lam(N) = \set{\si\cup\nu_N}{\si\in\LiftCong ST,\ \tau_N\sub\si},
\]
indexed by $N\normal G$.
Clearly, $\Lam(N)$ is a sublattice of $\Cong(T)$, and is isomorphic to the interval $[\tau_N,\nab_S]$ in the lattice $\LiftCong ST$; the bottom and top elements of $\Lam(N)$ are $\nu_N^\sharp=\tau_N\cup\nu_N$ (cf.~\eqref{eq:tau_N}) and $\nab_S\cup\nu_N=R_{S,N}^T$, respectively.  In particular, since $\tau_{\{\id_G\}}=\De_S$, we see that $\Lam\big(\{\id_G\}\big)$ is isomorphic to $\LiftCong ST$ itself. 

Given $N_1,N_2\in\NN(G)$, there is a comparison $\th_1\sub\th_2$ for some $\th_1\in\Lam(N_1)$ and $\th_2\in\Lam(N_2)$ if and only if $N_1\leq N_2$ (cf.~Lemma \ref{lem:nuN}\ref{nuN2}), in which case we must have $\th_1=\si_1\cup\nu_{N_1}$ and $\th_2=\si_2\cup\nu_{N_2}$ with $\si_1\sub\si_2$.  
Thus, if $N_1\lneq N_2$, then $\Lam(N_1)\cup\Lam(N_2)$ is a sublattice of $\Cong(T)$, and for each 
$\si\in[\tau_{N_2},\nab_S]$, 
we have $\si\cup\nu_{N_1}\subsetneq\si\cup\nu_{N_2}$.  
Figure~\ref{fig:E31} shows the Hasse diagram of the sublattice $\Lam(N_1)\cup\Lam(N_2)$, in the case that $N_1\lneq N_2$ and $\tau_{N_1}\subsetneq\tau_{N_2}$.  

Figure \ref{fig:E32} shows the entire lattice $\Cong(T)$, for $T$ as in Theorem \ref{thm:E6}\ref{E6ii}, in the case that $\NN(G)=\big\{\{\id_G\},N_1,N_2,G\big\}$ with $N_1$ and $N_2$ incomparable, with the associated congruences $\tau_{N_1}$ and $\tau_{N_2}$ incomparable as well, and with $\tau_G\not=\nab_S$.  
\end{rem}

\begin{figure}[ht]
\begin{center}
\begin{tikzpicture}[line width=0.2mm]
\nc\aaa7
\nc\bbb0
\nc\ccc{-.66}
\nc\ddd2
\small
\draw[fill=gray!40, opacity=0.7] (2,0)--(4,2)--(4,4)--(3.5,4.5)--(1,2)--(0,3)--(0,2)--(2,0);
\draw[fill=gray!70, opacity=0.7] (1,2)--(3.5,4.5)--(3,5)--(2,4)--(1,5)--(0,4)--(0,3)--(1,2);
\draw[fill=gray!100, opacity=0.7] (2,4)--(3,5)--(2,6)--(1,5)--(2,4);
\node () at (2.2,-.3) {$\De_S$};
\node () at (1.3,1.7) {$\tau_{N_1}$};
\node () at (1.7,3.8) {$\tau_{N_2}$};
\node () at (2.3,6.2) {$\nab_S$};
\draw (2,0)--(4,2)--(4,4)--(2,6)--(0,4)--(0,2)--(2,0);
\draw (1,2)--(3.5,4.5)--(2,6)--(0,4)--(0,3)--(1,2);
\draw (2,4)--(3,5)--(2,6)--(1,5)--(2,4);
\foreach \x/\y in {1/2,2/4} {\colvertx{\x}{\y}{white}{0.08}}
\foreach \x/\y in {2/0,2/6} {\colvertx{\x}{\y}{black}{0.08}}
\begin{scope}[shift={(\aaa,\bbb)}]
\draw[fill=gray!70, opacity=0.7] (1,2)--(3.5,4.5)--(3,5)--(2,4)--(1,5)--(0,4)--(0,3)--(1,2);
\draw[fill=gray!100, opacity=0.7] (2,4)--(3,5)--(2,6)--(1,5)--(2,4);
\foreach \x/\y in {2/4,3/5,1/5,2/6} {\draw(\x,\y)--(\x+\ccc,\y+\ddd);}
\draw (2,4)--(3,5)--(2,6)--(1,5)--(2,4);
\draw[fill=gray!100, opacity=0.7] (2+\ccc,4+\ddd)--(3+\ccc,5+\ddd)--(2+\ccc,6+\ddd)--(1+\ccc,5+\ddd)--(2+\ccc,4+\ddd);
\node () at (1.3,1.7) {$\tau_{N_1}\cup\nu_{N_1}$};
\node () at (1.6,3.75) {$\tau_{N_2}\cup\nu_{N_1}$};
\node () at (2.8+\ccc,6.2+\ddd) {$\nab_S\cup\nu_{N_2}=R_{S,N_2}^T$};
\node[rectangle,fill=white,inner sep=0] () at (3.5,6.15) {$\nab_S\cup\nu_{N_1}=R_{S,N_1}^T$};
\node[rectangle,fill=white,inner sep=0] () at (1.2+\ccc,3.9+\ddd) {$\tau_{N_2}\cup\nu_{N_2}$};
\draw (1,2)--(3.5,4.5)--(2,6)--(0,4)--(0,3)--(1,2);
\draw (2,4)--(3,5)--(2,6)--(1,5)--(2,4);
\draw (2+\ccc,4+\ddd)--(3+\ccc,5+\ddd)--(2+\ccc,6+\ddd)--(1+\ccc,5+\ddd)--(2+\ccc,4+\ddd);
\foreach \x/\y in {1/2,2/4,2/6,2+\ccc/6+\ddd,2+\ccc/4+\ddd} {\colvertx{\x}{\y}{white}{0.08}}
\end{scope}
\end{tikzpicture}
\caption{Left: the relative positions of $\tau_{N_1}\subsetneq\tau_{N_2}$ in the lattice $\LiftCong ST$.  Right: the corresponding layers $\Lam(N_1)$ and $\Lam(N_2)$ in the lattice $\Cong(T)$.  See Remark \ref{rem:E3}; cf.~Theorem \ref{thm:E6}\ref{E6ii}.}
\label{fig:E31}
\end{center}
\end{figure}
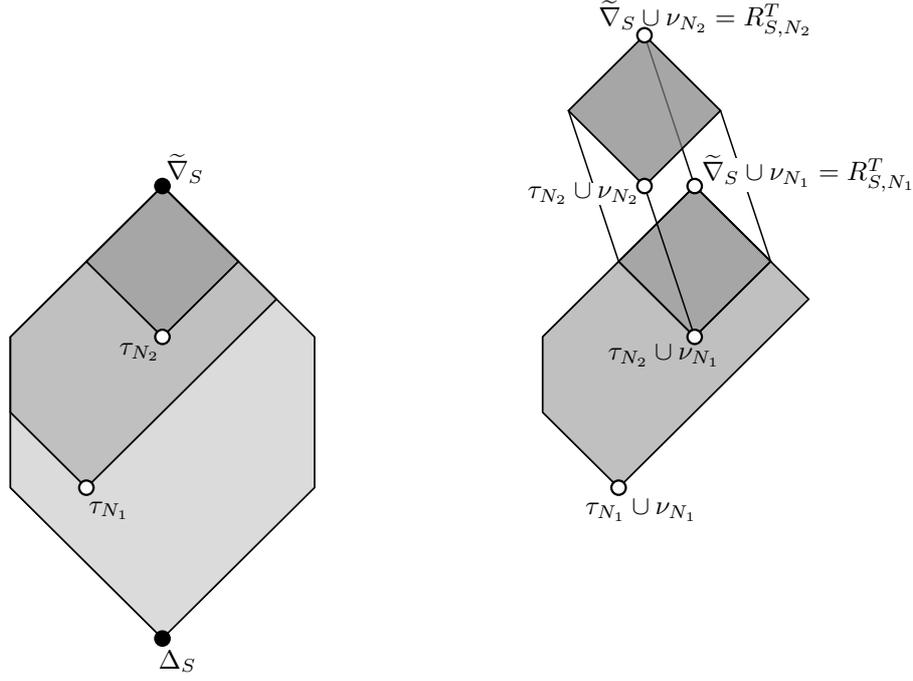

\section{Chains of ideals}\label{sect:chains}

We now use the theory from the previous section to describe the congruences of partial semigroups that are obtained by repeated extensions by single maximum $\J$-classes.

\begin{defn}\label{defn:chain}
A \emph{chain of ideals} of \emph{height} $\xi$, where $\xi$ is an ordinal from the set $\{1,2,3,\ldots,\om\}$, is a stable, regular partial semigroup $U$ whose $\J=\D$-classes form a chain of order type $\xi$; i.e., $U/\J = {\set{D_q}{0\leq q<\xi}}$ and  $D_0<D_1<\cdots$.
\end{defn}

The ideals of any partial semigroup~$S$ are in one-one correspondence with the non-empty downward-closed sets in $(S/\J,\leq)$.  Thus, the ideals of $U$ as in Definition \ref{defn:chain} are precisely the sets
\bit
\item $I_r = D_0\cup D_1\cup\cdots\cup D_r$ for each $0\leq r<\xi$, and
\item $I_\om = U = D_0\cup D_1\cup\cdots$ in the case that $\xi=\om$.
\eit
These ideals form a chain $I_0\subset I_1\subset\cdots$, which justifies the nomenclature.
Since each $D_q$ is regular and stable, each ideal $I_r$ is also regular and stable; hence Lemma \ref{lem:SD} applies, and tells us the following:
\bit
\item Green's relations on each $I_r$ are simply inherited from the corresponding relations on $U$, so we may dispense with superscripts when discussing these relations on $U$ and its ideals.
\item For finite $0\leq q\leq r$, $D_q$ is a (stable, regular) $\D$-class of $I_r$.
\item Every ideal of $I_r$ is an ideal of $U$, and hence the ideals of $I_r$ are precisely the $I_q$ with $0\leq q\leq r$.
\item
It follows that each $I_r$ is a chain of ideals in its own right; its height is $r+1$ when $r<\omega$, and $\omega$ when $r=\xi=\omega$.
\eit
Note that if $\xi=\om$, then the partial semigroup $U$ does not have a maximum (top) $\J$-class, even though every proper ideal of $U$ does.

This section develops the general theory of congruences of chains of ideals. 
When applied in the remaining sections to certain specific classes of chains of ideals, it will enable us to compute,
in one fell swoop, the congruence lattice of $U$ and of all its ideals $I_r$, for all $U$ belonging to the class.
Furthermore, these computations will be uniform across different classes, and the strategy underlying them is as follows.
There will be some (small) value $k$ such that, roughly speaking:
\bit
\item congruences on ideals $I_r$ ($r<k$) are characterised by great ``freedom'' (e.g., in some cases all possible equivalence relations are congruences; cf.~Remark~\ref{rem:T2} and Figure \ref{fig:T}); while
\item congruences on ideals $I_r$ ($r\geq k$) are more ``constrained'', and they all arise from one of the constructions discussed in Section~\ref{subsect:EMRT} (or a variation thereof; cf.~Section~\ref{sect:H});
\item
more specifically, for sufficiently large $r$, the lattice $\Cong(I_r)$ will be obtained 
by stacking on top of $\Cong(I_k)$ the normal subgroup lattices of maximal subgroups of the $\D$-classes
$D_{k+1},D_{k+2},\dots$; cf.~Figure \ref{fig:E2}.
\eit

The most significant contribution of this section is a description of this ``stacking process'' which we give in Section~\ref{subsect:chains}.
Sections \ref{subsect:PRB}--\ref{subsect:small} are devoted to techniques for computing the congruence lattices of ``short'' chains of ideals, which will be used to determined the lattices $\Cong(I_r)$ for $r\leq k$.  Section \ref{subsect:viscon} discusses visualisation conventions that will be used when depicting congruence lattices in specific categories and their ideals in the remainder of the paper.

\subsection{Congruences on chains of ideals}\label{subsect:chains}

The first main result of this section describes the congruences on a chain $U$ of ideals and all its sufficiently large ideals, modulo a description of the congruences on some fixed ideal $I_k$, and subject to certain assumptions as listed in Section~\ref{subsect:properties}.

\begin{thm}\label{thm:E2}
Let $U$ be a chain of ideals of height $\xi$, 
with $D_r$ and $I_r$ as in Definition \ref{defn:chain} and the discussion following it.
For each $0\leq q<\xi$ let $G_q$ be a group $\H$-class contained in $D_q$.  Suppose there exists $0\leq k<\xi$ such that the following hold:
\[
\LiftCong{I_k}U=\Cong(I_k)\COMMA \ngen{I_k}\COMMA \Sep{I_r}\ \text{for all}\ k<r<\xi.
\]
Then for any $r\geq k$ (including $r=\om$ if $\xi=\om$),
\[
\Cong(I_r) = \bigset{\si\cup\De_{I_r}}{\si\in\Cong(I_k)} \cup \bigset{R_{I_q,N}^{I_r}}{k\leq q<r,\ N\normal G_{q+1}} \cup \{\nab_{I_r}\}.
\]
\end{thm}

\pf
First note that all of the listed relations are congruences on~$I_r$.  Indeed, the relations of the form $\si\cup\De_{I_r}=(\si\cup\De_U)\restr_{I_r}$ are congruences because $\LiftCong{I_k}U=\Cong(I_k)$,
and those of the form $R_{I_q,N}^{I_r}$ because each $(I_q,N)$ is an IN-pair in $I_r$.  

To complete the proof, we must show that every congruence on $I_r$ has one of the forms listed in the theorem.

\pfstep{1}
We first prove the statement for every $r\in [k,\xi)$. The proof is by induction, and the anchor $r=k$ is obvious.
Now fix some $k<r<\xi$, and inductively assume that 
\[
\Cong(I_{r-1}) = \bigset{\si\cup\De_{I_{r-1}}}{\si\in\Cong(I_k)} \cup \bigset{R_{I_q,N}^{I_{r-1}}}{k\leq q<r-1,\ N\normal G_{q+1}} \cup \{\nab_{I_{r-1}}\}.
\]
By inspecting the congruences just listed (or by assumption if $r=k+1$), we see that $\ngen{I_{r-1}}$ holds.  
Since also $\Sep{I_r}$ holds by assumption, Theorem \ref{thm:E6}\ref{E6iii} applies to $T=I_r$ and $S=I_{r-1}$, and it tells us that
\[
\Cong(I_r) =  \bigset{\tau\cup\De_{D_r}}{\tau\in\LiftCong{I_{r-1}}{I_r}} \cup \bigset{R_{I_{r-1},N}^{I_r}}{N\normal G_r} \cup \{\nab_{I_r}\}.
\]
We also have $\LiftCong{I_{r-1}}{I_r}=\Cong(I_{r-1})$.  Indeed, every congruence on $I_{r-1}$ listed above lifts to $I_r$: those of the form $\si\cup\De_{I_{r-1}}$ because $\LiftCong{I_k}U=\Cong(I_k)$, and those of the form $R_{I_q,N}^{I_{r-1}}$ because of Lemma \ref{lem:E2}.
To complete the inductive step, it remains to check that $\tau\cup\De_{D_r}$ is one of the congruences listed in the theorem for each $\tau\in\Cong(I_{r-1})$.  But this is clear:
\bit
\item If $\tau=\si\cup\De_{I_{r-1}}$ for some $\si\in\Cong(I_k)$, then $\tau\cup\De_{D_r}=\si\cup\De_{I_r}$.
\item If $\tau=R_{I_q,N}^{I_{r-1}}$ for some $k\leq q<r-1$ and $N\normal G_{q+1}$, then $\tau\cup\De_{D_r}=R_{I_q,N}^{I_r}$.
\item If $\tau=\nab_{I_{r-1}}$, then $\tau\cup\De_{D_r}=R_{I_{r-1}}^{I_r}=R_{I_{r-1},\{\id_{G_r}\}}^{I_r}$.
\eit
This completes the inductive step, and hence Step 1 as well.

\pfstep{2}
We now consider the final case $r=\xi=\om$.  Let $\tau$ be an arbitrary congruence on $U=I_{\om}$.  We must show that $\tau$ has one of the forms listed in the theorem.  

\pfcase{1}  Suppose first that $\tau$ is contained in the Rees congruence $R_{I_t}^{U}$ for some $k\leq t<{\om}$.  Then $\tau=\tau\restr_{I_t}\cup\De_U$.  Since $t<\om$, we know the possible forms of the congruence $\tau\restr_{I_t}$, and it is easy to check (as above) that for each such form, $\tau=\tau\restr_{I_t}\cup\De_U$ has one of the desired forms.

\pfcase{2}  Now suppose $\tau$ is not contained in any of the Rees congruences $R_{I_t}^{U}$ ($k\leq t<\om$).  We claim that $\tau$ must contain each of these Rees congruences.  Before we prove the claim, note that it will then follow that $\tau$ contains $R_{I_k}^{U}\cup R_{I_{k+1}}^{U}\cup\cdots=\nab_U$, so that in fact $\tau=\nab_U$.

To prove the claim, fix some $k\leq t<\om$.  Since $\tau$ is not contained in $R_{I_t}^{U}$, we may fix some $(x,y)\in\tau\sm R_{I_t}^U$.  Without loss of generality, we may assume that $x\in D_u$ and $y\in D_v$ where~$u\geq v$.  So $x,y\in I_u$, and since $(x,y)\in\tau\sm R_{I_t}^U$, we have $u>t$.  Since $u<\om$, we have
\[
\Cong(I_{u}) = \bigset{\si\cup\De_{I_{u}}}{\si\in\Cong(I_k)} \cup \bigset{R_{I_q,N}^{I_{u}}}{k\leq q<u,\ N\normal G_{q+1}} \cup \{\nab_{I_{u}}\}.
\]
By the form of $(x,y)\in\tau\restr_{I_{u}}$, it follows that the congruence~$\tau\restr_{I_{u}}$ is either $\nab_{I_u}$ or else of the form~$R_{I_{u-1},N}^{I_u}$.  In either case, $\tau\restr_{I_{u}}$ contains $R_{I_{u-1}}^{I_{u}}$, and hence~$R_{I_t}^{I_{u}}$.  Thus, $\tau$ contains $R_{I_t}^{I_{u}}$, and hence also $R_{I_t}^{I_{u}}\cup\De_U=R_{I_t}^U$, as required.
\epf

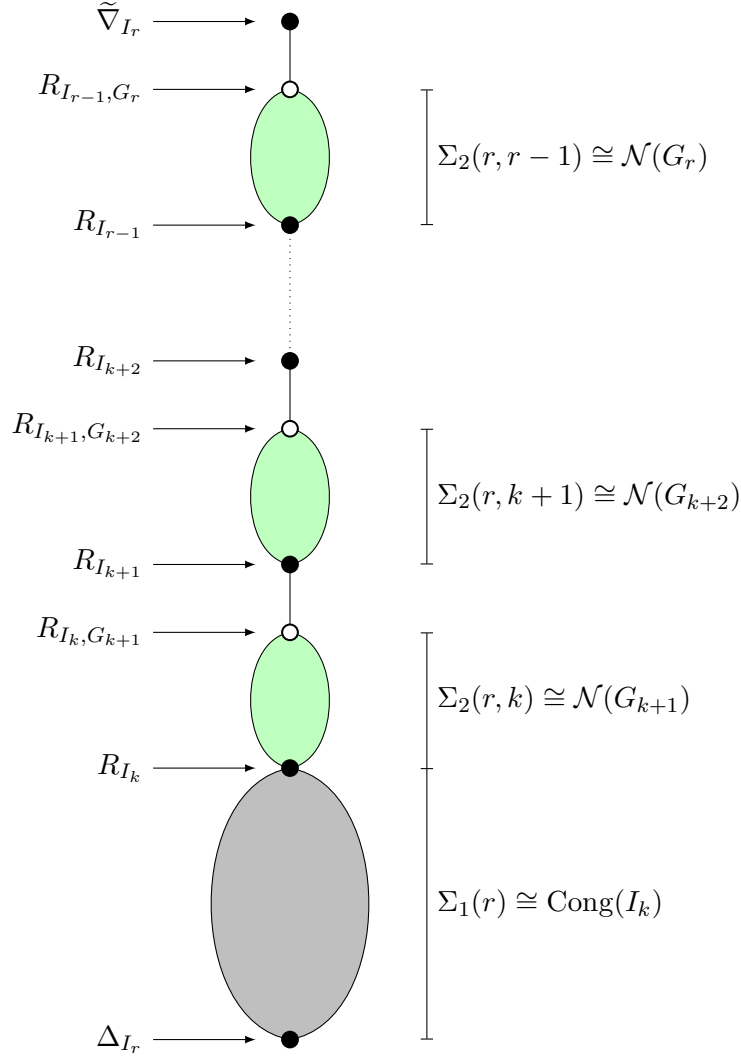
\begin{figure}[ht]
\begin{center}
\begin{tikzpicture}[scale=0.9]
\draw[fill=black!25] (0,0) to [bend left=80] (0,4) to [bend left=80] (0,0);
\draw[fill=green!25] (0,4) to [bend left=80] (0,6) to [bend left=80] (0,4);
\draw[fill=green!25] (0,7) to [bend left=80] (0,9) to [bend left=80] (0,7);
\draw[fill=green!25] (0,12) to [bend left=80] (0,14) to [bend left=80] (0,12);
\draw (0,6) -- (0,7) (0,9)--(0,10) (0,14)--(0,15);
\draw[dotted] (0,10)--(0,12);
\foreach \x in {0,4,6,7,9,10,12,14,15} {\draw[-{latex}] (-2,\x)--(-.5,\x);}
\foreach \x in {0,4,7,10,12,15} {\colvert{0}{\x}{black}}
\foreach \x in {6,9,14} {\colvert{0}{\x}{white}}
\foreach \x/\y in {
0/\De_{I_r},
4/R_{I_k},
6/R_{I_k,G_{k+1}},
7/R_{I_{k+1}},
9/R_{I_{k+1},G_{k+2}},
10/R_{I_{k+2}},
12/R_{I_{r-1}},
14/R_{I_{r-1},G_{r}},
15/\nab_{I_r},
} 
{\node[left] () at (-2,\x) {$\y$};}
\draw[|-|] (2,0)--(2,4);
\draw[-|] (2,4)--(2,6);
\draw[|-|] (2,7)--(2,9);
\draw[|-|] (2,12)--(2,14);
\node[right] () at (2,2) {$\Si_1(r)\cong\Cong(I_k)$};
\node[right] () at (2,5) {$\Si_2(r,k)\cong\NN(G_{k+1})$};
\node[right] () at (2,8) {$\Si_2(r,k+1)\cong\NN(G_{k+2})$};
\node[right] () at (2,13) {$\Si_2(r,r-1)\cong\NN(G_{r})$};
\end{tikzpicture}
\caption{Hasse diagram of $\Cong(I_r)$ for finite $r$, as in Theorem \ref{thm:E2}; see also Remark \ref{rem:E2}.  To avoid clutter, we write $R_{I_q}=R_{I_q}^{I_r}=R_{I_q,\{\id_{G_{q+1}}\}}^{I_r}$ and $R_{I_q,G_{q+1}}=R_{I_q,G_{q+1}}^{I_r}$.}
\label{fig:E2}
\end{center}
\end{figure}

\begin{rem}\label{rem:E2}
In the notation of Theorem \ref{thm:E2}, and for $r\geq k$ and $k\leq q<r$, define the sets
\[
\Si_1(r)=\bigset{\si\cup\De_{I_r}}{\si\in\Cong(I_k)} \AND \Si_2(r,q) = \bigset{R_{I_q,N}^{I_r}}{N\normal G_{q+1}}.
\]
Then, similarly to Remark \ref{rem:E1}, we have the following:
\bit
\item $\displaystyle\Cong(I_r) = \Si_1(r) \cup \bigcup_{k\leq q<r} \Si_2(r,q) \cup \{\nab_{I_r}\}$.
\item $\Si_1(r)$ is a sublattice of $\Cong(I_r)$ and is isomorphic to $\Cong(I_k)$.
\item $\Si_2(r,q)$ is a sublattice of $\Cong(I_r)$ and is isomorphic to $\NN(G_{q+1})$, the lattice of normal subgroups of $G_{q+1}$.
\item $\Si_1(r)$ has bottom element $\De_{I_r}$ and top element $R_{I_k}^{I_r}$.
\item $\Si_2(r,q)$ has bottom element $R_{I_q}^{I_r}$ and top element $R_{I_q,G_{q+1}}^{I_r}$.
\item The only remaining congruence on $I_r$ is $\nab_{I_r}$.  
If $r$ is finite, then $\nab_{I_r}$ sits directly above the top element of $\Si_2(r,r-1)$: i.e., $R_{I_{r-1},G_{r}}^{I_r}$.  
If $r=\om$, then $I_r=U$ and $\Cong(U)\sm\{\nab_U\}$ has no top element (though of course~$\nab_U$ still sits above this entire sublattice).
\eit
Informally, $\Cong(I_r)$ is obtained by successively stacking $\Cong(I_k),\NN(G_{k+1}),\NN(G_{k+2}),\dots$
on top of each other, and adding $\nab_{I_r}$ on top; this
is pictured in Figure \ref{fig:E2}, which gives the Hasse diagram of $\Cong(I_r)$ in the case that $r$ is finite.
Note that if the lattices $\Cong(I_k)$ and $\NN(G_{q+1})$, $k\leq q<r$, are all chains, then so too is $\Cong(I_r)$, including when $r=\om$; cf.~Figure \ref{fig:T}.
\end{rem}

\begin{rem}
Note that Theorem \ref{thm:E2} says nothing about the congruences on the proper ideals of $I_k$ (i.e., $I_0,\ldots,I_{k-1}$).  Nor does it tell us how to compute $\Cong(I_k)$ in the first place.  The results of Sections \ref{subsect:PRB}--\ref{subsect:small} will help us deal with these ideals.
\end{rem}

\subsection{Partial rectangular bands}\label{subsect:PRB}

In all of our applications in Sections \ref{sect:T}--\ref{sect:IB}, the minimal ideal will be a stable, regular, $\H$-trivial, $\D$-universal partial semigroup; these are also the $\H$-trivial chains of ideals of height~$1$.  This section shows how to construct the congruences on such partial semigroups.

Let $\C$ be a set, and suppose that for each $A\in\C$ we have non-empty sets $X_A$ and $Y_A$.  Suppose also that for distinct $A,B\in\C$, $X_A\cap X_B=\emptyset=Y_A\cap Y_B$, and put $X=\bigcup_{A\in\C}X_A$ and $Y=\bigcup_{B\in\C}Y_B$.  (We do not require $X$ and $Y$ to be disjoint.)  Put
\[
S =X\times Y = \bigcup_{A,B\in\C}(X_A\times Y_B).
\]
For $A,B\in\C$ and $x\in X_A$ and $y\in Y_B$, we write $\bd(x,y)=A$ and $\br(x,y)=B$.  
For $x_1,x_2\in X$ and $y_1,y_2\in Y$ such that $\br(x_1,y_1)=\bd(x_2,y_2)$, we define 
$(x_1,y_1)\cdot(x_2,y_2)=(x_1,y_2)$.  Then $(S,\C,\bd,\br,\cdot)$ is easily seen to be a partial semigroup.
We call such a partial semigroup $S$ a \emph{partial rectangular band}.  When $|\C|=1$, we obtain an ordinary rectangular band; cf.~\cite[p.~7]{Howie}. 

The hom-sets of the partial rectangular band~$S=X\times Y$, as above, are ${S_{A,B}=X_A\times Y_B}$ ($A,B\in\C$), and the endomorphism semigroups are $S_A=S_{A,A}=X_A\times Y_A$ ($A\in\C$); the latter are (ordinary) rectangular bands.  It is easy to see that the $\R$- and $\L$-classes in $S$ are, respectively, the sets
\[
\{x\}\times Y\text{ \ for each $x\in X$} \AND X\times\{y\} \text{ \ for each $y\in Y$.}
\]
It quickly follows from this that $S$ is $\H$-trivial and $\D$-universal (hence $\J=\D$), and also stable and regular.  In fact, the converse holds, as the next result shows.  For the statement, we say a partial semigroup $(T,\C,\bd,\br,\cdot)$ is \emph{without surplus objects} if every $A\in\C$ is in the image of at least one of the $\bd$ or $\br$ maps.

\begin{lemma}\label{lem:PRB}
Let $(T,\C,\bd,\br,\cdot)$ be an $\H$-trivial, $\D$-universal, stable, regular partial semigroup
without surplus objects.  For each $A\in\C$, let $X_A=T_A/\R^{T_A}$ and $Y_A=T_A/\L^{T_A}$ be the sets of all $\R^{T_A}$- and $\L^{T_A}$-classes of the endomorphism semigroup $T_A$.  Let $X=\bigcup_{A\in\C}X_A$ and $Y=\bigcup_{A\in\C}Y_A$, and let $S=X\times Y$ be the partial rectangular band as constructed above.  For $u\in T$, let $u\phi = (R_u\cap T_A,L_u\cap T_B)$, where $A=\bd(u)$ and $B=\br(u)$.  Then the mapping $u\mt u\phi$ determines an isomorphism $T\to S$.
\end{lemma}

\pf
We first show that $\phi$ is well defined.  To do so, let $u\in T$.  We must show that $R_u\cap T_A$ is an $\R^{T_A}$-class of $T_A$, and that $L_u\cap T_B$ an $\L^{T_B}$-class of $T_B$.  By symmetry, it suffices to do the former.  Since $u$ is regular, we have $u=uvu$ for some $v\in T$, and we put $w=uv$, noting that $w\in T_A$ and $(u,w)\in\R$; i.e., $R_u=R_w$.  We will show that $R_u\cap T_A=R_w^{T_A}$.  
Beginning with the backwards inclusion, let ~$x\in R_w^{T_A}$.  Then $x\in T_A$ and $(x,w)\in\R^{T_A}\sub\R$, and so $x\in R_w\cap T_A=R_u\cap T_A$.
Conversely, suppose $y\in R_u\cap T_A=R_w\cap T_A$.  Then $y\in T_A$, and we must show that $(y,w)\in\R^{T_A}$.  This is clear if $y=w$.  Otherwise, since $y\in R_w$, we have $y=wa$ and $w=yb$ for some $a,b\in T$.  By comparing domains and ranges, and keeping $y,w\in T_A$ in mind, it is clear that $a,b\in T_A$, which means that $(y,w)\in\R^{T_A}$, as required.

To show that $\phi$ is a homomorphism, suppose $u,v\in T$ are such that the product $uv$ exists.  Then we have $u\in T_{A,B}$ and $v\in T_{B,C}$ for some $A,B,C\in\C$, and $\bd(uv)=A$ and $\br(uv)=C$.  Since $T$ is $\D$-universal, we have $(uv,u),(uv,v)\in\D$, so stability gives $(uv,u)\in\R$ and $(uv,v)\in\L$; i.e., $R_{uv}=R_u$ and $L_{uv}=L_v$.  But then
\begin{align*}
(uv)\phi 
&= (R_{uv}\cap T_A,L_{uv}\cap T_C) \\
&= (R_u\cap T_A,L_v\cap T_C) = (R_u\cap T_A,L_u\cap T_B)\cdot(R_v\cap T_B,L_v\cap T_C) = (u\phi)\cdot(v\phi).
\end{align*}

To show that $\phi$ is surjective, let $(R,L)\in S$.  So $R=R_u\cap T_A$ and $L=L_v\cap T_B$ for some $A,B\in\C$, and for some $u,v\in T$ with $\bd(u)=A$ and $\br(v)=B$.  Let $a,b\in T$ be such that $u=uau$ and $v=vbv$.  Then $(u,ua)\in\R$ and $(v,bv)\in\L$, and we have $ua\in T_A$ and $bv\in T_B$.  Since $T$ is $\D$-universal and regular, we have $ua=sbvt$ for some $s,t\in T$, and we must have $s\in T_{A,B}$ (and $t\in T_{B,A}$).  Then the product $uasbv$ is defined, and we have $(uasbv,u)\in\R$ and $(uasbv,v)\in\L$ (again by $\D$-universality and stability), and also $\bd(uasbv)=A$ and $\br(uasbv)=B$.  But then
\[
(uasbv)\phi = (R_{uasbv}\cap T_A,L_{uasbv}\cap T_B) = (R_u\cap T_A,L_v\cap T_B) = (R,L).
\]

To show that $\phi$ is injective, suppose $u,v\in T$ and $u\phi=v\phi$; i.e., $(R_u\cap T_{\bd(u)},L_u\cap T_{\br(u)}) = (R_v\cap T_{\bd(v)},L_v\cap T_{\br(v)})$.  Note that $R_u\cap T_{\bd(u)}$ is a non-empty subset of $T_{\bd(u)}$ (it contains $ua$ for any $a\in T$ with $u=uau$), and similarly for the other three sets.  Since $R_u\cap T_{\bd(u)}=R_v\cap T_{\bd(v)}\not=\emptyset$, 
it follows that $R_u\cap R_v\not=\emptyset$, and so
$(u,v)\in\R$.  Similarly, 
$(u,v)\in\L$, and so $(u,v)\in\H$, and $u=v$ by $\H$-triviality.
\epf

The next result describes congruences on an (ordinary) rectangular band $X\times Y$;
they are entirely determined by pairs of equivalences on $X$ and $Y$.
 The proof is simple, and is omitted; it also follows from more general results on completely (0-)simple semigroups \cite[Section 10.7]{CPbook2}.

\begin{lemma}\label{lem:RB}
Let $S=X\times Y$ be a rectangular band.  For a pair of equivalences $\ve_1\in\Eq(X)$ and $\ve_2\in\Eq(Y)$, define the relation $\si_{\ve_1,\ve_2}$ on $S$ by
\[
\si_{\ve_1,\ve_2} = \bigset{ \big((x_1,y_1),(x_2,y_2)\big)}{(x_1,x_2)\in\ve_1,\ (y_1,y_2)\in\ve_2}.
\]
Then 
\ben
\item \label{RBi} $\Cong(S) = \bigset{\si_{\ve_1,\ve_2}}{\ve_1\in\Eq(X),\ \ve_2\in\Eq(Y)}$,
\item \label{RBii} $\Cong(S)$ is isomorphic to the lattice direct product $\Eq(X)\times\Eq(Y)$.  \epfres
\een
\end{lemma}

Here is the corresponding result for partial rectangular bands, where instead of pairs of equivalences, we need to consider sequences of pairs.

\begin{prop}\label{prop:PRB}
Let $S=X\times Y$ be a partial rectangular band as above.  For two $\C$-tuples of equivalences
\[
\E_1 = (\ve_1^A)_{A\in\C} \in \prod_{A\in\C}\Eq(X_A) \AND \E_2 = (\ve_2^A)_{A\in\C} \in \prod_{A\in\C}\Eq(Y_A),
\]
define the relation $\si_{\E_1,\E_2}$ by
\[
\si_{\E_1,\E_2} = \bigcup_{A,B\in\C} \bigset{\big((x_1,y_1),(x_2,y_2)\big)\in S_{A,B}\times S_{A,B}}{(x_1,x_2)\in\ve_1^A,\ (y_1,y_2)\in\ve_2^B}.
\]
Then
\ben
\item \label{PRBi} $\Cong(S) = \bigset{\si_{\E_1,\E_2}}{\E_1 \in \prod_{A\in\C}\Eq(X_A) ,\ \E_2 \in \prod_{A\in\C}\Eq(Y_A)}$,
\item \label{PRBii} $\Cong(S)$ is isomorphic to the lattice direct product $\prod_{A\in\C}\Eq(X_A) \times \prod_{A\in\C}\Eq(Y_A)$.
\een
\end{prop}

\pf
\firstpfitem{\ref{PRBi}}  It is routine to verify that such a relation $\si_{\E_1,\E_2}$ is a congruence.  Conversely, suppose $\si\in\Cong(S)$.  For each $A\in\C$, let $\si_A=\si\restr_{S_A}$ be the associated congruence on the rectangular band $S_A$.  By Lemma \ref{lem:RB} we have $\si_A=\si_{\ve_1^A,\ve_2^A}$ for some equivalences $\ve_1^A\in\Eq(X_A)$ and ${\ve_2^A\in\Eq(Y_A)}$.  Put $\E_1 = (\ve_1^A)_{A\in\C}$ and $\E_2 = (\ve_2^A)_{A\in\C}$.  We claim that $\si=\si_{\E_1,\E_2}$.  

Indeed, suppose first that $\big((x_1,y_1),(x_2,y_2)\big)\in\si$, say with $x_1,x_2\in X_A$ and $y_1,y_2\in Y_B$.  Then for any $x\in X_B$ and $y\in Y_A$, we have
\[
\big( (x_1,y),(x_2,y) \big) = \big( (x_1,y_1)(x,y),(x_2,y_2)(x,y) \big)\in\si_A,
\]
and so by definition $(x_1,x_2)\in\ve_1^A$.  Similarly, $(y_1,y_2)\in\ve_2^B$, and so $\big((x_1,y_1),(x_2,y_2)\big)\in\si_{\E_1,\E_2}$.

Conversely, suppose $\big((x_1,y_1),(x_2,y_2)\big)\in\si_{\E_1,\E_2}$, say with $x_1,x_2\in X_A$ and $y_1,y_2\in Y_B$, so that $(x_1,x_2)\in\ve_1^A$ and $(y_1,y_2)\in\ve_2^B$.  Then for any $x\in X_B$ and $y\in Y_A$, we have
\[
\big( (x_1,y),(x_2,y) \big)\in\si_{\ve_1^A,\ve_2^A}=\si_A\sub\si \ANDSIM \big( (x,y_1),(x,y_2)\big)\in\si.
\]
But then $\big((x_1,y_1),(x_2,y_2)\big) = \big((x_1,y)(x_1,y_1)(x,y_1),(x_2,y)(x_1,y_1)(x,y_2)\big)\in\si$.

\pfitem{\ref{PRBii}}  This follows quickly from \ref{PRBi} and the fact that 
\[
\si_{\E_1,\E_2}\sub\si_{\E_1',\E_2'}\quad \Leftrightarrow\quad
\ve_1^A\subseteq\ve_1^B \text{ and } \ve_2^A\subseteq \ve_2^B \text{ for all } A,B\in\C.\qedhere
\]
\epf

\subsection{Short retractable chains}\label{subsect:retract}

In many of the applications in Sections \ref{sect:P}--\ref{sect:J}, the minimal ideal is a partial rectangular band (cf.~Section \ref{subsect:PRB}), and the second smallest ideal is retractable.  It turns out that the congruences on such retractable ideals can be uniformly described, in terms of the congruences on the minimal ideal (cf.~Proposition \ref{prop:PRB}) and the normal subgroups of a group $\H$-class in the other $\D$-class of the ideal.  The next statement and proof refer to the congruences $\th_{I,\tau}$ and $\th_{I,N,\tau}$ discussed in Section \ref{subsect:EMRT}.  Recall that when $I$ is the minimal ideal, $\th_{I,N,\tau}=\tau\cup\nu_N\cup\De_S$.

\begin{prop}\label{prop:small_ideal}
Suppose $S$ is a chain of ideals of height 2, with $\D$-classes $D_0<D_1$, where~$D_0$ is $\H$-trivial.  
Suppose $S$ is retractable (onto $D_0$). Let~$G$ be a group $\H$-class in $D_1$.  Assume also that the following condition is satisfied:
\bit
\item For every $x\in D_1$ and $y\in D_1\sm H_x$ with $x\sim y$, there exists $(x',y')\in\cg xy$ with $x'\in D_1$ and $y'\in D_0$.
\eit
Then the following hold:
\ben
\item \label{small1} Every congruence on $D_0$ is liftable to~$S$, every IN-pair in $S$ is retractable, and we have
\[
\Cong(S) = \bigset{\tau\cup\nu_N}{\tau\in\Cong(D_0),\ N\normal G} \cup \bigset{\th_{S,\tau}}{\tau\in\Cong(D_0)}.
\]
\item \label{small2} The lattice ordering on $\Cong(S)$ is determined by:
\bit
\item $\tau\cup\nu_N\sub\tau'\cup\nu_{N'} \iff \tau\sub \tau'$ and $N\leq N'$,
\item $\th_{S,\tau}\sub\th_{S,\tau'} \iff \tau\sub\tau'$,
\item $\tau\cup\nu_N\sub\th_{S,\tau'} \iff \tau\sub\tau'$.
\eit
\item \label{small3} The lattice $\Cong(S)$ is isomorphic to the lattice direct product
\[
\Cong(D_0) \times \NN(G)^\top,
\]
where $\NN(G)^\top$ is the lattice obtained by adjoining a new top element $\top$ to~$\NN(G)$.
\een
\end{prop}

\pf
\firstpfitem{\ref{small1}}  Lemmas \ref{lem:Sretract} and \ref{lem:INretract} tell us that every congruence on $D_0$ is liftable, and that every IN-pair is retractable.  It follows from the results discussed in Section \ref{subsect:EMRT} that each of the stated relations is a congruence on $S$ (note that $\tau\cup\nu_N=\th_{D_0,N,\tau}$).

Conversely, let $\si$ be an arbitrary congruence on $S$, and put $\tau=\si\restr_{D_0}$.  Let $f:S\to D_0$ be the retraction.  If $\si\sub R_{D_0,G}$, then it follows from Theorem \ref{thm:E6}\ref{E6i} that $\si=\tau\cup\nu_N$ for some $N\normal G$.  This leaves us to consider the case in which $\si\not\sub R_{D_0,G}$.
Since $R_{D_0,G}=\nab_{D_0}\cup\H\restr_{D_1}$, there exists a pair $(x,y)\in\si$ with $x\in D_1$ and $(x,y)\not\in\H$.  By the assumption in the statement, we may assume without loss of generality that $y\in D_0$.
By Lemma~\ref{lem:mu_trick}, we have $\mu_S\sub\si$.
Now let $u,v\in S$ be arbitrary.  Then since $(u,uf),(v,vf)\in\mu_S\sub\si$, we have
\[
(u,v)\in\si \iff (uf,vf)\in\si \iff (uf,vf)\in\si\restr_{D_0}=\tau \iff (u,v)\in\th_{S,\tau},
\]
which shows that $\si=\th_{S,\tau}$ in this case.

\pfitem{\ref{small2}}  The first two items are easily checked.  For the third, first note that since $(\tau\cup\nu_N)\restr_{D_0}=\tau$ and $\th_{S,\tau'}\restr_{D_0}=\tau'$, we certainly have $\tau\cup\nu_N\sub\th_{S,\tau'} \implies \tau\sub\tau'$.  To prove the reverse implication, suppose $\tau\sub\tau'$.  Now, 
\[
\tau \cup \nu_N \sub \tau \cup \H\restr_{D_1} \AND \th_{S,\tau'} = \bigset{(x,y)\in S\times S}{(xf,yf)\in\tau'}.
\]
Thus, since $\tau\sub\tau'\sub\th_{S,\tau'}$, it suffices to show that $\H\restr_{D_1}\sub\th_{S,\tau'}$.  To do so, let $(x,y)\in\H\restr_{D_1}$.  Since $(x,y)\in\H$, \eqref{eq:axb} gives $(xf,yf)\in\H\restr_{D_0}=\De_{D_0}\sub\tau'$, which gives $(x,y)\in\th_{S,\tau'}$, as desired.

\pfitem{\ref{small3}}  It follows quickly from the previous parts that the mapping $\Cong(D_0)\times\NN(G)^\top\to\Cong(S)$ given by $(\tau,N)\mt \tau\cup\nu_N$ and $(\tau,\top)\mt\th_{S,\tau}$
is a lattice isomorphism. 
\epf

\begin{rem}
Figure \ref{fig:B4} gives the Hasse diagram of the congruence lattice of a retractable ideal of the form described in Proposition \ref{prop:small_ideal} in the specific case of the Brauer category $\B$, defined below.  (But note that the $\D$-classes in that example are denoted $D_0<D_2$, and the ideal itself by $I_2=D_0\cup D_2$, for reasons to be explained in Section \ref{sect:B}.)  As explained in Remark \ref{rem:Beven}, here $\Cong(I_0)$ is isomorphic to $\Eq_3\times\Eq_3$, and $\NN(G)$ is a two-element chain.
\end{rem}

\subsection{Short non-retractable chains}\label{subsect:small}

Recall that we are currently working towards a set-up within which we will be able to characterise the congruences for a range of specific categories and their ideals in Sections \ref{sect:T}--\ref{sect:IB}.
The main tool in doing this is Theorem \ref{thm:E2}, which highlights the importance of a certain (small) ideal $I_k$.
In all our applications in Sections \ref{sect:T}--\ref{sect:J}, $I_k$ will be the smallest non-retractable ideal in our category.
Furthermore, its height (as a chain of ideals) will always be $2$ or $3$.
The results of the last two sections will be used to determine $\Cong(I_r)$ for $r<k$.
In this section we turn our attention to $I_k$ itself.
We will consider separately the cases of height $2$ and $3$, beginning with the former, which is also the simpler of the two.
The statement refers to the relations $\wt\R=\set{(x,y)\in\R}{x\sim y}$ and $\wt\L=\set{(x,y)\in\L}{x\sim y}$.

\begin{prop}\label{prop:small01}
Suppose $T$ is a chain of ideals of height 2, with $\D$-classes $D_0<D_1$, where~$D_0$ is $\H$-trivial.  Let~$G$ be a group $\H$-class in $D_1$, and let $e$ be the identity of $G$.  Assume also that the following conditions are satisfied:
\ben
\item \label{01i} For every $x\in D_1$ and $y\in D_1\sm H_x$ with $x\sim y$, there exists $(x',y')\in\cg xy$ with $x'\in D_1$ and $y'\in D_0$.
\item \label{01ii} For every $x,y\in T$ with $x\sim y$ and $(x,y)\not\in\L$, we have $\rho_{D_0}\sub\cg xy$.
\item \label{01iii} For every $x,y\in T$ with $x\sim y$ and $(x,y)\not\in\R$, we have $\lam_{D_0}\sub\cg xy$.
\item \label{01iv} If $D_0$ is not $\wt\R$-trivial, then for every $x\in G\sm\{e\}$, there exists $a\in D_0$ such that $ae\not=ax$ (with both products defined).
\item \label{01v} If $D_0$ is not $\wt\L$-trivial, then for every $x\in G\sm\{e\}$, there exists $a\in D_0$ such that $ea\not=xa$ (with both products defined).
\een
Then
\[
\Cong(T) = \{\De_T,\ \lam_{D_0},\ \rho_{D_0},\ R_{D_0}\} \cup \set{R_{D_0,N}}{N\normal G} \cup \{\nab_T\}.
\]
\end{prop}

\pf
We first claim that for any congruence $\si$ on $T$, 
\begin{equation}
\label{eq:claim}
\si\not\sub\lam_{D_0}\quad \implies \quad \rho_{D_0}\sub\si . 
\end{equation} 
Indeed, this is obvious if $D_0$ is $\wt\R$-trivial, for then $\rho_{D_0}=\De_T$, so suppose instead that $D_0$ is not $\wt\R$-trivial.  Fix some $(x,y)\in\si\sm\lam_{D_0}$.  If $(x,y)\not\in\L$, then \ref{01ii} gives $\rho_{D_0}\sub\cg xy\sub\si$.  This leaves us to consider the case in which $(x,y)\in\L$.  Since $(x,y)\not\in\lam_{D_0}$, we must have $x,y\in D_1$.  If $y\not\in H_x$, then by \ref{01i} there exists $(x',y')\in\si$ such that $x'\in D_1$ and $y'\in D_0$, and \ref{01ii} then gives $\rho_{D_0}\sub\cg{x'}{y'}\sub\si$.  The remaining case to consider is when $y\in H_x\sm\{x\}$.  Here, by Green's Lemma \cite[Lemma 2.7]{DE2018}, there exist $u,v\in T^1$ such that $H_x\to G:h\mt uhv$ is a bijection, and such that $x\mt e$.  Then with $z=uyv$, we have $(e,z)=(uxv,uyv)\in\si$, and $z\in G\sm\{e\}$.  By \ref{01iv}, there exists $a\in D_0$ such that $ae\not=az$; since $(e,z)\in\si$, we also have $(ae,az)\in\si$.  Now,
\begin{equation}\label{eq:exN}
e,z\in G \implies (e,z)\in\H \implies (e,z)\in\R \implies (ae,az)\in\R,
\end{equation}
since $\R$ is a left congruence.  But $ae\not=az$, so since $D_0$ is $\H$-trivial, we must have $(ae,az)\not\in\L$.  But then \ref{01ii} gives $\rho_{D_0}\sub\cg{ae}{az}\sub\si$, and \eqref{eq:claim} is proved.

Next we show that $\ngen T$ holds.  To do so, fix some $x\in D_1$ and $y\in T\sm H_x$ with $x\sim y$, and put $\si=\cg xy$.  We must show that $\si=\nab_T$.  By \ref{01i} we may assume without loss of generality that $y\in D_0$.  We then have $\si\not\sub\lam_{D_0}$ and $\si\not\sub\rho_{D_0}$.  By \eqref{eq:claim} and its dual, it follows that $\si$ contains both $\rho_{D_0}$ and $\lam_{D_0}$, and hence also their join $\lam_{D_0}\vee\rho_{D_0}=R_{D_0}$.  Since also $(x,y)\in\si$ with $x\in D_1$ and $y\in D_0$, it follows from Lemma \ref{lem:D_trick} that $\nab_T=R_T\sub\si$, and so $\si=\nab_T$, as required.

Now that we have established $\ngen{T}$, it then follows from Theorem \ref{thm:E6}\ref{E6ii} that
\begin{equation}\label{eq:claimT}
\Cong(T) = \set{\si\cup\nu_N}{\si\in\LiftCong{D_0}T,\ N\normal G,\ \tau_N\sub\si} \cup \{\nab_T\}.
\end{equation}

Next we claim that
\begin{equation}\label{eq:claimI0}
\LiftCong{D_0}T = \{\De_{D_0},\ \lam_{D_0}^{D_0},\ \rho_{D_0}^{D_0},\ \nab_{D_0}\}.
\end{equation}
Indeed, each of the claimed congruences on $D_0$ is liftable by Lemma \ref{lem:E2}.  Conversely, consider some ${\si\in\LiftCong{D_0}T}$, and let $\tau=\si\cup\De_{D_1}\in\Cong(T)$.  
Suppose first that $\tau\sub\lam_{D_0}$, so that $\si=\tau\restr_{D_0}\sub\lam_{D_0}\restr_{D_0}=\lam_{D_0}^{D_0}$.  
\bit
\item If $\tau\sub\rho_{D_0}$, then $\tau\sub\lam_{D_0}\cap\rho_{D_0}=\De_T$ (as $D_0$ is $\H$-trivial), and so $\tau=\De_T$, and $\si=\De_{D_0}$.  
\item If $\tau\not\sub\rho_{D_0}$, then $\lam_{D_0}\sub\tau$ by the dual of~\eqref{eq:claim}, and so $\lam_{D_0}^{D_0}\sub\si$, and then $\si=\lam_{D_0}^{D_0}$.  
\eit
Similarly, if $\tau\sub\rho_{D_0}$, then $\si$ is either $\De_{D_0}$ or $\rho_{D_0}^{D_0}$.  
If $\tau$ is contained in neither $\lam_{D_0}$ nor $\rho_{D_0}$, then $\tau$ contains both $\rho_{D_0}$ and $\lam_{D_0}$, by \eqref{eq:claim} and its dual, and hence as above, $R_{D_0}\sub\tau$, so that $\nab_{D_0}\sub\si$, meaning that $\si=\nab_{D_0}$.
This completes the proof of \eqref{eq:claimI0}.

In light of \eqref{eq:claimT} and \eqref{eq:claimI0}, we may complete the proof by showing that for $N\normal G$,
\[
\tau_N = \begin{cases}
\De_{D_0} &\text{if $N=\{e\}$}\\
\nab_{D_0} &\text{if $N\not=\{e\}$.}
\end{cases}
\]
Again this is clear for $N=\{e\}$, so suppose otherwise, and let $x\in N\sm\{e\}$ be arbitrary.  Now, $(e,x)\in\nu_N\sub(\nu_N)_T^\sharp$ and $(e,x)$ belongs to neither $\lam_{D_0}$ nor $\rho_{D_0}$, so it follows that the congruence~$(\nu_N)_T^\sharp$ itself is contained in neither $\lam_{D_0}$ nor $\rho_{D_0}$.  As above, it follows that $R_{D_0}\sub(\nu_N)_T^\sharp$, and hence $\nab_{D_0}\sub(\nu_N)_T^\sharp\restr_{D_0}=\tau_N$, so that $\tau_N=\nab_{D_0}$.
\epf

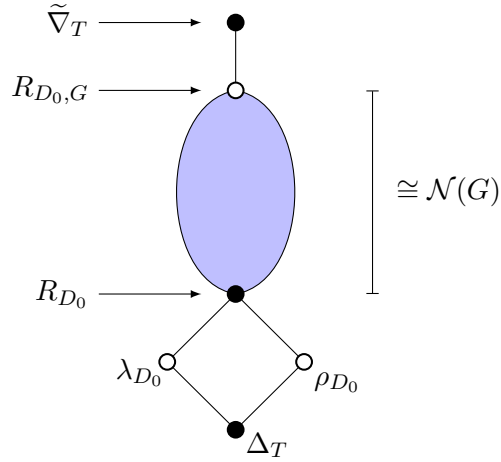
\begin{figure}[ht]
\begin{center}
\begin{tikzpicture}[scale=0.9]
\draw[fill=blue!25] (0,4) to [bend left=80] (0,7) to [bend left=80] (0,4);
\draw (0,7) -- (0,8) (0,2)--(1,3)--(0,4)--(-1,3)--(0,2);
\foreach \x in {4,7,8} {\draw[-{latex}] (-2,\x)--(-.5,\x);}
\foreach \x in {4,8} {\colvert{0}{\x}{black}}
\foreach \x in {7} {\colvert{0}{\x}{white}}
\colvert{-1}{3}{white}
\colvert{1}{3}{white}
\colvert{0}{2}{black}
\foreach \x/\y in {4/R_{D_0},7/R_{D_0,G},8/\nab_T} {\node[left] () at (-2,\x) {$\y$};}
\node[below right] () at (0,2.1) {$\De_T$};
\node[below left] () at (-.9,3.2) {$\lam_{D_0}$};
\node[below right] () at (1,3) {$\rho_{D_0}$};
\draw[|-|] (2,7)--(2,4);
\node[right] () at (2.2,5.5) {$\cong\NN(G)$};
\end{tikzpicture}
\caption{Hasse diagram of $\Cong(T)$, for $T$ as in Proposition \ref{prop:small01}; cf.~Remark \ref{rem:small01}.}
\label{fig:small01}
\end{center}
\end{figure}

\begin{rem}\label{rem:small01}
The lattice $\Cong(T)$ is easily visualised, for $T$ as in Proposition \ref{prop:small01}.  In general, the sublattice $\{\mu_{D_0}, \lam_{D_0}, \rho_{D_0}, R_{D_0}\}$ forms the usual diamond; the sublattice ${\set{R_{D_0,N}}{N\normal G}}$ is isomorphic to $\NN(G)$; and $\nab_T$ sits on top of the latter.  See Figure \ref{fig:small01}; cf.~Figures \ref{fig:Bodd}, \ref{fig:TL} and \ref{fig:J}.
\end{rem}

\begin{rem}\label{rem:small01_2}
Proposition \ref{prop:small01} simplifies if $D_0$ is $\wt\L$- and/or $\wt\R$-trivial.  For example, if $D_0$ is $\wt\L$-trivial, then $\lam_{D_0}=\De_T$, and also $R_{D_0}=\lam_{D_0}\vee\rho_{D_0}=\rho_{D_0}$, in which case the conclusion of the proposition is that
\[
\Cong(T) = \{\De_T\} \cup \set{R_{D_0,N}}{N\normal G} \cup \{\nab_T\}.
\]
Note that assumptions \ref{01iii} and \ref{01v} are vacuously satisfied when $D_0$ is $\wt\L$-trivial.  This simplification will occur in Theorem \ref{thm:T}, in the second smallest ideal of the transformation category~$\T$.
\end{rem}

Here is the result concerning partial semigroups with three $\D$-classes $D_0<D_1<D_2$.  Its technical assumptions are more detailed than those of Proposition \ref{prop:small01} because of the presence of a non-trivial retraction.  Assumption \ref{012ii} refers to a relation denoted $\L f^{-1}$; here $f$ is the retraction $S=D_0\cup D_1\to D_0$, and we write 
\[
\L f^{-1} = \bigset{(x,y)\in S\times S}{(xf,yf)\in\L}.
\]
Similar comments apply to the relation $\R f^{-1}$ in assumption \ref{012iii}.

\begin{prop}\label{prop:small012}
Suppose $T$ is a chain of ideals of height 3, with $\D$-classes ${D_0<D_1<D_2}$, where $D_0$ is $\H$-trivial.
Suppose the ideal $S=D_0\cup D_1$ is retractable, that $f:S\to D_0$ is the retraction, and that
property $\Sep{T}$ holds.
Let~$G_1$ and $G_2$ be group $\H$-classes in $D_1$ and $D_2$, respectively, and write $e$ for the identity of $G_2$.  Assume also that the following conditions are satisfied:
\ben
\item \label{012i} For every $x\in D_1$ and $y\in D_1\sm H_x$ with $x\sim y$, there exists $(x',y')\in\cg xy$ with $x'\in D_1$ and $y'\in D_0$.
\item \label{012ii} For every $x,y\in T$ with $x\sim y$ and $(x,y)\not\in\L\cup\L f^{-1}$, we have $\rho_{D_0}\sub\cg xy$.
\item \label{012iii} For every $x,y\in T$ with $x\sim y$ and $(x,y)\not\in\R\cup\R f^{-1}$, we have $\lam_{D_0}\sub\cg xy$.
\item \label{012iv} There exists a retractable IN-pair $(S,H)$ where $H\normal G_2$, and for every $x\in G_2\sm H$, there exist $a,b\in D_0$ such that $ae\not=ax$ and $eb\not=xb$ (with all four products defined).
\een
Then
\begin{align*}
\Cong(T) = \set{\mu_{D_0,N},\ \lam_{D_0,N},\ \rho_{D_0,N},\ R_{D_0,N}}{{}&N\normal G_1} \\[1truemm]
{} \cup \set{\mu_{S,N},\ \lam_{S,N},\ \rho_{S,N},\ R_{S,N}}{{}&N\normal G_2,\ N\leq H} \\[1truemm]
{} \cup \set{R_{S,N}}{{}&N\normal G_2,\ N\nleq H} 
 \cup \{\nab_T\}.
\end{align*}
\end{prop}

\pf
We begin with two claims that will be used in several of the calculations to follow.
There is an obvious dual to the first, but we will not need to state it.

\pfclaim{1}
For any $(x,y)\in R_S\sm\lam_S$, we have $\rho_{D_0}\sub\cg xy$.

Indeed, here we have $x,y\in S$ and $(xf,yf)\not\in\L$.  It follows from the latter that $(x,y)\not\in\L$ (consider the contrapositive; cf.~\eqref{eq:axb}), and that $(x,y)\not\in\L f^{-1}$.  So~\ref{012ii} gives $\rho_{D_0}\sub\cg xy$.

\pfclaim{2}
If $x,y\in T$ are distinct, and are such that $x\sim y$, $x\in D_1\cup D_2$, and $(x,y)\not\in\H\restr_{D_1}$, then $\mu_S\sub\cg xy$.

Indeed, swapping $x,y$ if necessary, we may assume that $x\geqJ y$.  Consider first the case in which $x\in D_1$.  Since $x\geqJ y$, we have $y\in S$.  By \ref{012i}, it suffices to assume that $y\in D_0$.  Lemma~\ref{lem:mu_trick} then gives $\mu_S\sub\cg xy$.

Next suppose $x\in D_2$.  Since $x\not=y$, it follows from $\Sep{T}$ that there exists $(x',y')\in\cg xy$ such that $x'\in D_1$ and $y'\in S\sm H_{x'}$.  This reduces to the previous case, and completes the proof of the claim.

\afterclaim
Next we show that $\ngen T$ holds.  
To do so, fix some $x\in D_2$ and $y\in T\sm H_x$ with $x\sim y$, and put $\si=\cg xy_T$.  We must show that $\si=\nab_T$.  By~\ref{S2} we may assume in fact that $y\in S$.  By Claim 2, we have $\mu_S\sub\si$.  Thus, $(y,yf)\in\si$, and so $(x,yf)\in\si$ by transitivity.  Since $x\in D_2$ and $y\in S$, \ref{012ii} and \ref{012iii} apply to the pair $(x,y)\in\si$, and tell us that $\si$ contains $\lam_{D_0}$ and $\rho_{D_0}$, and hence their join, $R_{D_0}$.  Since also $(x,yf)\in\si$ with $x\in D_2$ and $yf\in D_0$, Lemma \ref{lem:D_trick} gives $\nab_T=R_T\sub\si$, and so $\si=\nab_T$, as required.

Now that we know $\ngen{T}$ holds, it then follows from Theorem \ref{thm:E6}\ref{E6ii} that
\begin{equation}\label{eq:claimTT}
\Cong(T) = \set{\si\cup\nu_N}{\si\in\LiftCong ST,\ N\normal G_2,\ \tau_N\sub\si} \cup \{\nab_T\}.
\end{equation}

Next we claim that
\begin{equation}\label{eq:liftST}
\LiftCong ST = \set{\mu_{D_0,N}^S, \lam_{D_0,N}^S, \rho_{D_0,N}^S, R_{D_0,N}^S}{N\normal G_1} \cup \{\mu_{S}^S, \lam_{S}^S, \rho_{S}^S, R_{S}^S\}.
\end{equation}
To see this, first note that each of the stated relations is a congruence on $S$ by the results of Section \ref{subsect:EMRT}, and is liftable to $T$ by Lemma \ref{lem:E2}.  Conversely, suppose $\si\in\LiftCong ST$, and let $\tau=\si\cup\De_{D_2}\in\Cong(T)$ and $\th=\si\restr_{D_0}\in\Cong(D_0)$.  
We claim that
\begin{equation}\label{eq:claimlr}
\si\not\sub\lam_S^S \quad\implies\quad \rho_{D_0}^S\sub\si.
\end{equation}
Indeed, if there exists $(x,y)\in\si\sm\lam_S^S$, then since $\si\sm\lam_S^S\sub R_S\sm\lam_S$, it follows from Claim~1 that $\rho_{D_0}\sub\cg xy_T\sub\tau=\si\cup\De_{D_2}$, and so $\rho_{D_0}^S=\rho_{D_0}\restr_S\sub\si$.  As in the proof of \eqref{eq:claimI0}, it follows from~\eqref{eq:claimlr} and its dual that each of the following holds, concerning $\th=\si\restr_{D_0}$:
\bena
\item \label{a} If $\si\sub\lam_S^S$ and $\si\sub\rho_S^S$, then $\th=\De_{D_0}$.
\item \label{b} If $\si\sub\lam_S^S$ and $\si\not\sub\rho_S^S$, then $\th=\lam_{D_0}^{D_0}$.
\item \label{c} If $\si\not\sub\lam_S^S$ and $\si\sub\rho_S^S$, then $\th=\rho_{D_0}^{D_0}$.
\item \label{d} If $\si\not\sub\lam_S^S$ and $\si\not\sub\rho_S^S$, then $\th=\nab_{D_0}$.
\een
Next, note that if $\si\sub R_{D_0,G_1}^S=\nab_{D_0}\cup{\H}\restr_{D_1}$, then we have $\si=\si\restr_{D_0}\cup\si\restr_{D_1}=\th\cup\nu_N$ for some $N\normal G_1$, using Lemma \ref{lem:nu_trick} in the last step.  It follows that $\si=\mu_{D_0,N}^S$, $\lam_{D_0,N}^S$, $\rho_{D_0,N}^S$ or $R_{D_0,N}^S$, in cases \ref{a}--\ref{d}, respectively.

Now suppose $\si\not\sub R_{D_0,G_1}^S$, and let $(x,y)\in\si\sm R_{D_0,G_1}^S$.  Then by symmetry we have $x\in D_1$ and $y\in S\sm H_x$.  By Claim 2, it follows that $\mu_S\sub\cg xy_T\sub\tau=\si\cup\De_{D_2}$, and so $\mu_S^S=\mu_S\restr_S\sub\si$.  Now, $\si$ contains $\th\cup\De_{D_1}$, which is equal to $\De_S$, $\lam_{D_0}^S$, $\rho_{D_0}^S$ or~$R_{D_0}^S$, in cases \ref{a}--\ref{d}, respectively.  It then follows that $\si$ contains $\mu_S^S\vee(\th\cup\De_{D_1})=\mu_S^S$, $\lam_S^S$, $\rho_S^S$ or $R_S^S$, respectively (cf.~Figure \ref{fig:diamonds2}).  In cases \ref{b}--\ref{d}, it of course follows that $\si=\lam_S^S$, $\rho_S^S$ or $R_S^S$, respectively.  In case \ref{a}, we have $\si\sub\lam_S^S\cap\rho_S^S=\mu_S^S$, and so $\si=\mu_S^S$.  This completes the proof of \eqref{eq:liftST}.

In light of \eqref{eq:claimTT} and \eqref{eq:liftST}, we may complete the proof by showing that for $N\normal G_2$ we have
\[
\tau_N = \begin{cases}
\De_S &\text{if $N=\{e\}$}\\
\mu_S^S &\text{if $\{e\}\not=N\leq H$}\\
\nab_S &\text{if $N\nleq H$.}
\end{cases}
\]
As usual, this is clear for $N=\{e\}$.  For the other two cases, we write $\si=(\nu_N)_T^\sharp$, so that $\tau_N=\si\restr_S$.

Suppose first that $\{e\}\not=N\leq H$, and fix some $x\in N\sm\{e\}$, noting that $(e,x)\in\nu_N\sub\si$.  Claim~2 gives $\mu_S\sub\cg ex_T\sub\si$, and so $\mu_S^S=\mu_S\restr_S\sub\si\restr_S=\tau_N$.  Conversely, note that $(S,N)$ is retractable as $N\leq H$ and $(S,H)$ is retractable, so that $\mu_{S,N}$ is a congruence on $T$.  Since $\nu_N\sub\mu_{S,N}$, we have $\si=(\nu_N)_T^\sharp\sub\mu_{S,N}$, and it follows that $\tau_N=\si\restr_S\sub\mu_{S,N}\restr_S=\mu_S^S$.  Thus, $\tau_N=\mu_S^S$ in this case.

Finally, suppose $N\nleq H$, and fix some $x\in N\sm H$.  As in the previous case, we have $\mu_S\sub\si$.  By assumption \ref{012iv}, there exists $a\in D_0$ such that $ae\not=ax$.  As in \eqref{eq:exN} and the discussion immediately after it, we have $(ae,ax)\not\in\L$, so $(ae,ax)\in R_S\sm\lam_S$.  Claim 1 then gives $\rho_{D_0}\sub\cg{ae}{ax}_T\sub\si$.  Dually, $\lam_{D_0}\sub\si$.  It follows that $R_{S}=\lam_{D_0}\vee\rho_{D_0}\vee\mu_{S}\sub\si$, and so $\nab_S=R_S\restr_S\sub\si\restr_S=\tau_N$, whence $\tau_N=\nab_S$.
\epf

\begin{figure}[ht]
\begin{center}
\begin{tikzpicture}[line width=0.2mm]
\small


\begin{scope}[shift={(-1.5,1.5)}]
\draw[fill=red!60, opacity=0.7] (0,0)--(0.5,0.5)--(0.5,1.5)--(0,2)--(-0.5,1.5)--(-0.5,0.5)--(0,0);
\end{scope}

\begin{scope}[shift={(1.5,1)}]
\draw (0,0)--(-3,0.5); \draw (0.5,0.5)--(-2.5,1); \draw (0.5,1.5)--(-2.5,2); \draw (0,2)--(-3,2.5);
\draw (-0.5,1.5)--(-3.5,2); \draw (-0.5,0.5)--(-3.5,1);
\end{scope}

\begin{scope}[shift={(-3,0.5)}]
\draw (0,0)--(1.5,1); \draw (0.5,0.5)--(2,1.5); \draw (0.5,1.5)--(2,2.5); \draw (0,2)--(1.5,3);
\draw (-0.5,1.5)--(1,2.5); \draw (-0.5,0.5)--(1,1.5);
\end{scope}

\begin{scope}[shift={(1.5,1)}]
\draw[fill=red!60, opacity=0.7] (0,0)--(0.5,0.5)--(0.5,1.5)--(0,2)--(-0.5,1.5)--(-0.5,0.5)--(0,0);
\end{scope}

\draw (0,0)--(1.5,1); \draw (0.5,0.5)--(2,1.5); \draw (0.5,1.5)--(2,2.5); \draw (0,2)--(1.5,3);
\draw (-0.5,1.5)--(1,2.5); \draw (-0.5,0.5)--(1,1.5);

\begin{scope}[shift={(-3,0.5)}]
\draw[fill=red!60, opacity=0.7] (0,0)--(0.5,0.5)--(0.5,1.5)--(0,2)--(-0.5,1.5)--(-0.5,0.5)--(0,0);
\end{scope}

\draw (0,0)--(-3,0.5); \draw (0.5,0.5)--(-2.5,1); \draw (0.5,1.5)--(-2.5,2); \draw (0,2)--(-3,2.5);
\draw (-0.5,1.5)--(-3.5,2); \draw (-0.5,0.5)--(-3.5,1);

\draw[fill=red!60,opacity=0.7] (0,0)--(0.5,0.5)--(0.5,1.5)--(0,2)--(-0.5,1.5)--(-0.5,0.5)--(0,0);


\draw (0,2)--(0,4); \draw (-3,2.5)--(-3,4.5); \draw (1.5,3)--(1.5,5); \draw (-1.5,3.5)--(-1.5,5.5);


\begin{scope}[shift={(0,4)}]

\draw (-1.5,4.5)--(-1.5,5.5);

\begin{scope}[shift={(-1.5,1.5)}]
\draw[fill=blue!60, opacity=0.9] (0,0)--(0.5,0.5)--(0.5,1.5)--(0,2)--(-0.5,1.5)--(-0.5,0.5)--(0,0);
\draw[fill=blue!30, opacity=0.9] (0.5,1.5)--(0,2)--(-0.5,1.5)--(-0.5,2.5)--(0,3)--(0.5,2.5)--(0.5,1.5);
\end{scope}

\begin{scope}[shift={(1.5,1)}]
\draw (0,0)--(-3,0.5); \draw (0.5,0.5)--(-2.5,1); \draw (0.5,1.5)--(-2.5,2); \draw (0,2)--(-3,2.5);
\draw (-0.5,1.5)--(-3.5,2); \draw (-0.5,0.5)--(-3.5,1);
\end{scope}

\begin{scope}[shift={(-3,0.5)}]
\draw (0,0)--(1.5,1); \draw (0.5,0.5)--(2,1.5); \draw (0.5,1.5)--(2,2.5); \draw (0,2)--(1.5,3);
\draw (-0.5,1.5)--(1,2.5); \draw (-0.5,0.5)--(1,1.5);
\end{scope}

\begin{scope}[shift={(1.5,1)}]
\draw[fill=blue!60, opacity=0.9] (0,0)--(0.5,0.5)--(0.5,1.5)--(0,2)--(-0.5,1.5)--(-0.5,0.5)--(0,0);
\end{scope}

\draw (0,0)--(1.5,1); \draw (0.5,0.5)--(2,1.5); \draw (0.5,1.5)--(2,2.5); \draw (0,2)--(1.5,3);
\draw (-0.5,1.5)--(1,2.5); \draw (-0.5,0.5)--(1,1.5);

\begin{scope}[shift={(-3,0.5)}]
\draw[fill=blue!60, opacity=0.9] (0,0)--(0.5,0.5)--(0.5,1.5)--(0,2)--(-0.5,1.5)--(-0.5,0.5)--(0,0);
\end{scope}

\draw (0,0)--(-3,0.5); \draw (0.5,0.5)--(-2.5,1); \draw (0.5,1.5)--(-2.5,2); \draw (0,2)--(-3,2.5);
\draw (-0.5,1.5)--(-3.5,2); \draw (-0.5,0.5)--(-3.5,1);

\draw[fill=blue!60,opacity=0.9] (0,0)--(0.5,0.5)--(0.5,1.5)--(0,2)--(-0.5,1.5)--(-0.5,0.5)--(0,0);

\colvertx{-1.5}{5.5}{black}{0.08}

\colvertx{-1.5}{1.5}{black}{0.08}

\begin{scope}[shift={(-1.5,1.5)}]
\foreach \x/\y in {0/2,0/3} {\colvertx{\x}{\y}{white}{0.08}}
\end{scope}

\begin{scope}[shift={(1.5,1)}]
\foreach \x/\y in {0/0,0/2} {\colvertx{\x}{\y}{white}{0.08}}
\end{scope}

\begin{scope}[shift={(-3,0.5)}]
\foreach \x/\y in {0/0,0/2} {\colvertx{\x}{\y}{white}{0.08}}
\end{scope}

\foreach \x/\y in {0/0,0/2} {\colvertx{\x}{\y}{white}{0.08}}

\end{scope}

\begin{scope}[shift={(-1.5,1.5)}]
\foreach \x/\y in {0/2} {\colvertx{\x}{\y}{white}{0.08}}
\foreach \x/\y in {0/0} {\colvertx{\x}{\y}{black}{0.08}}
\end{scope}

\begin{scope}[shift={(1.5,1)}]
\foreach \x/\y in {0/0,0/2} {\colvertx{\x}{\y}{white}{0.08}}
\end{scope}

\begin{scope}[shift={(-3,0.5)}]
\foreach \x/\y in {0/0,0/2} {\colvertx{\x}{\y}{white}{0.08}}
\end{scope}

\foreach \x/\y in {0/2} {\colvertx{\x}{\y}{white}{0.08}}
\foreach \x/\y in {0/0} {\colvertx{\x}{\y}{black}{0.08}}


\node at (0,-0.3) {$\De_T$};
\node at (-3,0.2) {$\lambda_{D_0}$};
\node at (1.8,0.7) {$\rho_{D_0}$};
\node at (-1.4,1.2) {$R_{D_0}$};

\node[rectangle,fill=white,inner sep=0] at (0,2.3) {$\mu_{D_0,G_1}$};
\node[rectangle,fill=white,inner sep=0] at (-3.6,2.8) {$\lambda_{D_0,G_1}$};
\node[rectangle,fill=white,inner sep=0] at (2.1,3.2) {$\rho_{D_0,G_1}$};
\node[rectangle,fill=white,inner sep=0] at (-2.1,3.8) {$R_{D_0,G_1}$};

\node[rectangle,fill=white,inner sep=0] at (0.3,3.7) {$\mu_S$};
\node[rectangle,fill=white,inner sep=0] at (-3.3,4.2) {$\lambda_S$};
\node[rectangle,fill=white,inner sep=0] at (1.8,4.7) {$\rho_S$};
\node[rectangle,fill=white,inner sep=0] at (-1.2,5.2) {$R_S$};

\node[rectangle,fill=white,inner sep=0] at (-0.2,6.3) {$\mu_{S,H}$};
\node[rectangle,fill=white,inner sep=0] at (-3.2,6.8) {$\lambda_{S,H}$};
\node[rectangle,fill=white,inner sep=0] at (1.5,7.3) {$\rho_{S,H}$};
\node at (-1.4,7.8) {$R_{S,H}$};

\node[rectangle,fill=white,inner sep=0] at (-.9,8.6) {$R_{S,G_2}$};
\node[rectangle,fill=white,inner sep=0] at (-1.4,9.8) {$\nab_T$};

\begin{scope}[scale=1.2,shift={(-8,2.5)}]
\draw[fill=red!60,opacity=0.7] (0,0)--(0.5,0.5)--(0.5,1.5)--(0,2)--(-0.5,1.5)--(-0.5,0.5)--(0,0);
\node at (0,-0.3) {$\{\id_{G_1}\}$};
\node at (0, 2.3) {$G_1$};
\foreach \x/\y in {0/0,0/2} {\colvertx{\x}{\y}{white}{0.06}}
\end{scope}

\begin{scope}[scale=1.2,shift={(-5.25,2.5)}]
\draw[fill=blue!60, opacity=0.9] (0,0)--(0.5,0.5)--(0.5,1.5)--(0,2)--(-0.5,1.5)--(-0.5,0.5)--(0,0);
\draw[fill=blue!30, opacity=0.9] (0.5,1.5)--(0,2)--(-0.5,1.5)--(-0.5,2.5)--(0,3)--(0.5,2.5)--(0.5,1.5);
\node at (0,-0.3) {$\{\id_{G_2}\}$};
\node at (0, 2.3) {$H$};
\node at (0, 3.3) {$G_2$};
\foreach \x/\y in {0/0,0/2,0/3} {\colvertx{\x}{\y}{white}{0.06}}
\end{scope}

\end{tikzpicture}
\caption{Hasse diagram for $\Cong(T)$ (right), for $T$ as in Proposition \ref{prop:small012}, in terms of Hasse diagrams for $\NN(G_1)$ (left) and~$\NN(G_2)$ (middle); cf.~Remark \ref{rem:small012}.}
\label{fig:small012}
\end{center}
\end{figure}
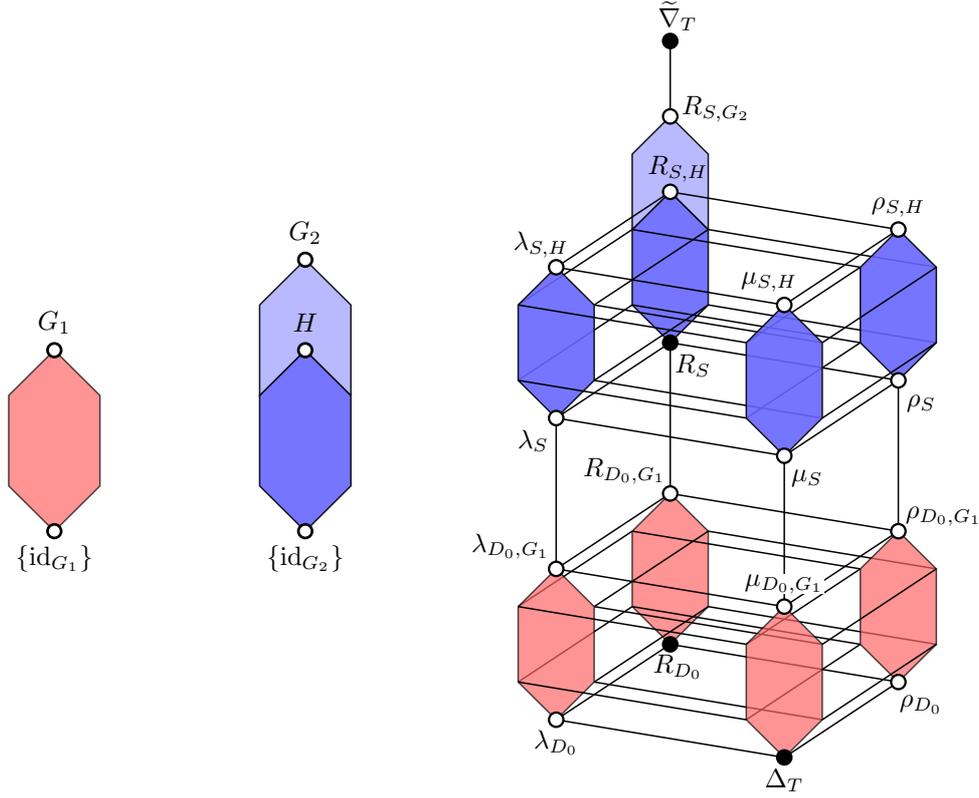

\begin{rem}\label{rem:small012}
The structure of $\Cong(T)$, for $T$ as in Proposition \ref{prop:small012}, is shown in Figure \ref{fig:small012}.
\bit
\item The sublattice $\Si_1=\set{\mu_{D_0,N},\ \lam_{D_0,N},\ \rho_{D_0,N},\ R_{D_0,N}}{N\normal G_1}$ is isomorphic to the direct product of $\NN(G_1)$ and the four-element diamond lattice.  See the red part of Figure \ref{fig:small012}.
\item The sublattice $\Si_2=\set{\mu_{S,N},\ \lam_{S,N},\ \rho_{S,N},\ R_{S,N}}{N\normal G_2,\ N\leq H}$ is isomorphic to the direct product of the interval $\big[\{\id_{G_2}\},H\big]$ in $\NN(G_2)$ with the diamond lattice.  See the dark blue part of Figure \ref{fig:small012}.
\item Congruences from the ``top diamond'' of $\Si_1$ are covered by those from the ``bottom diamond'' of $\Si_2$: 
$\ka_{D_0,G_1}\subset\ka_{S,\{\id_{G_2}\}}=\ka_S$, where $\ka$ denotes any of $\mu$, $\lam$, $\rho$ or $R$.
\item The sublattice $\Si_3=\set{R_{S,N}}{N\normal G_2}$ is isomorphic to $\NN(G_2)$, and contains the sublattice $\set{R_{S,N}}{N\normal G_2,\ N\leq H}$ of $\Si_2$.  See the light blue part of Figure \ref{fig:small012}.
\item The universal congruence $\nab_T$ sits above $R_{S,G_2}$, the top element of $\Si_3$.
\eit
See Figures \ref{fig:P}, \ref{fig:QP}, \ref{fig:Beven}, and especially \ref{fig:J}-\ref{fig:Jpmeven}, for some examples. 
\end{rem}

\subsection{Visualisation conventions}
\label{subsect:viscon}

In what follows we will present a series of applications of our theory.
Throughout we will
accompany the results with visual representations of the congruence lattices we compute.
Their Hasse diagrams will be additionally decorated to emphasise the link with the theory from this section.
Here are the main conventions that will be adhered to when representing the lattices $\Cong(I_r)$ in the generic case: i.e., for $r\geq k$, in the notation of this section:
\bit
\item
Each lattice will consist of two parts: the lower part isomorphic to $\Cong(I_k)$, and the upper part obtained by stacking 
the normal subgroup lattices $\NN(G_{k+1}),\NN(G_{k+2}),\dots$.
\item
In the upper part, the stacked copies of $\NN(G_{k+1}),\NN(G_{k+2}),\dots$ will be highlighted by means of green
background, linking back to Figure \ref{fig:E2}; see for example Figures \ref{fig:T} and \ref{fig:QT}.
\item
The lower part of the lattice is always one of two types, described either by Proposition \ref{prop:small01} or by Proposition \ref{prop:small012}.
\item
When the lower part is described by Proposition \ref{prop:small01}, then, referring to Figure \ref{fig:small01}, there will be no special colours for the bottom diamond, while the sublattice isomorphic to $\NN(G)$ will have blue background; for an example see Figure \ref{fig:Bodd}.
\item
When the lower part is described by Proposition \ref{prop:small012}, 
the copies of $\NN(G_1)$, of the sublattice of $\NN(G_2)$ lying below $H$, and of $\NN(G_2)$ itself will be indicated by the backgrounds of the same colours as in Figure \ref{fig:small012};
for some examples, see Figures \ref{fig:P} and \ref{fig:QP} (with $H=G_2$) and Figures~\ref{fig:Jeven} and~\ref{fig:Jpmeven} (with $H\not=G_2$).
\item
As further useful reference, the vertices corresponding to Rees congruences and the trivial congruence will be solid, while all the other vertices will be typeset with a white background.  
\eit

\section{Transformation categories}\label{sect:T}

Many of the most fundamental and well-studied categories consist of sets and various kinds of mappings, and our first applications of the results of Sections \ref{sect:IE} and \ref{sect:chains} are to several such categories and their ideals.  In Section \ref{subsect:defnT} we introduce the transformation category $\T=\T(\C)$, which consists of all mappings between sets from $\C$, a non-empty set of finite non-empty sets.  In Section \ref{subsect:prelimT} we describe Green's relations in $\T$ and its ideals, introduce a number of important subcategories, and prove a technical result (Lemma \ref{lem:P2T}) that will be used to show that $\T$ and its ideals satisfy various properties involved in the statements of the general results from Section~\ref{sect:chains} (cf.~Section \ref{subsect:properties}).  
These preparations done, in Section \ref{subsect:T} we come to the first application of the machinery developed in
Section \ref{sect:chains}: we rapidly and in one fell swoop classify the congruences on~$\T$ and on all its ideals,
and describe the lattices they form
(Theorem \ref{thm:T} and Figure~\ref{fig:T}).
In Section~\ref{subsect:QT} we describe even more rapidly the congruences on the ideals of certain subcategories of $\T$ that we call \emph{planar and annular reducts}.  Finally in Section \ref{subsect:OtherT} we describe how to do the same for some further categories of mappings.

\subsection{Definitions and preliminaries on $\T$}\label{subsect:defnT}

For the duration of Section \ref{sect:T}, we fix a non-empty set $\C$ of non-empty finite sets; note that $\C$ itself may be infinite.  We define the parameter
\[
\xi = \xi(\C) = {\LSUB}\bigset{|A|}{A\in\C}.
\]
Here, for a set $\Si$ of cardinals, $\LSUB(\Si)$ denotes the \emph{least strict upper bound} of $\Si$: i.e., the least cardinal that is strictly greater than all members of $\Si$.  Note that $\xi$ as above is at least $2$, and could be $\om$ if $\C$ contains arbitrarily large sets.

For $A,B\in\C$ we let
$\T_{A,B}$ stand for the set of all functions $\al:A\to B$, and define
\[
\T = \T(\C) = \bigset{(A,\al,B)}{A,B\in\C,\ \al\in\T_{A,B}}.
\]
For $A,B,C\in\C$, and for $\al\in\T_{A,B}$ and $\be\in\T_{B,C}$, we define
\[
\bd(A,\al,B)=A \COMMA \br(A,\al,B)=B \COMMA (A,\al,B)\cdot(B,\be,C)=(A,\al\be,C),
\]
where $\al\be\in\T_{A,C}$ denotes the composition of $\al$ and $\be$ (performing $\al$ first, then $\be$).  Then 
\[
 (\T,\C,\bd,\br,\cdot)
\]
is a partial semigroup, indeed a regular category \cite[Proposition 2.1]{Sandwich2}, which we call the \emph{transformation category (over $\C$)}.  For $A\in\C$, the endomorphism monoid $\T_A=\T_{A,A}$ is the \emph{full transformation semigroup} over $A$.  Note that when $\C=\{A\}$ consists of a single finite set, the category $\T$ is simply the full transformation semigroup $\T_A$ \cite{GMbook}.
Thus, every result we prove for the category $\T$ has a natural interpretation as a result concerning finite transformation semigroups.

In order to avoid cumbersome notation, we will identify an element ${(A,\al,B)\in\T}$ with the mapping $\al\in\T_{A,B}$, but we regard the sets $A$ and $B$ as ``encoded'' in $\al$ via $\bd(\al)=A$ and $\br(\al)=B$.

\subsection{Green's relations and multiplicative properties in $\T$}\label{subsect:prelimT}

For non-empty sets $A,B$, and for a function $\al:A\to B$, we write as usual
\[
\ker(\al)=\bigset{(x,y)\in A\times A}{x\al=y\al} \COMMA \im(\al) = \set{x\al}{x\in A} \COMMA \rank(\al)=|{\im(\al)}|
\]
for the \emph{kernel}, \emph{image} and \emph{rank} of $\al$.  
We remark that $\ker(\al)$ uniquely determines $\bd(\al)$, but that $\im(\al) $ does not generally determine $\br(\al)$.

 In what follows, we will often use the following result without explicit reference.  
 The proof can be found in \cite[Section 4]{Sandwich2}.

\begin{lemma}\label{lem:Green_T}
\ben
\item \label{GT1} For $\al,\be\in\T$, we have
\bit
\item $(\al,\be)\in\R \iff (\bd(\al)=\bd(\be)$ and) $\ker(\al)=\ker(\be)$,
\item $(\al,\be)\in\L \iff \br(\al)=\br(\be)$ and $\im(\al)=\im(\be)$,
\item $(\al,\be)\in\J\iff(\al,\be)\in\D\iff\rank(\al)=\rank(\be)$.
\eit
\item \label{GT2} The $\J=\D$-classes of $\T$ are the sets
\bit
\item $D_q = D_q(\T) = \set{\al\in\T}{\rank(\al)=q}$ for each $1\leq q<\xi$.
\eit
These are all regular and stable, and they form a chain: $D_1<D_2<\cdots$.
\item \label{GT3} The ideals of $\T$ are the sets
\bit
\item $I_r = I_r(\T) = \set{\al\in\T}{\rank(\al)\leq r}$ for each $1\leq r<\xi$, and
\item $I_\om = I_\om(\T) = \T $ in the case that $\xi=\om$.
\eit
These are all regular and stable, and they form a chain: $I_1\subset I_2\subset\cdots$.
\item \label{GT4} The $\H$-class of any idempotent from $D_q$ is isomorphic to $\S_q$, the symmetric group of degree~$q$.  \epfres
\een
\end{lemma}

Thus, $\T$ is a chain of ideals, as in Definition \ref{defn:chain}, though we note that the indexing on $\J=\D$-classes starts at $1$ instead of $0$.
Note that $I_1=D_1$ is the minimal ideal of $\T$.  By Lemma \ref{lem:SD}, the $\D$-classes and ideals of an ideal $I_r=I_r(\T)$ are precisely the $D_q$ and $I_q$ with $q\leq r$.

For any $\Om\sub\T$, and for any $1\leq q<\xi$, we write 
\bit
\item $D_q(\Om)=D_q\cap\Om$ for the set of transformations from $\Om$ of rank $q$,
\item $I_q(\Om)=I_q\cap\Om$ for the set of transformations from $\Om$ of rank at most $q$.
\eit
These notational conventions are for convenience only, and do not mean to imply that the sets $D_q(\Om)$ or $I_q(\Om)$ are themselves $\D$-classes or ideals of $\T$.

Our next goal is to prove a technical result (Lemma \ref{lem:P2T} below) that will be used to show, among other things, that certain ideals $I_r$ of $\T$ satisfy $\Mult{I_r}$, and hence also $\Sep{I_r}$; cf.~Definitions \ref{defn:Sep} and \ref{defn:Mult}.  Since we will also use this lemma in Section \ref{subsect:QT} when studying subcategories of $\T$, it will be convenient to first introduce some notation and terminology.

If $\al\in\T$, we write $\al = \partpermIII{A_1}\cdots{A_q}{b_1}\cdots{b_q}$
to indicate that $\rank(\al)=q$, $\im(\al)=\{b_1,\ldots,b_q\}$ and $b_i\al^{-1}=A_i$ ($1\leq i\leq q$).  Note that from this notation we may infer that $\bd(\al)=A_1\cup\cdots\cup A_q$, but the set $\br(\al)$ is not uniquely determined; thus, it will always be made clear which hom-set $\al$ belongs to.

In what follows, we assume that each set from $\C$ has a fixed total order; we will denote all these orders by $\leq$.  Let $\al\in\T_{A,B}$, where $A,B\in\C$, and write $A=\{a_1<\cdots<a_k\}$; by this we mean that $A=\{a_1,\ldots,a_k\}$ and $a_1<\cdots<a_k$.  We say $\al$ is
\bit
\item \emph{order-preserving} if $a_1\al\leq\cdots\leq a_k\al$,
\item \emph{order-reversing} if $a_1\al\geq\cdots\geq a_k\al$,
\item \emph{orientation-preserving} if $a_i\al\leq\cdots\leq a_k\al\leq a_1\al\leq\cdots\leq a_{i-1}\al$ for some $1\leq i\leq k$,
\item \emph{orientation-reversing} if $a_i\al\geq\cdots\geq a_k\al\geq a_1\al\geq\cdots\geq a_{i-1}\al$ for some $1\leq i\leq k$.
\eit
Note that order-preserving (or order-reversing) transformations are orientation-preserving (or orientation-reversing), respectively; take $i=1$ in the above definitions.  For any subset $\Om\sub\T$, we will write 
\bit
\item $\OrdP\Om$ for the set of all order-preserving transformations from $\Om$, 
\item $\OrdR\Om$ for the set of all order-reversing transformations from $\Om$, 
\item $\OrdPR\Om = \OrdP\Om\cup\OrdR\Om$,
\item $\OriP\Om$ for the set of all orientation-preserving transformations from $\Om$, 
\item $\OriR\Om$ for the set of all orientation-reversing transformations from $\Om$, 
\item $\OriPR\Om = \OriP\Om\cup\OriR\Om$.
\eit
This notation is not meant to imply that $\OrdP\Om$ and $\OrdR\Om$ are disjoint in general.  In fact, ${\OrdP\T\cap\OrdR\T}$ consists of all constant mappings, while ${\OriP\T\cap\OriR\T}$ consists of all orientation-preserving mappings of rank at most $2$.

We also call the elements of $\OrdP\T$ (respectively, $\OrdR\T$, $\OriP\T$ or $\OriR\T$) \emph{planar} (respectively, \emph{anti-planar}, \emph{annular} or \emph{anti-annular}), which explains the choice of the symbols $\mathscr P$ and~$\mathscr A$.  The reason for these alternative names is as follows.  Let $\al\in\T_{A,B}$, write $A=\{a_1<\cdots<a_k\}$ and $B=\{b_1<\cdots<b_l\}$, and let $m=\max(k,l)$.  We associate to $\al$ the graph with vertex set $A\cup B'$, where $B'=\{b_1'<\cdots<b_l'\}$ is assumed to be disjoint from $A$, and with an edge $\{x,(x\al)'\}$ for each $x\in A$.  The transformation $\al$ is then planar (i.e., order-preserving) if and only if this graph can be drawn in the plane $\RRR^2$ with:
\bit
\item vertex $a_i$ at $(i,1)$ for each $1\leq i\leq k$, 
\item vertex $b_i'$ at $(i,0)$ for each $1\leq i\leq l$, 
\item edges all in the rectangle $\bigset{(x,y)\in\mathbb R^2}{0\leq x\leq m+1,\ 0\leq y\leq1}$, and
\item with no edge crossings in the interior of this rectangle.
\eit
The transformation $\al$ is annular if and only if the above conditions are satisfied, but with the left- and right-hand borders of the rectangle identified to form a cylinder, and with edges allowed to ``wrap around'' this identified border.

For example, with $A=[9]=\{1<\cdots<9\}$ and $B=[7]=\{1<\cdots<7\}$, consider the transformations $\al,\be\in\T_{A,B}$ given by
\[
(1\al,\ldots,9\al)=(2,2,3,4,6,6,6,7,7) \AND (1\be,\ldots,9\be)=(6,7,1,1,1,4,5,6,6).
\]
By definition, we have $\al\in\OrdP{\T_{A,B}}$ and $\be\in\OriP{\T_{A,B}}$.  Graphs representing these transformations, and illustrating the above geometrical conditions, are pictured in Figure \ref{fig:PT_AT}.  Two diagrammatic representations of $\be$ are given; in the second, two lines ``wrap around'' the back of the cylinder.

\begin{figure}[ht]
\begin{center}
\begin{tikzpicture}[scale=.6,inner sep=0]
\begin{scope}[shift={(0,0)}]	
\uuvs{1,...,9}
\lvs{1,...,7}
\foreach \x/\y in {
1/2,
2/2,
3/3,
4/4,
5/6,
6/6,
7/6,
8/7,
9/7
}
{\draw(\x,4)--(\y,0);}
\foreach \x in {1,...,9} {\node () at (\x,4.5) {\tiny $\x$};}
\foreach \x in {1,...,7} {\node () at (\x,-.5) {\tiny $\phantom{'}\x'$};}
\draw(0.6,2)node[left]{$\al=$};
\end{scope}
\begin{scope}[shift={(14,0)}]	
\foreach \x/\y in {
3/1,
4/1,
5/1,
6/4,
7/5,
8/6,
9/6
}
{\draw(\x,4)--(\y,0);}
\draw[red,very thick](1,4)--(6,0);
\draw[blue,very thick](2,4)--(7,0);
\foreach \x in {1,...,9} {\node () at (\x,4.5) {\tiny $\x$};}
\foreach \x in {1,...,7} {\node () at (\x,-.5) {\tiny $\phantom{'}\x'$};}
\draw(0.6,2)node[left]{$\be=$};
\uuvs{1,...,9}
\lvs{1,...,7}
\end{scope}
\begin{scope}[shift={(19,-8)}]	
\nc\rr{4.6}
\draw[fill=gray!20] (0,0) ellipse (4.6 and 1);
\fill[gray!20](-\rr,0)--(-\rr,4)--(\rr,4)--(\rr,0)--(-\rr,0);
\draw(-\rr,0)--(-\rr,4) (\rr,4)--(\rr,0);
\draw[fill=gray!40] (0,4) ellipse (4.6 and 1);
\foreach \x/\t in {
1/30,
2/47,
3/62,
4/76,
5/90,
6/104,
7/118,
8/133,
9/150
}
{
\pgfmathsetmacro{\CosValue}{cos(180+\t)}
\pgfmathsetmacro{\SinValue}{sin(180+\t)}
\node (U\x) at (\rr*\CosValue,\SinValue+4) {};
\node (L\x) at (\rr*\CosValue,\SinValue) {};
}
\draw[red,very thick] (U1)--(-\rr,3) (\rr,2)--(L6);
\draw[blue,very thick] (U2)--(-\rr,2) (\rr,1)--(L7);
\draw (U3)--(L1);
\draw (U4)--(L1);
\draw (U5)--(L1);
\draw (U6)--(L4);
\draw (U7)--(L5);
\draw (U8)--(L6);
\draw (U9)--(L6);
\foreach \t in {
30,
47,
62,
76,
90,
104,
118,
133,
150
}
{
\pgfmathsetmacro{\CosValue}{cos(180+\t)}
\pgfmathsetmacro{\SinValue}{sin(180+\t)}
\fill (\rr*\CosValue,\SinValue+4)circle(.17);
}
\foreach \t in {
30,
47,
62,
76,
90,
104,
118
}
{
\pgfmathsetmacro{\CosValue}{cos(180+\t)}
\pgfmathsetmacro{\SinValue}{sin(180+\t)}
\fill (\rr*\CosValue,\SinValue)circle(.17);
}
\draw(0,6)node{$\mathbin{\rotatebox[origin=c]{-90}{$\equiv$}}$};
\end{scope}
\end{tikzpicture}
\caption{Left: an order-preserving (planar) transformation $\al\in\OrdP{\T_{A,B}}$, with ${A=[9]}$, and ${B=[7]}$.  Right: an orientation-preserving (annular) transformation $\be\in\OriP{\T_{A,B}}$.}
\label{fig:PT_AT}
\end{center}
\end{figure}

We will refer to both $\OrdP\T$ and $\OrdPR\T$ as \emph{planar reducts} of $\T$, and to $\OriP\T$ and $\OriPR\T$ as \emph{annular reducts} of $\T$.

As noted above, the main purpose of the next lemma is to ensure that certain ideals~$I_r$ of $\T$ (and of its planar and annular reducts) satisfy $\Mult{I_r}$.

\newpage

\begin{lemma}\label{lem:P2T}
Let $A,B\in\C$, and suppose $\al,\be\in\T_{A,B}$ with $r=\rank(\al)\geq\rank(\be)=q$ and $r\geq2$.
\ben
\item \label{P2Ti} If $q<r$, then there exists $\ga\in\OrdP{I_{r}}$ such that $\ga\al\in D_{r-1}$ and $\ga\be\in I_{r-1}\sm H_{\ga\al}$.
\item \label{P2Tii} If $q=r$ and $(\al,\be)\not\in\H$, then there exists $\ga\in\OrdP{I_{r}}$ such that, swapping $\al,\be$ if necessary,
\[
\text{$[\al\ga\in D_r$ and $\be\ga\in I_{r-1}]$ \ \ \ or \ \ \ $[\ga\al\in D_r$ and $\ga\be\in I_{r-1}]$.}
\]
\item \label{P2Tiii} If $\be\in H_\al\sm\{\al\}$, then there exists $\ga\in\OrdP{I_{r}}$ such that $\ga\al\in D_{r-1}$ and $\ga\be\in I_{r-1}\sm H_{\ga\al}$.
\een
\end{lemma}

\pf
For all parts of the proof, we write $\al=\partpermIII{A_1}\cdots{A_r}{b_1}\cdots{b_r}$, making no assumptions about the ordering on points unless specified.  For each $1\leq i\leq r$, we also fix some $a_i\in A_i$.

\pfitem{\ref{P2Ti}}  Since $q<r$, we may assume without loss of generality that $b_1\not\in\im(\be)$.  Let $\ga\in\OrdP{\T_A}$ be any (order-preserving) map with image $\{a_1,\ldots,a_{r-1}\}$.  Then $\im(\ga\al)=\{b_1,\ldots,b_{r-1}\}$, which gives $\ga\al\in D_{r-1}$.  But also $\im(\ga\be)\sub\im(\be)$ and $b_1\not\in\im(\be)$, so that $b_1\not\in\im(\ga\be)$.  Thus, $\ga\be\in I_{r-1}$ and $\im(\ga\al)\not=\im(\ga\be)$, and so $(\ga\al,\ga\be)\not\in\L$.  (This of course implies that $(\ga\al,\ga\be)\not\in\H$.)

\pfitem{\ref{P2Tii}}  There are two cases to consider.

\pfcase{1}  Suppose first that $(\al,\be)\not\in\L$, so that $\im(\al)\not=\im(\be)$.  This time we will assume that $b_1<\cdots<b_r$.  
Since $q=r$ and $\im(\al)\not=\im(\be)$, we have $b_i\not\in\im(\be)$ for some $1\leq i\leq r$.  

\pfcase{1.1}  If $1<i<r$, then take any $\ga\in\OrdP{D_r(\T_B)}$ extending the partial map $\partpermIII{b_1}\cdots{b_r}{b_1}\cdots{b_r}$ and for which $b_i\ga^{-1}=\{b_i\}$.  Then $\im(\al\ga)=\{b_1,\ldots,b_r\}$, which gives $\al\ga\in D_r$.  But also $\im(\be\ga)\sub\im(\ga)=\{b_1,\ldots,b_r\}$ and $b_i\not\in\im(\be\ga)$, which gives $\be\ga\in I_{r-1}$.

\pfcase{1.2}  Next suppose $i=1$.  Swapping $\al,\be$ if necessary, we may assume that ${b_1<\min(\im(\be))}$.  
Then we may take any $\ga\in\OrdP{D_r(\T_B)}$ extending the partial map $\partpermIII{b_1}\cdots{b_r}{b_1}\cdots{b_r}$ and such that $b_1\ga^{-1}=\set{x\in B}{x\leq b_1}$.  Since $b_1\ga^{-1}$ is disjoint from $\im(\be)$, the desired conclusion follows.

\pfcase{1.3}  The case in which $i=r$ is symmetrical to the previous case.

\pfcase{2}  Now suppose $(\al,\be)\not\in\R$, so that $\ker(\al)\not=\ker(\be)$.  We may assume without loss of generality that $(a_1,a_2)\in\ker(\be)$.  Then for any $\ga\in\OrdP{\T_A}$ with $\im(\ga)=\{a_1,\ldots,a_r\}$, we have $\ga\al\in D_r$ and $\ga\be\in I_{r-1}$.

\pfitem{\ref{P2Tiii}}  Here we have $\be=\partpermIII{A_1}\cdots{A_r}{b_{1\pi}}\cdots{b_{r\pi}}$ for some permutation $\pi\in\S_r$, and without loss of generality we may assume that $1\pi=r$.  We then take any $\ga\in\OrdP{\T_A}$ with $\im(\ga)=\{a_1,\ldots,a_{r-1}\}$; then $\ga\al,\ga\be\in D_{r-1}$ but $b_r\in\im(\ga\be)\sm\im(\ga\al)$, so that $(\ga\al,\ga\be)\not\in\L$.
\epf

\subsection{Congruences on ideals of $\T$}\label{subsect:T}

We now finally arrive at the first application of the theory developed in Section \ref{sect:chains}, to determine the congruences
of an arbitrary transformation category $\T=\T(\C)$ and of all its ideals.

Before we begin, we recall that historically the classification of congruences on $\T$ in the case when $|\C|=1$ (i.e., when $\T$ is the full transformation monoid) was obtained by Mal'cev \cite{Malcev1952} in~1952, the inaugural result in this area.
The congruences on ideals of full transformation monoids were determined by Klimov \cite{Klimov1977} some 15 years later.
It is worth pointing out that Klimov's classification does not subsume Mal'cev's, but instead relies on it.
By way of contrast, our methods yield a uniform description of congruences for $\T$ and all its ideals together.

For each $1\leq q<\xi$, fix some $X_q\in\C$ with $|X_q|\geq q$; such an $X_q$ exists by definition of $\xi$.  Note that it is possible to have $X_{q_1}=X_{q_2}$ with $q_1\not=q_2$; indeed, this is certainly the case if $\C$ contains only one set; moreover, if $\xi$ is finite, then we could take each $X_q$ to be a single fixed set of size $\xi-1$.  For each $q$, we write $n_q=|X_q|$, and swapping if necessary, we may in fact assume that $X_q=[n_q]=\{1,\ldots,n_q\}$.  For a permutation $\pi\in\S_q$, we define the transformation $\pi^\natural\in\T_{X_q}=\T_{X_q,X_q}$ by
\[
\pi^\natural = \partpermIV 1{q-1}{\{q,\ldots,n_q\}}{1\pi}{(q-1)\pi}{q\pi}.
\]
For $\Om\sub\S_q$, we define $\Om^\natural=\set{\pi^\natural}{\pi\in\Om}$.  In particular, we define $\G_q=\S_q^\natural$.  This is a group $\H$-class in $D_q$, and its identity is $\id_q^\natural = \partpermIV 1{q-1}{\{q,\ldots,n_q\}}{1}{q-1}{q}$.  Every normal subgroup of $\G_q$ is of the form $N^\natural$ for some $N\normal\S_q$.  So an ideal $I_r$ ($r\geq2$) contains the IN-pairs $(I_q,N^\natural)$ for each $1\leq q<r$ and $N\normal\S_{q+1}$, and these lead to the congruences
\[
R_{I_q,N^\natural}^{I_r} = R_{I_q}^{I_r}\cup\nu_{N^\natural}.
\]
We will soon see in Theorem~\ref{thm:T} that these, together with the universal congruence $\nab_{I_r}$, are \emph{all} the congruences on $I_r$ for $r\geq2$.  In order that Theorem \ref{thm:T} be as transparent as possible, we first provide in Lemma \ref{lem:nuN_T} a concrete description of the $\nu_{N^\natural}$ relations.

To do so, consider two $\H$-related elements $\al,\be\in D_q$.  We may then write 
\begin{equation}\label{eq:H_T}
\al = \partpermIII{A_1}\cdots{A_q}{b_1}\cdots{b_q} \AND \be = \partpermIII{A_1}\cdots{A_q}{b_{1\phi}}\cdots{b_{q\phi}} ,
\end{equation}
for some permutation $\phi\in\S_q$, and we define $\phi(\al,\be)=\phi$.  Note that $\phi$ depends on the ordering chosen on the $b_i$, so $\phi(\al,\be)$ is only well-defined up to conjugacy in $\S_q$.  But then it follows that for a normal subgroup $N$ of $\S_q$, the property of $\phi(\al,\be)$ belonging to $N$ is independent of the specific ordering.  Keeping this in mind we have the following:

\begin{lemma}\label{lem:nuN_T}
If $1\leq q<\xi$, and if $N\normal\S_q$, then $\nu_{N^\natural} = \bigset{(\al,\be)\in\H\restr_{D_q}}{\phi(\al,\be)\in N}$.
\end{lemma}

\pf
Let $(\al,\be)\in\H\restr_{D_q}$, say with $\al,\be\in\T_{A,B}$, and write $\al,\be$ as in \eqref{eq:H_T},
where $\phi=\phi(\al,\be)$.  
For each $1\leq i\leq q$, fix some $a_i\in A_i$.  Fix some partition $\{B_1,\ldots,B_q\}$ of $B$ such that $b_i\in B_i$ for each $i$.  Writing $Q=\{q,\ldots,n_q\}$, define
\begin{align*}
\ga_1 &= \partpermIV1{q-1}Q{a_1}{a_{q-1}}{a_q} \in \T_{X_q,A}, 
& \de_1 &= \partpermIV{A_1}{A_{q-1}}{A_q}1{q-1}q \in \T_{A,X_q}, \\
\ga_2 &= \partpermIV{B_1}{B_{q-1}}{B_q}1{q-1}q \in \T_{B,X_q}, 
& \de_2 &= \partpermIV1{q-1}Q{b_1}{b_{q-1}}{b_q} \in \T_{X_q,B}.
\end{align*}
Then since $(\al,\be) = (\de_1\id_q^\natural\de_2,\de_1\phi^\natural\de_2)$ and $(\id_q^\natural,\phi^\natural) = (\ga_1\al\ga_2,\ga_1\be\ga_2)$, Lemma \ref{lem:nuN}\ref{nuN1} gives
\[
(\al,\be)\in\nu_{N^\natural} \iff (\id_q^\natural,\phi^\natural)\in\nu_{N^\natural} \iff \phi^\natural\in N^\natural \iff \phi\in N. \qedhere
\]
\epf

We are now ready to state and prove the main result of this section.

\begin{thm}\label{thm:T}
Let $\C$ be a non-empty set of finite non-empty sets, let $\T=\T(\C)$ be the transformation category over $\C$, and keep the above notation.
\ben
\item \label{T1} The minimal ideal $I_1$ is a partial rectangular band, so its congruence lattice is described by Proposition \ref{prop:PRB}.  In particular:
\[
\Cong(I_1)\cong \prod_{A\in\C}\Eq(A).
\]
\item \label{T2} For $r\geq2$ (including $r=\om$ if $\xi=\om$),
\[
\Cong(I_r) = \{\De_{I_r}\}\cup
\set{R_{I_q,N^\natural}^{I_r}}{1\leq q<r,\ N\normal\S_{q+1}}\cup\{\nab_{I_r}\}
\]
forms a chain.
\een
\end{thm}

\pf
\firstpfitem{\ref{T1}} 
The minimal ideal of a stable, regular partial semigroup is always stable, regular and $\D$-universal; since $I_1$ is $\H$-trivial, it is a partial rectangular band, and so $\Cong(I_1)$ is indeed described by Proposition \ref{prop:PRB}.  For any $A\in\C$, $I_1(\T_A)$ contains only one $\R$-class (as there is only one kernel with one block), and the $\L$-classes of $I_1(\T_A)$ are indexed by the elements of $A$ (one for each image set of size $1$), so the isomorphism of $\Cong(I_1)$ with $\prod_{A\in\C}\Eq(A)$ is clear, after identifying $I_1$ with the abstract partial rectangular band as in Lemma \ref{lem:PRB}.

\pfitem{\ref{T2}}  We first consider the ideal $I_2=D_1\cup D_2$.  Here we must show that
\[
\Cong(I_2) = \{ \De_{I_2}, R_{I_1}^{I_2}, R_{I_1,\G_2}^{I_2}, \nab_{I_2}\}.
\]
But this follows quickly from Proposition \ref{prop:small01}.  Indeed, since~$D_1$ is $\wt\L$-trivial (but not $\wt\R$-trivial), it suffices to verify the following (cf.~Remark \ref{rem:small01_2}):
\bena
\item \label{Ta} For every $\al\in D_2$ and $\be\in D_2\sm H_\al$ with $\al\sim\be$, there exists $(\al',\be')\in\cg\al\be_{I_2}$ with $\al'\in D_2$ and $\be'\in D_1$.
\item \label{Tb} For every $\al,\be\in I_2$ with $\al\sim\be$ and $(\al,\be)\not\in\L$, we have $\rho_{D_1}^{I_2}\sub\cg\al\be_{I_2}$.
\item \label{Tc} Writing $\G_2=\{\al,\be\}$, there exists $\ga\in D_1$ such that $\ga\al\not=\ga\be$ (with both products defined).
\een
Item \ref{Ta} follows from Lemma \ref{lem:P2T}\ref{P2Tii}, and \ref{Tc} is easily checked (in fact, \emph{any} $\ga$ for which the products are defined has this property).  
For \ref{Tb}, let $\al,\be\in I_2$ with $\al\sim\be$ and $(\al,\be)\not\in\L$, say with $\al,\be\in\T_{A,B}$.
We may assume without loss of generality that there exists $x\in\im(\al)\sm\im(\be)$.  Fix some $u\in x\al^{-1}$, and write $y=u\be$, noting that $x\not=y$.  Now let~$C,D\in\C$, and let $\ga=\binom Ca$ and $\de=\binom Cb$ be arbitrary elements of $D_1(\T_{C,D})$.  To complete the proof of \ref{Tb}, it suffices to show that $(\ga,\de)\in\cg\al\be_{I_2}$.  Let $\ve_1=\binom Cu\in\T_{C,A}$, and let $\ve_2$ be any element of $I_2(\T_{B,D})$ such that $x\ve_2=a$ and $y\ve_2=b$.  Then $(\ga,\de)=(\ve_1\al\ve_2,\ve_1\be\ve_2)\in\cg\al\be_{I_2}$.

\aftercases
With the case of $r=2$ established, we now observe the following:
\bit
\item Every congruence on $I_2$ (as above) is of the form listed in Lemma \ref{lem:E2}, and hence we have $\LiftCong{I_2}{\T}=\Cong(I_2)$.
\item Property $\ngen{I_2}$ holds by inspection of the four congruences on $I_2$.
\item Property $\Sep{I_r}$ holds for each $2<r<\xi$; indeed, this follows from Lemma \ref{lem:P2T}, which as we have already noted implies the stronger property $\Mult{I_r}$ for all such $r$.
\eit
Thus, Theorem \ref{thm:E2} applies (with $k=2$, but again being aware of the different indexing), and the conclusion of that theorem gives precisely the desired description of $\Cong(I_r)$.  The fact that $\Cong(I_r)$ forms a chain follows from $\Cong(I_2)$ being a chain, and the normal subgroups of $\G_q\cong\S_q$ forming a chain for each $q$; cf.~Remark \ref{rem:E2} and Figure \ref{fig:E2}.
\epf

\begin{rem}\label{rem:T1}
In fact, Lemma \ref{lem:P2T} shows that $\Sep{I_2}$ holds as well. 
However, we may not use Theorem \ref{thm:E2} with $k=1$, since $\ngen{I_1}$ does not hold, and neither does $\LiftCong{I_1}\T=\Cong(I_1)$.
\end{rem}

\begin{figure}[ht]
\begin{center}
\begin{tikzpicture}[scale=0.9,inner sep=0pt]
\begin{scope}[shift={(-.5,0)}]
\foreach \x/\y/\L in {
0/0/0,
-2.5/1/12,
-1.5/1/34,
-0.5/1/13,
0.5/1/24,
1.5/1/14,
2.5/1/23,
-3/2/123,
-2/2/1234,
-1/2/134,
0/2/1324,
1/2/124,
2/2/1423,
3/2/234,
0/3/1
} {\node (\L) at (\x,\y) {};}
\foreach \x in {12,13,14,23,24,34} {\draw (0)--(\x);}
\foreach \x in {123,1234,134,1324,124,1423,234} {\draw (0,3)--(\x);}
\foreach \x in {123,124,1234} {\draw(12)--(\x);}
\foreach \x in {123,134,1324} {\draw(13)--(\x);}
\foreach \x in {124,134,1423} {\draw(14)--(\x);}
\foreach \x in {123,234,1423} {\draw(23)--(\x);}
\foreach \x in {124,234,1324} {\draw(24)--(\x);}
\foreach \x in {134,234,1234} {\draw(34)--(\x);}
\node[below right] () at (0.1,-0.1) {\small $\De_{I_1}$};
\node[above right] () at (0.1,3) {\small $\nab_{I_1}$};
\foreach \x/\y/\L in {
0/0/0,
-2.5/1/12,
-1.5/1/34,
-0.5/1/13,
0.5/1/24,
1.5/1/14,
2.5/1/23,
-3/2/123,
-2/2/1234,
-1/2/134,
0/2/1324,
1/2/124,
2/2/1423,
3/2/234,
0/3/1
}
{\colvertx{\x}{\y}{white}{0.08}}
\foreach \x in {0,3} {\colvertx0{\x}{black}{0.08}}

\end{scope}
\begin{scope}[shift={(5,0)}]
\draw[rounded corners=4,line width=10,blue!25,cap=round](0,1)--(0,2);
\draw(0,0)--(0,3);
\foreach \x in {0,...,3} {\colvertx0{\x}{white}{0.08}}
\node[right] () at (0.2,0) {\small $\De_{I_2}$};
\node[right] () at (0.2,1) {\small $R_{I_1}^{I_2}$};
\node[right] () at (0.2,3) {\small $\nab_{I_2}$};
\foreach \x in {0,1,3} {\colvertx0{\x}{black}{0.08}}
\end{scope}
\begin{scope}[shift={(8,0)}]
\draw[rounded corners=4,line width=10,blue!25,cap=round](0,1)--(0,2);
\draw[rounded corners=4,line width=10,green!25,cap=round](0,3)--(0,5);
\draw(0,0)--(0,6);
\foreach \x in {0,...,6} {\colvertx0{\x}{white}{0.08}}
\node[right] () at (0.2,0) {\small $\De_{I_3}$};
\node[right] () at (0.2,1) {\small $R_{I_1}^{I_3}$};
\node[right] () at (0.2,3) {\small $R_{I_2}^{I_3}$};
\node[right] () at (0.2,6) {\small $\nab_{I_3}$};
\foreach \x in {0,1,3,6} {\colvertx0{\x}{black}{0.08}}
\end{scope}
\begin{scope}[shift={(11,0)}]
\draw[rounded corners=4,line width=10,blue!25,cap=round](0,1)--(0,2);
\draw[rounded corners=4,line width=10,green!25,cap=round](0,3)--(0,5);
\draw[rounded corners=4,line width=10,green!25,cap=round](0,6)--(0,9);
\draw(0,0)--(0,10);
\foreach \x in {0,...,10} {\colvertx0{\x}{white}{0.08}}
\node[right] () at (0.2,0) {\small $\De_{\T}$};
\node[right] () at (0.2,1) {\small $R_{I_1}^{\T}$};
\node[right] () at (0.2,3) {\small $R_{I_2}^{\T}$};
\node[right] () at (0.2,6) {\small $R_{I_3}^{\T}$};
\node[right] () at (0.2,10) {\small $\nab_{\T}$};
\foreach \x in {0,1,3,6,10} {\colvertx0{\x}{black}{0.08}}
\end{scope}
\end{tikzpicture}
\caption{Hasse diagrams of $\Cong(I_r)$ for $r=1,2,3,4$ (left to right), in the category $\T=\T_4$ where~$\C$ consists of a single set of size $4$; cf.~Theorem \ref{thm:T} and Remark \ref{rem:T2}.}
\label{fig:T}
\end{center}
\end{figure}
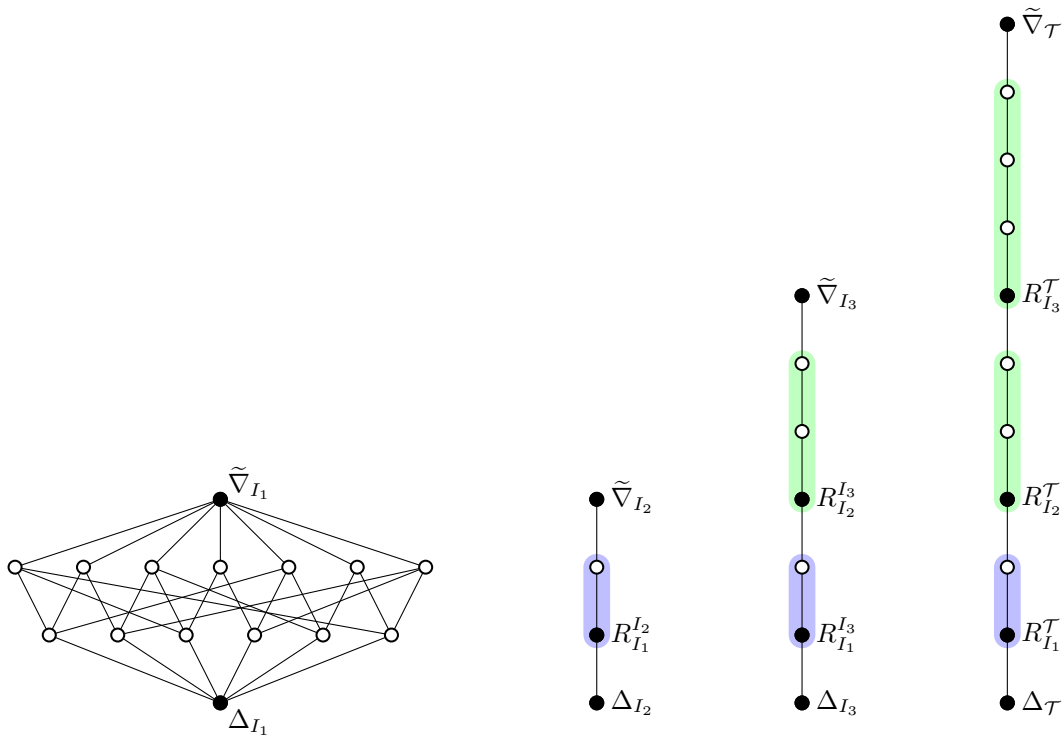

\begin{rem}\label{rem:T2}
Consider the case in which $\C=\{A\}$, where $A=\{1,2,3,4\}$, so that $\T(\C)=\T_4$ is the full transformation monoid of degree $4$.  The lattices $\Cong(I_r)$, $r=1,2,3,4$, are pictured in Figure \ref{fig:T}.  Note that $\Cong(I_1)=\Eq(I_1)$ is isomorphic to the lattice $\Eq_4$.
\end{rem}

\begin{rem}\label{rem:T3}
If we had two non-empty sets $\C_1$ and $\C_2$ of non-empty finite sets, then the ideals $I_r(\T(\C_1))$ and $I_r(\T(\C_2))$ with $r\geq2$ have isomorphic congruence lattices; when $r$ is finite, these are in fact isomorphic to the congruence lattice of the full transformation semigroup $\T_r$ \cite{Malcev1952,GMbook}.  Such an isomorphism would only exist for $r=1$ if there was a cardinality-preserving bijection $\C_1\to\C_2$.
\end{rem}

\subsection{Planar and annular reducts of $\T$}\label{subsect:QT}

We now consider the subcategories of $\T=\T(\C)$ consisting of all transformations preserving (or possibly reversing) order or orientation:
\[
\OrdP\T \COMMA \OrdPR\T \COMMA \OriP\T \COMMA \OriPR\T.
\]
Let $\Q$ denote any of the above categories.  First, it is easy to check that $\Q$ is regular; cf.~\cite[Sections 3 and~5]{CH1999} and \cite[Section 2]{GH1992}.  It follows from \cite[Lemma 2.9]{Sandwich1} that $\Q$ is stable.  It follows from regularity (and \cite[Lemma 2.8]{Sandwich1}) that Green's $\R$, $\L$ and $\H$ relations on $\Q$ are just the restrictions of the corresponding relations on $\T$; this is easily checked to be the case for the $\D=\J$ relation as well (again, cf.~\cite{CH1999,GH1992}).  Thus, parts \ref{GT1}--\ref{GT3} of Lemma \ref{lem:Green_T} remain true in $\Q$, so we have the (stable, regular) $\D=\J$-classes $D_q=D_q(\Q)$, $1\leq q<\xi$, and ideals $I_r=I_r(\Q)$, $1\leq r<\xi$, and also $I_\om=\Q$ if $\xi=\om$.  Part \ref{GT4} of Lemma \ref{lem:Green_T} is not true for the subcategory $\Q$, however; rather, the typical group $\H$-class $\G_q$ in $D_q$ is (with some trivial exceptions for small~$q$):
\bit
\item trivial when $\Q=\OrdP\T$,
\item cyclic of order $2$ when $\Q=\OrdPR\T$,
\item cyclic of order $q$ when $\Q=\OriP\T$,
\item dihedral of order $2q$ when $\Q=\OriPR\T$.
\eit

There is a version of Lemma \ref{lem:P2T} for $\Q$, which shows (among other things) that each ideal $I_r=I_r(\Q)$, $2\leq r<\xi$, has property $\Mult{I_r}$ and hence $\Sep{I_r}$.  The original proof of Lemma~\ref{lem:P2T} does not need to be modified at all, since every element whose existence is asserted belongs to~$\OrdP{\T}$, and hence to $\Q$.

Here is the classification of congruences on the ideals of the above categories:

\begin{thm}\label{thm:QT}
Let $\C$ be a non-empty set of finite non-empty sets, let $\T=\T(\C)$ be the transformation category over $\C$, let $\Q$ be one of $\OrdP\T$, $\OrdPR\T$, $\OriP\T$ or $\OriPR\T$, and keep the above notation.
\ben
\item \label{QT1} The minimal ideal $I_1$ is a partial rectangular band, so its congruence lattice is described by Proposition \ref{prop:PRB}.  In particular,
\[
\Cong(I_1)\cong \prod_{A\in\C}\Eq(A).
\]
\item \label{QT2} For $r\geq2$ (including $r=\om$ if $\xi=\om$),
\[
\Cong(I_r) = \{\De_{I_r}\} \cup \set{R_{I_q,N}^{I_r}}{1\leq q<r,\ N\normal\G_{q+1}}\cup\{\nab_{I_r}\}.
\]
\een
\end{thm}

\pf
\firstpfitem{\ref{QT1}}  This follows from Theorem \ref{thm:T}\ref{T1}, since $I_1(\Q)=I_1(\T)$.

\pfitem{\ref{QT2}}  The proof is essentially identical to that of Theorem \ref{thm:T}\ref{T2}.  (The only difference is that when $\Q=\OrdP\T$, we have $\Cong(I_2) = \{\De_{I_2},R_{I_1}^{I_2},\nab_{I_2}\}$, as $\OrdP\T$ is $\H$-trivial.)  
\epf

\begin{rem}\label{rem:QT}
For $r\geq2$, the lattice $\Cong(I_r(\Q))$ is best thought of as consisting of copies of the lattices $\Cong(I_2), \NN(\G_3), \NN(\G_4),\ldots$ stacked on top of each other.
Thus, it will always be a chain when~$\Q$ is $\OrdP\T$ or $\OrdPR\T$, where the groups are trivial or cyclic of order $2$;
but it will not in general be a chain when $\Q$ is  $\OriP\T$ or $\OriPR\T$, as the groups are cyclic of larger orders or dihedral.
 Figure~\ref{fig:QT} displays the lattice $\Cong(I_6)$ in the categories $\T$, $\OrdP\T$, $\OrdPR\T$, $\OriP\T$ and $\OriPR\T$.  The right-most lattice shown appears in \cite[p.~756]{FGJ2009}, although there it represents the congruence lattice of the monoid of orientation-preserving or -reversing \emph{partial} transformations of a set of size $6$.
\end{rem}

\begin{figure}[ht]
\begin{center}
\begin{tikzpicture}[scale=0.7,inner sep=0pt]
\begin{scope}[shift={(-4,0)}]
\draw[rounded corners=4,line width=10,blue!25,cap=round](0,1)--(0,2);
\draw[rounded corners=4,line width=10,green!25,cap=round]
(0,3)--(0,5)
(0,6)--(0,9)
(0,10)--(0,12)
(0,13)--(0,15)
;
\draw(0,0)--(0,16);
\foreach \x in {0,...,16} {\colvert0{\x}{white}}
\foreach \x in {0,1,3,6,10,13,16} {\colvert0{\x}{black}}
\end{scope}
\begin{scope}[shift={(0,0)}]
\draw[rounded corners=4,line width=10,blue!25,cap=round](0,1)--(0,1.001);
\draw[rounded corners=4,line width=10,green!25,cap=round]
(0,2)--(0,2.001)
(0,3)--(0,3.001)
(0,4)--(0,4.001)
(0,5)--(0,5.001)
;
\draw(0,0)--(0,6);
\foreach \x in {0,...,6} {\colvert0{\x}{black}}
\end{scope}
\begin{scope}[shift={(4,0)}]
\draw[rounded corners=4,line width=10,blue!25,cap=round](0,1)--(0,2);
\draw[rounded corners=4,line width=10,green!25,cap=round]
(0,3)--(0,4)
(0,5)--(0,6)
(0,7)--(0,8)
(0,9)--(0,10)
;
\draw(0,0)--(0,11);
\foreach \x in {0,...,11} {\colvert0{\x}{white}}
\foreach \x in {0,1,3,5,7,9,11} {\colvert0{\x}{black}}
\end{scope}
\begin{scope}[shift={(8,0)}]
\draw[rounded corners=4,line width=10,blue!25,cap=round](0,1)--(0,2);
\draw[rounded corners=4,line width=10,green!25,cap=round]
(0,3)--(0,4)
(0,5)--(0,7)
(0,8)--(0,9)
(0,10)--(-1,11) (-1,11)--(0,12) (0,12)--(1,11) (1,11)--(0,10)
;
\fill[green!25] 
(0,10)--(-1,11)--(0,12)--(1,11)--(0,10)
;
\draw(0,0)--(0,10)--(-1,11)--(0,12)--(0,13) (0,12)--(1,11)--(0,10);
\foreach \x in {0,...,10,12,13} {\colvert0{\x}{white}}
\foreach \x/\y in {-1/11,1/11} {\colvert{\x}{\y}{white}}
\foreach \x in {0,1,3,5,8,10,13} {\colvert0{\x}{black}}
\end{scope}
\begin{scope}[shift={(12,0)}]
\draw[rounded corners=4,line width=10,blue!25,cap=round](0,1)--(0,2);
\draw[rounded corners=4,line width=10,green!25,cap=round]
(0,3)--(0,5)
(0,6)--(0,7) (0,7)--(-1,8) (-1,8)--(0,9) (0,9)--(1,8) (1,8)--(0,7)
(0,10)--(0,12)
(0,13)--(-2,15) (-2,15)--(-1,16) (-1,16)--(1,14) (1,14)--(0,13)
;
\fill[green!25] 
(0,6)--(0,7)--(-1,8)--(0,9)--(1,8)--(0,7)
(0,13)--(-2,15)--(-1,16)--(1,14)--(0,13)
;
\draw(0,0)--(0,13)--(-2,15)--(-1,16)--(-1,17) (-1,14)--(-1,16)--(1,14)--(0,13) (-1,14)--(0,15) (0,7)--(-1,8)--(0,9)--(1,8)--(0,7);
\foreach \x in {0,...,13,15} {\colvert0{\x}{white}}
\foreach \x/\y in {-1/8,1/8,-1/14,-1/15,-1/16,-1/17,1/14,-2/15} {\colvert{\x}{\y}{white}}
\foreach \x in {0,1,3,6,10,13} {\colvert0{\x}{black}}
\foreach \x in {17} {\colvert{-1}{\x}{black}}
\end{scope}
\end{tikzpicture}
\caption{Left to right: Hasse diagrams of $\Cong(I_6)$ in the categories $\T$, $\OrdP\T$, $\OrdPR\T$, $\OriP\T$ and~$\OriPR\T$; cf.~Theorem \ref{thm:QT} and Remark \ref{rem:QT}. }
\label{fig:QT}
\end{center}
\end{figure}

\subsection{Other categories of transformations}\label{subsect:OtherT}

The methods of this section allow one to readily obtain analogous results for a number of other categories.  Specifically, let $\C$ be a non-empty set of finite sets (now $\C$ could include the empty set).  We denote by $\PT=\PT(\C)$ the category of all partial transformations between members of $\C$, and by $\I=\I(\C)$ the category of all injective partial transformations.  We also have planar and annular reducts such as $\OrdP{\PT}$, $\OriPR\I$, and so on.  Green's relations and ideals on all of these categories are also described in terms of domain, kernel, image and rank parameters, and all of the categories are stable and regular; cf.~\cite[Sections 3 and 5]{Sandwich1}.  The congruences on the ideals of these categories have the same basic kind of form, except for some minor differences, and one major simplification.  The main point is that the minimal ideal $I_0=D_0$ consists entirely of empty mappings (one for each pair of sets from $\C$), and therefore has only one congruence, $\De_{I_0}=\nab_{I_0}$.  We also have $\Cong(I_1)=\{\De_{I_1},\nab_{I_1}\}$; indeed, this follows from Proposition \ref{prop:small_ideal} and the fact that $D_1$ is $\H$-trivial.    We then apply Theorem \ref{thm:E2} for the larger ideals.
The details are routine and are omitted.

The category $\I$ can also be thought of as a ``rook category'': i.e., a category of $\{0,1\}$-matrices in which each row and each column has at most one non-zero entry; cf.~\cite{Solomon2002,FHH2009,BRY2012,JP2015}.

\section{Partition categories}\label{sect:P}

Sections \ref{sect:P}--\ref{sect:J} concern several diagram categories; for more details on such categories (and in particular their linear versions), see \cite{Martin2008}, which also gives an excellent overview of their origins in theoretical physics and representation theory.  We begin in this section with the largest of these categories, the so-called \emph{partition category} $\P$.  
It contains all the other diagram categories we will consider as subcategories.

Our goal in this section is to describe the congruences on $\P$ and its ideals.  Although there are some similarities with the transformation category $\T$ considered in Section \ref{sect:T}, there are also some significant differences.  
Like $\T$, the partition category $\P$ is a chain of ideals $I_0\subset I_1\subset\cdots$, and congruences of the form 
$R_{I,N}$ will play a key role.  However, the minimal ideal $I_0$ of $\P$ is neither $\wt\R$- nor~$\wt\L$-trivial, so there are non-trivial $\lam/\rho$ congruences (as defined in Section \ref{subsect:EMRT}), and hence ``diamonds'' in the lower parts of the congruence lattices; cf.~Figures \ref{fig:P} and \ref{fig:QP}.  Another major difference is that the ideal~$I_1$ of $\P$ is retractable (onto $I_0$).
This means that our stacking mechanism from Section \ref{subsect:chains} kicks in only from $I_2$ upwards.
Congruences of $I_0$ and $I_1$, of which there are many, will be handled by the theory from Sections \ref{subsect:PRB} and \ref{subsect:retract} respectively.

This section is structured in much the same way as Section \ref{sect:T}.  We begin in Section \ref{subsect:defnP} with the relevant definitions and notation.  We describe Green's relations in Section \ref{subsect:prelimP}, where we also prove a technical result (Lemma \ref{lem:P2P}) that will be used to establish the required multiplication/separation properties in ideals of $\P$. 
 In Section \ref{subsect:P} we apply the theory from Section \ref{sect:chains} and obtain a classification of the congruences on~$\P$ and its ideals, and describe the associated lattices; see Theorem \ref{thm:P} and Figure \ref{fig:P}.  We give the corresponding results for planar and annular reducts of $\P$ in Section \ref{subsect:QP}, and for various other subcategories such as the partial Brauer and Motzkin categories in Section \ref{subsect:PB}.  A number of further diagram categories will be treated in Sections \ref{sect:B}--\ref{sect:J}, including the Brauer, Temperley--Lieb and Jones categories.

\subsection{Definitions and preliminaries on $\P$}\label{subsect:defnP}

For the duration of Section \ref{sect:P}, we fix a non-empty set $\C$ of finite sets, possibly including the empty set.  Again, we write $\xi=\xi(\C)={\LSUB}\bigset{|A|}{A\in\C}$.  
In what follows, we will write $X=\bigcup_{A\in\C}A$ for the (possibly infinite) set of all points belonging to any of the members of $\C$.  We will need two additional disjoint copies $X'=\set{x'}{x\in X}$ and $X''=\set{x''}{x\in X}$ of this set.  For any $A\sub X$, we write $A'=\set{a'}{a\in A}$ and $A''=\set{a''}{a\in A}$.

For $A,B\in\C$, we write $\P_{A,B}$ for the set of all set partitions of $A\cup B'$.  Any such partition $\al\in\P_{A,B}$ will be called a \emph{partition from $A$ to $B$}, and we write $\bd(\al)=A$ and $\br(\al)=B$.  We define
\[
\P = \P(\C) = \bigcup_{A,B\in\C}\P_{A,B}
\]
to be the set of all partitions between members of $\C$.

A partial operation $\cdot$ is defined on $\P$ as follows.  Suppose $A,B,C\in\C$, and that $\al\in\P_{A,B}$ and $\be\in\P_{B,C}$.  First, we identify $\al$ and $\be$ with their corresponding equivalence relations on the sets $A\cup B'$ and $B\cup C'$.  Let $\al_\downarrow$ be the equivalence on $A\cup B''\cup C'$ obtained by renaming each element $b'\in B'$ by $b''$ and adding the diagonal relation on $C'$.  Similarly, let $\be^\uparrow$ be the equivalence on $A\cup B''\cup C'$ obtained by renaming each element $b\in B$ by $b''$ and adding the diagonal relation on $A$.  We now form the join $\al_\downarrow\vee\be^\uparrow$, which is again an equivalence on $A\cup B''\cup C'$.  The product $\al\be=\al\cdot\be\in\P_{A,C}$ is then defined to be the partition corresponding to the restriction $(\al_\downarrow\vee\be^\uparrow)\restr_{A\cup C'}$.  This partial product is associative, and 
\[
(\P,\C,\bd,\br,\cdot)
\]
is a partial semigroup, called the \emph{partition category (over $\C$)}.  
For $A\in\C$, the endomorphism monoid $\P_A=\P_{A,A}$ is the partition monoid over $A$ \cite{HR2005}, with identity $\bigset{\{a,a'\}}{a\in A}$.  In particular, $\P$ contains a copy of the symmetric group $\S_A$ for each $A\in\C$; as usual, we identify~$\S_A$ with the set of all partitions from $\P_A$ whose blocks are all of the form $\{a,b'\}$ for some $a,b\in A$.

In order to offer a
diagrammatic interpretation of this product we recall some terminology and notation.  A non-empty subset $A\sub X\cup X'$ is called
\bit
\item a \emph{transversal} if $A\cap X$ and $A\cap X'$ are both non-empty,
\item an \emph{upper non-transversal} if $A\sub X$,
\item a \emph{lower non-transversal} if $A\sub X'$.
\eit
If $\al\in\P$, we will write
\[
\al = \partI{A_1}{A_q}{C_1}{C_r}{B_1}{B_q}{D_1}{D_s}
\]
to indicate that $\al$ has the following blocks:  transversals $A_i\cup B_i'$ ($i=1,\ldots,q$), upper non-transversals $C_i$ ($i=1,\ldots,r$) and lower non-transversals $D_i'$ ($i=1,\ldots,s$).  From this notation, we may infer
\[
\bd(\al) = \bigcup_{i=1}^qA_i \cup \bigcup_{i=1}^rC_i \AND \br(\al) = \bigcup_{i=1}^qB_i \cup \bigcup_{i=1}^sD_i.
\]
Note that any of $q,r,s$ could be zero.

It will be convenient to sometimes use a variation of the above notation.  Specifically, we will often wish to omit some (but not necessarily all) singleton blocks from $\al\in\P$, and will indicate this by using square brackets:
\[
\al= \sqpartI{A_1}{A_r}{C_1}{C_s}{B_1}{B_r}{D_1}{D_t}.
\]
When using this notation, it will be important to specify which hom-set $\al$ belongs to: i.e., to specify $\bd(\al)=A$ and $\br(\al)=B$, in which case any element of $A$ or $B$ not listed in this square-bracket notation corresponds to a singleton block of $\al$.  

A partition $\al\in\P_{A,B}$ may also be identified with any graph on vertex set $A\cup B'$ whose connected components are the blocks of $\al$.  This leads to an alternative way to understand the product in $\P$.  Let $\al\in\P_{A,B}$ and $\be\in\P_{B,C}$, where $A,B,C\in\C$.  Take the graphs on vertex set $A\cup B''\cup C'$ corresponding to the partitions $\al_\downarrow$ and $\be^\uparrow$ defined above, and let $\Pi(\al,\be)$ be the graph on vertex set $A\cup B''\cup C'$ whose edge set is the union of the edge sets coming from~$\al_\downarrow$ and~$\be^\uparrow$.  This graph is called the \emph{product graph} of the pair $(\al,\be)$; note that $\Pi(\al,\be)$ is only defined for pairs with $\br(\al)=\bd(\be)$.  The product $\al\be\in\P_{A,C}$ is then the partition of $A\cup C'$ for which elements $x,y\in A\cup C'$ belong to the same block of $\al\be$ if and only if $x,y$ belong to the same component of $\Pi(\al,\be)$.  

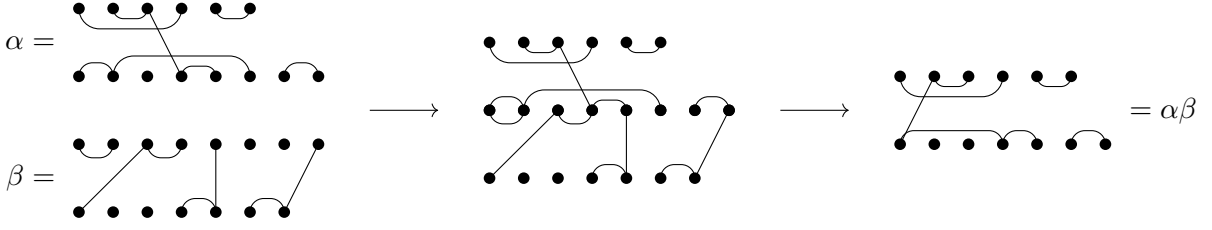
\begin{figure}[ht]
\begin{center}
\begin{tikzpicture}[scale=.45]
\begin{scope}[shift={(0,0)}]	
\uvs{1,...,6}
\lvs{1,...,8}
\uarcx14{.6}
\uarcx23{.3}
\uarcx56{.3}
\darc12
\darcx26{.6}
\darcx45{.3}
\darc78
\stline34
\draw(0.6,1)node[left]{$\al=$};
\draw[->](9.5,-1)--(11.5,-1);
\end{scope}
\begin{scope}[shift={(0,-4)}]	
\uvs{1,...,8}
\lvs{1,...,7}
\uarc12
\uarc34
\darc45
\darc67
\stline31
\stline55
\stline87
\draw(0.6,1)node[left]{$\be=$};
\end{scope}
\begin{scope}[shift={(12,-1)}]	
\uvs{1,...,6}
\lvs{1,...,8}
\uarcx14{.6}
\uarcx23{.3}
\uarcx56{.3}
\darc12
\darcx26{.6}
\darcx45{.3}
\darc78
\stline34
\draw[->](9.5,0)--(11.5,0);
\end{scope}
\begin{scope}[shift={(12,-3)}]	
\uvs{1,...,8}
\lvs{1,...,7}
\uarc12
\uarc34
\darc45
\darc67
\stline31
\stline55
\stline87
\end{scope}
\begin{scope}[shift={(24,-2)}]	
\uvs{1,...,6}
\lvs{1,...,7}
\uarcx14{.6}
\uarcx23{.3}
\uarcx56{.3}
\darc14
\darc45
\darc67
\stline21
\draw(7.4,1)node[right]{$=\al\be$};
\end{scope}
\end{tikzpicture}
\caption{Partitions $\al\in\P_{A,B}$ (top left), $\be\in\P_{B,C}$ (bottom left), the product graph $\Pi(\al,\be)$ (middle), and the product $\al\be\in\P_{A,C}$ (right).  Here $A=[6]$, $B=[8]$ and $C=[7]$, and vertices are drawn in the usual order, increasing from left to right.}
\label{fig:P_product}
\end{center}
\end{figure}

When depicting a graph representing a partition $\al\in\P_{A,B}$, we always draw the vertices from~$A$ on an upper row, and the vertices from $B'$ on a lower row; this explains the use of the terms \emph{upper} and \emph{lower} non-transversals.  When drawing a product graph, the double-dashed vertices are drawn on a middle row.  For example, let
\[
A = [6] = \{1,\ldots,6\} \COMMA
B = [8] = \{1,\ldots,8\} \COMMA
C = [7] = \{1,\ldots,7\} ,
\]
and
\begin{align*}
\al &= \big\{ \{1,4\},\{2,3,4',5'\},\{5,6\},\{1',2',6'\},\{3'\},\{7',8'\}\big\} \in \P_{A,B} ,\\
\be &= \big\{ \{1,2\}, \{3,4,1'\}, \{5,4',5'\}, \{6\}, \{7\}, \{8,6',7'\}, \{2'\}, \{3'\} \big\} \in \P_{B,C}.
\intertext{Figure \ref{fig:P_product} shows how to calculate the product}
\al\be &= \big\{ \{1,4\},\{2,3,1',4',5'\},\{5,6\},\{2'\},\{3'\},\{6',7'\}\big\} \in \P_{A,C} .
\end{align*}

\subsection{Green's relations and multiplicative properties in $\P$}\label{subsect:prelimP}

Let $\al\in\P_{A,B}$, and write $\al = \partI{A_1}{A_q}{C_1}{C_r}{B_1}{B_q}{D_1}{D_s}$.  The \emph{domain}, \emph{codomain}, \emph{kernel}, \emph{cokernel} and \emph{rank} of $\al$ are defined by
\bit
\item $\dom(\al) = \set{x\in A}{\text{$x$ belongs to a transversal of $\al$}} = \bigcup_{i=1}^qA_i$,
\item $\codom(\al) = \set{x\in B}{\text{$x'$ belongs to a transversal of $\al$}} = \bigcup_{i=1}^qB_i$,
\item $\ker(\al) = \bigset{(x,y)\in A\times A}{\text{$x$ and $y$ belong to the same block of $\al$}}$, the equivalence relation on $A$ 
associated with the partition $\set{A_i}{1\leq i\leq q}\cup\set{C_i}{1\leq i\leq r}$,
\item $\coker(\al) = \bigset{(x,y)\in B\times B}{\text{$x'$ and $y'$ belong to the same block of $\al$}}$, the equivalence relation on $B$ associated with the partition $\set{B_i}{1\leq i\leq q}\cup\set{D_i}{1\leq i\leq s}$,
\item $\rank(\al)=q$, the number of transversals of $\al$.
\eit
Note that $\rank(\al)$ can be anything from $0$ to $\min(|A|,|B|)$.  The next result is proved in the same fashion as \cite[Theorem 17]{Wilcox2007}, \cite[Theorem 3.3]{FL2011} and \cite[Lemma 4.3]{ER2018}, each of which concern partition \emph{monoids}.

\begin{lemma}\label{lem:Green_P}
\ben
\item \label{GP1} For $\al,\be\in\P$, we have
\bit
\item $(\al,\be)\in\R \iff \dom(\al)=\dom(\be)$ and $\ker(\al)=\ker(\be)$,
\item $(\al,\be)\in\L \iff \codom(\al)=\codom(\be)$ and $\coker(\al)=\coker(\be)$,
\item $(\al,\be)\in\J\iff(\al,\be)\in\D\iff\rank(\al)=\rank(\be)$.
\eit
\item \label{GP2} The $\J=\D$-classes of $\P$ are the sets
\bit
\item $D_q = D_q(\P) = \set{\al\in\P}{\rank(\al)=q}$ for each $0\leq q<\xi$.
\eit
These are all regular and stable, and they form a chain: $D_0<D_1<\cdots$.
\item \label{GP3} The ideals of $\P$ are the sets
\bit
\item $I_r = I_r(\P) = \set{\al\in\P}{\rank(\al)\leq r}$ for each $0\leq r<\xi$, and
\item $I_\om = I_\om(\P) = \P $ in the case that $\xi=\om$.
\eit
These are all regular and stable, and they form a chain: $I_0\subset I_1\subset\cdots$.
\item \label{GP4} The $\H$-class of any idempotent from $D_q$ is isomorphic to $\S_q$, the symmetric group of degree~$q$.  \epfres
\een
\end{lemma}

Thus, the category $\P$ is a chain of ideals, as in Definition \ref{defn:chain}.
In fact, $\P$ satisfies a property stronger than regularity, and analogous to the idea of a \emph{regular $*$-semigroup} \cite{NS1978}.  Specifically, for a partition $\al = \partI{A_1}{A_q}{C_1}{C_r}{B_1}{B_q}{D_1}{D_s}\in\P$, we define $\al^* = \partI{B_1}{B_q}{D_1}{D_s}{A_1}{A_q}{C_1}{C_r}$.  Then $\bd(\al^*)=\br(\al)$, $\br(\al^*)=\bd(\al)$ and $(\al^*)^*=\al=\al\al^*\al$, and if $\al,\be\in\P$ are composable, then $(\al\be)^*=\be^*\al^*$.  Accordingly, we call $\P$ a \emph{regular $*$-category}; cf.~\cite[Section 2]{DDE2019}.  We also have ${(\al,\be)\in\R\iff(\al^*,\be^*)\in\L}$, and other such identities.  This symmetry/duality will allow us to shorten many proofs.

By Lemmas \ref{lem:SD} and \ref{lem:Green_P}, the $\D$-classes and ideals of an ideal $I_r=I_r(\P)$ are precisely the $D_q$ and $I_q$ with $q\leq r$.

Here is the main technical result of this section; 
it will be used in checking the requisite multiplication and separation properties.

\begin{lemma}\label{lem:P2P}
Let $A,B\in\C$, and suppose $\al,\be\in\P_{A,B}$ with $r=\rank(\al)\geq\rank(\be)=q$.
\ben
\item \label{P2Pi} If $r\geq2$ and $q<r$, then there exists $\ga\in I_r$ such that $\ga\al\in D_{r-1}$ and $\ga\be\in I_{r-1}\sm H_{\ga\al}$.
\item \label{P2Pii} If $q=r\geq1$ and $(\al,\be)\not\in\H$, then there exists $\ga\in I_r$ such that, swapping $\al,\be$ if necessary,
\[
\text{$[\al\ga\in D_r$ and $\be\ga\in I_{r-1}]$ \ \ \ or \ \ \ $[\ga\al\in D_r$ and $\ga\be\in I_{r-1}]$.}
\]
\item \label{P2Piii} If $r\geq2$ and $\be\in H_\al\sm\{\al\}$, then there exists $\ga\in I_r$ such that $\ga\al\in D_{r-1}$ and $\ga\be\in I_{r-1}\sm H_{\ga\al}$.
\een
\end{lemma}

\pf
Throughout the proof, we write $\al = \partI{A_1}{A_r}{C_1}{C_s}{B_1}{B_r}{D_1}{D_t}$.  For each $1\leq i\leq r$, we fix some $a_i\in A_i$.

\pfitem{\ref{P2Pi}}  If $\dom(\al)\not\sub\dom(\be)$, then we may assume without loss of generality that $a_1\not\in\dom(\be)$.  If $\dom(\al)\sub\dom(\be)$, then by the pigeon-hole principle we may assume without loss of generality that $(a_1,a_2)\in\ker(\be)$.  In either case, we take $\ga=\sqpartpermIII{a_1}\cdots{a_{r-1}}{a_1}\cdots{a_{r-1}}\in\P_A$.  Then $\dom(\ga\al)=\{a_1,\ldots,a_{r-1}\}$ and $\ker(\ga\al)=\De_A$, so that $\ga\al\in D_{r-1}$.  Since $\rank(\ga\be)\leq\rank(\ga)=r-1$ we have $\ga\be\in I_{r-1}$.  In the first case we have $a_1\not\in\dom(\ga\be)$, and in the second case we have $(a_1,a_2)\in\ker(\ga\be)$.  Thus, $(\ga\al,\ga\be)\not\in\R$.

\pfitem{\ref{P2Pii}}  We assume that $(\al,\be)\not\in\R$, the case of $(\al,\be)\not\in\L$ being dual.  So either $\dom(\al)\not=\dom(\be)$ or $\ker(\al)\not=\ker(\be)$.

\pfcase{1}  Suppose first that $\dom(\al)\not=\dom(\be)$.  Without loss of generality, we may assume that $a_1\not\in\dom(\be)$.  We than take $\ga=\sqpartpermIII{a_1}\cdots{a_r}{a_1}\cdots{a_r}\in\P_A$.  Clearly $\ga\al\in D_r$.  On the other hand, $\dom(\ga\be)\sub\dom(\ga)=\{a_1,\ldots,a_r\}$ but $a_1\not\in\dom(\ga\be)$, which gives $\ga\be\in I_{r-1}$.

\pfcase{2}  Now suppose $\dom(\al)=\dom(\be)$ but $\ker(\al)\not=\ker(\be)$.  Without loss of generality, we may assume there exists $(x_1,x_2)\in\ker(\be)\sm\ker(\al)$.  Note then that $x_1$ and $x_2$ either both belong to $\dom(\be)=\dom(\al)$ or else both belong to $A\sm\dom(\al)$.

\pfcase{2.1}  Suppose first that $x_1,x_2\in\dom(\al)$.  Here we may assume that $x_1=a_1$ and $x_2=a_2$, and again we take $\ga=\sqpartpermIII{a_1}\cdots{a_r}{a_1}\cdots{a_r}\in\P_A$.  This time $\ga\al\in D_r$, while $\dom(\ga\be)\sub\{a_1,\ldots,a_r\}$ and $(a_1,a_2)\in\ker(\ga\be)$ together give $\ga\be\in I_{r-1}$.

\pfcase{2.2}  Now suppose $x_1,x_2\not\in\dom(\al)$.  We may also assume that $A_1,\ldots,A_r$ are the upper parts of the transversals of $\be$ (otherwise we would be in Case 2.1).  Here we take $\ga=\sqpartII{a_1}{a_{r-1}}{a_r}{}{a_1}{a_{r-1}}{x_1}{a_r,x_2}$.  This time we have $\ga\be\in D_r$ and $\ga\al\in I_{r-1}$; the latter is the case because $\dom(\ga\al)\sub\{a_1,\ldots,a_r\}$ yet $a_r\not\in\dom(\ga\al)$.

\pfitem{\ref{P2Piii}}  Here we have $\be = \partI{A_1}{A_r}{C_1}{C_s}{B_{1\pi}}{B_{r\pi}}{D_1}{D_t}$ for some permutation $\pi\in\S_r$, and without loss of generality we may assume that $1\pi=r$.  We then take $\ga=\sqpartpermIII{a_1}\cdots{a_{r-1}}{a_1}\cdots{a_{r-1}}\in\P_A$, noting that $\ga\al\in D_{r-1}$, $\ga\be\in I_{r-1}$ and $B_r\sub\codom(\ga\be)\sm\codom(\ga\al)$; the latter gives $(\ga\al,\ga\be)\not\in\L$.
\epf

\begin{rem}
It follows from Lemma \ref{lem:P2P} that the ideals $I_r=I_r(\P)$ have property $\Mult{I_r}$, and hence $\Sep{I_r}$, for $r\geq2$.  The fact that part \ref{P2Pii} of the lemma holds for $r\geq1$ (not just $r\geq2$) will be used to verify the technical assumption of Proposition \ref{prop:small_ideal} when treating the ideal $I_1$.
\end{rem}

\subsection{Congruences on ideals of $\P$}\label{subsect:P}

In what follows, for each $0\leq q<\xi$, it will be convenient to fix a specific group $\H$-class $\G_q$ in the (stable, regular) $\D$-class $D_q$.  To this end, as we did for the transformation category $\T$, for each such $q$ we fix some $X_q\in\C$ with $n_q=|X_q|\geq q$, and assume without loss of generality that $X_q=[n_q]=\{1,\ldots,n_q\}$.  Within $\P_{X_q}=\P_{X_q,X_q}$, we define
\[
\G_q = \bigset{\al\in\P_{X_q}}{\dom(\al)=\codom(\al)=[q],\ \ker(\al)=\coker(\al)=\De_{X_q}}.
\]
It is clear that $\G_q$ is isomorphic to the symmetric group $\S_q$ (cf.~Lemma \ref{lem:Green_P}\ref{GP4}).  
Every element~$\al$ of $\G_q$ has the form
\[
\al = \sqpartpermIII1\cdots q{1\pi}\cdots{q\pi} = \partI1q{q+1}{n_q}{1\pi}{q\pi}{q+1}{n_q} \qquad\text{for some permutation $\pi\in\S_q$,}
\]
and we write $\al=\pi^\natural$.  For any $\Om\sub\S_q$, we write $\Om^\natural=\set{\pi^\natural}{\pi\in\Om}$, so that in fact $\G_q=\S_q^\natural$.  Moreover, every normal subgroup of $\G_q$ is of the form $N^\natural$ for some $N\normal\S_q$.   It follows that an ideal $I_r$ ($r\geq2$) contains the IN-pairs $(I_q,N^\natural)$, for $0\leq q<r$ and $N\normal\S_{q+1}$.  These yield the congruences $R_{I_q,N^\natural}^{I_r}=R_{I_q}^{I_r}\cup\nu_{N^\natural}$, including the Rees congruence $R_{I_q}^{I_r}$ when $N$ is trivial.  In what follows, we will often use the abbreviation $R_N^{I_r}=R_{I_q,N^\natural}^{I_r}$, based on the fact that both $q$ and $N^\natural$ are ``encoded'' in $N$.

Since these congruences (and others to follow) involve the $\nu_{N^\natural}$ relations, we briefly give a concrete description of these relations.  Specifically, suppose $\al,\be\in\P$ are such that $(\al,\be)\in\H$.  Then we have
\[
\al = \partI{A_1}{A_q}{C_1}{C_r}{B_1}{B_q}{D_1}{D_s} \AND \be = \partI{A_1}{A_q}{C_1}{C_r}{B_{1\phi}}{B_{q\phi}}{D_1}{D_s}
\]
for some permutation $\phi\in\S_q$, and we define $\phi=\phi(\al,\be)$.  As in Section \ref{subsect:T}, note that $\phi$ depends on the ordering chosen on the transversals of $\al$, so $\phi(\al,\be)$ is only well-defined up to conjugacy in $\S_q$.  However, for any $N\normal\S_q$, we have (cf.~\cite[Lemma 5.6]{EMRT2018})
\[
\nu_{N^\natural} = \bigset{(\al,\be)\in\H\restr_{D_q}}{\phi(\al,\be)\in N}.
\]

As in \cite[Lemma 5.2]{EMRT2018}, the ideal $I_1=D_0\cup D_1$ is retractable.  
The retraction is given by
\begin{equation}\label{eq:retract_P}
f:I_1\to I_0: \alpha\mt\wh\al\qquad\text{where}\qquad
\wh\al=\begin{cases}
\partIII{A}{C_1}\cdots{C_r}{B}{D_1}\cdots{D_s} & \text{if }\al=\partIV{A}{C_1}\cdots{C_r}{B}{D_1}\cdots{D_s}\in D_1\\
\alpha &\text{if } \alpha\in D_0.
\end{cases}
\end{equation}

As in \cite[Lemma 5.3]{EMRT2018}, we have the following three retractable IN-pairs in $I_r$ ($r\geq2$):
\[
(I_0,\S_1^\natural) \COMMA
(I_1,\{\id_2^\natural\}) \COMMA
(I_1,\S_2^\natural).
\]
Each of the above retractable IN-pairs yields a family of $R/\lam/\rho/\mu$ congruences on $I_r$;
since $D_0$ is $\H$-trivial, the $\eta$ congruences coincide with the $\mu$ congruences, so we have no need to distinguish these, and will use the symbol $\mu$. 
Again, we will abbreviate the notation for these, by omitting the symbol $\natural$, and also either omiting the normal subgroup altogether if it is trivial, or else omitting the ideal $I_0/I_1$, as it can be deduced from the group.
Thus, for example, $\lam_{I_0}^{I_r}$ will stand for $\lam_{I_0,\S_1^\natural}^{I_r}$,
while $ \mu_{\S_2}^{I_r}$ will stand for $\mu_{I_1,\S_2^\natural}^{I_r}$.

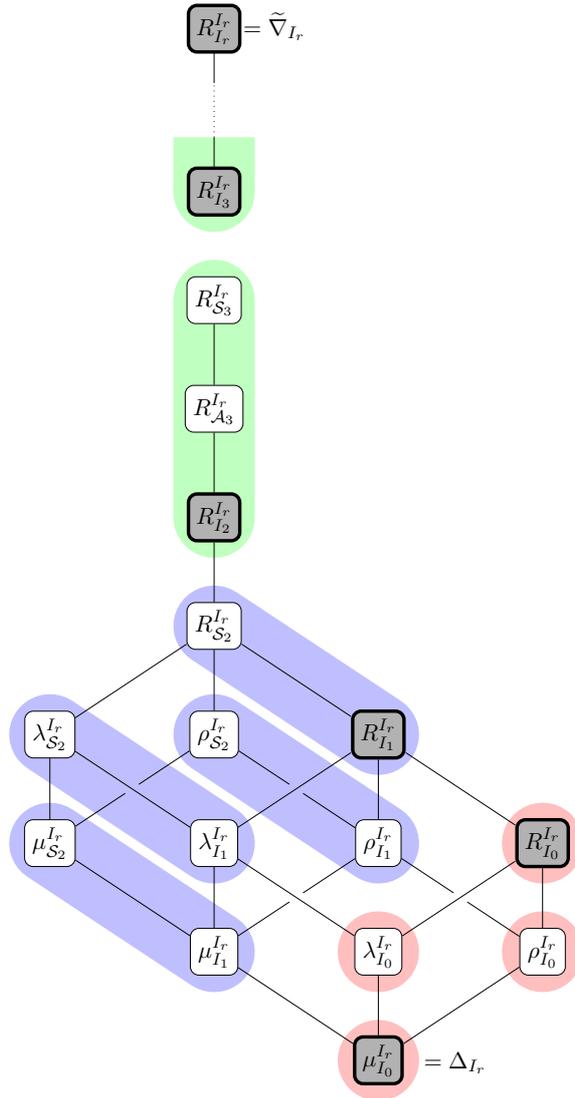
\begin{figure}[ht]
\begin{center}
\scalebox{0.8}{
\begin{tikzpicture}[scale=.9]
\foreach \x/\y in {0/0,0/2,3/2,3/4} {\draw[rounded corners=2em,line width=3.5em,red!25,cap=round] (3+\x,0+\y)--(3+\x,0.001+\y);}
\foreach \x/\y in {0/0,0/2,3/2,3/4} {\draw[rounded corners=2em,line width=3.5em,blue!25,cap=round] (0+\x,2+\y)--(-3+\x,4+\y);}
\draw[rounded corners=2em,line width=3.5em,green!25,cap=round] (0,10)--(0,14);
\draw[rounded corners=2em,line width=3.5em,green!25,cap=round] (0,16)--(0,16.01);
\draw[line width=3.5em,green!25] (0,16)--(0,17);
\dottedintervalx0{16}11{.43}
\node[rounded corners,rectangle,draw,fill=black!30,ultra thick] (N) at (0,19) {$R_{I_r}^{I_r}$};
   \draw(1.1,19) node {$=\nab_{I_r}$};
\node[rounded corners,rectangle,draw,fill=black!30,ultra thick] (R3) at (0,16) {$R_{I_3}^{I_r}$};
\node[rounded corners,rectangle,draw,fill=white] (S3) at (0,14) {$R_{\S_3}^{I_r}$};
\node[rounded corners,rectangle,draw,fill=white] (A3) at (0,12) {$R_{\A_3}^{I_r}$};
\node[rounded corners,rectangle,draw,fill=black!30,ultra thick] (R2) at (0,10) {$R_{I_2}^{I_r}$};
\node[rounded corners,rectangle,draw,fill=white] (krl) at (0,8) {$R_{\S_2}^{I_r}$};
\node[rounded corners,rectangle,draw,fill=white] (kr) at (-3,6) {$\lam_{\S_2}^{I_r}$};
\node[rounded corners,rectangle,draw,fill=white] (kl) at (0,6) {$\rho_{\S_2}^{I_r}$};
\node[rounded corners,rectangle,draw,fill=black!30,ultra thick] (mrl) at (3,6) {$R_{I_1}^{I_r}$};
\node[rounded corners,rectangle,draw,fill=white] (k) at (-3,4) {$\mu_{\S_2}^{I_r}$};
\node[rounded corners,rectangle,draw,fill=white] (mr) at (0,4) {$\lam_{I_1}^{I_r}$};
\node[rounded corners,rectangle,draw,fill=white] (ml) at (3,4) {$\rho_{I_1}^{I_r}$};
\node[rounded corners,rectangle,draw,fill=black!30,ultra thick] (rl) at (6,4) {$R_{I_0}^{I_r}$};
\node[rounded corners,rectangle,draw,fill=white] (m) at (0,2) {$\mu_{I_1}^{I_r}$};
\node[rounded corners,rectangle,draw,fill=white] (r) at (3,2) {$\lam_{I_0}^{I_r}$};
\node[rounded corners,rectangle,draw,fill=white] (l) at (6,2) {$\rho_{I_0}^{I_r}$};
\node[rounded corners,rectangle,draw,fill=black!30,ultra thick] (d) at (3,0) {$\mu_{I_0}^{I_r}$};
   \draw(4.4,-0.05) node {$=\De_{I_r}$};
\draw
(d)--(l) (m)--(ml) (k)--(kl)
(d)--(m)--(k) (l)--(ml)--(kl) (rl)--(mrl)--(krl)
(d)--(r) (m)--(mr) (k)--(kr)
(l)--(rl) (ml)--(mrl) (kl)--(krl)
(krl)--(R2)
(R2)--(A3)--(S3)
;
\fill[blue!25] (-1.5,5)circle(.15);
\fill[blue!25] (1.5,5)circle(.15);
\fill[white] (1.5,3)circle(.15);
\fill[white] (4.5,3)circle(.15);
\draw
(r)--(rl) (mr)--(mrl) (kr)--(krl) (r)--(mr)--(kr) 
;
\end{tikzpicture}
}
\caption{
Hasse diagram of $\Cong(I_r)$ for $r\geq2$, where $I_r=I_r(\P)$.  
See Theorem \ref{thm:P}.
}
\label{fig:P}
\end{center}
\end{figure}

As discussed in Section \ref{subsect:EMRT}, these congruences  together with the $R_N^{I_r}$ and the universal congruence~$\nab_{I_r}$, form the lattice shown in Figure \ref{fig:P}, where $\A_q\leq\S_q$ denotes the alternating group.  
Our goal is to show that these are \emph{all} the congruences on $I_r$ ($r\geq2$).  We will also have to describe separately the congruences on the two smallest ideals, $I_0=D_0$ and $I_1=D_0\cup D_1$, and these will be dealt with using Propositions \ref{prop:PRB} and~\ref{prop:small_ideal}.  As with Theorem~\ref{thm:T}, the main effort will be in describing $\Cong(I_2)$, and for this we will use Proposition \ref{prop:small012}; after that, we will use Theorem~\ref{thm:E2} to treat the larger ideals.

We begin with a lemma that will be used to verify the technical assumptions of Proposition~\ref{prop:small012}.  As in that proposition, the statement refers to the relation $\L f^{-1}$, where $f$ is the retraction \eqref{eq:retract_P}, so
\[
\L f^{-1} = \bigset{(\al,\be)\in I_1\times I_1}{(\wh\al,\wh\be)\in\L}.
\]

\begin{lemma}\label{lem:R_trick}
For any $\al,\be\in I_2=I_2(\P)$ with $\al\sim\be$ and $(\al,\be)\not\in\L\cup\L f^{-1}$, we have $\rho_{I_0}^{I_2}\sub\cg\al\be_{I_2}$.
\end{lemma}

\pf
Write $\si=\cg\al\be_{I_2}$, and suppose $\al,\be\in\P_{A,B}$.  Since $(\al,\be)\not\in\L$, it follows that $\codom(\al)\not=\codom(\be)$ or $\coker(\al)\not=\coker(\be)$.

\pfcase{1}
Suppose first that $\coker(\al)\not=\coker(\be)$.  
Without loss of generality, we may fix some ${(u,v)\in\coker(\al)\sm\coker(\be)}$.  Now let $(\ga,\de)\in\rho_{I_0}^{I_2}$ be arbitrary, say with $\ga,\de\in\P_{C,D}$; we must show that $(\ga,\de)\in\si$.  This is clear if $\ga=\de$, so suppose instead that $\ga,\de\in D_0$ and $(\ga,\de)\in\R$.  Then we may write ${\ga=\partII{A_1}\cdots{A_k}{B_1}\cdots{B_l}}$ and ${\de=\partII{A_1}\cdots{A_k}{C_1}\cdots{C_m}}$.  To show that $(\ga,\de)\in\si$, it suffices to show that $(\ga,\ga')\in\si$, where $\ga'=\sqpartXXIII{A_1}\cdots{A_k}{}{}{}\in\P_{C,D}$. Indeed, the same argument applied to $\de$ will also give $(\de,\ga')\in\si$, and so $(\ga,\de)\in\si$.

We prove that $(\ga,\ga')\in\si$ by descending induction on $l$.  If $l=|D|$, then $\ga=\ga'$, so suppose $l<|D|$.  Without loss of generality, we will assume that $|B_l|\geq2$; let $b\in B_l$ be arbitrary and let $C=B_l\sm\{b\}$.  Let $\ve_1\in\P_{C,A}$ and $\ve_2\in\P_{B,D}$ be such that:
\bit
\item $\ve_1\in D_0$, $\ker(\ve_1)=\ker(\ga)$ and $\coker(\ve_1)=\De_A$,
\item $\ve_2\in D_2$ has transversals $\{u\}\cup C'$ and $\{v,b'\}$ and lower non-transversals $B_1',\ldots,B_{l-1}'$, and $\ker(\ve_2)=\De_B$.
\eit
Then $(\ga,\ve_1\be\ve_2)=(\ve_1\al\ve_2,\ve_1\be\ve_2)\in\si$, and $\ve_1\be\ve_2=
\Big(  \hspace{-1.5 truemm}
{\scriptsize \renewcommand*{\arraystretch}{1} \begin{array} {\cstart|\c|\c|\c|\cend}
A_1 \:& \multicolumn{3}{c|}{\cdots\cdots\cdots} &\: A_k \\ \cline{1-5}
B_1 \:&\: \cdots \:&\: B_{l-1} \:&\: C \:&\: b
\rule[0mm]{0mm}{2.7mm}
\end{array}  }
\hspace{-1.5 truemm} \Big) 
$.  By induction, we have $(\ve_1\be\ve_2,\ga')\in\si$, and so $(\ga,\ga')\in\si$ by transitivity, as required.

\pfcase{2}
Now suppose $\coker(\al)=\coker(\be)$ but $\codom(\al)\not=\codom(\be)$.  We may assume without loss of generality that $\al\in D_2$.  Indeed, if $\al,\be\in I_1$, then we would have $\coker(\wh\al)=\coker(\al)=\coker(\be)=\coker(\wh\be)$, so that $(\wh\al,\wh\be)\in\L$; i.e., $(\al,\be)\in\L f^{-1}$, a contradiction.  Thus, we may write $\al=\partV{A_1}{A_2}{C_1}{C_s}{B_1}{B_2}{D_1}{D_t}$.  Choose any $x\in A_1$ and $y\in A_2$, and define $\ga=\sqpartXX{}{x,y}\in\P_A$ (so $\{x',y'\}$ is the only non-trivial block of $\ga$).  Then $B_1'\cup B_2'$ is a block of $\ga\al$, but at least one of $B_1'$ or $B_2'$ is a block of $\ga\be$ (since at least one of these is a block of $\be$).  Thus, $\coker(\ga\al)\not=\coker(\ga\be)$.  Since $(\ga\al,\ga\be)\in\si$, we have reduced to Case~1.
\epf

We may now prove the main result of this section.  The second part utilises the congruences~$\th_{I,\tau}^S$ defined in Section \ref{subsect:EMRT}.  For a positive integer $n$, we denote by $\bn$ the chain~$\{1<\cdots<n\}$.

\begin{thm}\label{thm:P}
Let $\C$ be a non-empty set of finite sets, let $\P=\P(\C)$ be the partition category over $\C$, and keep the above notation.
\ben
\item \label{cP1}  The minimal ideal $I_0$ is a partial rectangular band, so its congruence lattice is described by Proposition \ref{prop:PRB}.  In particular, $\Cong(I_0)$ is isomorphic to the lattice direct product
\[
\prod_{A\in\C}\Eq(\Eq(A))\times\prod_{A\in\C}\Eq(\Eq(A)).
\]
\item \label{cP2}  The ideal $I_1$ is retractable, every congruence on $I_0$ is liftable to~$I_1$, and
\[
\Cong(I_1)=\bigset{\tau\cup\De_{D_1}}{\tau\in\Cong(I_0)} \cup \bigset{\th_{I_0,\tau}^{I_1}}{\tau\in\Cong(I_0)}.
\]
In particular, $\Cong(I_1)$ is isomorphic to the lattice direct product $\Cong(I_0)\times\btwo$.
\newpage
\item \label{cP3} For $r\geq2$ (including $r=\om$ if $\xi=\om$), 
\begin{align*}
\Cong(I_r) = 
\set{\mu_{I_q,N^\natural}^{I_r},\ \lam_{I_q,N^\natural}^{I_r},\ \rho_{I_q,N^\natural}^{I_r},\ &R_{I_q,N^\natural}^{I_r}}{0\leq q\leq1,\ N\normal\S_{q+1}} \\[3truemm]
\qquad\qquad\qquad\qquad {}\cup\set{&R_{I_q,N^\natural}^{I_r}}{2\leq q<r,\ N\normal\S_{q+1}}\cup\{\nab_{I_r}\}
\end{align*}
forms the lattice pictured in Figure \ref{fig:P}.
\een
\end{thm}

\pf
\firstpfitem{\ref{cP1}}  The proof of this part is analogous to the corresponding part of Theorem \ref{thm:T}.  For any $A\in\C$, the $\R$- and $\L$-classes in $I_0(\P_A)$ are indexed by the elements of $\Eq(A)$: i.e., all possible kernels and cokernels.

\pfitem{\ref{cP2}}  This follows from Proposition \ref{prop:small_ideal}, and the fact that $\G_1$ is trivial, so that ${\Cong(\G_1)^\top\cong\btwo}$.  (The technical assumption in the statement of Proposition \ref{prop:small_ideal} follows from Lemma \ref{lem:P2P}\ref{P2Pii}.)

\pfitem{\ref{cP3}}  As in the proof of Theorem \ref{thm:T}\ref{T2}, it suffices to establish the claim in the case $r=2$.  Indeed, if we can do so, then it will follow that:
\bit
\item every congruence on $I_2$ is of the form listed in Lemma \ref{lem:E2}, and so $\LiftCong{I_2}{\P}=\Cong(I_2)$,
\item $I_2$ has property $\ngen{I_2}$,
\item $I_r$ has property $\Mult{I_r}$, and hence $\Sep{I_r}$, for each $2<r<\xi$ (cf.~Lemma \ref{lem:P2P}).
\eit
An application of Theorem \ref{thm:E2} will then complete the proof.
For the $r=2$ case, we apply Proposition \ref{prop:small012}, noting the following:
\bit
\item Property $\Sep{I_2}$ follows from Lemma \ref{lem:P2P} (which shows $\Mult{I_2}$).
\item Assumption \ref{012i} from the proposition follows from Lemma \ref{lem:P2P}\ref{P2Pii}.
\item Assumptions \ref{012ii} and \ref{012iii} follow from Lemma \ref{lem:R_trick} and its dual.
\item Assumption \ref{012iv} holds with $H=\G_2$, since $(I_1,\G_2)$ is retractable.  The elements $a,b$ are unneeded since $H=\G_2$.  \qedhere
\eit
\epf

\subsection{Planar and annular reducts of $\P$}\label{subsect:QP}

In this section we consider the partition analogues of the planar and annular reducts of the transformation category $\T$, introduced in Section \ref{subsect:QT}.
To this end we assume that there is a total order $\leq$ on each member of $\C$, and keep the notation from the previous three sections.

Consider some partition $\al\in\P_{A,B}$, where $A,B\in\C$.  Write $A=\{a_1<\cdots<a_k\}$ and ${B=\{b_1<\cdots<b_l\}}$, and put $m=\max(k,l)$.  We say that $\al$ is \emph{planar} if some graph representing~$\al$ can be drawn in the plane $\RRR^2$ with:
\bit
\item vertex $a_i$ at $(i,1)$ for each $1\leq i\leq k$, 
\item vertex $b_i'$ at $(i,0)$ for each $1\leq i\leq l$, 
\item edges all in the rectangle $\bigset{(x,y)\in\mathbb R^2}{0\leq x\leq m+1,\ 0\leq y\leq1}$, and
\item with no edge crossings in the interior of this rectangle.
\eit
We say that $\al$ is \emph{annular} if the above conditions are satisfied, but with the left- and right-hand borders of the rectangle identified to form a cylinder, and with edges allowed to ``wrap around''.  

\begin{figure}[ht]
\begin{center}
\begin{tikzpicture}[scale=.6]
\begin{scope}[shift={(0,0)}]	
\fill[gray!30] (0,0)--(9,0)--(9,2)--(0,2)--(0,0);
\uarcx23{.3}
\uarcx56{.3}
\darc12
\darcx45{.3}
\darc78
\stline34
\uarcxx14{.6}{red,ultra thick}
\darcxx26{.8}{blue,ultra thick}
\draw(-.3,1)node[left]{$\al=$};
\draw(11,1)node{$\equiv$};
\uvs{1,...,6}
\lvs{1,...,8}
\end{scope}
\begin{scope}[shift={(13,0)}]	
\fill[gray!30] (0,0)--(9,0)--(9,2)--(0,2)--(0,0);
\uarcx23{.3}
\uarcx56{.3}
\darc12
\darcx45{.3}
\darc78
\stline34
\draw[red,ultra thick] (0,1.4)--(.4,1.4) arc(270:360:.6);
\draw[red,ultra thick] (4,2)arc(180:270:.6) (4.6,1.4)--(9,1.4);
\draw[blue,ultra thick] (6,0)arc(180:90:.8) (6.8,.8)--(9,.8);
\draw[blue,ultra thick] (0,.8)--(1.2,.8) arc(90:0:.8);
\uvs{1,...,6}
\lvs{1,...,8}
\draw[dashed,line width=0.4mm] (0,0)--(0,2);
\draw[dashed,line width=0.4mm] (9,0)--(9,2);
\draw[white] (9.3,1)node[right]{$\al=$};
\end{scope}
\end{tikzpicture}
\caption{An annular representation of the partition $\al\in\P_{A,B}$ from Figure \ref{fig:P_product}.}
\label{fig:P_annular}
\end{center}
\end{figure}
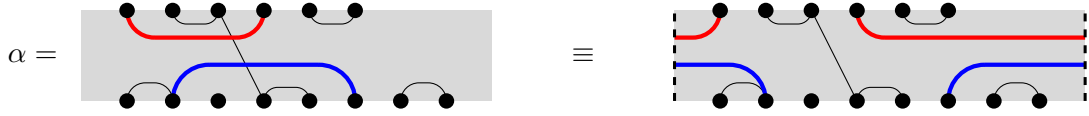

To illustrate these concepts, consider the partitions $\al$ and $\be$ pictured in Figure \ref{fig:P_product}.  It is clear that $\be$ is planar, and that $\al$ is not; indeed, in any graphical representation (as above) of $\al$, there will be a path joining vertices $1$ and $4$, and a path joining $3$ and $4'$, and these will always intersect.  On the other hand, $\al$ is annular, as we show in Figure \ref{fig:P_annular}.

We will introduce the anti-planar and anti-annular concepts shortly, but first make an important observation.  For $A=\{a_1<\cdots<a_k\}\in\C$, we define the following partitions from~$\P_A$:
\begin{equation}\label{eq:gade}
\ga_A = \partpermV{a_1}{a_2}{a_k}{a_k}{a_{k-1}}{a_1} \AND \de_A = \partpermIV{a_1}{a_{k-1}}{a_k}{a_2}{a_k}{a_1}.
\end{equation}
Note that in fact $\ga_A$ and $\de_A$ are permutations, and that $\ga_A$ is order-reversing, and $\de_A$ orientation-preserving.  By definition, it is clear that for any partition $\al\in\P_{A,B}$,
\begin{equation}\label{eq:annular}
\text{$\al$ is annular $\iff$ $\al=\de_A^m\be\de_B^n$ for some planar $\be$ and some $m,n\in\N$.}
\end{equation}
To see this, consider a representation of $\al$ as a cylindrical graph, and then ``twist'' the top and bottom of $\al$ as much as required to ensure that no edge needs to wrap around the back of the cylinder.
We say that 
\bit
\item a partition $\al\in\P_{A,B}$ is \emph{anti-planar} $\iff$ $\ga_A\al$ is planar $\iff$ $\al\ga_B$ is planar,
\item a partition $\al\in\P_{A,B}$ is \emph{anti-annular} $\iff$ $\ga_A\al$ is annular  $\iff$ $\al\ga_B$ is annular,
\eit
As we did with transformations, given any subset $\Om\sub\P$, we write
\bit
\item $\OrdP\Om$ for the set of all planar partitions from $\Om$, 
\item $\OrdR\Om$ for the set of all anti-planar partitions from $\Om$, 
\item $\OrdPR\Om = \OrdP\Om\cup\OrdR\Om$,
\item $\OriP\Om$ for the set of all annular partitions from $\Om$, 
\item $\OriR\Om$ for the set of all anti-annular partitions from $\Om$, 
\item $\OriPR\Om = \OriP\Om\cup\OriR\Om$.
\eit

If $\Q$ denotes any of $\OrdP\P$, $\OrdPR\P$, $\OriP\P$ or $\OriPR\P$, then $\Q$ is subcategory of $\P$.  Indeed, it is clear that $\OrdP\P$ and $\OriP\P$ are categories; if $\al,\be\in\P$ are composable and have planar (respectively, annular) graphs, then the product graph $\Pi(\al,\be)$ is also planar (respectively, annular), and so too therefore is the product $\al\be$.  
It quickly follows that $\OrdPR\P$ and $\OriPR\P$ are closed under well-defined composition.  For example, suppose $\al\in\P_{A,B}$ and $\be\in\P_{B,C}$.
\bit
\item If $\al\in\OrdP\P$ and $\be\in\OrdR{\P}$, then $\al\be = \al(\be\ga_C)\ga_C$, and since $\be\ga_C\in\OrdP\P$, we have $\al(\be\ga_C)\in\OrdP\P$ and hence $\al\be\in\OrdR{\P}$.
\item If $\al,\be\in\OriR{\P}$, then $\al\be=(\al\ga_B)(\ga_B\be)\in\OriP\P$, since both $\al\ga_B$ and $\ga_B\be$ belong to~$\OriP\P$.
\eit
Other combinations are checked analogously.

Clearly $\Q$ (still denoting any of the above subcategories of $\P$) is closed under the involution $\al\mt\al^*$; in particular, it is regular.  Thus, Green's $\R$, $\L$ and $\H$ relations on $\Q$ are simply inherited from those on $\P$.  This is easily seen to also be the case for the $\D=\J$ relation (cf.~\cite[Theorem 4.9]{DDE2019}), and hence the first three parts of Lemma \ref{lem:Green_P} hold for $\Q$ as well.  
As for the fourth part, ignoring the obvious exceptions for small $q$, the group $\H$-classes of $\Q$ contained in $D_q$ are:
\bit
\item trivial, if $\Q=\OrdP\P$,
\item cyclic of order $2$, if $\Q=\OrdPR\P$,
\item cyclic of order $q$, if $\Q=\OriP\P$,
\item dihedral of order $2q$, if $\Q=\OriPR\P$.
\eit
We also have a version of Lemma \ref{lem:P2P} for $\Q$.  Before we discuss this, we first consider the possible kernels (and cokernels) that elements of $\Q$ can have.  Let $A=\{a_1<\cdots<a_k\}\in\C$.  For subsets $U,V\sub A$, we say:
\bit
\item $U$ and $V$ are \emph{separated} if $u<v$ for all $u\in U$ and $v\in V$, or if $v<u$ for all such $u,v$,
\item $U$ is \emph{nested} by $V$ if there exist $v_1,v_2\in V$ such that $v_1<u<v_2$ for all $u\in U$, and no $v\in V$ satisfies $v_1<v<v_2$; here we say $U$ is nested by~$V$ \emph{via the pair $(v_1,v_2)$}.
\eit
We say an equivalence $\ve\in\Eq(A)$ is \emph{planar} if any pair of $\ve$-classes is either separated or nested.  We will write $\OrdP{\Eq(A)}$ for the set of all planar partitions on $A$.

\begin{lemma}\label{lem:PEqA}
If $\al\in\OriPR\P$, then $\ker(\al)$ and $\coker(\al)$ are both planar equivalences.
\end{lemma}

\pf
Write $A=\bd(\al)$ and $B=\br(\al)$.  By symmetry it suffices to prove that $\coker(\al)$ is planar.  In fact, it further suffices to assume that $\al\in\OriP\P$; indeed, if $\al\in\OriR\P$, then $\al=\ga_A(\ga_A\al)$ with $\coker(\al)=\coker(\ga_A\al)$ and $\ga_A\al\in\OriP\P$.

Thus, we assume for the rest of the proof that $\al\in\OriP\P$.  By \eqref{eq:annular}, we have $\al=\de_A^m\be\de_B^n$ for some $\be\in\OrdP\P$ and $m,n\in\N$, and we note that $\coker(\al)=\coker(\be\de_B^n)$.  By \cite[Lemma~7.1]{EMRT2018}, $\coker(\be)$ is a planar equivalence.  It therefore suffices to show that for any planar equivalence $\ve\in\OrdP{\Eq(B)}$, the equivalence
\[
\ve\cdot\de_B = \bigset{(x\de_B,y\de_B)}{(x,y)\in\ve}
\]
is also planar.  To do so, assume for convenience that $B=[k]=\{1,\ldots,k\}$, let the $\ve$-classes be $B_1,\ldots,B_r$, where $k\in B_r$, and write $\ve^+=\ve\cdot\de_B$.  The $\ve^+$-classes are $B_i^+=\set{x+1}{x\in B_i}$ for each $i$; here elements of $[k]$ are added modulo $k$ (so $1\in B_r^+$).  For $1\leq i,j<r$, the sets $B_i^+$ and~$B_j^+$ are nested (or separated) if and only if $B_i$ and $B_j$ are nested (or separated), respectively.  Write $B_r=\{b_1<\cdots<b_l<k\}$, where $l\geq0$.  Let $1\leq i<r$.  If $B_i$ is separated from $B_r$, then $\max(B_i)<\min(B_r)$, and so $B_i^+$ will be nested by $B_r^+$ via $(1,b_1+1)$ if $l>0$, or separated from~$B_r^+$ if $l=0$.  If $B_i$ is nested by $B_r$ via $(b_l,k)$, then $B_i^+$ and $B_r^+$ are separated.  If $B_i$ is nested by $B_r$ via $(b_p,b_{p+1})$ for some $1\leq p<l$, then $B_i^+$ is nested by $B_r^+$ via $(b_p+1,b_{p+1}+1)$.  Since $k\in B_r$, $B_r$ is not nested by any $B_i$.
\epf

As a partial converse to Lemma \ref{lem:PEqA}, it is clear that if $\ve_1\in\OrdP{\Eq(A)}$ and ${\ve_2\in\OrdP{\Eq(B)}}$, where $A,B\in\C$, then there exists $\al\in D_0(\OrdP{\P_{A,B}})$ such that $\ker(\al)=\ve_1$ and $\coker(\al)=\ve_2$.  
It follows that the set of all possible (co)kernels is the same for all four of the above categories.  In particular, all four categories have the same minimal ideal, which is the set of all partitions of rank $0$ with planar kernel and cokernel.

On the other hand, the categories $\OriP\P$ and $\OriPR\P$ have more $\R$- and $\L$-classes than the categories $\OrdP\P$ and $\OrdPR\P$.  For example, consider the planar equivalence on $\{1,2,3\}$ with equivalence classes $\{1,3\}$ and $\{2\}$.  There is no element from $\OrdPR\P$ with this kernel and domain $\{2\}$; however there is such an element in~$\OriP\P$.

Here is the version of Lemma \ref{lem:P2P} for $\Q$ (being any of the above four categories):

\begin{lemma}\label{lem:P2PQ}
Let $A,B\in\C$, and suppose $\al,\be\in\Q_{A,B}$ with $r=\rank(\al)\geq\rank(\be)=q$.
\ben
\item \label{P2PQi} If $r\geq2$ and $q<r$, then there exists $\ga\in I_r$ such that $\ga\al\in D_{r-1}$ and $\ga\be\in I_{r-1}\sm H_{\ga\al}$.
\item \label{P2PQii} If $q=r\geq1$ and $(\al,\be)\not\in\H$, then there exists $\ga\in I_r$ such that, swapping $\al,\be$ if necessary,
\[
\text{$[\al\ga\in D_r$ and $\be\ga\in I_{r-1}]$ \ \ \ or \ \ \ $[\ga\al\in D_r$ and $\ga\be\in I_{r-1}]$.}
\]
\item \label{P2PQiii} If $r\geq2$ and $\be\in H_\al\sm\{\al\}$, then there exists $\ga\in I_r$ such that $\ga\al\in D_{r-1}$ and $\ga\be\in I_{r-1}\sm H_{\ga\al}$.
\een
\end{lemma}

\pf
In the proof of Lemma \ref{lem:P2P}, all but one of the partitions $\ga$ constructed were of the form $\ga=\sqpartpermIII{a_1}\cdots{a_k}{a_1}\cdots{a_k}$, and hence belong to $\OrdP\P\sub\Q$.  The only exception is the partition $\ga=\sqpartII{a_1}{a_{r-1}}{a_r}{}{a_1}{a_{r-1}}{x_1}{a_r,x_2}$ constructed during Case 2.2 of the proof of \ref{P2Pii}.  So our task is to modify the definition of $\ga$ in this case to show that it may be taken from $\Q$ if $\al,\be$ themselves belong to $\Q$.  

\pfcase1
Suppose first that $\Q=\OrdP\P$.  Without loss of generality, we may assume that $a_i=\min(A_i)$ for all $i$.  Let $B$ be the upper non-transversal of $\be$ containing $x_1$ and $x_2$.  By planarity of $\be$, and since $a_r$ belongs to a transversal of $\be$, we have either $a_r<b$ for all $b\in B$, or else $a_r>b$ for all $b\in B$.  By symmetry, we will assume that $a_r<b$ for all $b$; then $B$ is either nested by $A_r$ or else $A_r$ and $B$ are separated, with $\max(A_r)<\min(B)$.  (The case in which $a_r>b$ for all $b\in B$ is symmetrical to the case of $\max(A_r)<\min(B)$.)  Permuting the transversals if necessary, we may assume that we do not have $a_r<a_i<b$ for any $b\in B$ (i.e., that $a_r$ is the closest of the $a_i$ to $B$ from the left).  We may also assume without loss of generality that $x_2<x_1$.  With these conditions satisfied, it is then clear that~$\ga$ as above is planar.  

It is important to note that the argument just given assumed only the planarity of~$\be$.

\pfcase2    
Next suppose $\Q=\OrdPR\P$.  Then $\al=\al'\ga_B^k$ and $\be=\be'\ga_B^l$ for some $k,l\in\{0,1\}$ and $\al',\be'\in\OrdP\P$.  We may then take $\ga\in\OrdP\P$ corresponding to $\al',\be'$, as in the previous case.

\pfcase3
Next suppose $\Q=\OriP\P$.  This time we have $\al=\de_A^k\al'\de_B^l$ and $\be=\de_A^m\be'\de_B^n$ for some $k,l,m,n\in\N$ and $\al',\be'\in\OrdP\P$.  Let $\al''=\de_A^{k-m}\al'\de_B^{l-n}$ (with possibly negative powers interpreted in the obvious way, as $\de_A\in\S_A$ and $\de_B\in\S_B$ are units).  Note that
\[
\al'' = \de_A^{-m} \al \de_B^{-n} \AND \be' = \de_A^{-m}\be\de_B^{-n}.
\]
Now, $\be'$ is planar (but $\al''$ might not be), and $\ker(\be')\sm\ker(\al'')$ is still non-empty, so by Case 1 there exists $\ga'\in D_r(\OrdP\P)$ such that $\ga'\be'\in D_r$ and $\ga'\al''\in I_{r-1}$.  Thus,
\bit
\item $\ga'\de_A^{-m}\be\de_B^{-n}=\ga'\be'\in D_r \implies \ga'\de_A^{-m}\cdot\be\in D_r$, and\\[-2truemm]
\item $\ga'\de_A^{-m}\al\de_B^{-n}
=\ga'\al''\in I_{r-1} \implies \ga'\de_A^{-m}\cdot\al\in I_{r-1}$, 
\eit
so we may take $\ga=\ga'\de_A^{-m}\in\OriP\P$ in this case.

\pfcase4  
When $\Q=\OriPR\P$, the assertion follows from Case 3 just as Case 2 followed from Case~1.
\epf

All other arguments of Sections \ref{subsect:prelimP} and \ref{subsect:P} may be modified in much more straightforward fashion to work for $\Q$.  In particular:
\bit
\item The category $\Q$ is closed under the mapping $\al\mt\wh\al$.  
\item The ideal $I_1=I_1(\Q)$ is retractable.
\item Choosing natural group $\H$-classes $\G_q\leq D_q$ for each $q$, the pairs $(I_q,N)$ are retractable for $q=0,1$ and $N\normal\G_{q+1}$.  If $\Q=\OrdP\P$, then there are only two such pairs since $\G_2$ is trivial.  Otherwise, there are three such pairs, since $\G_q(\Q)=\G_q(\P)=\S_q^\natural$ for $q\leq2$ if $\Q\not=\OrdP\P$.
\item In Case 1 of the proof of Lemma \ref{lem:R_trick}, the partition $\ve_1$ constructed will belong to $\Q$ (provided that~$\ga$ does); we can ensure that the partition $\ve_2$ belongs to $\OrdP\P$ by assuming that $u<v$, choosing~$B_l$ to be a lower non-transversal of $\ga$ that is not nested by any other blocks of~$\ga$, and choosing~${b=\max(B_l)}$.  The element $\ga$ constructed during Case 2 clearly belongs to $\OrdP\P$.
\eit
From this we may then deduce the following description of the congruences of arbitrary ideals $I_r=I_r(\Q)$.  For the first part, note that the $\R$- and $\L$-classes of the partial rectangular band $I_0=I_0(\Q)$ are indexed by all possible kernels/cokernels, which are the planar equivalences on the underlying sets $A\in\C$.  Recall that $\OrdP{\Eq(A)}$ is the set of all such planar equivalences.

\begin{thm}\label{thm:QP}
Let $\C$ be a non-empty set of finite sets, let $\P=\P(\C)$ be the category of all partitions between members of $\C$, let $\Q$ be one of $\OrdP\P$, $\OrdPR\P$, $\OriP\P$ or $\OriPR\P$, and keep the above notation.
\ben
\item \label{cQP1}  The minimal ideal $I_0$ is a partial rectangular band, so its congruence lattice is described by Proposition \ref{prop:PRB}.  In particular, $\Cong(I_0)$ is isomorphic to the lattice direct product
\[
\prod_{A\in\C}\Eq(\OrdP{\Eq(A)})\times\prod_{A\in\C}\Eq(\OrdP{\Eq(A)}).
\]
\item \label{cQP2}  The ideal $I_1$ is retractable, every congruence on $I_0$ is liftable to~$I_1$, and
\[
\Cong(I_1)=\set{\tau\cup\De_{D_1}}{\tau\in\Cong(I_0)} \cup \set{\th_{I_0,\tau}^{I_1}}{\tau\in\Cong(I_0)}.
\]
In particular, $\Cong(I_1)$ is isomorphic to the lattice direct product $\Cong(I_0)\times\btwo$.
\item \label{cQP3} For $r\geq2$ (including $r=\om$ if $\xi=\om$), 
\begin{align*}
\Cong(I_r) = 
\set{\mu_{I_q,N}^{I_r},\ \lam_{I_q,N}^{I_r},\ \rho_{I_q,N}^{I_r},\ &R_{I_q,N}^{I_r}}{0\leq q\leq1,\ N\normal\G_{q+1}} \\[3truemm]
\qquad\qquad\qquad\qquad {}\cup\set{&R_{I_q,N}^{I_r}}{2\leq q<r,\ N\normal\G_{q+1}}\cup\{\nab_{I_r}\}.
\epfreseq
\end{align*}
\een
\end{thm}

\begin{rem}\label{rem:QP}
Figure \ref{fig:QP} pictures the lattices $\Cong(I_4(\Q))$ for the categories $\Q=\P$, $\OrdP\P$, $\OrdPR\P$, $\OriP\P$ and $\OriPR\P$.

The upper parts of the lattices $\Cong(I_r(\Q))$, $r\geq2$, are always chains when $\Q$ is any of~$\P$,~$\OrdP\P$ or~$\OrdPR\P$, since the lattices $\NN(\G_q)$ are always chains for these categories (where~$\G_q$ is always a symmetric group, or of size at most $2$).  This is not the case when $\Q$ is $\OriP\P$ or $\OriPR\P$, however, even though Figure \ref{fig:QP} might suggest otherwise for $\OriP\P$; indeed, the lattices $\NN(\G_q)$ are not chains in general in these categories, as the groups $\G_q$ are generally cyclic or dihedral of orders $q$ and $2q$, respectively.
\end{rem}

\begin{figure}[ht]
\begin{center}
\begin{tikzpicture}[scale=0.7,inner sep=0pt]
\begin{scope}[shift={(-4,0)}]
\draw[rounded corners=4,line width=10,green!25,cap=round]
(1,5)--(1,7)
(1,8)--(1,11)
;
\foreach \x/\y in {0/0,0/-1,-1/-1,-1/-2} {\draw[rounded corners=4,line width=10,red!25,cap=round](3+\x,2+\y)--(3+\x,2.001+\y);}
\foreach \x/\y in {0/0,0/-1,-1/-1,-1/-2} {\draw[rounded corners=4,line width=10,blue!25,cap=round](2+\x,3+\y)--(1+\x,4+\y);}
\foreach \x in {0,1,2} {\draw(2-\x,0+\x)--(3-\x,1+\x)--(3-\x,2+\x)--(2-\x,1+\x)--(2-\x,0+\x);}
\draw(1,4)--(1,12);
\foreach \x/\y in {0/0,0/1,1/1,1/2} {\draw(2+\x,0+\y)--(0+\x,2+\y);}
\foreach \x in {1,...,12} {\colvert{1}{\x}{white}}
\foreach \x/\y in {0/2,0/3,2/0,2/1,2/2,2/3,3/1,3/2} {\colvert{\x}{\y}{white}}
\foreach \x in {5,8,12} {\colvert{1}{\x}{black}}
\foreach \x/\y in {3/2,2/3,2/0} {\colvert{\x}{\y}{black}}
\end{scope}
\begin{scope}[shift={(0,0)}]
\draw[rounded corners=4,line width=10,green!25,cap=round]
(2,4)--(2,4.001)
(2,5)--(2,5.001)
;
\foreach \x/\y in {0/0,0/-1,-1/-1,-1/-2} {\draw[rounded corners=4,line width=10,red!25,cap=round](3+\x,2+\y)--(3+\x,2.001+\y);}
\foreach \x/\y in {0/0,0/-1,-1/-1,-1/-2} {\draw[rounded corners=4,line width=10,blue!25,cap=round](2+\x,3+\y)--(2+\x,3.001+\y);}
\foreach \x in {0,1} {\draw(2-\x,0+\x)--(3-\x,1+\x)--(3-\x,2+\x)--(2-\x,1+\x)--(2-\x,0+\x);}
\draw(2,3)--(2,6);
\foreach \x/\y in {0/0,0/1,1/1,1/2} {\draw(2+\x,0+\y)--(1+\x,1+\y);}
\foreach \x/\y in {1/1,1/2,2/0,2/1,2/2,2/3,3/1,3/2} {\colvert{\x}{\y}{white}}
\foreach \x/\y in {3/2,2/3,2/4,2/5,2/6,2/0} {\colvert{\x}{\y}{black}}
\end{scope}
\begin{scope}[shift={(5,0)}]
\draw[rounded corners=4,line width=10,green!25,cap=round]
(1,5)--(1,6)
(1,7)--(1,8)
;
\foreach \x/\y in {0/0,0/-1,-1/-1,-1/-2} {\draw[rounded corners=4,line width=10,red!25,cap=round](3+\x,2+\y)--(3+\x,2.001+\y);}
\foreach \x/\y in {0/0,0/-1,-1/-1,-1/-2} {\draw[rounded corners=4,line width=10,blue!25,cap=round](2+\x,3+\y)--(1+\x,4+\y);}
\foreach \x in {0,1,2} {\draw(2-\x,0+\x)--(3-\x,1+\x)--(3-\x,2+\x)--(2-\x,1+\x)--(2-\x,0+\x);}
\draw(1,4)--(1,9);
\foreach \x/\y in {0/0,0/1,1/1,1/2} {\draw(2+\x,0+\y)--(0+\x,2+\y);}
\foreach \x in {1,...,9} {\colvert{1}{\x}{white}}
\foreach \x/\y in {0/2,0/3,2/0,2/1,2/2,2/3,3/1,3/2} {\colvert{\x}{\y}{white}}
\foreach \x in {5,7,9} {\colvert{1}{\x}{black}}
\foreach \x/\y in {3/2,2/3,2/0} {\colvert{\x}{\y}{black}}
\end{scope}
\begin{scope}[shift={(10,0)}]
\draw[rounded corners=4,line width=10,green!25,cap=round]
(1,5)--(1,6)
(1,7)--(1,9)
;
\foreach \x/\y in {0/0,0/-1,-1/-1,-1/-2} {\draw[rounded corners=4,line width=10,red!25,cap=round](3+\x,2+\y)--(3+\x,2.001+\y);}
\foreach \x/\y in {0/0,0/-1,-1/-1,-1/-2} {\draw[rounded corners=4,line width=10,blue!25,cap=round](2+\x,3+\y)--(1+\x,4+\y);}
\foreach \x in {0,1,2} {\draw(2-\x,0+\x)--(3-\x,1+\x)--(3-\x,2+\x)--(2-\x,1+\x)--(2-\x,0+\x);}
\draw(1,4)--(1,10);
\foreach \x/\y in {0/0,0/1,1/1,1/2} {\draw(2+\x,0+\y)--(0+\x,2+\y);}
\foreach \x in {1,...,9} {\colvert{1}{\x}{white}}
\foreach \x/\y in {0/2,0/3,2/0,2/1,2/2,2/3,3/1,3/2} {\colvert{\x}{\y}{white}}
\foreach \x in {5,7,10} {\colvert{1}{\x}{black}}
\foreach \x/\y in {3/2,2/3,2/0} {\colvert{\x}{\y}{black}}
\end{scope}
\begin{scope}[shift={(15,0)}]
\draw[rounded corners=4,line width=10,green!25,cap=round]
(1,5)--(1,7)
(1,8)--(1,9) (1,9)--(0,10) (0,10)--(1,11) (1,11)--(2,10) (2,10)--(1,9)
;
\foreach \x/\y in {0/0,0/-1,-1/-1,-1/-2} {\draw[rounded corners=4,line width=10,red!25,cap=round](3+\x,2+\y)--(3+\x,2.001+\y);}
\foreach \x/\y in {0/0,0/-1,-1/-1,-1/-2} {\draw[rounded corners=4,line width=10,blue!25,cap=round](2+\x,3+\y)--(1+\x,4+\y);}
\fill[green!25] 
(1,9)--(0,10)--(1,11)--(2,10)--(1,9)
;
\foreach \x in {0,1,2} {\draw(2-\x,0+\x)--(3-\x,1+\x)--(3-\x,2+\x)--(2-\x,1+\x)--(2-\x,0+\x);}
\draw(1,4)--(1,12);
\draw(1,9)--(0,10)--(1,11)--(2,10)--(1,9);
\foreach \x/\y in {0/0,0/1,1/1,1/2} {\draw(2+\x,0+\y)--(0+\x,2+\y);}
\foreach \x in {1,...,12} {\colvert{1}{\x}{white}}
\foreach \x/\y in {0/2,0/3,0/10,2/10,2/0,2/1,2/2,2/3,3/1,3/2} {\colvert{\x}{\y}{white}}
\foreach \x in {5,8,12} {\colvert{1}{\x}{black}}
\foreach \x/\y in {3/2,2/3,2/0} {\colvert{\x}{\y}{black}}
\end{scope}
\end{tikzpicture}
\caption{Left to right: Hasse diagrams of $\Cong(I_4)$ in the categories $\P$, $\OrdP\P$, $\OrdPR\P$, $\OriP\P$, and~$\OriPR\P$; cf.~Theorem \ref{thm:QT}. }
\label{fig:QP}
\end{center}
\end{figure}

\subsection{Partial Brauer categories}\label{subsect:PB}

A partition $\al\in\P$ is a \emph{partial Brauer partition} if each of its blocks has size at most $2$.  The subset $\PB=\PB(\C)$ of~$\P$ consisting of all partial Brauer partitions is a subcategory, called the \emph{partial Brauer category}; cf.~\cite{MM2014}.  
All partitions constructed during proofs in Sections \ref{subsect:prelimP}--\ref{subsect:QP} actually belonged to $\PB$.  It follows that the congruences on ideals of $\PB$ are exactly as described in Theorem \ref{thm:P}, but with one difference for the ideals $I_q=I_q(\PB)$ for $q=0,1$.  The kernels (and cokernels) of partial Brauer partitions are equivalences whose blocks have size at most $2$.  So, if we write $\Eq_{1,2}(A)$ for the set of all such equivalences on the set $A$, then in the category $\PB$, the lattice $\Cong(I_0)$ is isomorphic to
\[
\prod_{A\in\C}\Eq(\Eq_{1,2}(A))\times\prod_{A\in\C}\Eq(\Eq_{1,2}(A)),
\]
and as with $\P$, we have $\Cong(I_1)\cong\Cong(I_0)\times\btwo$.  

Similarly, the congruences on the ideals of the subcategories $\OrdP{\PB}$, $\OrdPR{\PB}$, $\OriP{\PB}$ and $\OriPR{\PB}$ are described in the same way as the corresponding subcategories of $\P$ in Theorem~\ref{thm:QP}, with the exception that $\Cong(I_0)$ is isomorphic to
\[
\prod_{A\in\C}\Eq(\OrdP{\Eq_{1,2}(A)})\times\prod_{A\in\C}\Eq(\OrdP{\Eq_{1,2}(A)}),
\]
and that $\Cong(I_1)\cong\Cong(I_0)\times\btwo$.  The category $\OrdP{\PB}$ could be called the \emph{Motzkin category}, given the obvious connections to the Motzkin algebras and monoids \cite{BH2014,DEG2017}.

It would be interesting to attempt to apply the methods of this paper to other natural diagram categories, such as rook partition categories; cf.~\cite{Grood2006,JErook}.

\section{Brauer categories}\label{sect:B}

The next three sections concern the \emph{Brauer categories} $\B=\B(\C)$, and their associated planar and annular reducts, which include the Temperley--Lieb and Jones categories.  Unlike the situation with the 
transformation and partition categories, $\B$ is not immediately amenable to applying our general machinery from Section \ref{sect:chains} because it is not a chain of ideals.  
However, it decomposes into ``even'' and ``odd'' subcategories, as we explain in Section \ref{subsect:defnB}, to each of which these general results do apply.  After describing Green's relations and establishing multiplicative properties in Section~\ref{subsect:prelimB}, we describe the congruences in Sections \ref{subsect:Bodd} and \ref{subsect:Beven} for the odd and even subcategories, respectively; see Theorems \ref{thm:Bodd} and \ref{thm:Beven}, and also Figures~\ref{fig:Bodd}--\ref{fig:B4}.  
As these figures indicate, the congruence lattices of ideals in the even subcategory are more complex than those for the odd subcategory; the reason for this is that the even subcategory has a non-trivial retractable ideal, and more retractable IN-pairs.
In Section~\ref{subsect:QB} we define the above-mentioned planar and annular reducts of~$\B$, which are then studied in Sections~\ref{sect:TL} and~\ref{sect:J}.

During Sections \ref{sect:B}--\ref{sect:J}, we say an equivalence $\ve$ on a finite set $A$ is a \emph{2-equivalence} if one of the following holds:
\bit
\item $|A|$ is even, and every $\ve$-class has size 2, or
\item $|A|$ is odd, and all but one $\ve$-class has size 2, with the other $\ve$-class being a singleton.
\eit
We write $\Eq_2(A)$ for the set of all such 2-equivalences on $A$.

\subsection{Definitions and preliminaries on $\B$}\label{subsect:defnB}

Throughout Section \ref{sect:B} we fix a non-empty set $\C$ of finite sets and  keep the conventions of Section~\ref{sect:P}.  This time it will also be convenient to define
\[
\C_\even=\set{A\in\C}{|A|\text{ is even}} \AND \C_\odd=\set{A\in\C}{|A|\text{ is odd}}.
\]
It is possible for one of these subsets of $\C$ to be empty.

A partition $\al$ from $\P=\P(\C)$ is a \emph{Brauer partition} if every block of $\al$ has size $2$.  The set $\B=\B(\C)$ of all Brauer partitions is a subcategory of $\P$, called the \emph{Brauer category (over~$\C$)}.  The endomorphism monoids $\B_A=\B_{A,A}$ are \emph{Brauer monoids} \cite{Maz1998}, and we note that again $\S_A\sub\B_A$ for each $A\in\C$.  Clearly $\B_{A,B}$ is empty if $|A|$ and $|B|$ are of opposite parities, so it follows that
\begin{equation}\label{eq:EO}
\B = \B(\C) = \B(\C_\even) \sqcup \B(\C_\odd),
\end{equation}
writing $\sqcup$ for disjoint union.  This decomposition has many obvious structural implications, but for our purposes, the most important are to do with ideals and congruences.  

An arbitrary ideal of $\B=\B(\C)$ is a disjoint union $I=U\sqcup V\neq\emptyset$ of
 an ideal $U$ of $\C_\even$ (or the empty set) and an ideal $V$ of $\C_\odd$ (or the empty set).
 In particular, the ideals of $\B$ do not form a chain, unless one of $\C_\even$ or $\C_\odd$ is empty.
Further, for $I=U\sqcup V$ we have
\[
\Cong(I) = \bigset{\si_1\sqcup\si_2}{\si_1\in\Cong(U),\ \si_2\in\Cong(V)},
\]
which means that when $U$ and $V$ are both non-empty, $\Cong(I)$ is isomorphic to the lattice direct product $\Cong(U)\times\Cong(V)$.  Thus, it suffices to describe the congruences, and the corresponding lattices, of ideals in the two subcategories $\B(\C_\even)$ and $\B(\C_\odd)$.  We do this in Sections \ref{subsect:Bodd} and~\ref{subsect:Beven}; see Theorems \ref{thm:Bodd} and~\ref{thm:Beven}.  But first we gather some properties relevant to both subcategories.

\subsection{Green's relations and multiplicative properties in $\B$}\label{subsect:prelimB}

Clearly $\B$ is closed under the involution $\al\mt\al^*$, so $\B$ is again a regular $*$-category.  A consequence of regularity is that Green's $\R$, $\L$ and $\H$ relations are simply inherited from those of~$\P$ (cf.~Lemma \ref{lem:Green_P}); again it is easy to see that this is the case for the $\D=\J$ relation as well.  Note also that the kernel of $\al\in\B$ uniquely determines the domain of $\al$ as well as $\bd(\al)$; indeed, $\dom(\al)$ is the union of the singleton $\ker(\al)$-classes, while $\bd(\al)$ is the union of all $\ker(\al)$-classes.  This (and the dual statement concerning codomains and cokernels) will allow us to simplify the following statement, for which we also define 
\begin{align*}
\xi_\even=\xi_\even(\C)={\LSUB}\bigset{|A|}{A\in\C_\even} \ANd \xi_\odd=\xi_\odd(\C)={\LSUB}\bigset{|A|}{A\in\C_\odd}.
\end{align*}
It is possible that one or both of $\xi_\even$ or $\xi_\odd$ could be $\om$, and that one (but not both) could be~$0$.

\begin{lemma}\label{lem:Green_B}
\ben
\item \label{GB1} For $\al,\be\in\B$, we have
\bit
\item $(\al,\be)\in\R \iff \ker(\al)=\ker(\be)$,
\item $(\al,\be)\in\L \iff \coker(\al)=\coker(\be)$,
\item $(\al,\be)\in\J\iff(\al,\be)\in\D\iff\rank(\al)=\rank(\be)$.
\eit
\item \label{GB2} The $\J=\D$-classes of $\B$ are the sets
\bit
\item $D_q = D_q(\B) = \set{\al\in\B}{\rank(\al)=q}$ for each even $0\leq q<\xi_\even$, and each odd $1\leq q<\xi_\odd$.
\eit
These are all regular and stable, and they form a pair of chains: $D_0<D_2<\cdots$ and $D_1<D_3<\cdots$.
\item \label{GB4} The $\H$-class of any idempotent from $D_q$ is isomorphic to $\S_q$, the symmetric group of degree~$q$.  \epfres
\een
\end{lemma}

\begin{rem}\label{rem:Green_B}
As in \cite{Maz1998} (cf.~\cite[Proposition 2.1]{DE4}), two Brauer elements $\al,\be\in\B_{A,B}$ are $\R$-related if and only if $\al=\be\pi$ for some permutation $\pi\in\S_B$.  Similar statements hold for the $\L$ and $\D=\J$ relations.
\end{rem}

As noted in Section \ref{subsect:defnB}, the ideals of $\B$ are essentially unions of ideals from $\B(\C_\even)$ and $\B(\C_\odd)$.  
It follows from Lemma~\ref{lem:Green_B}\ref{GB2} that the ideals of these subcategories are the sets
\bit
\item $I_r = D_0\cup D_2\cup\cdots\cup D_r = \set{\al\in\B}{\rank(\al)\leq r, \ \rank(\al)\text{ is even}}$ for even $0\leq r<\xi_\even$,
\item $I_r = D_1\cup D_3\cup\cdots\cup D_r = \set{\al\in\B}{\rank(\al)\leq r, \ \rank(\al)\text{ is odd}}$ for odd $1\leq r<\xi_\odd$,
\eit
as well as $\B(\C_\even)$ and $\B(\C_\odd)$ themselves, which are not included in the above lists if $\xi_\even=\om$ or $\xi_\odd=\om$. It follows that both   $\B(\C_\even)$ and $\B(\C_\odd)$ are chains of ideals, and therefore amenable to an application of our general theory from Section \ref{sect:chains}.

The next technical lemma is analogous to Lemma \ref{lem:P2P}, and will be used to establish the required multiplication/separation properties on various ideals of $\B$.  It requires its own proof, as the partitions constructed during the proof of Lemma \ref{lem:P2P} generally do not belong to $\B$.

\begin{lemma}\label{lem:P2B}
Let $A,B\in\C$, and suppose $\al,\be\in\B_{A,B}$ with $r=\rank(\al)\geq\rank(\be)=q$.
\ben
\item \label{P2Bi} If $r\geq3$ and $q<r$, then there exists $\ga\in I_r$ such that $\ga\al\in D_{r-2}$ and $\ga\be\in I_{r-2}\sm H_{\ga\al}$.
\item \label{P2Bii} If $q=r\geq2$ and $(\al,\be)\not\in\H$, then there exists $\ga\in I_r$ such that 
\[
\text{$[\al\ga\in D_r$ and $\be\ga\in I_{r-2}]$ \ \ \ or \ \ \ $[\ga\al\in D_r$ and $\ga\be\in I_{r-2}]$.}
\]
\item \label{P2Biii} If $r\geq3$ and $\be\in H_\al\sm\{\al\}$, then there exists $\ga\in I_r$ such that $\ga\al\in D_{r-2}$ and $\ga\be\in I_{r-2}\sm H_{\ga\al}$.
\een
\end{lemma}

\pf
Throughout the proof, we write $\al = \partI{a_1}{a_r}{C_1}{C_s}{b_1}{b_r}{D_1}{D_t}$.  

\pfitem{\ref{P2Bi}}  Since $q<r$, we may assume without loss of generality that $\{b_1',v'\}$ is a block of $\be$ for some ${v\in B}$.  We then let $\ga\in\B_A$ be any Brauer partition with codomain $\{a_1,\ldots,a_{r-2}\}$.  We clearly have $\ga\al\in D_{r-2}$ and $\ga\be\in I_{r-2}$.  Because $b_1\in\codom(\ga\al)\sm\codom(\ga\be)$, we also have $(\ga\al,\ga\be)\not\in\L$, so that $\ga\be\not\in H_{\ga\al}$.

\pfitem{\ref{P2Bii}}  We assume that $(\al,\be)\not\in\R$ (the case of $(\al,\be)\not\in\L$ being dual), so that $\ker(\al)\not=\ker(\be)$.  

\pfcase{1}  Suppose first that $\dom(\al)\not=\dom(\be)$.  Without loss of generality, we may assume that $\{a_1,v\}$ is a block of $\be$ for some ${v\in A}$.  We then let $\ga\in D_r$ be any Brauer partition with transversals $\{a_1,a_1'\},\ldots,\{a_{r-1},a_{r-1}'\},\{v,v'\}$.  Then
\[
\{a_1,\ldots,a_{r-1}\}\sub\dom(\ga\al) \implies \rank(\ga\al)\geq r-1 \implies \rank(\ga\al)=r,
\]
by parity, so $\ga\al\in D_r$.  On the other hand, $\dom(\ga\be)\sub\{a_1,\ldots,a_{r-1},v\}$ and $a_1,v\not\in\dom(\ga\be)$, so $\ga\be\in I_{r-2}$.

\pfcase{2}  Now suppose $\dom(\al)=\dom(\be)$.  Fix some $(u,v)\in\ker(\be)\sm\ker(\al)$.  Then $u$ and $v$ belong to distinct non-transversals of $\al$, say $\{u,x\}$ and $\{v,y\}$.  We then let $\ga\in D_r$ be any Brauer partition with codomain $\{a_1,\ldots,a_{r-2},u,v\}$ and the blocks $\{x',a_{r-1}'\}$ and $\{y',a_r'\}$.  It is again easy to check that $\ga\al\in D_r$ and $\ga\be\in I_{r-2}$.

\pfitem{\ref{P2Biii}}  Here we have $\be = \partI{a_1}{a_r}{C_1}{C_s}{b_{1\pi}}{b_{r\pi}}{D_1}{D_t}$ for some non-trivial permutation $\pi\in\S_r$.  Without loss of generality we may assume that $1\pi=r$.  Let $\ga\in\B_A$ be any Brauer diagram with codomain $\{a_1,\ldots,a_{r-2}\}$.  Then  $\codom(\ga\al)=\{b_1,\ldots,b_{r-2}\} \not= \{b_{1\pi},\ldots,b_{(r-2)\pi}\}=\codom(\ga\be)$, so we have $\ga\al\in D_{r-2}$ and $(\ga\al,\ga\be)\not\in\L$.
\epf

Following the development from Section \ref{subsect:P}, 
for each even $0\leq q<\xi_\even$, and for each odd $1\leq q<\xi_\odd$, we fix some set $X_q\in\C$ with $n_q=|X_q|\geq q$ and $n_q\equiv q\Mod 2$, and without loss of generality assume that $X_q=[n_q]=\{1,\ldots,n_q\}$.  For a permutation $\pi\in\S_q$ we will write
\[
\pi^\natural = \partI 1q{q+1,q+2}{n_q-1,n_q}{1\pi}{q\pi}{q+1,q+2}{n_q-1,n_q}\in\B_{X_q}.
\]
For $\Om\sub\S_q$ we write $\Om^\natural=\set{\pi^\natural}{\pi\in\Om}$.  Then $\S_q^\natural\cong\S_q$ is a group $\H$-class of $\B$ contained in $D_q$.  
As in Section \ref{subsect:P}, and using the same $\phi(\al,\be)$ notation, we have
\[
\nu_{N^\natural} = \bigset{(\al,\be)\in\H\restr_{D_q}}{\phi(\al,\be)\in N}
\]
for any normal subgroup $N\normal\S_q$.

We now split considerations into the odd and even parts of $\B$.

\subsection{Congruences on ideals of $\B$: the odd case}\label{subsect:Bodd}

Throughout Section \ref{subsect:Bodd}, we assume that every member of $\C$ has odd cardinality, and we define $\xi={\LSUB}\bigset{|A|}{A\in\C}$.  

As we have already noted, it follows from Lemma \ref{lem:Green_B} that the ideals of $\B=\B(\C)$ are the sets
\bit
\item $I_r=I_r(\B)=\set{\al\in\B}{\rank(\al)\leq r}$ for each odd $1\leq r<\xi$, and
\item $I_\om=I_\om(\B)=\B$ in the case that $\xi=\om$.
\eit

An ideal $I_r$ ($r\geq3$) contains the IN-pairs $(I_q,N^\natural)$ for odd $1\leq q<r$ and $N\normal\S_{q+2}$.
These yield the congruences $R_{I_q,N^\natural}^{I_r}=R_{I_q}^{I_r}\cup\nu_{N^\natural}$, including the Rees congruences $R_{I_q}^{I_r}$ when $N$ is the trivial subgroup.  Again, we will often use the abbreviation $R_{I_q,N^\natural}^{I_r}=R_N^{I_r}$.
Associated to the minimal ideal $I_1$, we have the congruences $R_{I_1}^{I_r}$, $\rho_{I_1}^{I_r}$, $\lam_{I_1}^{I_r}$ and $\mu_{I_1}^{I_r}=\De_{I_r}$.  These four congruences, together with the $R_N^{I_r}$ and the full congruence $\nab_{I_r}$, form the lattice shown in Figure \ref{fig:Bodd}, and we will see soon that these are \emph{all} the congruences on $I_r$ ($r\geq3$).

\begin{figure}[ht]
\begin{center}
\scalebox{0.8}{
\begin{tikzpicture}[scale=.9]
\draw[rounded corners=2em,line width=3.5em,blue!25,cap=round] (0,10)--(0,14);
\draw[rounded corners=2em,line width=3.5em,green!25,cap=round] (0,16)--(0,16.01);
\draw[line width=3.5em,green!25] (0,16)--(0,17);
\dottedintervalx0{16}11{.43}
\node[rounded corners,rectangle,draw,fill=black!30,ultra thick] (N) at (0,19) {$R_{I_r}^{I_r}$};
   \draw(1.15,19) node {$=\nab_{I_r}$};
\node[rounded corners,rectangle,draw,fill=black!30,ultra thick] (R3) at (0,16) {$R_{I_3}^{I_r}$};
\node[rounded corners,rectangle,draw,fill=white] (S3) at (0,14) {$R_{\S_3}^{I_r}$};
\node[rounded corners,rectangle,draw,fill=white] (A3) at (0,12) {$R_{\A_3}^{I_r}$};
\node[rounded corners,rectangle,draw,fill=black!30,ultra thick] (R2) at (0,10) {$R_{I_1}^{I_r}$};
\node[rounded corners,rectangle,draw,fill=white] (kr) at (-3,8) {$\lam_{I_1}^{I_r}$};
\node[rounded corners,rectangle,draw,fill=black!30,ultra thick] (kl) at (0,6) {$\mu_{I_1}^{I_r}$};
\node[rounded corners,rectangle,draw,fill=white] (mrl) at (3,8) {$\rho_{I_1}^{I_r}$};
   \draw(1.15,5.95) node {$=\De_{I_r}$};
\draw(R2)--(kr)--(kl)--(mrl)--(R2)--(A3)--(S3)--(R3);
\end{tikzpicture}
}
\caption{
Hasse diagram of $\Cong(I_r)$ for $r\geq3$, where $I_r=I_r(\B)$ in the case that all members of~$\C$ have odd cardinality.  
See Theorem \ref{thm:Bodd}.
}
\label{fig:Bodd}
\end{center}
\end{figure}

We begin with two lemmas that will help us deal with congruences on the ideal $I_3$.

\begin{lemma}\label{lem:R_trick_Bodd}
For any $\al,\be\in I_3=I_3(\B)$ with $\al\sim\be$ and $(\al,\be)\not\in\L$, we have $\rho_{I_1}^{I_3}\sub\cg\al\be_{I_3}$.
\end{lemma}

\pf
Write $\si=\cg\al\be_{I_3}$, and suppose $\al,\be\in\B_{A,B}$.  Without loss of generality, we may fix some $(u,v)\in\coker(\be)\sm\coker(\al)$.  

First we claim that there exists $(\al',\be')\in\si$ such that $(u,v)\in\coker(\be')$, and at least one of $u',v'$ belongs to a lower non-transversal of $\al'$ distinct from $\{u',v'\}$.  
Indeed, if one of $u,v$ does not belong to $\codom(\al)$, then we just take $(\al,\be)=(\al',\be')$.  So suppose instead that $u,v\in\codom(\al)$.  Let the transversal of $\al$ containing $u'$ be $\{x,u'\}$.  We then let $\ga\in\B_A$ be any Brauer partition with $\codom(\ga)=\{x\}$, and note that $(\al',\be')=(\ga\al,\ga\be)$ satisfies the desired property.

With the claim established, we may assume without loss of generality that $(\al,\be)$ is already such a pair, and that $\{v',w'\}$ is a block of $\al$ for some $w\in B\sm\{u,v\}$.  In fact, premultiplying this pair by an arbitrary element of $D_1(\B_A)$, we may assume that $\al,\be\in D_1$.

Now let $(\ga,\de)\in\rho_{I_1}^{I_3}$ be arbitrary, say with $\ga,\de\in\B_{C,D}$.  We must show that $(\ga,\de)\in\si$.  This is clear if $\ga=\de$, so suppose instead that $\ga,\de\in D_1$ with $(\ga,\de)\in\R$.  Then $\de=\ga\pi$ for some $\pi\in\S_D$; cf.~Remark \ref{rem:Green_B}.  (Note that $\pi$ itself may not belong to the ideal $I_3$.)  Since $\pi$ is a product of transpositions, it suffices inductively to assume that $\pi$ is a single transposition, say~$(x,y)$, so $\de=\ga(x,y)$.  Write $\ga=\partIV{c}{C_1}\cdots{C_k}{d}{D_1}\cdots{D_l}$.

\pfcase{1}  Suppose first that $d$ is one of $x,y$, say $d=x$.  Without loss of generality, we may assume that $D_1=\{y,z\}$ for some $z$.  Let $\ve_1\in\B_{C,A}$ and $\ve_2\in\B_{B,D}$ be any Brauer partitions such that $\ker(\ve_1)=\ker(\ga)$, and $\ve_2$ contains the blocks $\{u,x'\}$, $\{v,z'\}$, $\{w,y'\}$ and $D_2',\ldots,D_l'$.  Then $(\ga,\de)=(\ve_1\al\ve_2,\ve_1\be\ve_2)\in\si$.

\pfcase{2}  Now suppose $d\not\in\{x,y\}$.  Without loss of generality, we may assume that $x\in D_1$.  If also $y\in D_1$, then $\de=\ga$, so there is nothing to prove.  Thus, we may assume that $D_1=\{x,s\}$ and $D_2=\{y,t\}$ for some $s,t$.  By Case 1, we have
\[
\big(\ga,\ga(d,x)\big)\in\si \COMMA 
\big(\ga(d,x),\ga(d,x)(x,y) \big)\in\si \COMMA 
\big(\ga(d,x)(x,y),\ga(d,x)(x,y)(y,d) \big)\in\si.
\]
(Informally speaking, note that we are using the result of Case 1 to repeatedly ``move around'' the bottom part of the transversal until we have moved it back to its original position, and in the process switched $x,y$ as desired.)  By transitivity, it follows that $\big(\ga,\ga(d,x)(x,y)(y,d) \big)\in\si$, and we are done, since $\ga(d,x)(x,y)(y,d)=\de$.
\epf

\begin{lemma}\label{lem:N_Bodd}
For every $\pi\in\S_3\sm\{\id_3\}$ there exists $\al\in D_1$ such that $\al\pi^\natural\not=\al\id_3^\natural$ and $\pi^\natural\al\not=\id_3^\natural\al$, with all four products defined.
\end{lemma}

\pf
For any $\pi$ except the transposition $(2,3)$, we take $\al=\partIV1{2,3}\cdots{n_3-1,n_3}1{2,3}\cdots{n_3-1,n_3}$.  For $\pi=(2,3)$ we take $\al=\partVI3{1,2}{4,5}{n_3-1,n_3}3{1,2}{4,5}{n_3-1,n_3}$.
\epf

We are now ready to prove the main result of this section.  For the first part, note that the $\R$- and $\L$-classes in the minimal ideal $I_1=D_1$ are indexed by the possible (co)kernels, each of which is a 2-equivalence on some $A\in\C$.  Since $|A|$ is odd for $A\in\C$, each such 2-equivalence has a unique singleton block.

\begin{thm}\label{thm:Bodd}
Let $\C$ be a non-empty set of finite sets of odd cardinality, let $\B=\B(\C)$ be the Brauer category over $\C$, and keep the above notation.
\ben
\item \label{cBodd1}  The minimal ideal $I_1$ is a partial rectangular band, so its congruence lattice is described by Proposition \ref{prop:PRB}.  In particular, $\Cong(I_1)$ is isomorphic to the lattice direct product
\[
\prod_{A\in\C}\Eq(\Eq_2(A))\times\prod_{A\in\C}\Eq(\Eq_2(A)).
\]
\item \label{cBodd2} For odd $r\geq3$ (including $r=\om$ if $\xi=\om$), 
\[
\Cong(I_r) = 
\{\De_{I_r},\ \lam_{I_1}^{I_r},\ \rho_{I_1}^{I_r},\ R_{I_1}^{I_r}\} 
\cup\set{R_{I_q,N^\natural}^{I_r}}{1\leq q<r,\ q\text{ odd},\ N\normal\S_{q+2}}\cup\{\nab_{I_r}\}
\]
forms the lattice pictured in Figure \ref{fig:Bodd}.
\een
\end{thm}

\pf
\firstpfitem{\ref{cBodd1}}  This follows from Proposition \ref{prop:PRB}, and the above observations about 2-equivalences.

\pfitem{\ref{cBodd2}}  As usual, we just need to prove the result for $r=3$; Theorem \ref{thm:E2} and Lemma \ref{lem:P2B} will take care of the rest.  For $r=3$ we must show that
\[
\Cong(I_3) = \{\De_{I_3},\ \lam_{I_1}^{I_3},\ \rho_{I_1}^{I_3},\ R_{I_1}^{I_3},\ R_{I_1,\A_3^\natural}^{I_3},\ R_{I_1,\S_3^\natural}^{I_3},\ \nab_{I_3}\},
\]
and this in fact follows from Proposition \ref{prop:small01}.  The technical assumptions of that proposition follow from Lemmas \ref{lem:P2B}\ref{P2Bii}, \ref{lem:R_trick_Bodd} (and its dual), and \ref{lem:N_Bodd}.
\epf

\subsection{Congruences on ideals of $\B$: the even case}\label{subsect:Beven}

Throughout Section \ref{subsect:Beven}, we assume that each member of $\C$ has even cardinality, and we define $\xi={\LSUB}\bigset{|A|}{A\in\C}$.  This time the ideals of $\B=\B(\C)$ are the sets
\bit
\item $I_r=I_r(\B)=\set{\al\in\B}{\rank(\al)\leq r}$ for each even $0\leq r<\xi$, and
\item $I_\om=I_\om(\B)=\B$ in the case that $\xi=\om$.
\eit

In contrast with the odd case considered in Section \ref{subsect:Bodd}, this time the \emph{two} smallest ideals,~$I_0$ and $I_2$, of $\B$ will have different behaviour to the larger ideals; compare the first part of Theorem~\ref{thm:Bodd} with the first two parts of Theorem \ref{thm:Beven}, and see also Theorem \ref{thm:P}.  The reason for this is that, as in \cite[Lemma~8.2]{EMRT2018}, we have a retraction
\begin{equation}\label{eq:retract_B}
f:I_2\to I_0:  \al\mt\wh\al\qquad\text{where}\qquad
\wh\al=\begin{cases}
\partIII{a_1,a_2}{C_1}\cdots{C_k}{b_1,b_2}{D_1}\cdots{D_l} & \text{if } \al=\partV{a_1}{a_2}{C_1}{C_k}{b_1}{b_2}{D_1}{D_l}\in D_2\\
\al & \text{if } \al\in D_0.
\end{cases}
\end{equation}

An ideal $I_r$ ($r\geq4$) contains the IN-pairs $(I_q,N^\natural)$ for even $0\leq q<r$ and $N\normal\S_{q+2}$.
These yield the congruences
$R_N^{I_r} = R_{I_q,N^\natural}^{I_r}=R_{I_q}^{I_r}\cup\nu_{N^\natural}$,
including the Rees congruences $R_{I_q}^{I_r}$ for $N=\{\id_{q+2}\}$.

As in \cite[Lemma 8.3]{EMRT2018}, we also have the following \emph{four} retractable IN-pairs:
\[
(I_0,\{\id_2^\natural\})\COMMA
(I_0,\S_2^\natural)\COMMA
(I_2,\{\id_4^\natural\})\COMMA
(I_2,K^\natural),
\]
where
$K=\big\{{\id}_4,(1,2)(3,4),(1,3)(2,4),(1,4)(2,3)\big\}\normal \S_4$ is the \emph{Klein 4-group}.
Each of these retractable IN-pairs yields a family of $R/\lam/\rho/\mu$ congruences, 
Following the same convention as in Section \ref{subsect:P}, in the notation for these congruences we will
omit $\natural$ and also either
the normal subgroup if it is trivial, or else the ideal $I_0$ or $I_2$.

The lattice formed by the above congruences together with $\nab_{I_r}$
 is pictured in Figure~\ref{fig:Beven}.  Again, our goal is to show that these are \emph{all} the congruences on $I_r$ ($r\geq4$).
To achieve this goal, we first need two technical lemmas, both to help with the $r=4$ case.
The first refers to the relation $\L f^{-1}$, where $f$ is the retraction~\eqref{eq:retract_B}, so
\[
\L f^{-1} = \bigset{(\al,\be)\in I_2\times I_2}{(\wh\al,\wh\be)\in\L}.
\]

\begin{figure}[ht]
\begin{center}
\scalebox{0.78}{
\begin{tikzpicture}[scale=.95]
\draw[rounded corners=2em,line width=3.5em,blue!15,cap=round] (0,8)--(0,12);
\foreach \x/\y in {0/0,0/2,3/2,3/4} {\draw[rounded corners=2em,line width=3.5em,red!25,cap=round] (6+\x,-2+\y)--(3+\x,0+\y);}
\foreach \x/\y in {0/0,0/2,3/2,3/4} {\draw[rounded corners=2em,line width=3.5em,blue!35,cap=round] (0+\x,2+\y)--(-3+\x,4+\y);}
\draw[rounded corners=2em,line width=3.5em,green!25,cap=round] (0,14)--(0,14.01);
\draw[line width=3.5em,green!25] (0,14)--(0,15);
\dottedintervalx0{14}11{.43}
\node[rounded corners,rectangle,draw,fill=black!30,ultra thick] (N) at (0,17) {$R_{I_r}^{I_r}$};
  \draw(1.13,17.0) node {$=\nab_{I_r}$};
\node[rounded corners,rectangle,draw,fill=black!30,ultra thick] (R4) at (0,14) {$R_{I_4}^{I_r}$};
\node[rounded corners,rectangle,draw,fill=white] (S4) at (0,12) {$R_{\S_4}^{I_r}$};
\node[rounded corners,rectangle,draw,fill=white] (A4) at (0,10) {$R_{\A_4}^{I_r}$};
\node[rounded corners,rectangle,draw,fill=white] (krl) at (0,8) {$R_K^{I_r}$};
\node[rounded corners,rectangle,draw,fill=white] (kr) at (-3,6) {$\lambda_K^{I_r}$};
\node[rounded corners,rectangle,draw,fill=white] (kl) at (0,6) {$\rho_K^{I_r}$};
\node[rounded corners,rectangle,draw,fill=white] (k) at (-3,4) {$\mu_K^{I_r}$};
\node[rounded corners,rectangle,draw,fill=black!30,ultra thick] (erl) at (3,6) {$R_{I_2}^{I_r}$};
\node[rounded corners,rectangle,draw,fill=white] (er) at (0,4) {$\lambda_{I_2}^{I_r}$};
\node[rounded corners,rectangle,draw,fill=white] (el) at (3,4) {$\rho_{I_2}^{I_r}$};
\node[rounded corners,rectangle,draw,fill=white] (e) at (0,2) {$\mu_{I_2}^{I_r}$};
\node[rounded corners,rectangle,draw,fill=white] (nrl) at (6,4) {$R_{\S_2}^{I_r}$};
\node[rounded corners,rectangle,draw,fill=white] (nr) at (3,2) {$\lambda_{\S_2}^{I_r}$};
\node[rounded corners,rectangle,draw,fill=white] (nl) at (6,2) {$\rho_{\S_2}^{I_r}$};
\node[rounded corners,rectangle,draw,fill=white] (n) at (3,0) {$\mu_{\S_2}^{I_r}$};
\node[rounded corners,rectangle,draw,fill=black!30,ultra thick] (rl) at (9,2) {$R_{I_0}^{I_r}$};
\node[rounded corners,rectangle,draw,fill=white] (r) at (6,0) {$\lambda_{I_0}^{I_r}$};
\node[rounded corners,rectangle,draw,fill=white] (l) at (9,0) {$\rho_{I_0}^{I_r}$};
\node[rounded corners,rectangle,draw,fill=black!30,ultra thick] (D) at (6,-2) {$\mu_{I_0}^{I_r}$};
\draw(7.3,-2.05) node {$=\Delta_{I_r}$};
\draw
(D)--(n)--(e)--(k)
(l)--(nl)--(el)--(kl)
(k)--(kl) (e)--(el) (n)--(nl) (D)--(l) 
(k)--(kr) (e)--(er) (n)--(nr) (D)--(r) 
(l)--(rl) (nl)--(nrl) (el)--(erl) (kl)--(krl)
(rl)--(nrl)--(erl)--(krl)--(A4)--(S4)--(R4)
;
\fill[blue!35] (-1.5,5)circle(.15);
\fill[white] (1.5,3)circle(.15);
\fill[red!25] (4.5,1)circle(.15);
\fill[blue!35] (1.5,5)circle(.15);
\fill[white] (4.5,3)circle(.15);
\fill[red!25] (7.5,1)circle(.15);
\draw
(r)--(nr)--(er)--(kr)
(r)--(rl) (nr)--(nrl) (er)--(erl) (kr)--(krl) 
;
\end{tikzpicture}
}
\end{center}
\vspace{-7mm}
\caption{Hasse diagram of $\Cong(I_r)$ for $r\geq4$, where $I_r=I_r(\B)$ in the case that all members of~$\C$ have even cardinality.  
See Theorem \ref{thm:Beven}.}
\label{fig:Beven}
\end{figure}

\begin{lemma}\label{lem:R_trick_Beven}
For any $\al,\be\in I_4=I_4(\B)$ with $\al\sim\be$ and $(\al,\be)\not\in\L\cup\L f^{-1}$, we have ${\rho_{I_0}^{I_4}\sub\cg\al\be_{I_4}}$.
\end{lemma}

\pf
Write $\si=\cg\al\be_{I_4}$, and suppose $\al,\be\in\B_{A,B}$.  We consider two cases.

\pfcase{1}
Suppose first that there exist distinct $u,v,w\in B$ such that $\al$ and $\be$ contain the blocks $\{u',v'\}$ and $\{v',w'\}$, respectively.  

Suppose $\ga=\partII{C_1}\cdots{C_k}{D_1}\cdots{D_l}\in D_0$, say with $\ga\in\B_{C,D}$.  As in the proof of Lemma \ref{lem:R_trick_Bodd}, it suffices to show that $(\ga,\ga\pi)\in\si$ for any transposition $\pi\in\S_D$, say $\pi=(x,y)$.  If $x',y'$ belong to the same block of $\ga$, then $\ga=\ga\pi$, so there is nothing to show.  Thus, we assume that $D_1=\{x,s\}$ and $D_2=\{y,t\}$ for some $s,t$.  Let $\ve_1\in\B_{C,A}$ and $\ve_2\in\B_{B,D}$ be any Brauer partitions such that $\ker(\ve_1)=\ker(\ga)$, and $\ve_2$ contains the blocks $\{u,x'\}$, $\{v,s'\}$, $\{w,y'\}$ and $D_3,\ldots,D_l$.  Then $(\ga,\ga\pi)=(\ve_1\al\ve_2,\ve_1\be\ve_2)\in\si$.

\pfcase{2}
Now suppose we are not in Case 1.  This means that:
\begin{align}
\nonumber &\text{for every lower non-transversal $\{u',v'\}$ of $\al$, }\\
\label{eq:R_trick_Beven} &\text{either $\{u',v'\}$ is also a block of $\be$ or else ${u,v\in\codom(\be)}$,}
\end{align} 
with a symmetrical statement for the lower non-transversals of $\be$.

If we had $\al,\be\in I_2$, then \eqref{eq:R_trick_Beven} would imply $\coker(\wh\al)=\coker(\wh\be)$, and so $(\wh\al,\wh\be)\in\L$, contrary to our assumption that $(\al,\be)\not\in\L f^{-1}$.  So without loss of generality, we assume that~${\al\in D_4}$.  Let the transversals of $\al$ be $\{a_i,b_i'\}$, $i=1,2,3,4$.  If $\codom(\al)=\codom(\be)$, then~\eqref{eq:R_trick_Beven} would give $\coker(\al)=\coker(\be)$, whence $(\al,\be)\in\L$, a contradiction.  So we may assume that ${b_1\not\in\codom(\be)}$; together with \eqref{eq:R_trick_Beven}, we may further assume that~$\be$ contains the block $\{b_1',b_2'\}$.  Choose any $\ga\in D_0(\B_A)$ containing the block $\{a_2',a_3'\}$, and put $(\al',\be')=(\ga\al,\ga\be)\in\si$.  Then $\{b_1',b_2'\}$ and $\{b_2',b_3'\}$ are blocks of $\be'$ and $\al'$, respectively, so we have reduced to Case 1.
\epf

\begin{lemma}\label{lem:N_Beven}
For every $\pi\in\S_4\sm K$ there exists $\al\in D_0$ such that $\al\pi^\natural\not=\al\id_4^\natural$ and $\pi^\natural\al\not=\id_4^\natural\al$, with all four products defined.
\end{lemma}

\pf
If $\pi$ is one of $(1,2)$, $(3,4)$, $(1,3,2,4)$ or $(1,4,2,3)$,
we take $\al=\partVII{1,4}{2,3}{4,5}{n_4-1,4_3}{1,4}{2,3}{4,5}{n_4-1,4_3}$.
Otherwise we take $\al=\partXXIV{1,2}\cdots{n_4-1,n_4}{1,2}\cdots{n_4-1,n_4}$.
\epf

We may now prove the main result of this section.  The first part again involves 2-equivalences.

\begin{thm}\label{thm:Beven}
Let $\C$ be a non-empty set of finite sets of even cardinality, let $\B=\B(\C)$ be the Brauer category over $\C$, and keep the above notation.
\ben
\item \label{cBeven1}  The minimal ideal $I_0$ is a partial rectangular band, so its congruence lattice is described by Proposition \ref{prop:PRB}.  In particular, $\Cong(I_0)$ is isomorphic to the lattice direct product
\[
\prod_{A\in\C}\Eq(\Eq_2(A))\times\prod_{A\in\C}\Eq(\Eq_2(A)).
\]
\item \label{cBeven2}  The ideal $I_2$ is retractable, every congruence on $I_0$ is liftable to~$I_2$, and
\[
\Cong(I_2)=\bigset{\tau\cup\nu_N^\natural}{\tau\in\Cong(I_0),\ N\normal\S_2} \cup \bigset{\th_{I_0,\tau}^{I_1}}{\tau\in\Cong(I_0)}.
\]
In particular, $\Cong(I_2)$ is isomorphic to the lattice direct product $\Cong(I_0)\times\bthree$.
\item \label{cBeven3} For even $r\geq4$ (including $r=\om$ if $\xi=\om$), 
\begin{align*}
\Cong(I_r) = 
\set{\mu_{I_0,N^\natural}^{I_r},\ \lam_{I_0,N^\natural}^{I_r},\ \rho_{I_0,N^\natural}^{I_r},\ &R_{I_0,N^\natural}^{I_r}}{N\normal\S_2} \\[3truemm]
{}\cup\set{\mu_{I_2,N^\natural}^{I_r},\ \lam_{I_2,N^\natural}^{I_r},\ \rho_{I_2,N^\natural}^{I_r},\ &R_{I_2,N^\natural}^{I_r}}{N=\{\id_4\}\text{ or }N=K} \\[3truemm]
\qquad\qquad\qquad\qquad {}\cup\set{&R_{I_q,N^\natural}^{I_r}}{2\leq q<r,\ q\text{ even},\ N\normal\S_{q+2}}\cup\{\nab_{I_r}\}
\end{align*}
forms the lattice pictured in Figure \ref{fig:Beven}.
\een
\end{thm}

\pf
Parts \ref{cBeven1} and \ref{cBeven2} again follow from Propositions \ref{prop:PRB} and \ref{prop:small_ideal}.  This leaves us to prove~\ref{cBeven3}, and as usual it suffices (by Theorem \ref{thm:E2} and Lemma \ref{lem:P2B}) to consider the case $r=4$.  But this follows from Proposition \ref{prop:small012}, whose technical assumptions follow from Lemmas~\ref{lem:P2B},~\ref{lem:R_trick_Beven} (and its dual) and~\ref{lem:N_Beven}.
\epf

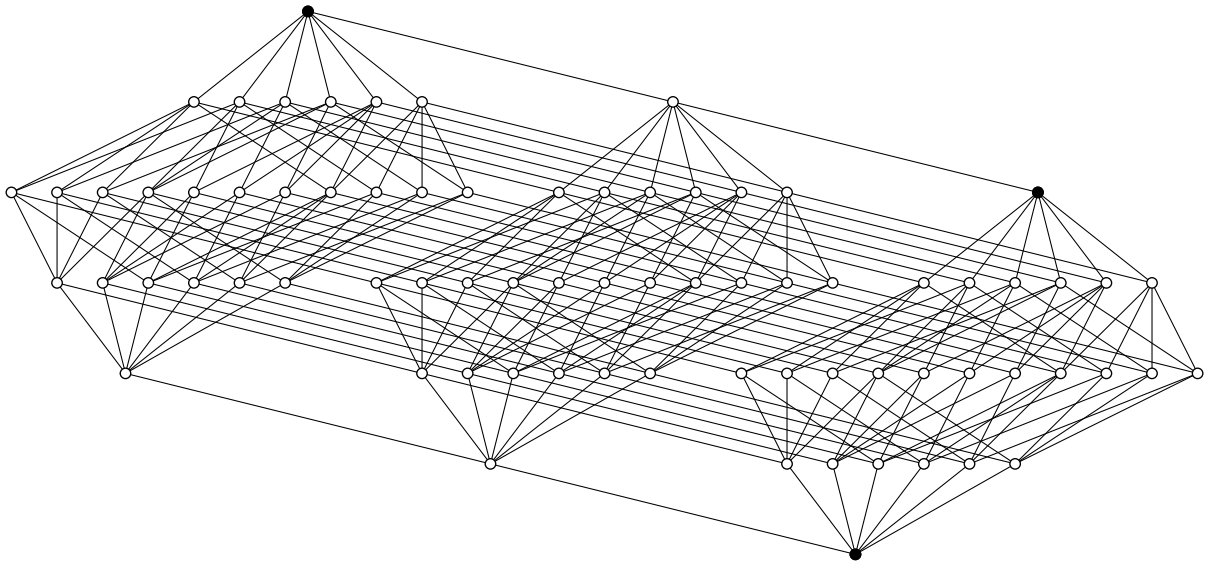
\begin{figure}[ht]
\begin{center}
\begin{tikzpicture}[scale=0.6,inner sep=0pt]
\begin{scope}[shift={(0,-15)}]
\foreach \L/\x/\y in {
a/2.5/0,
b/1/1,
c/2/1,
d/3/1,
e/4/1,
f/5/1,
g/6/1,
h/0/2,
i/1/2,
j/2/2,
k/3/2,
l/4/2,
m/5/2,
n/6/2,
o/7/2,
p/8/2,
q/9/2,
r/10/2,
s/4/3,
t/5/3,
u/6/3,
v/7/3,
w/8/3,
x/9/3,
y/6.5/4
} {\node (\L) at (\x,2*\y) {};}
\foreach \x in {a,b,c,d,e,f,g,h,i,j,k,l,m,n,o,p,q,r,s,t,u,v,w,x,y} {\draw(\x)-++(-16,4);}
\foreach \y in {b,c,d,e,f,g} {\draw(a)--(\y);}
\foreach \y in {h,i,j,k} {\draw(b)--(\y);}
\foreach \y in {k,l,m,n} {\draw(c)--(\y);}
\foreach \y in {h,l,o,p} {\draw(d)--(\y);}
\foreach \y in {i,m,o,q} {\draw(e)--(\y);}
\foreach \y in {j,n,o,r} {\draw(f)--(\y);}
\foreach \y in {k,p,q,r} {\draw(g)--(\y);}
\foreach \y in {s,t} {\draw(h)--(\y);}
\foreach \y in {s,u} {\draw(i)--(\y);}
\foreach \y in {s,v} {\draw(j)--(\y);}
\foreach \y in {t,u,v} {\draw(k)--(\y);}
\foreach \y in {t,w} {\draw(l)--(\y);}
\foreach \y in {u,w} {\draw(m)--(\y);}
\foreach \y in {v,w} {\draw(n)--(\y);}
\foreach \y in {s,w,x} {\draw(o)--(\y);}
\foreach \y in {t,x} {\draw(p)--(\y);}
\foreach \y in {u,x} {\draw(q)--(\y);}
\foreach \y in {v,x} {\draw(r)--(\y);}
\foreach \y in {s,t,u,v,w,x} {\draw(y)--(\y);}
\foreach \L/\x/\y in {
a/2.5/0,
b/1/1,
c/2/1,
d/3/1,
e/4/1,
f/5/1,
g/6/1,
h/0/2,
i/1/2,
j/2/2,
k/3/2,
l/4/2,
m/5/2,
n/6/2,
o/7/2,
p/8/2,
q/9/2,
r/10/2,
s/4/3,
t/5/3,
u/6/3,
v/7/3,
w/8/3,
x/9/3,
y/6.5/4
} {\colvert{\x}{2*\y}{white}}
\foreach \x/\y in {2.5/0,6.5/4} {\colvert{\x}{2*\y}{black}}
\end{scope}
\begin{scope}[shift={(-8,-13)}]
\foreach \L/\x/\y in {
a/2.5/0,
b/1/1,
c/2/1,
d/3/1,
e/4/1,
f/5/1,
g/6/1,
h/0/2,
i/1/2,
j/2/2,
k/3/2,
l/4/2,
m/5/2,
n/6/2,
o/7/2,
p/8/2,
q/9/2,
r/10/2,
s/4/3,
t/5/3,
u/6/3,
v/7/3,
w/8/3,
x/9/3,
y/6.5/4
} {\node (\L) at (\x,2*\y) {};}
\foreach \y in {b,c,d,e,f,g} {\draw(a)--(\y);}
\foreach \y in {h,i,j,k} {\draw(b)--(\y);}
\foreach \y in {k,l,m,n} {\draw(c)--(\y);}
\foreach \y in {h,l,o,p} {\draw(d)--(\y);}
\foreach \y in {i,m,o,q} {\draw(e)--(\y);}
\foreach \y in {j,n,o,r} {\draw(f)--(\y);}
\foreach \y in {k,p,q,r} {\draw(g)--(\y);}
\foreach \y in {s,t} {\draw(h)--(\y);}
\foreach \y in {s,u} {\draw(i)--(\y);}
\foreach \y in {s,v} {\draw(j)--(\y);}
\foreach \y in {t,u,v} {\draw(k)--(\y);}
\foreach \y in {t,w} {\draw(l)--(\y);}
\foreach \y in {u,w} {\draw(m)--(\y);}
\foreach \y in {v,w} {\draw(n)--(\y);}
\foreach \y in {s,w,x} {\draw(o)--(\y);}
\foreach \y in {t,x} {\draw(p)--(\y);}
\foreach \y in {u,x} {\draw(q)--(\y);}
\foreach \y in {v,x} {\draw(r)--(\y);}
\foreach \y in {s,t,u,v,w,x} {\draw(y)--(\y);}
\foreach \L/\x/\y in {
a/2.5/0,
b/1/1,
c/2/1,
d/3/1,
e/4/1,
f/5/1,
g/6/1,
h/0/2,
i/1/2,
j/2/2,
k/3/2,
l/4/2,
m/5/2,
n/6/2,
o/7/2,
p/8/2,
q/9/2,
r/10/2,
s/4/3,
t/5/3,
u/6/3,
v/7/3,
w/8/3,
x/9/3,
y/6.5/4
} {\colvert{\x}{2*\y}{white}}
\end{scope}
\begin{scope}[shift={(-16,-11)}]
\foreach \L/\x/\y in {
a/2.5/0,
b/1/1,
c/2/1,
d/3/1,
e/4/1,
f/5/1,
g/6/1,
h/0/2,
i/1/2,
j/2/2,
k/3/2,
l/4/2,
m/5/2,
n/6/2,
o/7/2,
p/8/2,
q/9/2,
r/10/2,
s/4/3,
t/5/3,
u/6/3,
v/7/3,
w/8/3,
x/9/3,
y/6.5/4
} {\node (\L) at (\x,2*\y) {};}
\foreach \y in {b,c,d,e,f,g} {\draw(a)--(\y);}
\foreach \y in {h,i,j,k} {\draw(b)--(\y);}
\foreach \y in {k,l,m,n} {\draw(c)--(\y);}
\foreach \y in {h,l,o,p} {\draw(d)--(\y);}
\foreach \y in {i,m,o,q} {\draw(e)--(\y);}
\foreach \y in {j,n,o,r} {\draw(f)--(\y);}
\foreach \y in {k,p,q,r} {\draw(g)--(\y);}
\foreach \y in {s,t} {\draw(h)--(\y);}
\foreach \y in {s,u} {\draw(i)--(\y);}
\foreach \y in {s,v} {\draw(j)--(\y);}
\foreach \y in {t,u,v} {\draw(k)--(\y);}
\foreach \y in {t,w} {\draw(l)--(\y);}
\foreach \y in {u,w} {\draw(m)--(\y);}
\foreach \y in {v,w} {\draw(n)--(\y);}
\foreach \y in {s,w,x} {\draw(o)--(\y);}
\foreach \y in {t,x} {\draw(p)--(\y);}
\foreach \y in {u,x} {\draw(q)--(\y);}
\foreach \y in {v,x} {\draw(r)--(\y);}
\foreach \y in {s,t,u,v,w,x} {\draw(y)--(\y);}
\foreach \L/\x/\y in {
a/2.5/0,
b/1/1,
c/2/1,
d/3/1,
e/4/1,
f/5/1,
g/6/1,
h/0/2,
i/1/2,
j/2/2,
k/3/2,
l/4/2,
m/5/2,
n/6/2,
o/7/2,
p/8/2,
q/9/2,
r/10/2,
s/4/3,
t/5/3,
u/6/3,
v/7/3,
w/8/3,
x/9/3,
y/6.5/4
} {\colvert{\x}{2*\y}{white}}
\foreach \x/\y in {6.5/4} {\colvert{\x}{2*\y}{black}}
\end{scope}
\end{tikzpicture}
\caption{Hasse diagram of $\Cong(I_2)$ in $\B_4$; cf.~Theorem~\ref{thm:Beven} and Remark~\ref{rem:Beven}.
}
\label{fig:B4}
\end{center}
\end{figure}

\begin{rem}\label{rem:Beven}
Figure \ref{fig:B4} pictures the lattice $\Cong(I_2)$ in the case that $\C$ consists of a single set of size $4$, produced with the aid of GAP \cite{GAP4,GAP}.  Here $\B=\B_4$ is just the Brauer monoid of degree~$4$.  By Theorem \ref{thm:Beven}\ref{cBeven2}, this lattice is essentially three copies of $\Cong(I_0)$.  By Theorem~\ref{thm:Beven}\ref{cBeven1}, the latter is isomorphic to $\Eq_3\times\Eq_3$, since there are precisely three 2-equivalences on a set of size $4$.
\end{rem}

\subsection{Planar and annular reducts of $\B$}\label{subsect:QB}

In Section~\ref{subsect:QP} we considered subcategories of $\P$ consisting of (anti-)planar and (anti-)annular partitions.  We will soon consider the corresponding subcategories of the Brauer category ${\B=\B(\C)}$, namely:
\bit
\item $\OrdP\B$, the category of planar Brauer partitions,
\item $\OrdPR\B$, the category of planar and anti-planar Brauer partitions,
\item $\OriP\B$, the category of annular Brauer partitions, 
\item $\OriPR\B$, the category of annular and anti-annular Brauer partitions.
\eit
The category $\OrdP\B$ is the \emph{Temperley--Lieb category}, which is usually denoted $\TL$.  We call the category $\OriP\B$ the \emph{Jones category}, and denote it $\JJ$, in light of the corresponding algebras named after Jones \cite{Jones2001,Jones1994_a}; cf.~\cite{GL1996,ADV2012_2}.  
(But note that the Temperley--Lieb monoids $\TL_n$ are often called Jones monoids and denoted $\JJ_n$; see for example \cite[Section 9]{EMRT2018}, and also \cite{LF2006,BDP2002} for a discussion of naming conventions.)
We call $\OrdPR\B$ and $\OriPR\B$ the \emph{anti-Temperley--Lieb} and \emph{anti-Jones categories}, respectively, and we will denote them by $\TLpm$ and~$\Jpm$.

Unlike the situation with $\P$ in Section \ref{subsect:QP}, it is not simply a matter of modifying the above arguments for $\B$ to immediately describe the congruences on the ideals of its planar and annular reducts.  
The reason for this is that the customary verification of the multiplication and separation properties requires a fairly fine geometrical analysis.
Thus, we treat the planar reducts in Section~\ref{sect:TL}, and the annular reducts in Section \ref{sect:J}.  Before we do, however, note that these subcategories are all closed under the $\al\mt\al^*$ map, so again Green's relations are simply inherited from~$\B$ (cf.~Lemma \ref{lem:Green_B}), although group $\H$-classes need to be described separately:

\begin{lemma}\label{lem:Green_QB}
Let $\Q$ be any of $\TL$, $\TLpm$, $\JJ$ or $\Jpm$.
\ben
\item \label{GQB1} For $\al,\be\in\Q$, we have
\bit
\item $(\al,\be)\in\R \iff \ker(\al)=\ker(\be)$,
\item $(\al,\be)\in\L \iff \coker(\al)=\coker(\be)$,
\item $(\al,\be)\in\J\iff(\al,\be)\in\D\iff\rank(\al)=\rank(\be)$.
\eit
\item \label{GQB2} The $\J=\D$-classes of $\Q$ are the sets
\bit
\item $D_q = D_q(\Q) = \set{\al\in\Q}{\rank(\al)=q}$ for each even $0\leq q<\xi_\even$, and each odd $1\leq q<\xi_\odd$.
\eit
These are all regular and stable, and they form a pair of chains: $D_0<D_2<\cdots$ and $D_1<D_3<\cdots$.
\item \label{GQB3} Apart from trivial exceptions for small $q$, the $\H$-class of any idempotent from $D_q$ is:
\bit
\item trivial when $\Q=\TL$,
\item cyclic of order $2$ when $\Q=\TLpm$,
\item cyclic of order $q$ when $\Q=\JJ$,
\item dihedral of order $2q$ when $\Q=\Jpm$.  \epfres 
\eit
\een
\end{lemma}

If $\Q$ is any of the above subcategories of $\B$, we define the sets
\bit
\item $I_r = I_r(\Q) = \set{\al\in\Q}{\rank(\al)\leq r, \ \rank(\al)\text{ is even}}$ for even $0\leq r<\xi_\even$,
\item $I_r = I_r(\Q) = \set{\al\in\Q}{\rank(\al)\leq r, \ \rank(\al)\text{ is odd}}$ for odd $1\leq r<\xi_\odd$.
\eit
Again, these may be used to describe all the ideals of $\Q$, as in Section \ref{subsect:defnB}.

\section{Temperley--Lieb and anti-Temperley--Lieb categories}\label{sect:TL}

Throughout Section \ref{sect:TL}, we fix a non-empty set $\C$ of finite sets, and denote by $\TL=\TL(\C)$ and $\TLpm=\TLpm(\C)$ the Temperley--Lieb and anti-Temperley--Lieb categories, as defined in Section~\ref{subsect:QB}.  We often let $\Q$ stand for either of these two categories.  
We keep the notational conventions of Sections~\ref{sect:P} and~\ref{sect:B}.
This includes the permutations $\ga_A\in\S_A$ ($A\in\C$), as in~\eqref{eq:gade}.

Just as was the case with the Brauer category $\B$, we have $\Q=\Q(\C_\even)\sqcup\Q(\C_\odd)$, and the ideals of $\Q$ similarly decompose into even and odd parts.  It is therefore again sufficient to
consider separate cases where the members of $\C$ all have even size, or all odd size.
In the odd case we treat $\TL$ and $\TLpm$ simultaneously (Section \ref{subsect:TLodd}), but in the even case we 
have to separate them (Sections \ref{subsect:TLeven} and \ref{subsect:TLpmeven}).

Since we will often be proving results for $\TL$ and $\TLpm$ at the same time, we will need to be careful that any partitions constructed during proofs belong to $\TL$.  As noted in \cite[Section~9]{EMRT2018}, the planarity constraint leads to some important points regarding parity.  Let $A\in\C$ and write $A=\{a_1<\cdots<a_k\}$.  By the \emph{parity} of $a_i$, we mean the parity of $i$; so $a_1$ is odd, $a_2$ is even, and so on.  

Consider some $\al\in\TL$, and write $\al=\partI{a_1}{a_r}{C_1}{C_s}{b_1}{b_r}{D_1}{D_t}$ with $a_1<\cdots<a_r$.  Then we also have $b_1<\cdots<b_r$, and moreover the parities of $a_1,a_2,a_3,\ldots$ and $b_1,b_2,b_3,\ldots$ alternate: odd, even, odd, and so on.  For any non-transversal $\{x,y\}$ or $\{x',y'\}$ of $\al$, one of $x,y$ is odd and the other even.  If $A,B\in\C$ are such that $|A|$ and $|B|$ have the same parity, and if we have some collection of pairwise-disjoint edges between vertices from $A\cup B'$ that may be drawn in planar fashion (within the rectangle determined by the vertices), and obey the above parity conditions, then it is always possible to extend these edges to a Temperley--Lieb partition from $\TL_{A,B}$.

If $\al=\partI{a_1}{a_r}{C_1}{C_s}{b_1}{b_r}{D_1}{D_t}\in\TL^-$ with $a_1<\cdots<a_r$, then we have $b_1>\cdots>b_r$.  The parities of $b_1,b_2,b_3,\ldots$ still alternate, but $b_1$ can be even or odd, as $b_1\equiv r\Mod2$.

\subsection{Multiplicative properties in $\TL$ and $\TLpm$}

Throughout this section, $\Q$ denotes either of $\TL$ or $\TLpm$.  Before we begin, we require the following lemma from \cite{EMRT2018}, concerning Temperley--Lieb \emph{monoids} $\TL_A$.  For the statement, and for later use, recall from \cite{NS1978} that an element $x$ of a regular $*$-semigroup is a \emph{projection} if $x^2=x=x^*$.  It is easy to see that the projections of $\TL$ are of the form $\partI{a_1}{a_q}{B_1}{B_r}{a_1}{a_q}{B_1}{B_r}$.  Also, since the product of two (anti-)planar partitions is planar, it quickly follows that every idempotent (and hence every projection) from $\TLpm$ belongs to $\TL$.

\begin{lemma}\label{lem:TL_proj}
Let $\al,\be$ be distinct projections from $\TL_A$ with $r=\rank(\al)=\rank(\be)\geq2$.  Then there exists $\ga\in D_r(\TL_A)$ such that, swapping $\al,\be$ if necessary, $\ga\al\in D_r$ and $\ga\be\in I_{r-2}$.
\end{lemma}

\pf
Actually, \cite[Lemma 9.19]{EMRT2018} demonstrates the existence of a projection $\ga\in D_r(\TL_A)$ such that, swapping if necessary, $\al\ga\in D_r$ and $\be\ga\in I_{r-2}$.  But $\ga\al = \ga^*\al^* = (\al\ga)^* \mathrel\J \al\ga$, so that $\ga\al\in D_r$, and similarly $\ga\be\in I_{r-2}$.
\epf

Following the established pattern, we prove the following technical lemma, which will
be used to prove the relevant multiplication/separation properties.
Note that part \ref{P2TLiii} is vacuously true when $\Q=\TL$.

\newpage

\begin{lemma}\label{lem:P2TL}
Let $A,B\in\C$, and suppose $\al,\be\in\Q_{A,B}$ with $r=\rank(\al)\geq\rank(\be)=q$.
\ben
\item \label{P2TLi} If $r\geq3$ and $q<r$, then there exists $\ga\in I_r$ such that $\ga\al\in D_{r-2}$ and $\ga\be\in I_{r-2}\sm H_{\ga\al}$.
\item \label{P2TLii} If $q=r\geq2$ and $(\al,\be)\not\in\H$, then there exists $\ga\in I_r$ such that, swapping $\al,\be$ if necessary,
\[
\text{$[\al\ga\in D_r$ and $\be\ga\in I_{r-2}]$ \ \ \ or \ \ \ $[\ga\al\in D_r$ and $\ga\be\in I_{r-2}]$.}
\]
\item \label{P2TLiii} If $r\geq3$ and $\be\in H_\al\sm\{\al\}$, then there exists $\ga\in I_r$ such that $\ga\al\in D_{r-2}$ and $\ga\be\in I_{r-2}\sm H_{\ga\al}$.
\een
\end{lemma}

\pf
\firstpfitem{\ref{P2TLi}}
Write $\al = \partI{a_1}{a_r}{C_1}{C_s}{b_1}{b_r}{D_1}{D_t}$, where $a_1<\cdots<a_r$.  Put $\be'=\al\al^*\be\in\Q_{A,B}$, and note that $\be'$ belongs to $I_{r-2}$ (since $\be$ does).  If $\be'\in I_{r-4}$, then we may take $\ga=\de\al\al^*$, for any $\de\in\TL_A$ with codomain $\{a_1,\ldots,a_{r-2}\}$.  (Since $\al\al^*\in\TL$, we have $\ga\in\TL$.)

Suppose from now on that $\be'\in D_{r-2}$.  Define the set $Q=\codom(\al)\cap\codom(\be')$, noting that $|Q|\leq|{\codom(\be')}|=\rank(\be')=r-2$.  If $|Q|=r-2$, then $Q=\codom(\be')$, so we cannot have $Q=\{b_2,\ldots,b_{r-1}\}$, because $b_2$ is even.  It follows that, regardless of the size of $Q$, there exists $2\leq i\leq r-1$ such that $b_i\not\in Q$.  Let the block of~$\be'$ containing $b_i'$ be $\{b_i',x'\}$.  Let $j\in\{i+1,i-1\}$ be such that $x\not=b_j$.  Let $\ga=\de\al\al^*$, for any $\de\in\TL_A$ with codomain $\{a_1,\ldots,a_r\}\sm\{a_i,a_j\}$ and containing the block $\{a_i',a_j'\}$.  Then $\ga\be=\de\be'\in I_{r-2}$ has the block $\{b_i',x'\}$, but $\ga\al=\de\al\in D_{r-2}$ has the block $\{b_i',b_j'\}$; since $b_j\not=x$, it follows that $\coker(\ga\al)\not=\coker(\ga\be)$, whence $(\ga\al,\ga\be)\not\in\L$, and so $\ga\be\not\in H_{\ga\al}$.

\pfitem{\ref{P2TLii}}  By symmetry we may assume that $(\al,\be)\not\in\R$.  Let $\al'=\al\al^*$ and $\be'=\be\be^*$.  Note that~$\al'$ and $\be'$ are projections in the Temperley--Lieb monoid $\TL_A$ (whether $\al,\be$ are planar or anti-planar).  Since $(\al,\al')\in\R$ and $(\be,\be')\in\R$, but $(\al,\be)\not\in\R$, we must have $\al'\not=\be'$.  We also have $\rank(\al')=\rank(\be')=r\geq2$.  It follows from Lemma \ref{lem:TL_proj} that there exists $\ga\in D_r(\TL_A)$ such that, swapping $\al',\be'$ if necessary, $\ga\al'\in D_r$ and $\ga\be'\in I_{r-2}$.  Since $(\ga\al,\ga\al')\in\R\sub\J$, we have $\ga\al\in D_r$.  Similarly, $\ga\be\in I_{r-2}$.

\pfitem{\ref{P2TLiii}}  The proof here is virtually identical to that of Lemma \ref{lem:P2B}\ref{P2Biii}.  Here we assume that ${a_1<\cdots<a_r}$; the permutation $\pi\in\S_r$ then must in fact be $\ga_{[r]}=(1,r)(2,r-1)\cdots$, so that we \emph{must} have $1\pi=r$.  We then make sure we choose $\ga$ to belong to $\TL$.
\epf

\subsection{Congruences on ideals of $\TL$ and $\TLpm$: the odd case}\label{subsect:TLodd}

For the duration of Section \ref{subsect:TLodd}, we assume that every member of $\C$ has odd cardinality, write $\xi={\LSUB}\bigset{|A|}{A\in\C}$, and  let $\Q$ stand for either $\TL$ or $\TLpm$.
We follow the established procedure to find the congruences on $\Q$ and its ideals.
The brunt of the work goes into enabling an application of  Proposition \ref{prop:small01}, which is achieved in
Lemma \ref{lem:R_trick_TLodd}.
The proof of Lemma \ref{lem:R_trick_TLodd} is modelled on that of Lemma \ref{lem:R_trick_Bodd} concerning the Brauer category $\B$.  Several steps of the proof require a much more subtle argument due to planarity constraints, hence the need for the next two technical lemmas.

\begin{lemma}\label{lem:R_trick_TLodd2}
Suppose $\si$ is a congruence on $I_3=I_3(\Q)$.  If there exists a pair $(\al,\be)\in\si$ with $\coker(\al)\not=\coker(\be)$, then there exists $(\al',\be')\in\si$ with $\al',\be'\in D_1$ and $\codom(\al')\not=\codom(\be')$.
\end{lemma}

\pf
Let $(\al,\be)\in\si$ be such a pair, say with $\al,\be\in\Q_{A,B}$.  For convenience, we may assume without loss of generality that $B=[k]$ for some $k\geq 3$.  

\pfcase{1}  Suppose there exists some odd $u\in\codom(\al)\sm\codom(\be)$.  Then $\al$ has a block $\{x,u'\}$ for some $x\in A$, and $\be$ has a block $\{u',y'\}$ for some $y\in B\sm\{u\}$.  We then take $(\al',\be')=(\ga\al,\ga\be)$ for any $\ga\in\TL_A$ with $\codom(\ga)=\{x\}$.  (Here we have $u\in\codom(\al')\sm\codom(\be')$, and clearly $\al',\be'\in D_1$.)

\pfcase{2}  Suppose we are not in Case 1, but that there exists some even $v\in\codom(\al)\sm\codom(\be)$.  Since there is an even element of $\codom(\al)$, and since $\rank(\al)\leq3$, we must in fact have $\rank(\al)=3$, so we may write $\codom(\al)=\{u<v<w\}$, where $u$ and $w$ are odd.  Since we are not in Case 1, we have $u,w\in\codom(\be)$ but, by assumption, $v\not\in\codom(\be)$.  So $\be$ contains a block of the form $\{v',t'\}$ for some $t\in B\sm\{u,v,w\}$; by planarity (of $\be$ or $\ga_A\be$) we have $u<t<w$, and we note that $t$ is odd (as $v$ is even).

\pfcase{2.1}  If $v<t$, we take $(\al'',\be'')=(\al\ga,\be\ga)$ for any $\ga\in I_3(\TL_B)$ with transversals $\{u,u'\}$, $\{v,v'\}$ and $\{t,t'\}$.  Then $\codom(\al'')=\{u,v,t\}$ yet $\{v',t'\}$ is a block of $\be''$.  So $(\al'',\be'')\in\si$ satisfies the assumptions of Case 1 ($t\in\codom(\al'')\sm\codom(\be'')$ is odd), and we have reduced to that case.

\pfcase{2.2}  If $v>t$, we take $(\al'',\be'')=(\al\ga,\be\ga)$ for any $\ga\in I_3(\TL_B)$ with transversals $\{t,t'\}$, $\{v,v'\}$ and $\{w,w'\}$.  Again, $(\al'',\be'')\in\si$ satisfies the assumptions of Case 1.

\aftercases  We have now dealt with the case in which $\codom(\al)\not\sub\codom(\be)$.  The $\codom(\be)\not\sub\codom(\al)$ case is dual.

\pfcase{3} Now suppose $\codom(\al)=\codom(\be)$.  If necessary (and without causing cokernels to become equal), we first pre-multiply by an arbitrary element of $D_1(\TL_A)$, so we may assume that $\al,\be\in D_1$, say with ${\codom(\al)=\codom(\be)=\{x\}}$.  By assumption, some lower non-transversal of~$\al$ is not a block of $\be$.  By planarity, this non-transversal of $\al$ is contained in either ${\{1',\ldots,(x-1)'\}}$ or $\{(x+1)',\ldots,k'\}$.  These cases have left-right symmetry, so without loss of generality we assume the former is the case.  

Let $u$ be the minimum element of $B$ such that the block of $\al$ containing $u'$ is not equal to the block of $\be$ containing $u'$.  By assumption, $u<x$.  Let the blocks of $\al$ and $\be$ containing $u'$ be $\{u',v'\}$ and $\{u',w'\}$, respectively.  By minimality of $u$ we have $u<v$ and $u<w$; by symmetry, we may assume that $v<w$.  

\pfcase{3.1}  If $u$ is odd, then we take $(\al',\be')=(\al\ga,\be\ga)$ for any $\ga\in D_3(\TL_B)$ containing the blocks $\{u,u'\}$, $\{v,v'\}$, $\{v+1,(v+1)'\}$ and $\{w,x\}$.  Then $u\in\codom(\be')\sm\codom(\al')$.

\pfcase{3.2}  Now suppose $u$ is even.  Let $U=\{1,\ldots,u-1\}$, and let $A_1',\ldots,A_p'$ be all of the lower non-transversals of $\al$ contained in $U'$.  These involve $p$ even and $p$ odd elements of $U$.  Since $|U|=u-1$ is odd, $U$ contains more odd elements than even elements, so there exists some odd $s\in U\sm(A_1\cup\cdots\cup A_p)$.  Then $\al$ contains the block $\{s',t'\}$ for some $t$, and by definition of the~$A_i$, we have $t>u$.  Since $s<u$, minimality of $u$ ensures that $\{s',t'\}$ is a block of $\be$ as well; since $\{u',w'\}$ is a block of $\be$, planarity gives $t>w$.  We then take $(\al',\be')=(\al\ga,\be\ga)$ for any $\ga\in D_3(\TL_B)$ containing the blocks $\{v,v'\}$, $\{v+1,(v+1)'\}$, $\{w,w'\}$, $\{s,u\}$ and $\{t,x\}$.  Then $\codom(\al')=\{v\}\not=\{w\}=\codom(\be')$.
\epf

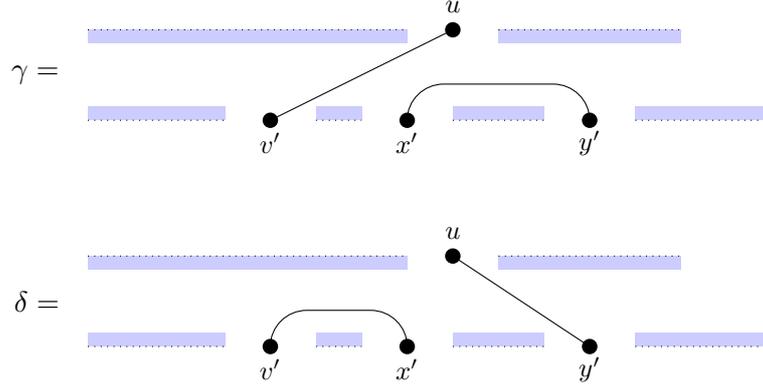
\begin{figure}[ht]
\begin{center}
\begin{tikzpicture}[scale=.6]
\begin{scope}[shift={(0,0)}]	
\darcx8{12}{.8}
\stline95
\ur18{.3}
\ur{10}{14}{.3}
\lr14{.3}
\lr67{.3}
\lr9{11}{.3}
\lr{13}{16}{.3}
\draw(.6,1)node[left]{$\ga=$};
\uvs{9}
\lvs{5,8,12}
\draw(9,2.5)node{\small $u$};
\draw(5,-.5)node{\small $v'$};
\draw(8,-.5)node{\small $x'$};
\draw(12,-.5)node{\small $y'$};
\end{scope}
\begin{scope}[shift={(0,-5)}]	
\darcx58{.8}
\stline9{12}
\ur18{.3}
\ur{10}{14}{.3}
\lr14{.3}
\lr67{.3}
\lr9{11}{.3}
\lr{13}{16}{.3}
\draw(.6,1)node[left]{$\de=$};
\uvs{9}
\lvs{5,8,12}
\draw(9,2.5)node{\small $u$};
\draw(5,-.5)node{\small $v'$};
\draw(8,-.5)node{\small $x'$};
\draw(12,-.5)node{\small $y'$};
\end{scope}
\end{tikzpicture}
\caption{The partitions $\ga,\de\in D_1(\TL_{A,B})$ from Lemma \ref{lem:R_trick_TLodd4}.  The shaded regions denote arbitrary (non-intersecting) non-transversals, with the regions in $\ga$ the same as the corresponding regions in~$\de$.}
\label{fig:R_trick_TLodd2}
\end{center}
\end{figure}

The next lemma will be used frequently in the proof of Lemma \ref{lem:R_trick_TLodd}.  The partitions $\ga$ and $\de$ from the statement are shown in Figure \ref{fig:R_trick_TLodd2}.  Recall the definition of nested blocks from just before Lemma \ref{lem:PEqA}.

\begin{lemma}\label{lem:R_trick_TLodd4}
Suppose $\si$ is a congruence on $I_3=I_3(\Q)$, and suppose there exists a pair $(\al,\be)\in\si$ with $\coker(\al)\not=\coker(\be)$.  Suppose $\ga,\de\in D_1(\TL_{A,B})$ are given by
\[
\ga = \partVI u{A_1}{A_2}{A_s}v{x,y}{B_2}{B_t} \AND \de = \partVI u{A_1}{A_2}{A_s}y{v,x}{B_2}{B_t},
\]
where $v<x<y$ and $x$ is even, and where $\{x',y'\}$and $\{v',x'\}$ are not nested in $\ga$ and $\de$, respectively.  Then $(\ga,\de)\in\si$.
\end{lemma}

\pf
By Lemma \ref{lem:R_trick_TLodd2} we may assume that $\al,\be\in D_1$ and $\{a\}=\codom(\al)\not=\codom(\be)$, say with $\al,\be\in\TL_{C,D}$.  As before, we may assume without loss of generality that $D=[k]$ for some~$k$.  So $\be$ has a block $\{a',b'\}$ for some $b$; since $a$ is odd, $b$ is even; the latter gives $1<b<k$ (since $k=|D|$ is odd).
\bit
\item If $a<b$, then let $\ve_2\in D_3(\TL_{D,B})$ be any Temperley--Lieb partition with the blocks $\{a,v'\}$, $\{b,x'\}$, $\{k,y'\}$ and $B_2',\ldots,B_t'$, and put $(\al',\be')=(\al\ve_2,\be\ve_2)\in\si$.
\item If $a>b$, then let $\ve_2\in D_3(\TL_{D,B})$ be any Temperley--Lieb partition with the blocks $\{1,v'\}$, $\{b,x'\}$, $\{a,y'\}$  and $B_2',\ldots,B_t'$, and put $(\al',\be')=(\be\ve_2,\al\ve_2)\in\si$.
\eit
In either case, $(\ga,\de)=(\ve_1\al',\ve_1\be')\in\si$ for any $\ve_1\in\TL_{A,C}$ with $\ker(\ve_1)=\ker(\ga)$.  
\epf

\begin{lemma}\label{lem:R_trick_TLodd}
For any $\al,\be\in I_3=I_3(\Q)$ with $\al\sim\be$ and $(\al,\be)\not\in\L$, we have $\rho_{I_1}^{I_3}\sub\cg\al\be_{I_3}$.
\end{lemma}

\pf
Write $\si=\cg\al\be_{I_3}$.  Since $\coker(\al)\not=\coker(\be)$, note that Lemma \ref{lem:R_trick_TLodd4} applies to $\si$.  For each positive integer $i$, write $H_i=\{2i,2i+1\}$.

Let $(\ga,\de)\in\rho_{I_1}^{I_3}$ be arbitrary, say with $\ga,\de\in\TL_{A,B}$.  We must show that $(\ga,\de)\in\si$.  This is clear if $\ga=\de$, so suppose instead that $\ga,\de\in D_1$ and $(\ga,\de)\in\R$.  Then we may write ${\ga=\partIV u{A_1}\cdots{A_s}v{B_1}\cdots{B_t}}$ and $\de=\partIV u{A_1}\cdots{A_s}w{C_1}\cdots{C_t}$.  To prove the lemma, it suffices to show that $(\ga,\ve)\in\si$, where~${\ve=\partIV u{A_1}\cdots{A_s}1{H_1}\cdots{H_t}}$, because by symmetry we will also have $(\de,\ve)\in\si$.

To show that $(\ga,\ve)\in\si$, we show by induction that for each $0\leq i\leq t$, we have $(\ga,\ve_i)\in\si$ for some $\ve_i$ with codomain~$\{1\}$ and the blocks $A_1,\ldots,A_s$ and $H_1',\ldots,H_i'$.  Taking $i=t$ will give the desired result.  

For $i=0$, we take $\ve_0=\ga$ if $v=1$.  Otherwise, we may assume that $B_1=\{1,p\}$, and we note that $1<p<v$ and $p$ is even.  By Lemma~\ref{lem:R_trick_TLodd4}, we have
\[
\ga = \partVI u{A_1}{A_2}{A_s}v{1,p}{B_2}{B_t} \mathrel\si \partVI u{A_1}{A_2}{A_s}1{p,v}{B_2}{B_t},
\]
so we take $\ve_0$ to be the latter partition.

\begin{figure}[ht]
\begin{center}
\begin{tikzpicture}[scale=.5]
\begin{scope}[shift={(0,0)}]	
\darcx6{10}{.6}
\darcx5{15}{1.2}
\uustline7{-1}
\uuur{-1}6{.3}
\uuur8{18}{.3}
\lrr04{.3}
\lr79{.3}
\lr{11}{14}{.3}
\lr{16}{23}{.3}
\draw(-2,1.5)node[left]{$\ve_{i-1}\ \ \ =$};
\uuuvs{7}
\lvs{-1,5,6,10,15}
\end{scope}
\begin{scope}[shift={(0,-6)}]	
\darcx6{10}{.8}
\darcx{-1}5{.8}
\uustline7{15}
\uuur{-1}6{.3}
\uuur8{18}{.3}
\lrr04{.3}
\lr79{.3}
\lr{11}{14}{.3}
\lr{16}{23}{.3}
\draw(-2,1.5)node[left]{$\si$};
\uuuvs{7}
\lvs{-1,5,6,10,15}
\end{scope}
\begin{scope}[shift={(0,-12)}]	
\darcx{10}{15}{.8}
\darcx{-1}5{.8}
\uustline76
\uuur{-1}6{.3}
\uuur8{18}{.3}
\lrr04{.3}
\lr79{.3}
\lr{11}{14}{.3}
\lr{16}{23}{.3}
\draw(-2,1.5)node[left]{$\si$};
\uuuvs{7}
\lvs{-1,5,6,10,15}
\end{scope}
\begin{scope}[shift={(0,-18)}]	
\darcx{10}{15}{.8}
\darcx56{.5}
\uustline7{-1}
\uuur{-1}6{.3}
\uuur8{18}{.3}
\lrr04{.3}
\lr79{.3}
\lr{11}{14}{.3}
\lr{16}{23}{.3}
\draw(-2,1.5)node[left]{$\si$};
\uuuvs{7}
\lvs{-1,5,6,10,15}
\draw(5.9,.9)node{\small$H_i'$};
\draw(24,1.5)node[right]{$=\ \ \ \ve_i$};
\end{scope}
\end{tikzpicture}
\caption{Diagrammatic verification of the calculation \eqref{eq:ve'} in the proof of Lemma \ref{lem:R_trick_TLodd}.  The only vertices pictured are $u$, and $1'<(2i)'<(2i+1)'<p'<q'$.  Shading conventions are similar to Figure \ref{fig:R_trick_TLodd2}, but the red shaded region denotes the blocks $H_1',\ldots,H_{i-1}'$.  See the text for more details.}
\label{fig:ve'}
\end{center}
\end{figure}
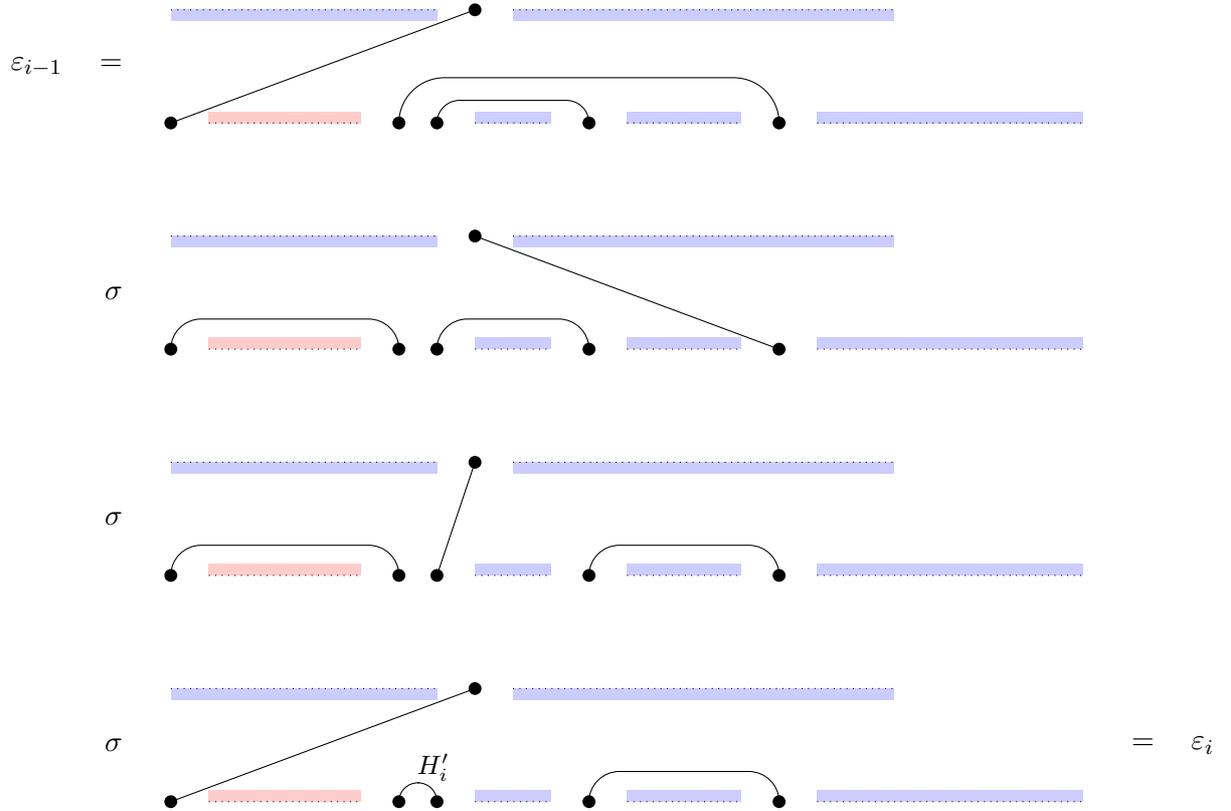

Now let $1\leq i\leq t$, and suppose $(\ga,\ve_{i-1})\in\si$ for some $\ve_{i-1}$ with codomain $\{1\}$ and the blocks $A_1,\ldots,A_s$ and $H_1',\ldots,H_{i-1}'$.  Let the other $\coker(\ve_{i-1})$-classes be $K_i,\ldots,K_t$, where $2i\in K_i$.  If $K_i=H_i$, then we take $\ve_i=\ve_{i-1}$, so suppose otherwise.  Then $K_i=\{2i,q\}$ for some $q>2i+1$, and we may assume that $K_{i+1}=\{2i+1,p\}$ for some $p<q$.  Note that $1<2i<2i+1<p<q$, and that $p$ is even and $q$ odd.  Several applications of Lemma \ref{lem:R_trick_TLodd4} gives the following, with each step shown in Figure \ref{fig:ve'}: 
\newpage
\begin{align}
\nonumber \ga \mathrel\si \ve_{i-1} &= \Big( 
{ \scriptsize \renewcommand*{\arraystretch}{1}
\begin{array} {\c|\c|\c|\c|\c|\c|\c|\c|\cend}
u \:&\: A_1 \:& \multicolumn{6}{c|}{\cdots\cdots\cdots\cdots\cdots\cdots\cdots\cdots\cdots\cdots} &\: A_s \\ \cline{2-9}
1 \:&\: 2i,q \:&\: 2i+1,p \:&\: H_1 \:&\: \cdots \:&\: H_{i-1} \:&\: K_{i+2} \:&\: \cdots \:&\: K_t
\rule[0mm]{0mm}{2.7mm}
\end{array} 
}
\hspace{-1.5 truemm} \Big)\\[2truemm]
\nonumber &\hspace{0.7mm}\mathrel\si \Big( 
{ \scriptsize \renewcommand*{\arraystretch}{1}
\begin{array} {\c|\c|\c|\c|\c|\c|\c|\c|\cend}
u \:&\: A_1 \:& \multicolumn{6}{c|}{\cdots\cdots\cdots\cdots\cdots\cdots\cdots\cdots\cdots\cdots} &\: A_s \\ \cline{2-9}
q \:&\: 1,2i \:&\: 2i+1,p \:&\: H_1 \:&\: \cdots \:&\: H_{i-1} \:&\: K_{i+2} \:&\: \cdots \:&\: K_t
\rule[0mm]{0mm}{2.7mm}
\end{array} 
}
\hspace{-1.5 truemm} \Big)\\[2truemm]
\nonumber &\hspace{0.7mm}\mathrel\si \Big( 
{ \scriptsize \renewcommand*{\arraystretch}{1}
\begin{array} {\c|\c|\c|\c|\c|\c|\c|\c|\cend}
u \:&\: A_1 \:& \multicolumn{6}{c|}{\cdots\cdots\cdots\cdots\cdots\cdots\cdots\cdots} &\: A_s \\ \cline{2-9}
2i+1 \:&\: 1,2i \:&\: p,q \:&\: H_1 \:&\: \cdots \:&\: H_{i-1} \:&\: K_{i+2} \:&\: \cdots \:&\: K_t
\rule[0mm]{0mm}{2.7mm}
\end{array} 
}
\hspace{-1.5 truemm} \Big)\\[2truemm]
\nonumber &\hspace{0.7mm}\mathrel\si \Big( 
{ \scriptsize \renewcommand*{\arraystretch}{1}
\begin{array} {\c|\c|\c|\c|\c|\c|\c|\c|\cend}
u \:&\: A_1 \:& \multicolumn{6}{c|}{\cdots\cdots\cdots\cdots\cdots\cdots\cdots\cdots} &\: A_s \\ \cline{2-9}
1 \:&\: 2i,2i+1 \:&\: p,q \:&\: H_1 \:&\: \cdots \:&\: H_{i-1} \:&\: K_{i+2} \:&\: \cdots \:&\: K_t
\rule[0mm]{0mm}{2.7mm}
\end{array} 
}
\hspace{-1.5 truemm} \Big) \\[2truemm]
\label{eq:ve'}  
&= \Big( 
{ \scriptsize \renewcommand*{\arraystretch}{1}
\begin{array} {\c|\c|\c|\c|\c|\c|\c|\c|\cend}
u \:&\: A_1 \:& \multicolumn{6}{c|}{\cdots\cdots\cdots\cdots\cdots\cdots\cdots\cdots} &\: A_s \\ \cline{2-9}
1 \:&\: H_1 \:&\: \cdots \:&\: H_{i-1} \:&\: H_i \:&\: p,q \:&\: K_{i+2} \:&\: \cdots \:&\: K_t
\rule[0mm]{0mm}{2.7mm}
\end{array} 
}
\hspace{-1.5 truemm} \Big).
\end{align}
Taking $\ve_i$ to be this last partition, this completes the inductive step, and indeed of the lemma.
\epf

We may now prove the main result of this section.  For part \ref{cTLodd1}, note that ${I_1(\TLpm)=I_1(\TL)}$.  Also, the $\R$- and $\L$-classes in $I_1$ are indexed by the allowable (co)kernels; by Lemma \ref{lem:PEqA}, such (co)kernels are planar, and they must in fact be 2-equivalences (with a unique singleton block).  However, not every planar 2-equivalence arises as a (co)kernel of an element of~$I_1$; indeed, the singleton block corresponds to a transversal, and so must not be nested by any of the blocks of size $2$.  For the statement, if $A\in\C$, we will write $\OrdPP{\Eq_2(A)}$ for the set of all planar 2-equivalences on $A$ in which the unique singleton block is not nested.  The isomorphism types of the groups $\G_{r+2}$ in part \ref{cTLodd2} are determined by Lemma \ref{lem:Green_QB}\ref{GQB3}.

\newpage

\begin{thm}\label{thm:TLodd}
Let $\C$ be a non-empty set of finite sets of odd cardinality, let $\Q$ be either the Temperley--Lieb category $\TL=\TL(\C)$ or the anti-Temperley--Lieb category $\TLpm=\TLpm(\C)$, and keep the above notation.
\ben
\item \label{cTLodd1}  The minimal ideal $I_1$ is a partial rectangular band, so its congruence lattice is described by Proposition \ref{prop:PRB}.  In particular, $\Cong(I_1)$ is isomorphic to the lattice direct product
\[
\prod_{A\in\C}\Eq(\OrdPP{\Eq_2(A)})\times\prod_{A\in\C}\Eq(\OrdPP{\Eq_2(A)}).
\]
\item \label{cTLodd2} For odd $r\geq3$ (including $r=\om$ if $\xi=\om$), 
\begin{align*}
\Cong(I_r) = \{\De_{I_r},\ \lam_{I_1}^{I_r},\ \rho_{I_1}^{I_r},\ R_{I_1}^{I_r}\}
\cup\set{R_{I_q,N}^{I_r}}{1\leq q<r,\ q\text{ odd},\ N\normal\G_{q+2}}\cup\{\nab_{I_r}\}
\end{align*}
forms the lattice pictured in Figure \ref{fig:TL}.  
\een
\end{thm}

\pf
Part \ref{cTLodd1} follows from Proposition \ref{prop:PRB}.  For part \ref{cTLodd2}, it suffices as usual (by Theorem~\ref{thm:E2} and Lemma \ref{lem:P2TL}) to show that the result holds for $r=3$.  This follows from Proposition~\ref{prop:small01}, the technical assumptions of which follow from Lemmas \ref{lem:P2TL}\ref{P2TLii}, \ref{lem:R_trick_TLodd} (and its dual) and \ref{lem:N_Bodd}.  (It should be noted that Lemma \ref{lem:N_Bodd} concerned the Brauer category $\B$, but all partitions constructed during its proof were from $\TL$.)
\epf

\begin{figure}[ht]
\begin{center}
\begin{tikzpicture}[scale=0.7,inner sep=0pt]
\begin{scope}[shift={(0,0)}]
\draw[rounded corners=4,line width=10,blue!25,cap=round](0,2)--(0,2.001);
\draw[rounded corners=4,line width=10,green!25,cap=round]
(0,3)--(0,3.001)
(0,4)--(0,4.001)
(0,6)--(0,6.001)
;
\draw(0,2)--(1,1)--(0,0)--(-1,1)--(0,2)--(0,4.5) (0,5.5)--(0,7);
\draw[dotted](0,4.5)--(0,5.5);
\foreach \x/\y in {1/1,-1/1} {\colvert{\x}{\y}{white}}
\foreach \y in {0,2,3,4,6,7} {\colvert{0}{\y}{black}}
\end{scope}
\begin{scope}[shift={(5,0)}]
\draw[rounded corners=4,line width=10,blue!25,cap=round](0,2)--(0,3);
\draw[rounded corners=4,line width=10,green!25,cap=round]
(0,4)--(0,5)
(0,6)--(0,7)
(0,9)--(0,10)
;
\draw(0,2)--(1,1)--(0,0)--(-1,1)--(0,2)--(0,7.5) (0,8.5)--(0,11);
\draw[dotted](0,7.5)--(0,8.5);
\foreach \x/\y in {1/1,-1/1} {\colvert{\x}{\y}{white}}
\foreach \y in {3,5,7,10} {\colvert{0}{\y}{white}}
\foreach \y in {0,2,4,6,9,11} {\colvert{0}{\y}{black}}
\end{scope}
\begin{scope}[shift={(10,1)}]
\foreach \x/\y in {0/0,0/-1,-1/-1,-1/-2} {\draw[rounded corners=4,line width=10,red!25,cap=round](1+\x,1+\y)--(1+\x,1.001+\y);}
\foreach \x/\y in {0/0,0/-1,-1/-1,-1/-2} {\draw[rounded corners=4,line width=10,blue!25,cap=round](0+\x,2+\y)--(0+\x,2.001+\y);}
\draw[rounded corners=4,line width=10,green!25,cap=round]
(0,3)--(0,3.001)
(0,5)--(0,5.001)
;
\draw(0,2)--(1,1)--(0,0)--(-1,1)--(0,2)--(0,3.5) (0,4.5)--(0,6);
\draw[dotted](0,3.5)--(0,4.5);
\draw(0,-1)--(1,0)--(0,1)--(-1,0)--(0,-1);
\draw(0,-1)--(0,0) (0,1)--(0,2) (1,0)--(1,1) (-1,0)--(-1,1);
\foreach \x/\y in {-1/0,-1/1,0/0,0/1,1/0} {\colvert{\x}{\y}{white}}
\foreach \y in {-1,2,3,5,6} {\colvert{0}{\y}{black}}
\foreach \x/\y in {1/1} {\colvert{\x}{\y}{black}}
\end{scope}
\begin{scope}[shift={(15,0)}]
\foreach \x/\y in {0/0,0/-1,-1/-1,-1/-2} {\draw[rounded corners=4,line width=10,red!25,cap=round](3+\x,2+\y)--(2+\x,3+\y);}
\foreach \x/\y in {0/0,0/-1,-1/-1,-1/-2} {\draw[rounded corners=4,line width=10,blue!25,cap=round](1+\x,4+\y)--(0+\x,5+\y);}
\draw[rounded corners=4,line width=10,green!25,cap=round]
(0,6)--(0,7)
(0,9)--(0,10)
;
\foreach \x in {0,1,2,3} {\draw(2-\x,0+\x)--(3-\x,1+\x)--(3-\x,2+\x)--(2-\x,1+\x)--(2-\x,0+\x);}
\draw(0,5)--(0,7.5) (0,8.5)--(0,11);
\draw[dotted](0,7.5)--(0,8.5);
\foreach \x/\y in {0/0,0/1,1/1,1/2} {\draw(2+\x,0+\y)--(-1+\x,3+\y);}
\foreach \x in {1,2,3} {\colvert{1}{\x}{white}}
\foreach \x in {7,10} {\colvert{0}{\x}{white}}
\foreach \x/\y in {0/2,0/3,2/0,2/1,2/2,2/3,3/1,3/2,-1/3,-1/4,0/4,0/5} {\colvert{\x}{\y}{white}}
\foreach \x in {6,9,11} {\colvert0{\x}{black}}
\foreach \x/\y in {3/2,1/4,2/0} {\colvert{\x}{\y}{black}}
\end{scope}
\end{tikzpicture}
\end{center}
\vspace{-7mm}
\caption{Left to right:  Hasse diagrams of $\Cong(I_r(\TL))$, $\Cong(I_r(\TLpm))$, $\Cong(I_s(\TL))$ and $\Cong(I_s(\TLpm))$, where $r\geq3$ is odd and $s\geq4$ is even.  See Theorems~\ref{thm:TLodd},~\ref{thm:TLeven} and~\ref{thm:TLpmeven}.}
\label{fig:TL}
\end{figure}

\subsection{Congruences on ideals of $\TL$: the even case}\label{subsect:TLeven}

We may very quickly describe the congruences on ideals of $\TL$ in the case in which every member of $\C$ has even size.  The reason for this is that $\TL(\C)$ is then isomorphic to an ordinary planar partition category $\OrdP{\P(\C')}$ over a different collection $\C'$; cf.~\cite[Section~1]{HR2005}.

In this paragraph, $\C$ will be an arbitrary set of finite sets (of even and/or odd size).  For $A\in\C$, we define $2A=A\times\{1,2\}$, and we also write $2\C=\set{2A}{A\in\C}$.  Given a planar partition $\al$ from $\OrdP{\P_{A,B}}$, we ``trace around'' the edges of a canonical graph representing $\al$ (see \cite[Lemma 7.1]{EMRT2018}) to construct a Temperley--Lieb partition $\alt$ from $\TL_{2A,2B}$, as shown in Figure \ref{fig:TLiso}.  The map $\al\mt\alt$ is then an isomorphism $\OrdP{\P(\C)}\to\TL(2\C)$.

\begin{figure}[ht]
\begin{center}
\begin{tikzpicture}[scale=1]
\uarcxx{1.25}{3.75}{.65}{gray!60}
\uarcxx{1.75}{3.25}{.5}{gray!60}
\uarcxx{2.25}{2.75}{.2}{gray!60}
\uarcxx{4.75}{5.25}{.2}{gray!60}
\uarcx14{.8}
\uarcx23{.35}
\darcxx{2.75}{3.25}{.2}{gray!60}
\darcxx{1.25}{1.75}{.2}{gray!60}
\darcxx{5.25}{5.75}{.2}{gray!60}
\darcxx{6.25}{6.75}{.2}{gray!60}
\darcxx{2.25}{3.75}{.4}{gray!60}
\darcxx{4.75}{7.25}{.55}{gray!60}
\stlinex{.75}{.75}{gray!60}
\stlinex{4.25}{4.25}{gray!60}
\stlinex{5.75}{7.75}{gray!60}
\stlinex{6.25}{8.25}{gray!60}
\darcx12{.35}
\darcx24{.6}
\darcx56{.35}
\darcx67{.35}
\stlines{1/1,4/4,6/8}
\foreach \x in {1,...,8} {
\fill (\x,0)circle(.1);
\fill[gray!60] (\x+.25,0)circle(.1);
\fill[gray!60] (\x-.25,0)circle(.1);
}
\foreach \x in {1,...,6} {
\fill (\x,2)circle(.1);
\fill[gray!60] (\x+.25,2)circle(.1);
\fill[gray!60] (\x-.25,2)circle(.1);
}
\end{tikzpicture}
\end{center}
\vspace{-5mm}
\caption{A planar partition $\al$ from $\OrdP{\P_{A,B}}$ (black), with its corresponding Temperley--Lieb partition $\alt$ from $\TL_{2A,2B}$ (grey).}
\label{fig:TLiso}
\end{figure}

As a result of this isomorphism, we may easily describe the congruences on ideals of $\TL$ (in the case that all members of $\C$ are of even size) by suitably translating Theorem \ref{thm:QP}.  Recall that $\OrdP{\Eq_2(A)}$ denotes the set of all planar 2-equivalences on $A$.  As usual, we write  ${\xi={\LSUB}\bigset{|A|}{A\in\C}}$.

\begin{thm}\label{thm:TLeven}
Let $\C$ be a non-empty set of finite sets of even cardinality, let $\TL=\TL(\C)$ be the Temperley--Lieb category over $\C$, and keep the above notation.
\ben
\item \label{cTLeven1}  The minimal ideal $I_0$ is a partial rectangular band, so its congruence lattice is described by Proposition \ref{prop:PRB}.  In particular, $\Cong(I_0)$ is isomorphic to the lattice direct product
\[
\prod_{A\in\C}\Eq(\OrdP{\Eq_2(A)})\times\prod_{A\in\C}\Eq(\OrdP{\Eq_2(A)}).
\]
\item \label{cTLeven2}  The ideal $I_2$ is retractable, every congruence on $I_0$ is liftable to~$I_2$, and
\[
\Cong(I_2)=\set{\tau\cup\De_{D_2}}{\tau\in\Cong(I_0)} \cup \set{\th_{I_0,\tau}^{I_2}}{\tau\in\Cong(I_0)}.
\]
In particular, $\Cong(I_2)$ is isomorphic to the lattice direct product $\Cong(I_0)\times\btwo$.
\item \label{cTLeven3} For even $r\geq4$ (including $r=\om$ if $\xi=\om$), 
\begin{align*}
\Cong(I_r) = \set{\mu_{I_q}^{I_r},\ \lam_{I_q}^{I_r},\ \rho_{I_q}^{I_r},\ R_{I_q}^{I_r}}{q=0,2}
\cup\set{R_{I_q}^{I_r}}{4\leq q<r,\ q\text{ even}}\cup\{\nab_{I_r}\}
\end{align*}
forms the lattice pictured in Figure \ref{fig:TL}.  \epfres
\een
\end{thm}

\subsection{Congruences on ideals of $\TLpm$: the even case}\label{subsect:TLpmeven}

Now we turn to the anti-Temperley--Lieb category $\TLpm=\TLpm(\C)$ where every member of $\C$ has even cardinality.
The trick from the previous section no longer applies: even the monoids~$\P^\pm_n$ and~$\TLpm_{2n}$ are not isomorphic as they have different sizes; we will also see that they have non-isomorphic congruence lattices (cf.~Figures \ref{fig:QP} and \ref{fig:TL}).
As usual let ${\xi={\LSUB}\bigset{|A|}{A\in\C}}$.

We continue to use the notation introduced earlier in this section.
And again we fix a specific group $\H$-class $\G_q$ in each $\D$-class~$D_q$.  For $q\geq2$, such a group 
is cyclic of order $2$, specifically $\G_q=\{\ve_q,\ga_q\}$ where
\[
\ve_q=\id_q^\natural = \partI1q{q+1,q+2}{n_q-1,n_q}1q{q+1,q+2}{n_q-1,n_q} \AND 
\ga_q=\ga_{[q]}^\natural = \partI1q{q+1,q+2}{n_q-1,n_q}q1{q+1,q+2}{n_q-1,n_q};
\]
and, of course, $\G_0=\{\ve_0\}$ is trivial.

As with $\B$ and $\TL$, the congruences on the bottom two ideals $I_r=I_r(\TLpm)$, $r=0,2$, are described by Propositions \ref{prop:PRB} and \ref{prop:small_ideal}.  Thus, the main issue is with the ideals $I_r$ for $r\geq4$.  

An ideal $I_r$ ($r\geq4$) contains the IN-pairs $(I_q,\{\ve_{q+2}\})$ and $(I_q,\G_{q+2})$, for even ${0\leq q<r}$.  The retraction in $\B$ persists in $\TLpm$, since for any $\al\in I_2(\TLpm)$, we clearly have $\wh\al\in\TLpm$ (indeed,~$\wh\al\in\TL$).  So the ideal $I_2=I_2(\TLpm)$ is retractable, and we continue to denote the retraction by
\[
f:I_2\to I_0 : \al\mt \wh\al.
\]
Further, we have the following four retractable IN-pairs in~$I_r$:
\[
(I_0,\{\ve_2\}) \COMMA
(I_0,\G_2) \COMMA
(I_2,\{\ve_4\}) \COMMA
(I_2,\G_4).
\]
It is clear that the first and third of the above pairs are retractable.  The second is retractable because $\G_2=\S_2^\natural$ and $(I_0(\B),\S_2^\natural)$ is a retractable IN-pair in $\B$ (it also follows from Lemma \ref{lem:INretract}).  The fourth is retractable because $\G_4=\big\{{\id_4},(1,4)(2,3)\big\}^\natural$ is a subgroup of $K^\natural$ (where $K$ is the Klein 4-group), and since $(I_2(\B),K^\natural)$ is a retractable IN-pair in $\B$.  As usual, these retractable IN-pairs each yield a family of $\mu/\lam/\rho/R$ congruences.

\begin{lemma}\label{lem:R_trick_TLpm}
For any $\al,\be\in I_4=I_4(\TLpm)$ with $\al\sim\be$ and $(\al,\be)\not\in\L\cup\L f^{-1}$, we have $\rho_{I_0}^{I_4}\sub\cg\al\be_{I_4}$.
\end{lemma}

\pf
Write $\si=\cg\al\be_{I_4}$, and suppose $\al,\be\in\TLpm_{A,B}$.  We consider two cases.  

\pfcase{1}  
Suppose first that there exist distinct $u,v,w\in B$ such that $\al$ and $\be$ contain the blocks $\{u',v'\}$ and $\{v',w'\}$, respectively.  Premultiplying the pair $(\al,\be)$ by an arbitrary element of $D_0(\TL_A)$, we may in fact assume that $\al,\be\in D_0$.

Since $D_0(\TLpm)=D_0(\TL)$, we have $\al,\be\in D_0(\TL)$.  By Theorem \ref{thm:TLeven}, the smallest congruence on $I_4(\TL)$ containing the pair $(\al,\be)$ is $\rho_{I_0(\TL)}^{I_4(\TL)}$; i.e., $\cg\al\be_{I_4(\TL)}=\rho_{I_0(\TL)}^{I_4(\TL)}$.  But then
\[
\wt\R\restr_{D_0(\TLpm)} = \wt\R\restr_{D_0(\TL)} \sub \rho_{I_0(\TL)}^{I_4(\TL)} = \cg\al\be_{I_4(\TL)} \sub \cg\al\be_{I_4(\TLpm)} = \si,
\]
so it follows that $\rho_{I_0(\TLpm)}^{I_4(\TLpm)} = \wt\R\restr_{D_0(\TLpm)} \cup \De_{I_4(\TLpm)} \sub \si$.

\pfcase{2}
If we are not in Case 1, then we follow the proof of Lemma \ref{lem:R_trick_Beven} to reduce to Case~1.  We just need to be careful that the element $\ga$ constructed is planar, but we can do this by noting that $a_1$ and $a_2$ have different parities (since $\{b_1',b_2'\}$ is a block of $\be$), so we may choose $a_3$ of different parity to $a_2$.
\epf

\begin{thm}\label{thm:TLpmeven}
Let $\C$ be a non-empty set of finite sets of even cardinality, let $\TLpm=\TLpm(\C)$ be the anti-Temperley--Lieb category, and keep the above notation.
\ben
\item \label{cTLpmeven1}  The minimal ideal $I_0$ is a partial rectangular band, so its congruence lattice is described by Proposition \ref{prop:PRB}.  In particular, $\Cong(I_0)$ is isomorphic to the lattice direct product
\[
\prod_{A\in\C}\Eq(\OrdP{\Eq_2(A)})\times\prod_{A\in\C}\Eq(\OrdP{\Eq_2(A)}).
\]
\item \label{cTLpmeven2}  The ideal $I_2$ is retractable, every congruence on $I_0$ is liftable to~$I_2$, and
\[
\Cong(I_2)=\set{\tau\cup\nu_N}{\tau\in\Cong(I_0),\ N\normal\G_2} \cup \set{\th_{I_0,\tau}^{I_2}}{\tau\in\Cong(I_0)}.
\]
In particular, $\Cong(I_2)$ is isomorphic to the lattice direct product $\Cong(I_0)\times\bthree$.
\item \label{cTLpmeven3} For even $r\geq4$ (including $r=\om$ if $\xi=\om$), 
\begin{align*}
\Cong(I_r) = 
\set{\mu_{I_q,N}^{I_r},\ \lam_{I_q,N}^{I_r},\ \rho_{I_q,N}^{I_r},\ &R_{I_q,N}^{I_r}}{q=0,2,\ N\normal\G_{q+2}} \\[3truemm]
{}\cup\set{&R_{I_q,N}^{I_r}}{4\leq q<r,\ q\text{ even},\ N\normal\G_{q+2}}\cup\{\nab_{I_r}\}
\end{align*}
forms the lattice pictured in Figure \ref{fig:TL}.  
\een
\end{thm}

\pf
As usual, Propositions \ref{prop:PRB} and \ref{prop:small_ideal} give the first two parts.  For the third, it suffices (by Theorem \ref{thm:E2} and Lemma \ref{lem:P2TL}) to prove the $r=4$ case.  This follows from Proposition~\ref{prop:small012}, together with Lemmas \ref{lem:P2TL}\ref{P2TLii} and \ref{lem:R_trick_TLpm}, and the dual of the latter.  (In assumption \ref{012iv} of Proposition~\ref{prop:small012} we have $H=\G_2$.)
\epf

\section{Jones and anti-Jones categories}\label{sect:J}

This is the last section of the paper concerning diagram categories.  
It is also the last where we use the most straightforward form of the stacking mechanism from Section \ref{subsect:chains}.
We describe the congruences on the ideals of the Jones and anti-Jones categories, $\JJ=\JJ(\C)$ and $\Jpm=\Jpm(\C)$, as defined in Section \ref{subsect:QB}.  Here, as usual, $\C$ is a fixed non-empty set of finite sets.  Throughout Section \ref{sect:J}, $\Q$ will stand for either of $\JJ$ or $\Jpm$.  As in Sections \ref{sect:B} and \ref{sect:TL}, we have decompositions
\[
\Q = \Q(\C) = \Q(\C_\even) \sqcup \Q(\C_\odd),
\]
which will allow us to assume without loss of generality that the sizes of all members of $\C$ have the same parity; we consider the odd case in Section \ref{subsect:Jodd} and the even case in Section \ref{subsect:Jeven}.  Technical results relevant to both cases are in Section \ref{subsect:prelim_J}.

Before we begin, we make some observations regarding an element $\al=\partI{a_1}{a_r}{C_1}{C_s}{b_1}{b_r}{D_1}{D_s}$ of $\Jpm$, say with $\bd(\al)=A$ and $\br(\al)=B$.  First, concerning transversals, if $a_1<\cdots<a_r$, then the parities $a_1,a_2,\ldots,a_r$ must alternate, but it is no longer necessarily the case that $a_1$ is odd.

Next, if $a_1<\cdots<a_r$, then we do not necessarily have $b_1<\cdots<b_r$ or $b_1>\cdots>b_r$, though the possible orderings on the $b_i$ are restricted to cyclic orderings and/or reversals, depending on whether $\al$ is annular or anti-annular.  We also do not necessarily have $a_i\equiv b_i\Mod2$, even if $a_1<\cdots<a_r$ and $b_1<\cdots<b_r$ both hold.  In an upper non-transversal $C_i=\{u,v\}$, $u$ and~$v$ do not necessarily have opposite parities, though these parities can only be equal if the 
edge connecting $u$ and $v$ wraps around the identified border in an annular representation of $\al$;
for example, consider~the annular partition~$\partn{2}{1,3}{1}{2,3}$.

Of crucial importance is the fact that if $a_1=\min(A)$ and $b_1=\min(B)$, then $\al$ as above belongs to $\TL_{A,B}$.  Indeed, consider some annular representation of $\al$ (as a planar graph drawn in the usual way on the cylinder); rotating the bottom circle of the cylinder a whole number of times if necessary, we can assume that the edge $\{a_1,b_1'\}$ is drawn as a vertical line; all the remaining edges of the graph must not cross this line, and hence do not wrap around the identified border, so that $\al$ is actually planar.

The above considerations mean that there will sometimes be more cases to consider in proofs; but on the other hand, they will sometimes allow us more freedom when we wish to define an (anti-)Jones partition with a desired set of blocks.

We are now almost ready to begin, but first we introduce one more piece of notation.  For an integer $n\geq1$, we write $\JJ_n=\JJ_{[n]}=\JJ_{[n],[n]}$, where as usual $[n]=\{1,\ldots,n\}$, and with similar meanings for $\Jpm_n$, $\TL_n$, $\P_n$ and so on.  We also write $\ga_n=\ga_{[n]}$ and $\de_n=\de_{[n]}$, where the permutations $\ga_A,\de_A\in\S_A$ are defined in \eqref{eq:gade}.

\subsection{Multiplicative properties in $\JJ$ and $\Jpm$}\label{subsect:prelim_J}

Once again, we need to prove the usual technical result (Lemma \ref{lem:P2J} below)
that is used for establishing the required multiplication and separation properties.
As preparation we need an ``annular version'' of Lemma~\ref{lem:TL_proj}.
In the proof, we will use the following construction: for partitions $\al\in\P_k$ and $\be\in\P_l$, we form $\al\oplus\be\in\P_{k+l}$ by placing a translated copy of $\be$ to the right of~$\al$. 
When $\C=\bigset{[n]}{n\in\N}$, this operation (extended in the natural way to arbitrary hom-sets) gives the partition category $\P=\P(\C)$ the structure of a \emph{strict monoidal category} in the sense of \cite{JS1993}; cf.~\cite{LZ2015,Martin2008,JEpresDC}.
Not all subcategories of $\P$ are closed under the~$\oplus$ operation; $\B$ and~$\TL$ are, for example, but $\TLpm$, $\JJ$ and $\Jpm$ are not.
But note that, for example, if $\al,\be\in\JJ$ are such that $\rank(\al)\leq1$ and $\rank(\be)=0$, then $\al\oplus\be\in\JJ$.

\begin{lemma}\label{lem:J_proj}
Let $\al,\be$ be distinct projections from $\JJ_n$ with $r=\rank(\al)=\rank(\be)\geq2$.  Then there exists $\ga\in D_r(\JJ_n)$ such that, swapping $\al,\be$ if necessary, $\ga\al\in D_r$ and $\ga\be\in I_{r-2}$.
\end{lemma}

\pf
Write $\al=\partI{a_1}{a_r}{B_1}{B_s}{a_1}{a_r}{B_1}{B_s}$, where $a_1<\cdots<a_r$.  Let $\al'=\de_n^{a_1-1}\al\de_n^{-a_1+1}$ and $\be'=\de_n^{a_1-1}\be\de_n^{-a_1+1}$, noting that $\al',\be'$ are still distinct projections of rank $r$, and that $\{1,1'\}$ is a transversal of $\al$.  
If we can find $\ga'\in D_r$ such that $\ga'\al'\in D_r$ and $\ga'\be'\in I_{r-2}$, then we may take $\ga=\ga'\de_n^{a_1-1}$.

Thus, without loss of generality, we may assume that $\al=\al'$ and $\be=\be'$, so in particular $a_1=1$.  From this it follows that $\al$ in fact belongs to $\TL_n$.  If $1\in\dom(\be)$, then $\be$ also belongs to $\TL_n$, so the lemma follows from Lemma \ref{lem:TL_proj}.  Thus, for the remainder of the proof, we will assume that $1\not\in\dom(\be)$, and denote by $\{1,x\}$ the block of $\be$ containing $1$.

Since $\al\in\TL$, the parities of $1=a_1,a_2,a_3,\ldots$ are odd, even, odd, and so on.  On the other hand, $x$ could be even or odd; this depends on the parity of $n$, and whether or not the edge~$\{1,x\}$ wraps around the back of the cylinder.

\pfcase{1}
Suppose first that $x\leq a_r$.  Since $x>1=a_1$, it follows that $a_i<x\leq a_{i+1}$ for some $1\leq i\leq r-1$.  Let $j\in\{i,i+1\}$ be such that $x$ and $a_j$ have the same parity.  

\pfcase{1.1}
If $j\geq2$, then we take $\ga\in D_r(\TL_n)$ to be any Temperley--Lieb projection with domain $\big(\{a_1,\ldots,a_r\}\sm\{a_j\}\big)\cup\{x\}$.  (Note that we might in fact have $a_j=x$.)  Then $\ga\al\in D_r$, since $\{a_1,\ldots,a_r\}\sm\{a_j\}\sub\dom(\ga\al)$; but $\ga\be\in I_{r-2}$, since $\dom(\ga\be)\sub\dom(\ga)$ yet $\{1,x\}=\{a_1,x\}$ is a non-transversal of $\ga\be$.

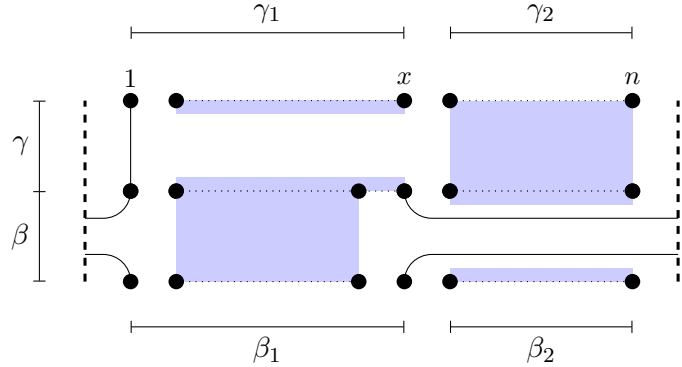
\begin{figure}[ht]
\begin{center}
\begin{tikzpicture}[scale=.6]
\begin{scope}[shift={(0,0)}]	
\stline11
\ur27{.3}
\lr27{.3}
\ur8{12}{1.5}
\lr8{12}{1.5}
\draw[|-] (-1,2)--(-1,0);
\draw(-1,1)node[left]{$\ga$};
\uvs{1,2,7,8,12}
\lvs{1,2,7,8,12}
\draw(1,2.5)node{\small $1$};
\draw(7,2.5)node{\small $x$};
\draw(12,2.5)node{\small $n$};
\draw[|-|] (1,3.5)--(7,3.5);
\draw(4,3.5)node[above]{$\ga_1$};
\draw[|-|] (12,3.5)--(8,3.5);
\draw(10,3.5)node[above]{$\ga_2$};
\end{scope}
\begin{scope}[shift={(0,-2)}]	
\draw (0,1.4)--(.4,1.4) arc(270:360:.6);
\draw (7,2)arc(180:270:.6) (7.6,1.4)--(13,1.4);
\draw (7,0)arc(180:90:.6) (7.6,.6)--(13,.6);
\draw (0,.6)--(.4,.6) arc(90:0:.6);
\ur26{1.3}
\lr26{1.3}
\ur8{12}{.3}
\lr8{12}{.3}
\draw[|-|] (-1,2)--(-1,0);
\draw(-1,1)node[left]{$\be$};
\uvs{1,2,6,7,8,12}
\lvs{1,2,6,7,8,12}
\draw[dashed,line width=0.4mm] (0,0)--(0,4);
\draw[dashed,line width=0.4mm] (13,0)--(13,4);
\draw[|-|] (1,-1)--(7,-1);
\draw(4,-1)node[below]{$\be_1$};
\draw[|-|] (12,-1)--(8,-1);
\draw(10,-1)node[below]{$\be_2$};
\end{scope}
\end{tikzpicture}
\caption{The partitions $\be=\be_1\oplus\be_2$ and $\ga=\ga_1\oplus\ga_2$ from Case 1.2 of the proof of Lemma \ref{lem:J_proj}.}
\label{fig:J_proj}
\end{center}
\end{figure}

\pfcase{1.2}
Now suppose $j=1$, and note then that $i=1$ as well.  Since $x$ has the same parity as $a_j=a_1=1$, it follows that $x$ is odd.  This then implies that the edge $\{1,x\}$ of $\be$ wraps around the back of the cylinder.  By planarity, each vertex from $\{x+1,\ldots,n\}$ is contained in an upper non-transversal of $\be$, and this in turn tells us that $n$ is odd.  It also follows that $\be=\be_1\oplus\be_2$ for some $\be_1\in D_r(\JJ_x)$ and $\be_2\in D_0(\JJ_{n-x})$.  The partitions $\be_1$ and $\be_2$ are shown in Figure \ref{fig:J_proj}, which also shows further partitions constructed during this case of the proof.  Further, by the definition of $x$, and since $a_2$ is even, we have $a_1<x<a_2<a_3<\cdots<a_r$, so that each of $a_2,\ldots,a_r$ belong to $\{x+1,\ldots,n\}$.  Let~$\ga_1$ be any projection from $\TL_x$ with (co)domain $\{1\}$, and let $\ga_2$ be any projection from $\TL_{n-x}$ with (co)domain $\{a_2-x,\ldots,a_r-x\}$.  (Note that $\ga_2$ is well defined since $a_2-x,a_3-x,\ldots$ are odd, even, and so on, since $x$ is odd.)  Now put $\ga=\ga_1\oplus\ga_2\in\TL_n\sub\JJ_n$, and note that $\codom(\ga)=\dom(\al)$ by construction.  It follows that $\ga$ and $\ga\al$ both belong to $D_r$.  On the other hand,
\begin{equation}\label{eq:D1D0}
\ga\be = (\ga_1\be_1)\oplus(\ga_2\be_2) \in D_1(\JJ_x)\oplus D_0(\JJ_{n-x}) \sub D_1 \sub I_{r-2}.
\end{equation}

\pfcase{2}
Now suppose $x>a_r$.  

\pfcase{2.1}
Suppose first that a graph of $\be$ may be drawn in which $\{1,x\}$ does not wrap around the back of the cylinder.  By planarity, $x$ is even, and each vertex from $\{2,\ldots,x-1\}$ belongs to a non-transversal of $x$.  Thus, keeping in mind that $\be$ is a projection, we have $\be=\be_1\oplus\be_2$ for some $\be_1\in D_0(\TL_x)$ and some $\be_2\in D_r(\JJ_{n-x})$.  (A diagram similar to Figure \ref{fig:J_proj} could be drawn.)  Since $x$ is even, $n$, $r$ and $n-x$ all have the same parity.

Suppose first that $n$, $r$ and $n-x$ are all even.  Noting that $a_r<x$, let $\ga_1$ be any projection from $\TL_x$ with domain $\{a_1,\ldots,a_r\}$, let $\ga_2$ be an arbitrary projection from $D_0(\TL_{n-x})$, and put $\ga=\ga_1\oplus\ga_2$.  Then certainly $\ga,\ga\al\in D_r$, while a calculation similar to \eqref{eq:D1D0} gives $\ga\be\in D_0\sub I_{r-2}$.

If $n$, $r$ and $n-x$ are all odd, then we let $\ga=\ga_1\oplus\ga_2$ where $\ga_1\in\TL_x$ is a projection with domain $\{a_1,\ldots,a_{r-1}\}$, and $\ga_2\in \TL_{n-x}$ a projection with rank $1$.  This time, $\ga,\ga\al\in D_r$ and $\ga\be\in D_1\sub I_{r-2}$.

\pfcase{2.2}
Now suppose a graph of $\be$ may be drawn in which $\{1,x\}$ does wrap around the back of the cylinder.  By planarity, the $n-x$ vertices from $\{x+1,\ldots,n\}$ all belong to non-transversals of $\be$ nested by this edge, so this set must have even size; i.e., $x\equiv n\Mod2$.  But also $a_r\equiv n\Mod2$, since $\al\in\TL_n$, so it follows that $x$ and $a_r$ have the same parity.  Keeping in mind that $a_r<x$, we take $\ga\in D_r(\TL_n)$ to be any Temperley--Lieb projection with (co)domain $\{a_1,\ldots,a_{r-1},x\}$.  Then $\ga,\ga\al\in D_r$, and $\ga\be\in I_{r-2}$ since $\dom(\ga\be)\sub\dom(\ga)$ and $\ga\be$ contains the non-transversal $\{1,x\}=\{a_1,x\}$.
\epf

\begin{lemma}\label{lem:P2J}
Let $A,B\in\C$, and suppose $\al,\be\in\Q_{A,B}$ with $r=\rank(\al)\geq\rank(\be)=q$.
\ben
\item \label{P2Ji} If $r\geq3$ and $q<r$, then there exists $\ga\in I_r$ such that $\ga\al\in D_{r-2}$ and $\ga\be\in I_{r-2}\sm H_{\ga\al}$.
\item \label{P2Jii} If $q=r\geq2$ and $(\al,\be)\not\in\H$, then there exists $\ga\in I_r$ such that, swapping $\al,\be$ if necessary,
\[
\text{$[\al\ga\in D_r$ and $\be\ga\in I_{r-2}]$ \ \ \ or \ \ \ $[\ga\al\in D_r$ and $\ga\be\in I_{r-2}]$.}
\]
\item \label{P2Jiii} If $r\geq3$ and $\be\in H_\al\sm\{\al\}$, then there exists $\ga\in I_r$ such that $\ga\al\in D_{r-2}$ and $\ga\be\in I_{r-2}\sm H_{\ga\al}$.
\een
\end{lemma}

\pf
\firstpfitem{\ref{P2Ji}}  The proof is the same as that of Lemma \ref{lem:P2TL}\ref{P2TLi}, except that (in the case $\be'\in D_{r-2}$), we \emph{could} have $Q=\{b_2,\ldots,b_{r-1}\}$.  In this case, we simply take $i=1$ and continue as before, but interpret $i-1=r$ in the definition of $j$.

\pfitem{\ref{P2Jii}}  This follows from Lemma \ref{lem:J_proj} in the same way that Lemma \ref{lem:P2TL}\ref{P2TLii} followed from Lemma~\ref{lem:TL_proj}.  (Note that all projections from $\Jpm$ belong to $\JJ$.)

\pfitem{\ref{P2Jiii}}  We follow the proof of Lemma \ref{lem:P2B}\ref{P2Biii}, assuming the ordering $a_1<\cdots<a_r$.  If $1\pi=r$, then we define $\ga$ as in the above-mentioned proof, but making sure that $\ga\in\JJ$.  Otherwise, we have $1\pi=i$ for some $1<i<r$.  Here we have $2\pi=i\pm1$ (depending on whether~$\be$ is annular or anti-annular), but since $r\geq3$ we do not have $\{1\pi,2\pi\}=\{1,2\}$.  In this case we take $\ga\in\TL_A$ to be any Temperley--Lieb partition with codomain $\{a_3,\ldots,a_r\}$, noting that $\codom(\ga\al)=\{b_3,\ldots,b_r\} \not= \{b_1,\ldots,b_r\}\sm\{b_{1\pi},b_{2\pi}\}=\codom(\ga\be)$.
\epf

\subsection{Congruences on ideals of $\JJ$ and $\Jpm$: the odd case}\label{subsect:Jodd}

For the duration of Section \ref{subsect:Jodd}, we assume that every member of $\C$ has odd cardinality.  
Most of the work goes into describing the congruences on the ideal $I_3=I_3(\Q)$, which largely amounts to verifying the technical assumptions of Proposition \ref{prop:small01}; see Lemma \ref{lem:R_trick_Jodd}.  The rest is taken care of by our stacking mechanism.

\begin{lemma}\label{lem:R_trick_Jodd1}
Suppose $\si$ is a congruence on $I_3=I_3(\Q)$.  If there exists a pair $(\al,\be)\in\si$ with $\coker(\al)\not=\coker(\be)$, then there exists $(\al',\be')\in\si$ with $\al',\be'\in D_1$ and $\codom(\al')\not=\codom(\be')$.
\end{lemma}

\pf
Let $(\al,\be)\in\si$ be such a pair, say with $\al,\be\in\Q_{A,B}$, and assume for convenience that $B=[k]$ for some $k$.  If $\codom(\al)\not\sub\codom(\be)$, then we fix $u\in \codom(\al)\sm\codom(\be)$, and (regardless of the parity of $u$) we follow Case 1 in the proof of Lemma \ref{lem:R_trick_TLodd2} to find such a pair~$(\al',\be')$.  By symmetry, the only remaining case is where $\codom(\al)=\codom(\be)$, and again premultiplying by an arbitrary element of $D_1(\TL_A)$, we may assume that $\al,\be\in D_1$, say with $\codom(\al)=\codom(\be)=\{x\}$.  Further premultiplying by an element of $\JJ_{B,A}$ with domain~$\{x\}$, we may assume that $\al,\be\in\Q_B$ both contain the transversal $\{x,x'\}$.

Now let $\al''=\de_B^{x-1}\al\de_B^{1-x}$ and $\be''=\de_B^{x-1}\be\de_B^{1-x}$.  
(Since $\de_B$ may not belong to $I_3$, we are not able to assume that $(\al'',\be'')\in\si$.)
Since $\coker(\al)\not=\coker(\be)$, and since $\de_B$ is a unit, it follows that $\coker(\al'')\not=\coker(\be'')$.  Since $\al''$ and $\be''$ both contain the transversal~$\{1,1'\}$, they belong to $\TL_B$, indeed to $D_1(\TL_B)$.  We then apply the argument of Case 3 in the proof of Lemma \ref{lem:R_trick_TLodd2} to the pair $(\al'',\be'')$, which tells us that there exists $\ga\in D_3(\TL_B)$ such that $\codom(\al''\ga)\not=\codom(\be''\ga)$; i.e., $\codom(\de_B^{x-1}\al\de_B^{1-x}\ga)\not=\codom(\de_B^{x-1}\be\de_B^{1-x}\ga)$.  Since $\de_B$ is a unit, it follows from this that $\codom(\al\cdot\de_B^{1-x}\ga)\not=\codom(\be\cdot\de_B^{1-x}\ga)$.  Since $\de_B^{1-x}\ga\in I_3(\Q)$, we may take $(\al',\be')=(\al\cdot\de_B^{1-x}\ga,\be\cdot\de_B^{1-x}\ga)\in\si$.
\epf

The next lemma improves on the previous one by showing that the pair $(\al',\be')$ can be taken with $\al',\be'\in\TL$.  

\begin{lemma}\label{lem:R_trick_Jodd2}
Suppose $\si$ is a congruence on $I_3=I_3(\Q)$.  If there exists a pair $(\al,\be)\in\si$ with $\coker(\al)\not=\coker(\be)$, then 
there exists $(\al',\be')\in\si$ with $\al',\be'\in D_1(\TL)$ and ${\codom(\al')\not=\codom(\be')}$.
\end{lemma}

\pf
Again we assume that $\al,\be\in\Q_{A,B}$, with $B=[k]$ for some odd $k=2l+1$.
By Lemma~\ref{lem:R_trick_Jodd1} we may assume that $\{u\}=\codom(\al)\not=\codom(\be)=\{v\}$, and by symmetry that $u<v$.  For $1\leq i\leq l$, let $H_i=\{2i,2i+1\}$, and put  
\[
\al' = \partVI1{H_1}{H_2}{H_l}1{H_1}{H_2}{H_l} \AND \be' = \partVI1{H_1}{H_2}{H_l}3{1,2}{H_2}{H_l},
\]
both from $\TL_B$.  Also let $\ve_1$ be any element of $D_1(\TL_{B,A})$ with $\ker(\ve_1)=\ker(\al')$.

\pfcase{1}
If $u$ and $v$ have the same parity, then $(\al',\be')=(\ve_1\al\ve_2,\ve_1\be\ve_2)\in\si$ for any Jones partition $\ve_2\in\JJ_B$ containing the blocks $\{u,1'\}$, $\{u+1,2'\}$, $\{v,3'\}$ and $H_2',\ldots,H_l'$.

\pfcase{2}
Next suppose $u$ is odd and $v$ is even.  Since $k$ is odd, we have $v<k$.  In this case, $(\al',\be')=(\ve_1\be\ve_2,\ve_1\al\ve_2)\in\si$ for any Jones partition $\ve_2\in\JJ_B$ containing the blocks $\{u,3'\}$, $\{v,1'\}$, $\{k,2'\}$ and $H_2',\ldots,H_l'$.

\pfcase{3}
Similarly, if $u$ is even and $v$ is odd, then $(\al',\be')=(\ve_1\be\ve_2,\ve_1\al\ve_2)\in\si$ for any Jones partition $\ve_2\in\JJ_B$ containing the blocks $\{1,2'\}$, $\{u,3'\}$, $\{v,1'\}$ and $H_2',\ldots,H_l'$.
\epf

Here is the analogue of Lemmas \ref{lem:R_trick_Bodd} and \ref{lem:R_trick_TLodd} for $\JJ$ and $\Jpm$:

\begin{lemma}\label{lem:R_trick_Jodd}
For any $\al,\be\in I_3=I_3(\Q)$ with $\al\sim\be$ and $(\al,\be)\not\in\L$, we have $\rho_{I_1}^{I_3}\sub\cg\al\be_{I_3}$.
\end{lemma}

\pf
By the proof of Lemma \ref{lem:R_trick_Jodd2}, we may assume that
\[
\al = \partVI1{H_1}{H_2}{H_l}1{H_1}{H_2}{H_l} \AND \be = \partVI1{H_1}{H_2}{H_l}3{1,2}{H_2}{H_l},
\]
where $\al,\be\in\TL_{2l+1}$ for some $l\geq1$, and where $H_i=\{2i,2i+1\}$ for all $i\geq1$.  Write $C=[2l+1]$.

Write $\tau = \rho_{I_1(\TL)}^{I_3(\TL)}$ for the $\rho_{I_1}$ congruence on the ideal $I_3(\TL)$ of the Temperley--Lieb category $\TL=\TL(\C)$.  Note that~$\tau$ is not a congruence on $I_3(\Q)$, but by Lemma \ref{lem:R_trick_TLodd}, we have
\[
\tau\sub \cg\al\be_{I_3(\TL)} \sub \cg\al\be_{I_3(\Q)} \sub \si.
\]

Now let $\ga\in D_1(\Q)$ be arbitrary, say with $\ga\in\Q_{A,B}$, and write $\ga=\partIV u{A_1}\cdots{A_s}v{B_1}\cdots{B_t}$.  Without loss of generality we may assume that $B=[k]$ where $k=1+2t$ (we may or may not have $t=l$).  
As usual, we may prove the lemma by showing that $(\ga,\de)\in\si$, where $\de=\partIV u{A_1}\cdots{A_s}1{H_1}\cdots{H_t}$.

First we claim that
\begin{equation}\label{eq:gaga'}
(\ga,\ga')\in\si \qquad\text{for some $\ga'\in\Q_{A,B}$ with $\ker(\ga')=\ker(\ga)$ and ${\codom(\ga')=\{1\}}$. }
\end{equation}
We prove \eqref{eq:gaga'} by descending induction on $v$.  If $v=1$ we just take $\ga'=\ga$, so assume instead that $v>1$.
Fix any $\ve_1\in\Q_{A,C}$ with $\ker(\ve_1)=\ker(\ga)$.
Without loss of generality, we may assume that $B_1=\{p,v-1\}$, and we note that the edge~$B_1'$ is not nested in (an annular drawing of) $\ga$.  We either have $p<v-1$ or else $p>v$.

\pfcase{1}
Suppose first that $p<v-1$.  Since $B_1'$ does not wrap around the back of the cylinder,~$p$ and $v-1$ have opposite parities.  Thus, $p<v-1<v$ are of alternating parities.  Put ${\ve_2 = 
\Big( 
{ \scriptsize \renewcommand*{\arraystretch}{1}
\begin{array} {\c|\c|\c|\c|\c|\cend}
 1 \:&\: 2 \:&\: 3 \:&\: H_2 \:&\: \cdots \:&\: H_l \\ \cline{4-6}
 p \:&\: v-1 \:&\: v \:&\: B_2 \:&\: \cdots \:&\: B_t 
\rule[0mm]{0mm}{2.7mm}
\end{array} 
}
\hspace{-1.5 truemm} \Big)}$, noting that $\ve_2\in\JJ_{C,B}$ (as $B_1'$ is not nested in $\ga$).  Then with $\ga''=\ve_1\al\ve_2$, we have $(\ga,\ga'') = (\ve_1\be\ve_2,\ve_1\al\ve_2)\in\si$.  But $\ga'' = \Big( 
{ \scriptsize \renewcommand*{\arraystretch}{1}
\begin{array} {\c|\c|\c|\c|\cend}
 u \:&\: A_1 \:& \multicolumn{2}{c|}{\cdots\cdots} &\: A_s \\ \cline{2-5}
 p \:&\: v-1,v \:&\: B_2 \:&\: \cdots \:&\: B_t
\rule[0mm]{0mm}{2.7mm}
\end{array} 
}
\hspace{-1.5 truemm} \Big)$.  Since $p<v$, we have inductively that $(\ga'',\ga')\in\si$ for some $\ga'$ with $\codom(\ga')=\{1\}$ and $\ker(\ga')=\ker(\ga'')=\ker(\ga)$.  By transitivity, we also have $(\ga,\ga')\in\si$, completing the proof of \eqref{eq:gaga'} in this case.

\pfcase{2}
Now suppose $p>v$.  Since $B_1'$ wraps around the back of the cylinder (and since $k$ is odd), $p$ and $v-1$ have the same parity.  Thus, $v-1<v<p$ have alternating parities.  
This time put ${\ve_2 = 
\Big( 
{ \scriptsize \renewcommand*{\arraystretch}{1}
\begin{array} {\c|\c|\c|\c|\c|\cend}
 1 \:&\: 2 \:&\: 3 \:&\: H_2 \:&\: \cdots \:&\: H_l \\ \cline{4-6}
 v \:&\: p \:&\: v-1 \:&\: B_2 \:&\: \cdots \:&\: B_t 
\rule[0mm]{0mm}{2.7mm}
\end{array} 
}
\hspace{-1.5 truemm} \Big)}$, again noting that $\ve_2\in\JJ_{C,B}$ (as $B_1'$ is not nested in $\ga$).  Then with $\ga''=\ve_1\be\ve_2$, we have $(\ga,\ga'') = (\ve_1\al\ve_2,\ve_1\be\ve_2)\in\si$.  But $\ga'' = \Big( 
{ \scriptsize \renewcommand*{\arraystretch}{1}
\begin{array} {\c|\c|\c|\c|\cend}
 u \:&\: A_1 \:& \multicolumn{2}{c|}{\cdots\cdots} &\: A_s \\ \cline{2-5}
 v-1 \:&\: v,p \:&\: B_2 \:&\: \cdots \:&\: B_t
\rule[0mm]{0mm}{2.7mm}
\end{array} 
}
\hspace{-1.5 truemm} \Big)$, and again the inductive assumption completes the proof of \eqref{eq:gaga'} in this case.

\aftercases
Now that we have established \eqref{eq:gaga'}, write $\ga'=\partIV u{A_1}\cdots{A_s}1{C_1}\cdots{C_t}$, so $(\ga,\ga')\in\si$.  Let
\[
\de_1 = \partIV 1{H_1}\cdots{H_t}1{C_1}\cdots{C_t} \AND \de_2 = \partIV 1{H_1}\cdots{H_t}1{H_1}\cdots{H_t},
\]
both from $\Q_B$.
Since $\de_1$ and $\de_2$ both contain the transversal $\{1,1'\}$, they both in fact belong to~$\TL_B$.
But $(\de_1,\de_2)\in \rho_{I_1(\TL)}^{I_3(\TL)}=\tau$, so since $\tau\sub\si$ (as shown near the beginning of the proof), we have $(\de_1,\de_2)\in\si$.  But then $(\ga',\de) = (\ga\de_1,\ga\de_2)\in\si$.  (The element $\de$ was defined near the start of the proof.)  Thus, $(\ga,\de)\in\si$ by transitivity.  As noted above, this completes the proof.
\epf

We are now ready to prove the main result of this section.  For the statement, note that ${I_1(\Jpm)=I_1(\JJ)}$.  Also, the $\R$- and $\L$-classes in $I_1$ are indexed by the allowable (co)kernels; by Lemma \ref{lem:PEqA}, such (co)kernels are planar, and they must in fact be 2-equivalences.  Unlike the case with $\TL$ and $\TLpm$, however, \emph{any} planar 2-equivalence arises as a (co)kernel of an element of $I_1$; indeed, although the unique block of size $1$ corresponds to a transversal, if this block is nested by other blocks, then these blocks can wrap around the back of the cylinder.  Thus, the $\R$- and $\L$-classes in $D_1(\Q_A)$ are indexed by $\OrdP{\Eq_2(A)}$, the set of all planar 2-equivalences on $A$.

\begin{thm}\label{thm:Jodd}
Let $\C$ be a non-empty set of finite sets of odd cardinality, let $\Q$ be either the Jones category $\JJ=\JJ(\C)$ or the anti-Jones category $\Jpm=\Jpm(\C)$, and keep the above notation. 
\ben
\item \label{cTLodd1}  The minimal ideal $I_1$ is a partial rectangular band, so its congruence lattice is described by Proposition \ref{prop:PRB}.  In particular, $\Cong(I_1)$ is isomorphic to the lattice direct product
\[
\prod_{A\in\C}\Eq(\OrdP{\Eq_2(A)})\times\prod_{A\in\C}\Eq(\OrdP{\Eq_2(A)}).
\]
\item \label{cTLodd2} For odd $r\geq3$ (including $r=\om$ if $\xi=\om$), 
\begin{align*}
\Cong(I_r) = \{\De_{I_r},\ \lam_{I_1}^{I_r},\ \rho_{I_1}^{I_r},\ R_{I_1}^{I_r}\}
\cup\set{R_{I_q,N}^{I_r}}{1\leq q<r,\ q\text{ odd},\ N\normal\G_{q+2}}\cup\{\nab_{I_r}\}.
\end{align*}
\een
\end{thm}

\pf
The proof is virtually identical to that of Theorems \ref{thm:Bodd} and \ref{thm:TLodd}.  We use Lemmas~\ref{lem:P2J},~\ref{lem:R_trick_Jodd} and \ref{lem:N_Bodd} to support Proposition \ref{prop:small01} in the $r=3$ case.
\epf

\begin{rem}
Figure \ref{fig:J} shows the lattices $\Cong(I_9(\JJ))$ and $\Cong(I_9(\Jpm))$.  While these diagrams may give the impression that the upper parts of the lattices $\Cong(I_r(\JJ))$ and $\Cong(I_r(\Jpm))$ are always chains when $r$ is odd, this is not actually the case.  Indeed, as in Remark \ref{rem:E2} (cf.~Figure~\ref{fig:E2}), these upper parts consist of copies of $\NN(\CC_q)$ or $\NN(\DD_q)$ for $q=3,5,\ldots,r$ stacked on top of each other.  These normal subgroup lattices are chains for $q\leq 13$, but, of course, not in general.
\end{rem}

\begin{figure}[ht]
\begin{center}
\begin{tikzpicture}[scale=0.7,inner sep=0pt]
\begin{scope}[shift={(0,0)}]
\draw[rounded corners=4,line width=10,blue!25,cap=round](0,2)--(0,3);
\draw[rounded corners=4,line width=10,green!25,cap=round]
(0,4)--(0,5)
(0,6)--(0,7)
(0,8)--(0,10)
;
\draw(0,2)--(1,1)--(0,0)--(-1,1)--(0,2)--(0,11);
\foreach \x/\y in {1/1,-1/1,0/3,0/5,0/7,0/9,0/10} {\colvert{\x}{\y}{white}}
\foreach \y in {0,2,4,6,8,11} {\colvert{0}{\y}{black}}
\end{scope}
\begin{scope}[shift={(5,0)}]
\draw[rounded corners=4,line width=10,blue!25,cap=round](0,2)--(0,4);
\draw[rounded corners=4,line width=10,green!25,cap=round]
(0,5)--(0,7)
(0,8)--(0,10)
(0,11)--(0,14)
;
\draw(0,2)--(1,1)--(0,0)--(-1,1)--(0,2)--(0,15);
\foreach \x/\y in {1/1,-1/1} {\colvert{\x}{\y}{white}}
\foreach \y in {3,4,6,7,9,10,11,12,13,14} {\colvert{0}{\y}{white}}
\foreach \y in {0,2,5,8,11,15} {\colvert{0}{\y}{black}}
\end{scope}
\begin{scope}[shift={(10,0)}]
\draw[rounded corners=4,line width=10,green!25,cap=round]
(0,7)--(-1,8) (-1,8)--(0,9) (0,9)--(1,8) (1,8)--(0,7)
(0,10)--(0,13)
;
\draw[rounded corners=4,line width=10,blue!13,cap=round] (0,5)--(0,6);
\foreach \x/\y in {0/0,0/-1,-1/-1,-1/-2} {\draw[rounded corners=4,line width=10,red!25,cap=round](3+\x,2+\y)--(2+\x,3+\y);}
\foreach \x/\y in {0/0,0/-1,-1/-1,-1/-2} {\draw[rounded corners=4,line width=10,blue!25,cap=round](1+\x,4+\y)--(0+\x,5+\y);}
\fill[green!25] (0,7)--(-1,8)--(0,9)--(1,8)--(0,7);
\foreach \x in {0,1,2,3} {\draw(2-\x,0+\x)--(3-\x,1+\x)--(3-\x,2+\x)--(2-\x,1+\x)--(2-\x,0+\x);}
\draw(0,5)--(0,7)--(-1,8)--(0,9)--(0,14) (0,7)--(1,8)--(0,9);
\foreach \x/\y in {0/0,0/1,1/1,1/2} {\draw(2+\x,0+\y)--(-1+\x,3+\y);}
\foreach \x in {1,2,3} {\colvert{1}{\x}{white}}
\foreach \x in {7,9,11,12,13} {\colvert{0}{\x}{white}}
\foreach \x/\y in {0/2,0/3,2/0,2/1,2/2,2/3,3/1,3/2,-1/3,-1/4,0/4,0/5,0/6,-1/8,1/8} {\colvert{\x}{\y}{white}}
\foreach \x in {7,10,14} {\colvert0{\x}{black}}
\foreach \x/\y in {3/2,1/4,2/0} {\colvert{\x}{\y}{black}}
\end{scope}
\begin{scope}[shift={(17,0)}]
\fill[blue!13] 
(0,5)--(-1,6)--(0,7)--(1,6)--(0,5)
;
\draw[rounded corners=4,line width=10,blue!13,cap=round]
(-1,6)--(0,7) (0,7)--(1,6) (1,6)--(0,5)
;
\draw[rounded corners=4,line width=10,green!25,cap=round]
(0,8)--(1,9) (1,9)--(-1,11) (-1,11)--(-2,10) (-2,10)--(0,8)
(-1,12)--(-1,14) (-1,14)--(-2,15) (-2,15)--(-1,16) (-1,16)--(0,15) (0,15)--(-1,14)
;
\foreach \x/\y in {0/0,0/-1,-1/-1,-1/-2} {\draw[rounded corners=4,line width=10,red!25,cap=round](3+\x,2+\y)--(2+\x,3+\y);}
\foreach \x/\y in {0/0,0/-1,-1/-1,-1/-2} {\draw[rounded corners=4,line width=10,blue!25,cap=round](1+\x,4+\y)--(-1+\x,6+\y);}
\fill[green!25] 
(0,8)--(-2,10)--(-1,11)--(1,9)--(0,8)
(-1,14)--(-2,15)--(-1,16)--(0,15)--(-1,14)
;
\foreach \x in {0,1,2,3,4} {\draw(2-\x,0+\x)--(3-\x,1+\x)--(3-\x,2+\x)--(2-\x,1+\x)--(2-\x,0+\x);}
\draw(0,5)--(0,7)--(-1,6) (0,5)--(1,6)--(0,7)--(0,8)--(1,9)--(-1,11)--(-2,10)--(0,8) (-1,9)--(-1,17) (-1,14)--(0,15)--(-1,16)--(-2,15)--(-1,14) (-1,9)--(0,10);
\foreach \x/\y in {0/0,0/1,1/1,1/2} {\draw(2+\x,0+\y)--(-2+\x,4+\y);}
\foreach \x in {1,2,3} {\colvert{1}{\x}{white}}
\foreach \x in {7,10} {\colvert{0}{\x}{white}}
\foreach \x in {9,10,11,13,14,15,16} {\colvert{-1}{\x}{white}}
\foreach \x/\y in {0/2,0/3,2/0,2/1,2/2,2/3,3/1,3/2,-1/3,-1/4,0/4,0/5,0/6,-1/9,1/9,-2/15,0/15,-2/10,-2/4,-2/5,-1/5,-1/6,1/6} {\colvert{\x}{\y}{white}}
\foreach \x in {8} {\colvert0{\x}{black}}
\foreach \x/\y in {3/2,1/4,2/0,-1/12,-1/17} {\colvert{\x}{\y}{black}}
\end{scope}
\end{tikzpicture}
\end{center}
\vspace{-7mm}
\caption{Left to right:  Hasse diagrams of $\Cong(I_9(\JJ))$, $\Cong(I_9(\Jpm))$, $\Cong(I_8(\JJ))$ and $\Cong(I_8(\Jpm))$.  See Theorems and \ref{thm:Jodd} and~\ref{thm:Jeven}.}
\label{fig:J}
\end{figure}

\subsection{Congruences on ideals of $\JJ$ and $\Jpm$: the even case}\label{subsect:Jeven}

We now assume every member of $\C$ has even cardinality, and define $\xi={\LSUB}\bigset{|A|}{A\in\C}$.  Throughout this section, $\Q$ denotes either $\JJ=\JJ(\C)$ or $\Jpm=\Jpm(\C)$.  Even though the congruence lattices of $\Q$ and of its ideals are more involved than in the odd case (see Figure~\ref{fig:J}), the actual proof of the classification result (Theorem \ref{thm:Jeven}) is easier than its odd counterpart (Theorem~\ref{thm:Jodd}), for the following reason.
As we showed in Section \ref{subsect:QP}, the (co)kernel of any $\al\in\Jpm$ is a planar equivalence.  It follows from this that the minimal ideal $I_0(\Q)=D_0(\Q)$ is in fact equal to $I_0(\TL)=D_0(\TL)$.  The proof of the next lemma is therefore virtually identical to that of Lemma~\ref{lem:R_trick_TLpm}.  In the statement, we write $f:I_2\to I_0:\al\mt\wh\al$ for the retraction, which is the appropriate restriction of \eqref{eq:retract_B}.

\begin{lemma}\label{lem:R_trick_Jeven}
For any $\al,\be\in I_4=I_4(\Q)$ with $\al\sim\be$ and $(\al,\be)\not\in\L\cup\L f^{-1}$, we have $\rho_{I_0}^{I_4}\sub\cg\al\be_{I_4}$.  \epfres
\end{lemma}

In order to state the main result of this section, we first fix some group $\H$-classes in $\Q$, and discuss the (retractable) IN-pairs.  For $n\geq2$, recall that we write $\ga_n=\ga_{[n]}=(1,2,\ldots,n)$ and $\de_n=\de_{[n]}=(1,n)(2,n-1)\cdots$, as in \eqref{eq:gade}.  We then define
\[
\CC_n = \la\ga_n\ra \AND \DD_n = \la\ga_n,\de_n\ra,
\]
the cyclic group of order $n$ and the dihedral group of order $2n$, respectively, for (even) $n\geq4$; since $\ga_2=\de_2$, we have $\CC_2=\DD_2$.  For any even $2\leq q<\xi$, and continuing to use the $\pi^\natural$ notation, we define
\[
\G_q = \begin{cases}
\CC_q^\natural &\text{if $\Q=\JJ$}\\
\DD_q^\natural &\text{if $\Q=\Jpm$.}
\end{cases}
\]
Then $\G_q$ is a group $\H$-class contained in $D_q$.  In what follows, we will also write $\ve_q=\id_q^\natural$, so that $\ve_q$ is the identity element of the group $\G_q$.

An ideal $I_r$ ($r\geq4$) contains the IN-pairs $(I_q,N)$ for even $0\leq q<r$ and $N\normal\G_{q+2}$.  These yield the congruences 
\[
R_N^{I_r}=R_{I_q,N}^{I_r}=R_{I_q}^{I_r}\cup\nu_{N},
\]
each involving the $\nu_{N}$ relation, membership of which is again determined by the $\phi(\al,\be)$ parameter, as defined in Section \ref{subsect:P}.

We also have retractable IN-pairs in two forms:
\bit
\item $(I_0,N)$ for $N\normal\G_2$, and 
\item $(I_2,N)$ for $N\normal\G_4$ satisfying $N\leq K^\natural$, where $K\normal \S_4$ is the Klein 4-group.
\eit
When $\Q=\JJ$, the only normal subgroups of $\G_4$ contained in $K^\natural$ are $\{\id_4^\natural\}$ and $C^\natural$, where $C=\la\ga_4^2\ra = \big\{{\id}_4,(1,3)(2,4)\big\}$.  So the above list yields four retractable IN-pairs in $I_r(\JJ)$ for even~$r\geq4$:
\begin{align*}
&(I_0,\{\ve_2\}) \COMMA
(I_0,\G_2) \COMMA
(I_2,\{\ve_4\}) \COMMA
(I_2,C^\natural).
\intertext{When $\Q=\Jpm$, however, there are \emph{three} normal subgroups of $\G_4$ contained in $K^\natural$:  $\{\id_4^\natural\}$, $C^\natural$ and $K^\natural=\la\ga_4^2,\de_4\ra^\natural$ itself.  So the above list yields five retractable IN-pairs in $I_r(\Jpm)$ for even $r\geq4$:}
&(I_0,\{\ve_2\}) \COMMA
(I_0,\G_2) \COMMA
(I_2,\{\ve_4\}) \COMMA
(I_2,C^\natural) \COMMA
(I_2,K^\natural).
\end{align*}
These retractable IN-pairs yield the usual families of $R/\lam/\rho/\mu$ congruences.  But note that for both categories, there are normal subgroups $N\normal\G_4$ for which $(I_2,N)$ is not retractable.  When $\Q=\JJ$, the only such $N$ is $\G_4=\CC_4^\natural$ itself.  When $\Q=\Jpm$, there is $\G_4=\DD_4^\natural$ itself, and two others: $\CC_4^\natural$, and a normal subgroup we will denote by $\ol K^\natural$, where
\[
\overline{K} = \big\{{\id}_4,(1,3),(2,4),(1,3)(2,4)\big\} \normal\DD_4.
\]
(Note that $\overline{K}$ is isomorphic to $K$.)  

Figures \ref{fig:Jeven} and \ref{fig:Jpmeven} show the congruences coming from the (retractable) IN-pairs just discussed, in the case $r=4$, as well as the universal congruence $\nab_{I_4}$, with the usual abbreviations.  
Note that the lattices in Figures \ref{fig:Jeven} and \ref{fig:Jpmeven} are instances of the general case pictured in Figure \ref{fig:small012}.
For $r=8$, but without labels on the congruences, see Figure~\ref{fig:J}.  The ``non-chainlike'' structure of the upper part of the lattices is due to the complexity of the normal subgroup lattices $\NN(\G_q)$.

\begin{figure}[ht]
\begin{center}
\scalebox{0.75}{
\begin{tikzpicture}[scale=.95]
\draw[rounded corners=2em,line width=3.5em,blue!15,cap=round] (0,8)--(0,10);
\foreach \x/\y in {0/0,0/2,3/2,3/4} {\draw[rounded corners=2em,line width=3.5em,red!25,cap=round] (6+\x,-2+\y)--(3+\x,0+\y);}
\foreach \x/\y in {0/0,0/2,3/2,3/4} {\draw[rounded corners=2em,line width=3.5em,blue!35,cap=round] (0+\x,2+\y)--(-3+\x,4+\y);}
\draw(1.13,12.0) node {$=\nab_{I_4}$};
\node[rounded corners,rectangle,draw,fill=black!30,ultra thick] (S4) at (0,12) {$R_{I_4}^{I_4}$};
\node[rounded corners,rectangle,draw,fill=white] (A4) at (0,10) {$R_{\CC_4}^{I_4}$};
\node[rounded corners,rectangle,draw,fill=white] (krl) at (0,8) {$R_{C}^{I_4}$};
\node[rounded corners,rectangle,draw,fill=white] (kr) at (-3,6) {$\lambda_{C}^{I_4}$};
\node[rounded corners,rectangle,draw,fill=white] (kl) at (0,6) {$\rho_{C}^{I_4}$};
\node[rounded corners,rectangle,draw,fill=white] (k) at (-3,4) {$\mu_{C}^{I_4}$};
\node[rounded corners,rectangle,draw,fill=black!30,ultra thick] (erl) at (3,6) {$R_{I_2}^{I_4}$};
\node[rounded corners,rectangle,draw,fill=white] (er) at (0,4) {$\lambda_{I_2}^{I_4}$};
\node[rounded corners,rectangle,draw,fill=white] (el) at (3,4) {$\rho_{I_2}^{I_4}$};
\node[rounded corners,rectangle,draw,fill=white] (e) at (0,2) {$\mu_{I_2}^{I_4}$};
\node[rounded corners,rectangle,draw,fill=white] (nrl) at (6,4) {$R_{\CC_2}^{I_4}$};
\node[rounded corners,rectangle,draw,fill=white] (nr) at (3,2) {$\lambda_{\CC_2}^{I_4}$};
\node[rounded corners,rectangle,draw,fill=white] (nl) at (6,2) {$\rho_{\CC_2}^{I_4}$};
\node[rounded corners,rectangle,draw,fill=white] (n) at (3,0) {$\mu_{\CC_2}^{I_4}$};
\node[rounded corners,rectangle,draw,fill=black!30,ultra thick] (rl) at (9,2) {$R_{I_0}^{I_4}$};
\node[rounded corners,rectangle,draw,fill=white] (r) at (6,0) {$\lambda_{I_0}^{I_4}$};
\node[rounded corners,rectangle,draw,fill=white] (l) at (9,0) {$\rho_{I_0}^{I_4}$};
\node[rounded corners,rectangle,draw,fill=black!30,ultra thick] (D) at (6,-2) {$\mu_{I_0}^{I_4}$};
\draw(7.3,-2.05) node {$=\Delta_{I_4}$};
\draw
(D)--(n)--(e)--(k)
(l)--(nl)--(el)--(kl)
(k)--(kl) (e)--(el) (n)--(nl) (D)--(l) 
(k)--(kr) (e)--(er) (n)--(nr) (D)--(r) 
(l)--(rl) (nl)--(nrl) (el)--(erl) (kl)--(krl)
(rl)--(nrl)--(erl)--(krl)--(A4)--(S4)
;
\fill[blue!35] (-1.5,5)circle(.15);
\fill[white] (1.5,3)circle(.15);
\fill[red!25] (4.5,1)circle(.15);
\fill[blue!35] (1.5,5)circle(.15);
\fill[white] (4.5,3)circle(.15);
\fill[red!25] (7.5,1)circle(.15);
\draw
(r)--(nr)--(er)--(kr)
(r)--(rl) (nr)--(nrl) (er)--(erl) (kr)--(krl) 
;
\end{tikzpicture}
}
\caption{Hasse diagram of $\Cong(I_4)$, where $I_4=I_4(\JJ)$ in the case that all members of~$\C$ have even cardinality.  We abbreviate $R_{I_{q-2},N^\natural}^{I_4}=R_{N}^{I_4}$, and similarly for $\lam$, $\rho$ and $\mu$. See Theorem~\ref{thm:Jeven}.}
\label{fig:Jeven}
\end{center}
\end{figure}

\begin{figure}[ht]
\begin{center}
\scalebox{0.75}{
\begin{tikzpicture}[scale=.95]
\draw[rounded corners=2em,line width=3.5em,blue!15,cap=round] (0,8)--(3,10) (3,10)--(0,12) (0,12)--(-3,10) (-3,10)--(0,8);
\fill[blue!15] (0,8)--(3,10)--(0,12)--(-3,10)--(0,8);
\foreach \x/\y in {0/0,0/2,3/2,3/4} {\draw[rounded corners=2em,line width=3.5em,red!25,cap=round] (6+\x,-2+\y)--(3+\x,0+\y);}
\foreach \x/\y in {0/0,0/2,3/2,3/4} {\draw[rounded corners=2em,line width=3.5em,blue!35,cap=round] (0+\x,2+\y)--(-6+\x,6+\y);}
\draw(1.13,14.0) node {$=\nab_{I_4}$};
\node[rounded corners,rectangle,draw,fill=black!30,ultra thick] (i) at (0,14) {$R_{I_4}^{I_4}$};
\node[rounded corners,rectangle,draw,fill=white] (d) at (0,12) {$R_{\DD_4}^{I_4}$};
\node[rounded corners,rectangle,draw,fill=white] (n1) at (0,10) {$R_{\overline{K}}^{I_4}$};
\node[rounded corners,rectangle,draw,fill=white] (n2) at (3,10) {$R_{\CC_4}^{I_4}$};
\node[rounded corners,rectangle,draw,fill=white] (kkrl) at (-3,10) {$R_{K}^{I_4}$};
\node[rounded corners,rectangle,draw,fill=white] (kkr) at (-6,8) {$\lambda_{K}^{I_4}$};
\node[rounded corners,rectangle,draw,fill=white] (kkl) at (-3,8) {$\rho_{K}^{I_4}$};
\node[rounded corners,rectangle,draw,fill=white] (kk) at (-6,6) {$\mu_{K}^{I_4}$};
\node[rounded corners,rectangle,draw,fill=white] (krl) at (0,8) {$R_{C}^{I_4}$};
\node[rounded corners,rectangle,draw,fill=white] (kr) at (-3,6) {$\lambda_{C}^{I_4}$};
\node[rounded corners,rectangle,draw,fill=white] (kl) at (0,6) {$\rho_{C}^{I_4}$};
\node[rounded corners,rectangle,draw,fill=white] (k) at (-3,4) {$\mu_{C}^{I_4}$};
\node[rounded corners,rectangle,draw,fill=black!30,ultra thick] (erl) at (3,6) {$R_{I_2}^{I_4}$};
\node[rounded corners,rectangle,draw,fill=white] (er) at (0,4) {$\lambda_{I_2}^{I_4}$};
\node[rounded corners,rectangle,draw,fill=white] (el) at (3,4) {$\rho_{I_2}^{I_4}$};
\node[rounded corners,rectangle,draw,fill=white] (e) at (0,2) {$\mu_{I_2}^{I_4}$};
\node[rounded corners,rectangle,draw,fill=white] (nrl) at (6,4) {$R_{\DD_2}^{I_4}$};
\node[rounded corners,rectangle,draw,fill=white] (nr) at (3,2) {$\lambda_{\DD_2}^{I_4}$};
\node[rounded corners,rectangle,draw,fill=white] (nl) at (6,2) {$\rho_{\DD_2}^{I_4}$};
\node[rounded corners,rectangle,draw,fill=white] (n) at (3,0) {$\mu_{\DD_2}^{I_4}$};
\node[rounded corners,rectangle,draw,fill=black!30,ultra thick] (rl) at (9,2) {$R_{I_0}^{I_4}$};
\node[rounded corners,rectangle,draw,fill=white] (r) at (6,0) {$\lambda_{I_0}^{I_4}$};
\node[rounded corners,rectangle,draw,fill=white] (l) at (9,0) {$\rho_{I_0}^{I_4}$};
\node[rounded corners,rectangle,draw,fill=black!30,ultra thick] (D) at (6,-2) {$\mu_{I_0}^{I_4}$};
\draw
(D)--(n)--(e)--(k)
(l)--(nl)--(el)--(kl)
(k)--(kl) (e)--(el) (n)--(nl) (D)--(l) 
(k)--(kr) (e)--(er) (n)--(nr) (D)--(r) 
(l)--(rl) (nl)--(nrl) (el)--(erl) (kl)--(krl)
(rl)--(nrl)--(erl)--(krl) 
(k)--(kk) (kl)--(kkl) (kr)--(kkr) (krl)--(kkrl)
(kk)--(kkl) (kk)--(kkr)--(kkrl)--(kkl)
(krl)--(n2)--(d)--(n1)--(krl)
(kkrl)--(d)--(i)
;
\fill[blue!35] (-4.5,7)circle(.15);
\fill[blue!35] (-1.5,7)circle(.15);
\fill[blue!35] (-1.5,5)circle(.15);
\fill[blue!35] (1.5,5)circle(.15);
\fill[white] (1.5,3)circle(.15);
\fill[white] (4.5,3)circle(.15);
\fill[red!25] (4.5,1)circle(.15);
\fill[red!25] (7.5,1)circle(.15);
\draw
(r)--(nr)--(er)--(kr)
(r)--(rl) (nr)--(nrl) (er)--(erl) (kr)--(krl) 
(kr)--(kkr)
;
\draw(7.3,-2.05) node {$=\Delta_{I_4}$};
\end{tikzpicture}
}
\caption{Hasse diagram of $\Cong(I_4)$, where $I_4=I_4(\Jpm)$ in the case that all members of~$\C$ have even cardinality.  We abbreviate $R_{I_{q-2},N^\natural}^{I_4}=R_{N}^{I_4}$, and similarly for $\lam$, $\rho$ and $\mu$.  See Theorem~\ref{thm:Jeven}.}
\label{fig:Jpmeven}
\end{center}
\end{figure}

For the first part of the following statement, keep in mind that $D_0(\JJ)=D_0(\Jpm)=D_0(\TL)$, so that the $\R$- and $\L$-classes in $D_0$ are again indexed by the planar 2-equivalences on the sets from $\C$.

\begin{thm}\label{thm:Jeven}
Let $\C$ be a non-empty set of finite sets of even cardinality, let $\Q$ be either the Jones category $\JJ=\JJ(\C)$ or the anti-Jones category $\Jpm=\Jpm(\C)$, and keep the above notation.
\ben
\item \label{cTLeven1}  The minimal ideal $I_0$ is a partial rectangular band, so its congruence lattice is described by Proposition \ref{prop:PRB}.  In particular, $\Cong(I_0)$ is isomorphic to the lattice direct product
\[
\prod_{A\in\C}\Eq(\OrdP{\Eq_2(A)})\times\prod_{A\in\C}\Eq(\OrdP{\Eq_2(A)}).
\]
\item \label{cTLeven2}  The ideal $I_2$ is retractable, every congruence on $I_0$ is liftable to~$I_2$, and
\[
\Cong(I_2)=\set{\tau\cup\nu_N}{\tau\in\Cong(I_0),\ N\normal\G_2} \cup \set{\th_{I_0,\tau}^{I_2}}{\tau\in\Cong(I_0)}.
\]
In particular, $\Cong(I_2)$ is isomorphic to the lattice direct product $\Cong(I_0)\times\bthree$.
\item \label{cTLeven3} For even $r\geq4$ (including $r=\om$ if $\xi=\om$), 
\begin{align*}
\Cong(I_r) = 
\set{\mu_{I_0,N}^{I_r},\ \lam_{I_0,N}^{I_r},\ \rho_{I_0,N}^{I_r},\ &R_{I_0,N}^{I_r}}{N\normal\G_2} \\[3truemm]
{}\cup\set{\mu_{I_2,N}^{I_r},\ \lam_{I_2,N}^{I_r},\ \rho_{I_2,N}^{I_r},\ &R_{I_2,N}^{I_r}}{N\normal\G_4,\ N\leq K^\natural} \\[3truemm]
\qquad\qquad\qquad\qquad {}\cup\set{&R_{I_q,N}^{I_r}}{2\leq q<r,\ q\text{ even},\ N\normal\G_{q+2}}\cup\{\nab_{I_r}\}.
\end{align*}
\een
\end{thm}

\pf
The proof is virtually identical to that of Theorems \ref{thm:Beven} and \ref{thm:TLeven}.  We use Lemmas~\ref{lem:P2J},~\ref{lem:R_trick_Jeven} and \ref{lem:N_Beven} to support Proposition \ref{prop:small012} in the $r=4$ case.
\epf

\section{Categories and partial semigroups with $\H$-congruences}\label{sect:H}

In Sections \ref{sect:T}--\ref{sect:J}, we have seen many applications of the stacking mechanism established by Theorem \ref{thm:E2}.
However this theorem does not apply to two of our remaining examples:  the linear category~$\LL$ (Section \ref{sect:L}) and the partial braid category~$\IB$ (Section~\ref{sect:IB}).  Although these categories are both chains of ideals (as in Definition \ref{defn:chain}), their ideals $I_r$ do not satisfy the property $\Sep{I_r}$.
This, in turn, is a consequence of the presence of \emph{$\H$-congruences}: i.e., congruences contained in Green's $\H$-relation.

We note that we have already encountered some $\H$-congruences:
$\mu_{I_0,\S_2}$ is the unique non-trivial $\H$-congruence
in the categories~$\B$,~$\TLpm$,~$\JJ$ and~$\Jpm$ in the case that all underlying sets are of even size.
However, this congruence occurs sufficiently ``low down'' in the chain of ideals that it is subsumed in the computation
of the bottom portion of the lattice, $\Cong(I_k)$ in the notation of Theorem \ref{thm:E2}.

In this section we undertake another theoretical development that will enable us to deal with categories such as~$\LL$ and~$\IB$.
The development again builds on Section \ref{sect:IE}, and yields extensions of the results from Section \ref{sect:chains}
with $\Sep{I_r}$ replaced by a weaker condition we call $\Sepz{I_r}$.
We begin in Section \ref{subsect:H} by discussing $\H$-congruences, the main result being Proposition \ref{prop:ThN}, which describes all such congruences on an arbitrary stable, regular partial semigroup.  
Then in Section~\ref{subsect:chains2}, we return to chains of ideals and prove the above-mentioned extension  of Theorem \ref{thm:E2},
which is Theorem \ref{thm:E4}.
We also establish a certain specialisation of Theorem \ref{thm:E4} in Theorem~\ref{thm:E5}, which
is in fact sufficient to deal with our intended applications.

\subsection{$\H$-congruences}\label{subsect:H}

Throughout Section \ref{subsect:H}, $S$ denotes a fixed partial semigroup.  
By an \emph{$\H$-congruence} on $S$ we mean a congruence that is contained in $\H$.  

A well-known result of Lallement \cite{Lallement1967} states that a congruence $\si$ on a regular semigroup is contained in $\H$ if and only if it is \emph{idempotent-separating}, meaning that distinct idempotents are never $\si$-related.  Although it is not essential for our purposes, we give a short argument, communicated to us by John Meakin, that these concepts are equivalent for partial semigroups as well.  As with ordinary semigroups, an \emph{inverse} of an element $x$ of the partial semigroup $S$ is an element $a\in S$ for which $x=xax$ and $a=axa$.  

\begin{lemma}\label{lem:id-sep}
Let $\si$ be a congruence on a regular partial semigroup $S$.  Then $\si$ is an $\H$-congruence if and only if it is idempotent-separating.
\end{lemma}

\pf
Since an $\H$-class contains at most one idempotent, the forwards implication is trivial.  Conversely, suppose $\si$ is idempotent-separating, and let $(x,y)\in\si$.  We must show that ${(x,y)\in\H}$.  We just show that $(x,y)\in\R$, and $(x,y)\in\L$ will follow by symmetry.  Let $a$ be an inverse of $x$, and note that $(xa,ya)\in\si$ with $xa$ an idempotent.  Then $ya\mathrel\si xa =(xa)^2 \mathrel\si (ya)^2$.  Now let $b$ be an inverse of $(ya)^2$.  It is easy to check that $(ya)b(ya)$ is an idempotent.  But then
\[
xa \mathrel\si ya \mathrel\si (ya)^2 = (ya)^2b(ya)^2 \mathrel\si (ya)b(ya).
\]
Since $xa$ and $(ya)b(ya)$ are both idempotents, it follows that $xa=(ya)b(ya)$, and from this that ${x=(xa)x=(ya)b(ya)x\in yS}$.  A symmetrical argument shows that $y\in xS$, and so $(x,y)\in\R$, as required.
\epf

Since $\H$ is an equivalence, the join and meet of any family of $\H$-congruences is an $\H$-congruence, so the set $\CongH(S)$ of all $\H$-congruences is a (complete) sublattice of $\Cong(S)$.  The join of all $\H$-congruences is the maximum $\H$-congruence, and this will be denoted by $\ze=\ze_S$.  (The maximum idempotent-separating congruence on a regular semigroup is typically denoted $\mu=\mu_S$ \cite{Hall1973}, but we use $\ze$ to avoid confusion with our other frequent use of the $\mu$ symbol.)  By analogy with semigroup terminology, we say the partial semigroup $S$ is \emph{fundamental} if $\ze_S=\De_S$.  In general, the quotient $S/\ze_S$ is fundamental.

In the next statement, recall that an element $x$ of a partial semigroup $S$ is self-divisible if $x=ux=xv$ for some $u,v\in S$.  If $S$ is regular or a category, then every element of $S$ is self-divisible.

\begin{lemma}\label{lem:zeST}
If $S$ is a partial semigroup in which every element is self-divisible, and if~$T$ is an ideal extension of~$S$, then 
\ben
\item \label{zeST1} every $\H^S$-congruence on $S$ is liftable to $T$,
\item \label{zeST2} $\ze_S=\ze_T\restr_S$.
\een
\end{lemma}

\pf
\firstpfitem{\ref{zeST1}}
Let $\si$ be an $\H^S$-congruence on $S$, and let $(x,y)\in\si$ and $a\in T$.  We need to show that
\[
\br(a)=\bd(x) \implies (ax,ay)\in\si \AND \bd(a)=\br(x) \implies (xa,ya)\in\si.
\]
It is sufficient to establish the first, so suppose $\br(a)=\bd(x)$.  As $x\in S$ is self-divisible, we have $x=ux$ for some $u\in S$.  Since $(x,y)\in\si\sub\H\sub\R$, it follows that $y=uy$ also (indeed, $y=xb$ for some $b\in S^1$, and then $uy=uxb=xb=y$).  But then $(ax,ay)=((au)x,(au)y)\in\si$, since $au\in S$ (as $S$ is an ideal) and $\si$ is a congruence on $S$.  

\pfitem{\ref{zeST2}}
First note that $\ze_T\restr_S$ is a congruence on $S$.  We also have ${\ze_T\restr_S\sub\H^T\restr_S=\H^S}$ by Lemma~\ref{lem:SD}\ref{SD1}.  This shows that $\ze_T\restr_S$ is an $\H^S$-congruence on $S$, so it follows that $\ze_T\restr_S\sub\ze_S$.

Conversely, part \ref{zeST1} tells us that $\ze_S\cup\De_T$ is a congruence on $T$, and it is in fact an $\H^T$-congruence.  It follows that $\ze_S\cup\De_T\sub\ze_T$, and so $\ze_S\sub\ze_T\restr_S$.
\epf

In what follows, it will be important to understand the $\H$-congruences on stable, regular partial semigroups.  Thus, for the remainder of Section \ref{subsect:H}, we assume $S$ is stable and regular.  Let the $\D=\J$-classes of $S$ be $\set{D_q}{q\in Q}$.  Define a partial order $\leq$ on $Q$ by $p\leq q \iff D_p\leq D_q$ in the usual ordering on $\J$-classes.  For each $q\in Q$, fix some group $\H$-class $G_q$ in $D_q$.

Since  $\ze\restr_{D_q}\sub\H$ for each $q$, Lemma \ref{lem:nu_trick} gives $\ze\restr_{D_q} = \nu_{Z_q}$ for some $Z_q\normal G_q$.  Thus,
\[
\ze=\bigcup_{q\in Q}\nu_{Z_q}.
\]
Similarly, for any $\H$-congruence $\si$, noting that $\si\sub\ze$, we have 
\begin{equation}\label{eq:Hcong}
\si=\bigcup_{q\in Q}\nu_{N_q} \quad\text{for some $Q$-tuple $\bN=(N_q)_{q\in Q}$, where $N_q\normal G_q$ and $N_q\leq Z_q$ for each $q$.}
\end{equation}
Not all such $Q$-tuples $\bN$ give rise to congruences of the form \eqref{eq:Hcong}, however.  To classify those that do, we require the following concept.

Fix some $p,q\in Q$ with $p\leq q$, and some $N\normal G_q$ with $N\leq Z_q$.  Consider the relation $\nu_N^\sharp\restr_{D_p}$: i.e., the restriction of the congruence $\nu_N^\sharp\in\Cong(S)$ to the $\J=\D$-class $D_p$.  Note that
\[
N\leq Z_q \ \ \implies\ \  \nu_N\sub\nu_{Z_q}\sub\ze \ \ \implies\ \  \nu_N^\sharp\sub\ze \ \ \implies\ \   
\nu_N^\sharp\restr_{D_p}\subseteq \ze\restr_{D_p}\sub\H.
\]
By Lemma \ref{lem:nu_trick},
there exists a normal subgroup $N^{\downarrow_p}\unlhd G_p$ such that
$N^{\downarrow_p}\leq Z_p$ and
\[
\nu_N^\sharp\restr_{D_p} = \nu_{N^{\downarrow_p}}.
\]
For $r\in Q$ with $r\not\leq q$, we have $\nu_N^\sharp\restr_{D_r}=\De_{D_r}$ because 
$\nu_N^\sharp\sub R_I$ for the ideal $I=\bigcup_{p\leq q}D_p$, and so
\begin{equation}\label{eq:N_down}
\nu_N^\sharp = \De_S \cup \bigcup_{p\leq q} \nu_{N^{\downarrow_p}}.
\end{equation}
Note that with $N\normal G_q$ and $N\leq Z_q$ as above, we have $N^{\downarrow_q}=N$; indeed, if $I=\bigcup_{p<q}D_p$ denotes the ideal of all elements strictly $\leqJ$-below  $D_q$, then $\nu_N\sub\nu_N^\sharp\restr_{D_q}\sub R_{I,N}^\sharp\restr_{D_q}=R_{I,N}\restr_{D_q}=\nu_N$, so that $\nu_N=\nu_N^\sharp\restr_{D_q}=\nu_{N^{\downarrow_q}}$.
(As in the proof of Lemma \ref{lem:nuN}, we have $I=\emptyset$ if $D_q$ is the minimal ideal of $S$; in this case, the previous calculation simplifies, as then $\nu_N$ is a congruence on $D_q$, so that $\nu_N^\sharp\restr_{D_q}=\nu_N\restr_{D_q}=\nu_N$.)
It is also clear that $\{\id_{G_q}\}^{\downarrow_p}=\{\id_{G_p}\}$ for any $p,q\in Q$.  

We call a $Q$-tuple $\bN=(N_q)_{q\in Q}$ an \emph{NZ-tuple} if
\bit
\item $N_q\normal G_q$ and $N_q\leq Z_q$ for all $q\in Q$, and
\item $N_q^{\downarrow_p}\leq N_p$ for all $p,q\in Q$ with $p\leq q$.
\eit
To such an NZ-tuple $\bN$, we associate the relation 
\[
\Th(\bN) = \Th^S(\bN) = \bigcup_{q\in Q}\nu_{N_q}.
\]

\begin{prop}\label{prop:ThN}
If $S$ is a stable, regular partial semigroup, then
\[
\CongH(S) = \set{\Th^S(\bN)}{\bN \text{ is an NZ-tuple}}.
\]
\end{prop}

\pf
We first show that $\Th(\bN)$ is a congruence for any NZ-tuple $\bN=(N_q)_{q\in Q}$.  Since $\Th(\bN)$ is a union of equivalences on disjoint sets, it suffices to show that for all $(x,y)\in\Th(\bN)$ and $a\in S$,
\[
\br(a)=\bd(x) \implies (ax,ay)\in\Th(\bN) \AND \bd(a)=\br(x) \implies (xa,ya)\in\Th(\bN).
\]
As usual, we just do the first, so suppose $\br(a)=\bd(x)$, and also suppose $(x,y)\in\nu_{N_q}$.  Now,
\[
N_q\leq Z_q \ \implies\  (x,y)\in\nu_{Z_q}\sub\ze \ \implies\  (ax,ay)\in\ze\sub\H\sub\D \ \implies\  ax,ay\in D_p \text{ for some $p\leq q$.}
\]
But then $(ax,ay)\in\cg xy\restr_{D_p}\sub\nu_{N_q}^\sharp\restr_{D_p}=\nu_{N_q^{\downarrow_p}}\sub\nu_{N_p}\sub\Th(\bN)$, as required.

Conversely, we have already seen that any $\H$-congruence $\si$ has the form \eqref{eq:Hcong}.  Suppose $p,q\in Q$ are such that $p\leq q$.  Then $\nu_{N_q^{\downarrow_p}}=\nu_{N_q}^\sharp\restr_{D_p}\sub\si^\sharp\restr_{D_p}=\si\restr_{D_p}=\nu_{N_p}$, which gives $N_q^{\downarrow_p}\leq N_p$.
\epf

So if $\th$ is an $\H$-congruence on the stable, regular partial semigroup $S$, then $\th=\Th(\bN)$ for some NZ-tuple $\bN$; we denote this tuple by~$\bN=\bN(\th)$.  If $\bN(\th)=(N_q)_{q\in Q}$, then we will write $\bN(\th)_q=N_q$ for any $q\in Q$.

The ordering in the lattice $\CongH(S)$ is given by coordinate-wise inclusion of (normal) subgroups:  for two NZ-tuples $\bN=(N_q)_{q\in Q}$ and $\bN'=(N_q')_{q\in Q}$, we have 
\[
\Th(\bN)\sub\Th(\bN')\quad  \Leftrightarrow \quad N_q\leq N_q' \text{ for all } q\in Q.
\]

\subsection{More on chains of ideals}\label{subsect:chains2}

We now wish to prove an analogue of Theorem \ref{thm:E2} for partial semigroups that may have non-trivial $\H$-congruences.  The proof of Theorem \ref{thm:E2} relied on iterating Theorem \ref{thm:E6}\ref{E6iii}, and this iteration required suitably large ideals to have property $\Sep{I_r}$.  In the absence of this property, we iterate the more general Theorem \ref{thm:E6}\ref{E6ii}, and look to replace $\Sep{I_r}$ by a weaker property.

An immediate consequence of item \ref{S3} in $\Sep{T}$ is that the principal congruence generated by a non-trivial pair of $\H$-related elements from $D_T$ is never an $\H$-congruence.  Thus, if there are non-trivial $\ze$-related pairs in $D_T$, then $\Sep{T}$ cannot hold.
The minimum weakening of $\Sep T$ to partial semigroups that may have non-trivial $\H$-congruences, which is simultaneously the maximum strengthening of $\Sepb T$ to such partial semigroups, is to exclude the pairs belonging to $\ze\restr_{D_T}$ from \ref{S3}.  Thus, we arrive at the following:

\begin{defn}
\label{defn:Sepz}
A partial semigroup $T$ satisfies the property $\Sepz{T}$ if it satisfies
$\Dmax{T}$, the ideal $S=T\setminus D_T$ satisfies $\Dmax{S}$, and
the following hold:
\begin{enumerate}[label=\textup{\textsf{(S\arabic*)}$^{\white \ze}$},leftmargin=11.7mm]
\item \label{Sepz1} for any $x\in D_T$ and $y\in S$ with $x\sim y$, there exists $(x',y')\in\cg xy_T$ such that $x'\in D_S$ and $y'\in S\sm H_{x'}^S$, 
\item \label{Sepz2} for any $x\in D_T$ and $y\in D_T\setminus H_x^T$ with $x\sim y$, there exists $(x',y')\in\cg xy_T$ such that $x'\in D_T$ and $y'\in S$,
\end{enumerate}
\begin{enumerate}[label=\textup{\textsf{(S\arabic*)}$^\ze$},leftmargin=11.7mm]\addtocounter{enumi}{2}
\item \label{Sepz3} for every $x,y\in D_T$ with $(x,y)\in\H^T\sm\ze_T$, there exists $(x',y')\in\cg xy_T$ such that $x'\in D_S$ and $y'\in S\sm H_{x'}^S$,
\end{enumerate}
\end{defn}

We also have the corresponding multiplicative property, which we denote by $\Multz T$, and which clearly implies $\Sepz T$; cf.~Definitions \ref{defn:Sep} and \ref{defn:Mult}.

\newpage

\begin{defn}
\label{defn:Multz}
A partial semigroup $T$ satisfies the property $\Multz{T}$ if it satisfies $\Dmax{T}$, the ideal $S=T\setminus D_T$ satisfies $\Dmax{S}$, and the following hold:
\begin{enumerate}[label=\textup{\textsf{(M\arabic*)}$^{\white \ze}$},leftmargin=11.7mm]
\item \label{Multz1} for any $x\in D_T$ and $y\in S$ with $x\sim y$, there exist $a,b\in T^1$ such that (swapping $x,y$ if necessary) $axb\in D_S$ and $ayb\in S\sm H_{axb}^S$, 
\item \label{Multz2} for any $x\in D_T$ and $y\in D_T\setminus H_x^T$ with $x\sim y$, there exists $a,b\in T^1$ such that (swapping $x,y$ if necessary) $axb\in D_T$ and $ayb\in S$,
\end{enumerate}
\begin{enumerate}[label=\textup{\textsf{(M\arabic*)}$^\ze$},leftmargin=11.7mm]\addtocounter{enumi}{2}
\item \label{Multz3} for every $x,y\in D_T$ with $(x,y)\in\H^T\sm\ze_T$, there exists $a,b\in T^1$ such that (swapping $x,y$ if necessary) $axb\in D_S $ and $ayb\in S\sm H_{axb}^S$.
\end{enumerate}
\end{defn}

The next result concerns a chain of ideals $U$, as in Definition~\ref{defn:chain}.  
Thus, $U$ is a stable, regular partial semigroup, and $U/\J = {\set{D_q}{0\leq q<\xi}}$ with $D_0<D_1<\cdots$, where $\xi$ is some ordinal from the set $\{1,2,3,\ldots,\om\}$.
As in the discussion following Definition \ref{defn:chain}, we have the ideals $I_r=D_0\cup D_1\cup\cdots\cup D_r$ for $0\leq r<\xi$, and also $I_\om=U$ in the case $\xi=\om$.  We also fix group $\H$-classes $G_q$ in $D_q$ ($0\leq q<\xi$), and the normal subgroups $Z_q\normal G_q$ for which $\ze_U=\nu_{Z_0}\cup\nu_{Z_1}\cup\cdots$.  Lemma \ref{lem:zeST}\ref{zeST2} gives $\ze_{I_r}=\nu_{Z_0}\cup\nu_{Z_1}\cup\cdots\cup\nu_{Z_r}$ for $0\leq r<\xi$; thus, in what follows, we will not use subscripts to distinguish the $\ze$ relations on $U$ and its ideals.  It also follows that NZ-tuples in ideals of $U$ are just truncations of NZ-tuples in $U$.

\begin{thm}\label{thm:E4}
Let $U$ be a chain of ideals of height $\xi$, and keep the notation of Definition \ref{defn:chain} and the above discussion.  Suppose there exists $0\leq k<\xi$ such that the following hold:
\[
\LiftCong{I_k}U=\Cong(I_k)\COMMA \ngen{I_k}\COMMA \Sepz{I_r}\ \text{for all}\ k<r<\xi.
\]
Then for any $r\geq k$ (including $r=\om$ if $\xi=\om$),
\ben
\item \label{E41} $\CongH(I_r) = \bigset{\Th^{I_r}(\bN)}{\bN\text{ is an NZ-tuple in }I_r}$,
\item \label{E42} $\Cong(I_r) = \bigset{\si\cup\th}{\si\in\Cong(I_k),\ \th\in\CongH(I_r),\ \th\restr_{I_k}\sub\si}$
\item[] \hspace{0.751cm} ${} \cup \bigset{R_{I_q,N}^{I_r}\cup\th}{k\leq q<r,\ N\normal G_{q+1},\ \th\in\CongH(I_r),\ \bN(\th)_{q+1}\leq N} \cup \{\nab_{I_r}\}$.
\een
\end{thm}

\pf
Part \ref{E41} is 
Proposition \ref{prop:ThN},  so we just need to prove \ref{E42}.  
We begin by showing that the stated relations are congruences.  
First suppose $\tau=\si\cup\th$, where $\si\in\Cong(I_k)$ and ${\th\in\CongH(I_r)}$ are such that $\th\restr_{I_k}\sub\si$, and write $\bN=\bN(\th)=(N_0,N_1,\ldots)$.  
Since $\theta\subseteq \H\subseteq\D$ and $\theta\restr_{I_k}\subseteq \sigma$,
we have $\tau=\sigma\cup\theta\restr_{I_r\setminus I_k}$, the union of two equivalence relations on disjoint sets;
hence,~$\tau$ is an equivalence relation.
Then from $\tau=(\sigma\cup\Delta_{I_r})\cup\theta$ it follows that
$\tau=(\sigma\cup\Delta_{I_r})\vee\theta$, the join of two congruences, and so $\tau$ is a congruence.
(Note that $\si$ is liftable to $I_r$ as it is liftable to $U$.)
The proof that relations of the form $\tau=R_{I_q,N}^{I_r}\cup\th$ are congruences is analogous.

The rest of the proof is devoted to showing
that every congruence on $I_r$ has one of the forms listed. 
The proof follows that of Theorem \ref{thm:E2}, but with the view to using part \ref{E6ii} of Theorem~\ref{thm:E6}, and not part
\ref{E6iii}.

\pfstep{1}
We first show by induction that this is true for $k\leq r<\xi$.  The $r=k$ case being clear, suppose $k<r<\xi$, and inductively assume that
\begin{align}
\nonumber \Cong(I_{r-1}) = \bigset{\si\cup\th}{{}&\si\in\Cong(I_k),\ \th\in\CongH(I_{r-1}),\ \th\restr_{I_k}\sub\si} \\[2truemm]
\label{eq:E4_1}  {} \cup \bigset{R_{I_q,N}^{I_{r-1}}\cup\th}{{}&k\leq q<{r-1},\ N\normal G_{q+1},\ \th\in\CongH(I_{r-1}),\ \bN(\th)_{q+1}\leq N} \cup \{\nab_{I_{r-1}}\}.
\end{align}
It follows from this (or by assumption, if $r=k+1$) that $\ngen{I_{r-1}}$ holds.  Since $\Sepz{I_r}$ holds by assumption, certainly $\Sepb{I_r}$ holds.  Thus, $\ngen{I_r}$ holds, by Lemma \ref{lem:sepb_ngen}, and so Theorem~\ref{thm:E6}\ref{E6ii} applies.  To facilitate this application note that
\begin{equation}\label{eq:Ir-1}
\LiftCong{I_{r-1}}{I_r}=\Cong(I_{r-1}).
\end{equation}
Indeed, if~$\ka$ is any congruence listed in \eqref{eq:E4_1}, then $\ka\cup\De_{I_r}$ is one of the relations listed in the theorem: 
for $\ka\not=\nab_{I_{r-1}}$, we simply append $N_r=\{\id_{G_r}\}$ to $\bN(\th)$.  Thus, by Theorem~\ref{thm:E6}\ref{E6ii}, we have
\begin{equation}\label{eq:E4_2} 
\Cong(I_r) = \set{\ka\cup\nu_{N_r}}{\ka\in\Cong(I_{r-1}),\ N_r\normal G_r,\ \tau_{N_r}\sub\ka} \cup \{\nab_{I_r}\},
\end{equation}
where $\tau_{N_r}=(\nu_{N_r})_{I_r}^\sharp\restr_{I_{r-1}}$.
It is now a matter of showing that each congruence in \eqref{eq:E4_2} is one of those listed in the theorem.  To do so, we first note that for any $N_r\normal G_r$,
\begin{numcases}{\tau_{N_r} = }
\label{E4N1} \Th^{I_{r-1}}(N_r^{\downarrow_0},N_r^{\downarrow_1},\ldots,N_r^{\downarrow_{r-1}}) &\text{if $N_r\leq Z_r$}\\
\label{E4N2} \nab_{I_{r-1}} &\text{if $N_r\nleq Z_r$.}
\end{numcases}
Indeed, \eqref{E4N1} is true by definition (cf.~\eqref{eq:N_down}).  For \eqref{E4N2}, we apply \ref{Sepz3} to a pair ${(e,x)\in\nu_{N_r}}$ with $x\in N_r\sm Z_r$, and then apply~$\ngen{I_{r-1}}$.

Returning now to the main proof, we need to show that any congruence $\ka\cup\nu_{N_r}$ on $I_r$, as in~\eqref{eq:E4_2}, is among those listed in the theorem.

\pfcase{1}
 If ${N_r}\nleq Z_r$, then by \eqref{E4N2} we have $\ka\cup\nu_{N_r} = \nab_{I_{r-1}}\cup\nu_{N_r} = R_{I_{r-1},N_r}^{I_r} = R_{I_{r-1},N_r}^{I_r}\cup\De_{I_r}$.

\pfcase{2}
If ${N_r}\leq Z_r$ and $\ka=\nab_{I_{r-1}}$, then again we have $\ka\cup\nu_{N_r}=R_{I_{r-1},N_r}^{I_r}$.

\pfcase{3}
Next, suppose $N_r\leq Z_r$, and that $\ka$ is of the form $\si\cup\th$, as in \eqref{eq:E4_1}.  Keeping Proposition~\ref{prop:ThN} and Lemma~\ref{lem:zeST} in mind, we may write
\begin{align*}
\th &= \nu_{N_0} \cup \nu_{N_1} \cup \cdots \cup \nu_{N_k}\cup\cdots\cup\nu_{N_{r-1}} &&\text{where $(N_0,\ldots,N_{r-1})$ is an NZ-tuple,}\\
\si\cap\ze &= \nu_{N_0'} \cup \nu_{N_1'} \cup \cdots \cup \nu_{N_k'} &&\text{where $(N_0',\ldots,N_k')$ is an NZ-tuple.}
\end{align*}
Since $\th\restr_{I_k}\sub\si$ and also $\th\restr_{I_k}\sub\ze$, it follows that $\th\restr_{I_k}\sub\si\cap\ze$, and so $N_t\leq N_t'$ for all $0\leq t\leq k$.  In particular, $(N_0',N_1',\ldots,N_k',N_{k+1},\ldots,N_{r-1})$ is an NZ-tuple.  
By \eqref{E4N1}, and since $\tau_{N_r} \sub \ka$ and $\tau_{N_r}\sub\ze$ (the latter because $N_r\leq Z_r$), we have
\[
\nu_{N_r^{\downarrow_0}} \cup \cdots \cup \nu_{N_r^{\downarrow_{r-1}}} = \tau_{N_r} \sub \ka\cap\ze = (\si\cup\th)\cap\ze = (\si\cap\ze)\cup\th = \nu_{N_0'}\cup\cdots\cup\nu_{N_k'}\cup\nu_{N_{k+1}}\cup\cdots\cup\nu_{N_{r-1}}.
\]
It follows that $\bN=(N_0',\ldots,N_k',N_{k+1},\ldots,N_{r-1},N_r)$ is an NZ-tuple.  But then
\[
\ka\cup\nu_{N_r} = \si\cup\th\cup\nu_{N_r} = \si \cup \Th^{I_r}(\bN),
\]
with $\Th^{I_r}(\bN)\restr_{I_k}=\si\cap\ze\sub\si$.

\pfcase{4}
 Finally, suppose $N_r\leq Z_r$, and that $\ka$ is of the form $R_{I_q,N}^{I_{r-1}}\cup\th$, as in \eqref{eq:E4_1}.  Again, write $\th= \nu_{N_0} \cup \nu_{N_1} \cup \cdots \cup \nu_{N_{r-1}}$, where $(N_0,\ldots,N_{r-1})$ is an NZ-tuple, noting that $N_{q+1}\leq N$ (cf.~\eqref{eq:E4_1}).  This time, 
\[
\nu_{N_r^{\downarrow_0}} \cup \nu_{N_r^{\downarrow_1}} \cup \cdots \cup \nu_{N_r^{\downarrow_{r-1}}} = \tau_{N_r}
\sub\ka
=R_{I_q,N}^{I_{r-1}}\cup\th
=\nab_{I_q}\cup\nu_N\cup\nu_{N_{q+2}}\cup\cdots\cup\nu_{N_{r-1}}
\]
gives $N_r^{\downarrow_{q+1}}\leq N$ and $N_r^{\downarrow_t}\leq N_t$ for all $q+2\leq t\leq r-1$.  Put $N_{q+1}'=N_{q+1}N_r^{\downarrow_{q+1}}$, noting that $N_{q+1}'\normal G_{q+1}$, $N_{q+1}'\leq Z_{q+1}$ and $N_{q+1}'\leq N$.  Then
\[
\ka\cup\nu_{N_r} = \nab_{I_q}\cup\nu_N\cup\nu_{N_{q+2}}\cup\cdots\cup\nu_{N_{r-1}}\cup\nu_{N_r} = R_{I_q,N}^{I_r}\cup\Th^{I_r}(\bN),
\]
for the NZ-tuple $\bN = (Z_0,Z_1,\ldots,Z_q,N_{q+1}',N_{q+2},\ldots,N_{r-1},N_r)$, and we have already noted that $N_{q+1}'\leq N$.

\aftercases
This completes the inductive step, and hence Step 1 as well.

\pfstep{2}
We now consider the case in which $r=\xi=\om$.  Suppose $\ka$ is an arbitrary congruence on $U=I_\om$.  
The proof will be complete if we can show that $\ka$ is one of those listed in the theorem.
For each $k\leq q<\om$, let $\ka_q = \ka\restr_{I_q}$.  If $\ka_q=\nab_{I_q}$ for all such $q$, then $\ka=\nab_U$.  So suppose otherwise, and let
\begin{equation}\label{eq:E4_3}
p=\min\set{q}{k< q<\om,\ \ka_q\not=\nab_{I_q}},
\end{equation}
noting that $p>k$.
For any $p\leq t<\om$, the finite case of the theorem (proved above) tells us that $\ka_t$ $(\not=\nab_{I_t})$ is contained in $R_{I_{t-1},G_t}^{I_t} = \nab_{I_{t-1}} \cup \H\restr_{D_t}$, and so $\ka\restr_{D_t} = \ka_t\restr_{D_t}\sub\H$; Lemma~\ref{lem:nu_trick} then gives $\ka\restr_{D_t}=\nu_{N_t}$ for some $N_t\normal G_t$.  Thus,
\begin{equation}\label{eq:ka_p-1}
\ka = \ka_{p-1} \cup \nu_{N_p} \cup \nu_{N_{p+1}} \cup \cdots.
\end{equation}
If $N_t\not\leq Z_t$ for some $t>p$, then \eqref{E4N2} would give $\nab_{I_{t-1}} = (\nu_{N_t})_{I_t}^\sharp\restr_{I_{t-1}} \sub (\nu_{N_t})_U^\sharp \sub \ka$, and hence $\nab_{I_p}\sub\ka$, contradicting $\ka_p\not=\nab_{I_p}$.  Thus, $N_t\leq Z_t$ for all $t>p$.  We also note that for any $t>s\geq p$,
\[
\nu_{N_t^{\downarrow_s}} = (\nu_{N_t})_U^\sharp\restr_{D_s} \sub \ka\restr_{D_s} = \nu_{N_s},
\]
which gives $N_t^{\downarrow_s}\leq N_s$ for all such $t,s$.

\pfcase{1}
If $\ka_{p-1}=\nab_{I_{p-1}}$, then \eqref{eq:ka_p-1} gives $\ka = \nab_{I_{p-1}}\cup\nu_{N_p}\cup\nu_{N_{p+1}}\cup\cdots = R_{I_{p-1},N_p}^U\cup\Th(\bN)$, for the NZ-tuple $\bN=(Z_0,Z_1,\ldots,Z_{p-1},N_p\cap Z_p,N_{p+1},N_{p+2},\ldots)$.  (Note that for $t>p$, we have $\nu_{N_t^{\downarrow_p}}=(\nu_{N_t})_U^\sharp\restr_{D_p}\sub\ze\restr_{D_p}=\nu_{Z_p}$, which gives $N_t^{\downarrow_p}\leq N_p\cap Z_p$.)

\pfcase{2}
Now suppose $\ka_{p-1}\not=\nab_{I_{p-1}}$.  As above, we must have $N_p\leq Z_p$ (or else $\nab_{I_{p-1}}\sub\ka$, a contradiction).  By the definition of $p$ in \eqref{eq:E4_3}, we must in fact have $p=k+1$.  Thus, writing $\si=\ka_{p-1}=\ka_k\in\Cong(I_k)$, \eqref{eq:ka_p-1} gives 
\[
\ka = \si \cup \nu_{N_{k+1}} \cup \nu_{N_{k+2}} \cup \cdots,
\]
and we have seen that $N_t\normal G_t$ and $N_t\leq Z_t$ for all $t>k$.  
Now, $\si\cap\ze=\si\cap\ze_{I_k}$ is an $\H$-congruence on $I_k$, so by Proposition \ref{prop:ThN} (cf.~Lemma \ref{lem:zeST}), $\si\cap\ze=\nu_{N_0}\cup\nu_{N_1}\cup\cdots\cup\nu_{N_k}$ for some NZ-tuple $(N_0,N_1,\ldots,N_k)$.  Now, for any $0\leq s\leq k<t$, we have
\[
\nu_{N_t^{\downarrow_s}} = (\nu_{N_t})_U^\sharp\restr_{D_s} \sub \ka\restr_{D_s} = \si\restr_{D_s} \AND \nu_{N_t^{\downarrow_s}} = (\nu_{N_t})_U^\sharp\restr_{D_s} \sub (\nu_{Z_t})_U^\sharp \sub \ze,
\]
from which it follows that $\nu_{N_t^{\downarrow_s}} \sub (\si\cap\ze)\restr_{D_s} = \nu_{N_s}$, and so $N_t^{\downarrow_s}\leq N_s$.  All of this shows that $\bN=(N_0,N_1,\ldots,N_k,N_{k+1},\ldots)$ is an NZ-tuple in $U$.  But then $\ka = \si \cup \Th(\bN)$, and also $\Th(\bN)\restr_{I_k} = \nu_{N_0}\cup\nu_{N_1}\cup\cdots\cup\nu_{N_k} = \si\cap\ze \sub \si$.
\epf

\begin{rem}
Theorem \ref{thm:E2} may be viewed as a special case of Theorem \ref{thm:E4}\ref{E42}.  
This follows from the observation that
 $\Sep{I_r}$ implies $Z_r=\{\id_{G_r}\}$ for all $r>k$.  
Thus, every NZ-tuple $(N_0,N_1,\ldots)$ satisfies $N_r=\{\id_{G_r}\}$ for $r>k$.  
Hence the congruences listed in Theorem \ref{thm:E4}\ref{E42} take the simpler forms $\si\cup\th=\si\cup\De_{I_r}$ and $R_{I_q,N}^{I_r}\cup\th=R_{I_q,N}^{I_r}$, exactly as in Theorem \ref{thm:E2}.
\end{rem}

A particular simplification of Theorem \ref{thm:E4} will be of use in our forthcoming applications.  
Here we are dealing with the case in which $k=1$, and the minimal ideal $I_0=D_0$ of $U$ (as above) is \emph{trivial} in the sense that each hom-set $U_{A,B}$ contains a unique element of $D_0$.  In this case, every ideal of~$U$ is trivially retractable.

For the next statement, if $I$ is an ideal of a partial semigroup $S$, and if $\bN=(N_q)_{q\in Q}$ where the $N_q$ are normal subgroups of group $\H$-classes contained in distinct stable, regular $\J$-classes contained in $S\sm I$, then we define the relation
\[
R_{I,\bN}^S = R_I^S \cup \bigcup_{q\in Q}\nu_{N_q}.
\]
Note that here we do not assume the $\J$-classes containing the $N_q$ are minimal in $(S\sm I)/\J$, or that these are \emph{all} the $\J$-classes from $S\sm I$.

\begin{thm}\label{thm:E5}
Let $U$ be a chain of ideals of height $\xi$, and keep the notation of Definition \ref{defn:chain} and the discussion before Theorem \ref{thm:E4}.  Suppose also that
\bit
\item the minimal ideal $I_0=D_0$ is trivial,
\item for every $x\in D_1$ and $y\in D_1\sm H_x$ with $x\sim y$, there exists $(x',y')\in\cg xy_{I_1}$ with $x'\in D_1$ and $y'\in D_0$,
\item $\Sepz{I_r}$ holds for each $2\leq r<\xi$.
\eit
Then for any $r\geq 0$ (including $r=\om$ if $\xi=\om$), the congruences on the ideal $I_r$ are precisely $\nab_{I_r}$ and those of the form
\bit
\item $R_{I_q,\bN}^{I_r}$ for some $0\leq q<r$ and some $\bN=(N_{q+1},N_{q+2},\ldots)$,
\eit
where
\bit
\item $N_p\normal G_p$ for all $p\geq q+1$,
\item $N_p\leq Z_p$ for all $p\geq q+2$,
\item $N_t^{\downarrow_p}\leq N_p$ for all $t>p\geq q+1$.
\eit
\end{thm}

\pf
This is clear for $r=0$; indeed, the trivial $I_0$ has only one congruence: $\De_{I_0}=\nab_{I_0}$.  
Proposition \ref{prop:small_ideal} gives the $r=1$ case, noting that $I_1$ is retractable.  In the notation of that proposition, and keeping in mind that $\Cong(I_0)=\{\De_{I_0}=\nab_{I_0}\}$, the ``$\tau\cup\nu_N$-type'' congruences are all of the form $\nab_{I_0}\cup\nu_N=R_{I_0,N}^{I_1}$ ($N\normal G_1$), and the only ``$\th_{S,\tau}$-type'' congruence is ${\th_{I_1,\nab_{I_0}}=\nab_{I_1}}$.  In particular,~$\ngen{I_1}$ holds, and we also have $\LiftCong{I_1}U=\Cong(I_1)$; cf.~Lemma~\ref{lem:E2}.  Thus, Theorem~\ref{thm:E4} applies with $k=1$, and it is easy to check that it gives precisely the desired descriptions of $\Cong(I_r)$ for $r\geq2$.
\epf

\section{Linear and projective linear categories}\label{sect:L}

We now turn our attention to categories of linear transformations, specifically the linear and projective linear categories, $\LL$ and $\PL$, both defined in Section~\ref{subsect:defnL}.  
After discussing Green's relations and establishing the required multiplicative properties in Section \ref{subsect:prelimL}, 
we describe the congruences on arbitrary ideals of $\LL$ in Section~\ref{subsect:L}, and on $\PL$ in Section~\ref{subsect:PL}. 
It turns out that $\LL$ has non-trivial $\H$-congruences, and so we will be applying the results from Section~\ref{sect:H};
for $\PL$, however, we will return to those of Section \ref{sect:chains}.
In Section \ref{subsect:visual_L} we discuss how to visualise the congruence lattices of ideals in~$\LL$.
The congruences of the ideals of the linear \emph{monoid} were determined by Mal'cev \cite{Malcev1953}, one of the early papers motivating the present work.

\subsection{Definitions and preliminaries on $\LL$ and $\PL$}\label{subsect:defnL}

Throughout Section \ref{sect:L}, $\F$ will be a fixed field, and $\Fx=\F\sm\{0\}$ will denote the multiplicative group of units of $\F$.  We denote by $\GL_q=\GL_q(\F)$ and $\SL_q=\SL_q(\F)$, respectively, the general and special linear groups of $q\times q$ matrices over $\F$.  We also write $\PGL_q=\PGL_q(\F)=\GL_q/Z(\GL_q)$ for the
projective linear group; here $Z(\GL_q)=\set{cE_q}{c\in\Fx}$ is the centre of $\GL_q$, consisting of all non-zero scalar multiples of the  
$q\times q$ identity matrix $E_q$.

Fix some non-empty set $\C$ of finite-dimensional vector spaces over $\F$, and define the parameter
\[
\xi=\xi(\C)={\LSUB}\bigset{{\dim}(V)}{V\in\C},
\]
where as usual $\dim(V)$ is the dimension of $V$.  
For vector spaces $U$ and $V$ we write $\LL_{U,V}$ for the set of all linear transformations from $U\to V$.  We now define
\[
\LL = \LL(\C) 
= \bigset{(U,\al,V)}{U,V\in\C,\ \al\in\LL_{U,V}}.
\]
For $U,V,W\in\C$, and for $\al\in\LL_{U,V}$ and $\be\in\LL_{V,W}$, we define
\[
\bd(U,\al,V)=U \COMMA \br(U,\al,V)=V \COMMA (U,\al,V)\cdot(V,\be,W) = (U,\al\be,W),
\]
where $\al\be\in\LL_{U,W}$ denotes the composition (performing $\al$ first, then $\be$).  Then
\[
 (\LL,\C,\bd,\br,\cdot)
\]
is a partial semigroup, indeed a stable, regular category (cf.~\cite[Lemmas 3.1 and 3.2]{DE2018}), which we will call the \emph{linear category (over $\C$)}.  For $V\in\C$, the endomorphism monoid $\LL_V=\LL_{V,V}$ is the \emph{full linear monoid} over $V$ \cite{Okninski1998}.  Note that when $\C=\{V\}$ consists of a single finite-dimensional vector space, the category $\LL$ is simply the monoid $\LL_V$.

As with the transformation category $\T$, we will avoid cumbersome notation by identifying $(U,\al,V)\in\LL$ with the transformation $\al\in\LL_{U,V}$, regarding the vector spaces $U,V$ as ``encoded'' in $\al$, and writing $\bd(\al)=U$ and $\br(\al)=V$.

Of importance in all that follows will be the relation $\ze$ on $\LL$ defined by
\[
\ze = \bigset{(\al,\be)\in\LL\times\LL}{\al=c\be \text{ for some }c\in\Fx}.
\]
It is a routine matter to verify that $\ze$ is a congruence, and that~$\ze\sub\H$.  
We have taken the liberty of using the suggestive notation $\ze$ for this congruence because it is in fact the maximum $\H$-congruence on $\LL$; cf.~Corollary \ref{cor:zeL}.

For $\al\in\LL$, we write $\ol\al$ for the $\ze$-class of~$\al$.  
We call the quotient category
\[
\PL = \PL(\C) = \LL/\ze = \set{\ol\al}{\al\in\LL}
\]
the \emph{projective linear category (over $\C$)}.  Endomorphism monoids in $\PL$ are projective linear monoids $\PL_U$ ($U\in\C$).
It is worth observing that $\ze=\De_\LL$ if $|\F|=2$, so that $\PL=\LL$ in this case.

\subsection{Green's relations and multiplicative properties in $\LL$ and $\PL$}\label{subsect:prelimL}

For a linear transformation $\al:U\to V$, we write as usual
\[
\Ker(\al) = \set{\vu\in U}{\vu\al=\vzero} \COMMA \im(\al) = \set{\vu\al}{\vu\in U} \COMMA \rank(\al) = \dim(\im(\al))
\]
for the \emph{kernel}, \emph{image} and \emph{rank} of $\al$.  Throughout $\vzero$ will  denote the zero of any vector space.  Note that $\Ker(\al)$ is a \emph{subspace} of $U$ (not an equivalence), $\im(\al)$ is a subspace of $V$, and $\rank(\al)$ is a cardinal between~$0$ and $\min(\dim(U),\dim(V))$, inclusive.

For a proof of the next result, see \cite[Lemma 3.1]{DE2018}, which is formulated equivalently in terms of matrices and row/column spaces; see also \cite[Lemma 2.1]{Okninski1998}.  In what follows, it will occasionally be convenient to move between the viewpoints of linear transformations and matrices.

\begin{lemma}\label{lem:Green_L}
\ben
\item \label{GL1} For $\al,\be\in\LL$, we have
\bit
\item $(\al,\be)\in\R \iff \bd(\al)=\bd(\be)$ and $\Ker(\al)=\Ker(\be)$,
\item $(\al,\be)\in\L \iff \br(\al)=\br(\be)$ and $\im(\al)=\im(\be)$,
\item $(\al,\be)\in\J\iff(\al,\be)\in\D\iff\rank(\al)=\rank(\be)$.
\eit
\item \label{GL2} The $\J=\D$-classes of $\LL$ are the sets
\bit
\item $D_q = D_q(\LL) = \set{\al\in\LL}{\rank(\al)=q}$ for each $0\leq q<\xi$.
\eit
These are all regular and stable, and they form a chain: $D_0<D_1<\cdots$.
\item \label{GL3} The ideals of $\LL$ are the sets
\bit
\item $I_r = I_r(\LL) = \set{\al\in\LL}{\rank(\al)\leq r}$ for each $0\leq r<\xi$, and
\item $I_\om = I_\om(\LL) = \LL $ in the case that $\xi=\om$.
\eit
These are all regular and stable, and they form a chain: $I_0\subset I_1\subset\cdots$.
\item \label{GL4} The $\H$-class of any idempotent from $D_q$ is isomorphic to $\GL_q(\F)$.
 \epfres
\een
\end{lemma}

Thus, $\LL$ is a chain of ideals, as in Definition \ref{defn:chain}.  Note that the minimal ideal $I_0=D_0$ consists entirely of zero-mappings, and is hence \emph{trivial} in the sense described before Theorem~\ref{thm:E5}.

Since the congruence $\ze$ is contained in $\H$, it follows quickly that Green's relations on the projective linear category $\PL=\LL/\ze$ are induced by those on $\LL$ in the following sense.  If $\al,\be\in\LL$, and if~$\K$ is any of Green's relations, then $\ol\al$ and $\ol\be$ are $\K$-related in $\PL$ if and only if $\al$ and $\be$ are $\K$-related in $\LL$.  Thus, the (stable, regular) $\J=\D$-classes and ideals of $\PL$ are the sets 
\[
D_q(\PL) = \set{\ol\al}{\al\in D_q(\LL)} \AND I_r(\PL) = \set{\ol\al}{\al\in I_r(\LL)} \qquad\text{for appropriate $q,r$.}
\]
These again form chains:  $D_0(\PL)<D_1(\PL)<\cdots$ and $I_0(\PL)\subset I_1(\PL)\subset\cdots$.  The group $\H$-classes in $D_q(\PL)$ are isomorphic to $\PGL_q(\F)$.

As usual, the following technical lemma will underpin the investigations to follow.

\begin{lemma}\label{lem:P2L}
Let $U,V\in\C$, and suppose $\al,\be\in\LL_{U,V}$ with $r=\rank(\al)\geq\rank(\be)=q$.  
\ben
\item \label{P2Li} If $r\geq2$ and $q<r$, then there exists $\ga\in I_r$ such that $\al\ga\in D_{r-1}$ and $\be\ga\in I_{r-1}\sm H_{\al\ga}$.
\item \label{P2Lii} If $q=r\geq1$ and $(\al,\be)\not\in\H$, then there exists $\ga\in I_r$ such that 
\[
\text{$[\al\ga\in D_r$ and $\be\ga\in I_{r-1}]$ \ \ \ or \ \ \ $[\ga\al\in D_r$ and $\ga\be\in I_{r-1}]$.}
\]
\item \label{P2Liii} If $r\geq2$ and $(\al,\be)\in\H\sm\ze$, then there exists $\ga\in I_r$ such that $\al\ga\in D_{r-1}$ and $\be\ga\in I_{r-1}\sm H_{\al\ga}$.
\een
\end{lemma}

\pf
\firstpfitem{\ref{P2Li}}  
Since $\dim(\Ker(\be))=\dim(U)-\rank(\be)>\dim(U)-\rank(\al)=\dim(\Ker(\al))$, we may fix some $\vu_1\in U$ such that $\vv_1=\vu_1\al\not=\vzero=\vu_1\be$.  Let $\{\vv_1,\ldots,\vv_d\}$ be a basis of $V$ such that $\{\vv_1,\ldots,\vv_r\}$ is a basis of $\im(\al)$.  Define $\ga\in I_r(\LL_V)$ by
\[
\vv_i\ga = \begin{cases}
\vv_i &\text{if $1\leq i\leq r-1$}\\
\vzero &\text{otherwise.}
\end{cases}
\]
Since $\{\vv_1,\ldots,\vv_{r-1}\}$ is a basis of $\im(\al\ga)$, we have $\al\ga\in D_{r-1}$.  Since $\vu_1\al\ga=\vv_1\not=\vzero=\vu_1\be\ga$, we have $\Ker(\al\ga)\not=\Ker(\be\ga)$, and so $(\al\ga,\be\ga)\not\in\R$, whence $\be\ga\in I_{r-1}\sm H_{\al\ga}$.

\pfitem{\ref{P2Lii}}  There are two cases to consider.

\pfcase{1}
Suppose first that $(\al,\be)\not\in\R$, so that $\vu_1\be=\vzero\not=\vu_1\al$ for some $\vu_1\in U$.  Put $\vv_1=\vu_1\al$, and let $\{\vv_1,\ldots,\vv_r\}$ be a basis of $\im(\al)$.  Choose some $\vu_i\in\vv_i\al^{-1}$ for each $2\leq i\leq r$.  We extend these vectors to a basis $\{\vu_1,\ldots,\vu_d\}$ of $U$.  Define $\ga\in I_r(\LL_U)$ by
\[
\vu_i\ga = \begin{cases}
\vu_i &\text{if $1\leq i\leq r$}\\
\vzero &\text{otherwise.}
\end{cases}
\]
Then $\ga\al\in D_r$ and $\ga\be\in I_{r-1}$.

\pfcase{2}  The case in which $(\al,\be)\not\in\L$ is dual.  Fix bases of $U$ and $V$, and let the matrices associated to the transformations $\al$ and $\be$ (and the fixed bases) be $A$ and $B$, respectively.  Since~$A$ and $B$ are $\L$-unrelated, their transposes $A^\tr$ and $B^\tr$ are $\R$-unrelated.  Thus, by the previous case there is a matrix $C$ of rank $r$ such that $\rank(CA^\tr)=r>\rank(CB^\tr)$.  We then take $\ga\in I_r(\LL)$ to be the transformation associated to $C^\tr$.

\pfitem{\ref{P2Liii}}  
Let $\{\vu_1,\ldots,\vu_d\}$ be a basis of $U$ such that $\{\vu_{r+1},\ldots,\vu_d\}$ is a basis of $\Ker(\al)=\Ker(\be)$.  For $1\leq i\leq r$, let $\vv_i=\vu_i\al$, so that $\{\vv_1,\ldots,\vv_r\}$ is a basis of $\im(\al)=\im(\be)$.  Extend this to a basis $\{\vv_1,\ldots,\vv_e\}$ of $V$.

\pfcase{1}  Suppose first that there exists some $1\leq i\leq r$ such that $\vu_i\be$ is not a scalar multiple of $\vv_i$.  Without loss of generality, we may assume that $i=r$, so that $\vu_r\be=\sum_{j=1}^rc_j\vv_j$ for some $c_j\in\F$ and such that $c_j\not=0$ for some $1\leq j\leq r-1$.  Define $\ga\in I_r(\LL_V)$ by
\[
\vv_i\ga = \begin{cases}
\vv_i &\text{if $1\leq i\leq r-1$}\\
\vzero &\text{otherwise.}
\end{cases}
\]
Then $\al\ga\in D_{r-1}$ and $\be\ga\in I_{r-1}$.  But also $\vu_r(\al\ga)=\vzero$, yet $\vu_r(\be\ga)=\sum_{j=1}^{r-1}c_j\vv_j\not=\vzero$, so that $(\al\ga,\be\ga)\not\in\R$.

\pfcase{2}  Now suppose that for each $1\leq i\leq r$ we have $\vu_i\be=c_i\vv_i$ for some $c_i\in\F$.
If $c_1=\cdots=c_r$, then (since also $\vu_i\al=\vu_i\be=\vzero$ for $i>r$) we would have $\be=c_1\al$; since $c_1\not=0$ (as $\rank(\be)=r>0$), it would follow that $(\al,\be)\in\ze$, a contradiction.  Thus, we may assume without loss of generality that $c_{r-1}\not=c_r$.  Define $\ga\in I_r(\LL_V)$ by
\[
\vv_i\ga = \begin{cases}
\vv_i &\text{if $1\leq i\leq r-1$}\\
-\vv_{r-1} &\text{otherwise.}
\end{cases}
\]
Then $\al\ga,\be\ga\in D_{r-1}$, as $\{\vv_1,\ldots,\vv_{r-1}\}$ is clearly a basis of $\im(\al\ga)=\im(\be\ga)$.  But since ${(\vu_{r-1}+\vu_r)\al\ga=\vzero}$ and $(\vu_{r-1}+\vu_r)\be\ga=(c_{r-1}-c_r)\vv_{r-1}\not=\vzero$,
we have ${(\al\ga,\be\ga)\not\in\R}$.
\epf

\begin{rem}\label{rem:P2L}
It follows from Lemma \ref{lem:P2L} that the ideal $I_r$ of $\LL$ ($2\leq r<\xi$) satisfies $\Multz{I_r}$, and hence $\Sepz{I_r}$; cf.~Definitions \ref{defn:Sepz} and \ref{defn:Multz}.

On the other hand, since $\ze\sub\H^{\LL}$, it follows that the ideal $I_r$ of $\PL=\LL/\ze$ ($2\leq r<\xi$) satisfies the stronger $\Mult{I_r}$, and hence $\Sep{I_r}$; cf.~Definitions~\ref{defn:Sep} and~\ref{defn:Mult}.  
\end{rem}

\begin{cor}\label{cor:zeL}
The congruence $\ze=\bigset{(\al,\be)\in\LL\times\LL}{\al=c\be \text{ for some }c\in\Fx}$ is the maximum $\H$-congruence on $\LL$.
\end{cor}

\pf
We first claim that $\ze\restr_{D_q}=\H\restr_{D_q}$ for $q=0,1$.  Indeed, this is clear for $q=0$.  Next suppose $(\al,\be)\in\H\restr_{D_1}$, say with $\al,\be\in\LL_{U,V}$.  Fix a basis $\{\vu_1,\ldots,\vu_d\}$ of $U$ so that $\{\vu_1,\ldots,\vu_{d-1}\}$ is a basis of $\Ker(\al)=\Ker(\be)$, and let $\vv=\vu_d\al$.  Since $\im(\al)=\im(\be)$, we have $\vu_d\be=c\vv$ for some $c\in\Fx$, and so $\be=c\al$, which means that $(\al,\be)\in\ze$.  This completes the proof of the claim.

Now, we have already observed that $\ze$ is an $\H$-congruence.  If $\ze$ was strictly contained in some other $\H$-congruence $\si$, then for any $(\al,\be)\in\si\sm\ze$, we would have $(\al,\be)\in\H\sm\ze$, and so $\rank(\al)=\rank(\be)\geq2$ (by the claim).  Lemma \ref{lem:P2L}\ref{P2Liii} then tells us that $\si\sm\H$ is non-empty, a contradiction.
\epf

\subsection{Congruences on ideals of $\LL$}\label{subsect:L}

In this section we use Theorem \ref{thm:E5} to describe the congruences on arbitrary ideals of the linear category 
$\LL=\LL(\C)$.  We keep all the notation as in the previous section.

We have already noted that $I_0=D_0$ is trivial, and the remaining technical assumptions from Theorem \ref{thm:E5} follow from Lemma \ref{lem:P2L}.
Therefore all congruences on $\LL$ are of the form $R_{I_q,\bN}^{I_r}$, as in Theorem \ref{thm:E5}.
However, we can derive more concrete descriptions of these congruences by being more explicit about the permitted sequences $\bN$.

For each $0\leq q<\xi$, fix some vector space $V_q$ from $\C$ with $n_q = \dim(V_q)\geq q$, and fix a basis $\BB_q=\{\vv_{q,1},\ldots,\vv_{q,n_q}\}$ of $V_q$.  We define the idempotent $\ve_q\in D_q(\LL_{V_q})$ by
\[
\vv_{q,i} \ve_q = \begin{cases}
\vv_{q,i} &\text{if $1\leq i\leq q$}\\
\vzero &\text{otherwise.}
\end{cases}
\]
We then let $\G_q$ be the $\H$-class of the idempotent $\ve_q$.  

We have already noted that $\G_q\cong\GL_q=\GL_q(\F)$, and it is easy to see this in matrix form.
With respect to the basis $\BB_q$, the linear transformation $\ve_q$ corresponds to the $n_q\times n_q$ matrix given in block form by $\tmat EOOO$, where $E=E_q$ is the $q\times q$ identity matrix, and where the $O$'s denote zero matrices of appropriate sizes.
Any matrix $\H$-related to $\tmat EOOO$ is of the form $\tmat AOOO$ for some invertible $q\times q$ matrix $A\in\GL_q$ (cf.~\cite[Lemma 3.1]{DE2018}).  The element of~$\G_q$ corresponding to the matrix $\tmat AOOO$ will be denoted $A^\natural$.  For $\Om\sub\GL_q$, we write $\Om^\natural=\set{A^\natural}{A\in\Om}$.  In particular,~$\G_q=\GL_q^\natural$.

The next step is to give a concrete description of the $\nu_N$ relations.  We do this in Lemma~\ref{lem:nuN_L}, for which we first introduce some terminology analogous to that used for transformations and partitions.  With this in mind, fix some $\H$-related elements $\al,\be$ from $D_q(\LL_{U,V})$, where $U,V\in\C$.  Let $\BB=\{\vu_1,\ldots,\vu_d\}$ be a basis of $U$ such that $\{\vu_{q+1},\ldots,\vu_d\}$ is a basis of $\Ker(\al)$.  Let $\BB'=\{\vv_1,\ldots,\vv_e\}$ be a basis of $V$ such that $\vv_i=\vu_i\al$ for all $1\leq i\leq q$.  Since $\Ker(\al)=\Ker(\be)$ and $\im(\al)=\im(\be)$, the matrix of $\be$ with respect to the bases~$\BB$ and $\BB'$ has the form $\tmat AOOO$ for some invertible $q\times q$ matrix $A\in\GL_q$ (and where again the $O$'s denote zero matrices of appropriate sizes).  We then define $\phi(\al,\be)=A\in\GL_q$, again noting that this is defined only up to conjugacy (change of basis) in $\GL_q$.  

\begin{lemma}\label{lem:nuN_L}
If $0\leq q<\xi$, and if $N\normal\GL_q$, then $\nu_{N^\natural} = \bigset{(\al,\be)\in\H\restr_{D_q}}{\phi(\al,\be)\in N}$.  
\end{lemma}

\pf
Let $(\al,\be)\in\H\restr_{D_q}$, say with $\al,\be\in\LL_{U,V}$, and fix the bases $\BB$ and $\BB'$ as above.  Define $\ga_1\in\LL_{V_q,U}$, $\ga_2\in\LL_{V,V_q}$, $\de_1\in\LL_{U,V_q}$ and $\de_2\in\LL_{V_q,V}$ by
\begin{align*}
\vv_{q,i}\ga_1 &= \begin{cases}
\vu_i\phantom{\vv_{q,i}} &\text{for $1\leq i\leq q$}\\
\vzero &\text{for $q<i\leq n_q$,}
\end{cases}
&
\vu_{i}\de_1 &= \begin{cases}
\vv_{q,i}\phantom{\vv_i} &\text{for $1\leq i\leq q$}\\
\vzero &\text{for $q<i\leq d$,}
\end{cases}
\\[1truemm]
\vv_{i}\ga_2 &= \begin{cases}
\vv_{q,i}\phantom{\vu_i} &\text{for $1\leq i\leq q$}\\
\vzero &\text{for $q<i\leq e$,}
\end{cases}
&
\vv_{q,i}\de_2 &= \begin{cases}
\vv_{i}\phantom{\vv_{q,i}} &\text{for $1\leq i\leq q$}\\
\vzero &\text{for $q<i\leq n_q$.}
\end{cases}
\end{align*}
Since $(E_q^\natural,A^\natural) = (\ga_1\al\ga_2,\ga_1\be\ga_2)$ and $(\al,\be) = (\de_1E_q^\natural\de_2,\de_1A^\natural\de_2)$, we have (cf.~Lemma \ref{lem:nuN}\ref{nuN1})
\[
(\al,\be)\in\nu_{N^\natural} \iff (E_q^\natural,A^\natural)\in\nu_{N^\natural} \iff A^\natural\in N^\natural \iff \phi(\al,\be)=A\in N. \qedhere
\]
\epf

With the above definition of $\G_q=\GL_q^\natural$ in place, we may also calculate the groups $Z_q$.  Indeed, since $\nu_{Z_q}=\ze\restr_{D_q}$, Lemma \ref{lem:nu_trick} tells us that
\[
Z_q = \set{\al\in\G_q}{(\ve_q,\al)\in\ze} = \set{c\ve_q}{c\in\Fx} = Z(\GL_q)^\natural.
\]
Thus,
\[
\ze = \nu_{Z(\GL_0)^\natural} \cup \nu_{Z(\GL_1)^\natural} \cup  \nu_{Z(\GL_2)^\natural} \cup \cdots.
\]

The subgroups of $Z(\GL_q)$ have the form 
\[
\Z_q(H) = \set{cE_q}{c\in H} \qquad\text{where $H\leq\Fx$,}
\]
and these are all normal in $\GL_q$.
The NZ-tuples in $\LL$ are of the form $\bN=(N_0,N_1,\ldots)$, where each $N_q\normal \G_q$ and $N_q\leq Z_q$, and where $N_q^{\downarrow_p}\leq N_p$ for all $0\leq p<q$.  The next result explicates the effect of $\downarrow_p$.

\begin{lemma}\label{lem:Zqp_L}
For any $0\leq p\leq q<\xi$, and for any $H\leq\Fx$, we have $(\Z_q(H)^\natural)^{\downarrow_p}=\Z_p(H)^\natural$.
\end{lemma}

\pf
For convenience, we write $N_q=\Z_q(H)^\natural$ and $N_p=\Z_p(H)^\natural$; we must show that $N_q^{\downarrow_p}=N_p$.  Note that $N_q=\set{c\ve_q}{c\in H}$ and similarly for $N_p$.

First, since $\ve_p\leqJ\ve_q$, we have $\ve_p=\ga_1\ve_q\ga_2$ for some $\ga_1,\ga_2\in\LL$.  For any $c\in H$, we have $(\ve_q,c\ve_q)\in\nu_{N_q}$, and so $(\ve_p,c\ve_p) = (\ga_1\ve_q\ga_2,\ga_1(c\ve_q)\ga_2) \in \nu_{N_q}^\sharp\restr_{D_p}=\nu_{N_q^{\downarrow_p}}$.  It follows from Lemma~\ref{lem:nuN}\ref{nuN1} that $c\ve_p\in N_q^{\downarrow_p}$.  Since $c\in H$ was arbitrary, it follows that $N_p \leq N_q^{\downarrow_p}$.  

We will show that $N_q^{\downarrow_p}\leq N_p$ by showing that $\nu_{N_q^{\downarrow_p}}\sub\nu_{N_p}$.  Now, $\nu_{N_q^{\downarrow_p}}=\nu_{N_q}^\sharp\restr_{D_p}$.  As in Remark \ref{rem:tau_N}, the congruence $\nu_{N_q}^\sharp$ is the union of $\De_\LL$ with the transitive closure of the set
\[
\Om = \bigset{(\ga_1\al\ga_2,\ga_1\be\ga_2)}{\al,\be\in N_q,\ \ga_1,\ga_2\in \LL,\ \br(\ga_1)=\bd(\ga_2)=V_q}.
\]
Any pair from $\Om$ satisfies $(\ga_1\al\ga_2,\ga_1\be\ga_2)\in\H\sub\J$, since $(\al,\be)\in\nu_{N_q}\sub\nu_{Z_q}\sub\ze$, and since~$\ze$ is an $\H$-congruence.  It follows that $\nu_{N_q^{\downarrow_p}}=\nu_{N_q}^\sharp\restr_{D_p}$ is in fact the transitive closure of the set~$\Om\restr_{D_p}$.  Since $\nu_{N_p}$ is an equivalence, we can therefore show that $\nu_{N_q^{\downarrow_p}}\sub\nu_{N_p}$ by showing that $\Om\restr_{D_p}\sub\nu_{N_p}$.  With this in mind, consider some pair $(\ga_1\al\ga_2,\ga_1\be\ga_2)\in\Om\restr_{D_p}$.  Since $\al,\be\in N_q$, we have $\be=c\al$ for some $c\in H$ (since $\al=c_1\ve_q$ and $\be=c_2\ve_q$ for some $c_1,c_2\in H$).  
Since $\ga_1\al\ga_2\in D_p$, we have $(\ga_1\al\ga_2,\ve_p)\in\D$, and so $\ga_1\al\ga_2=\de_1\ve_p\de_2$ for some $\de_1,\de_2\in\LL$.
But then $(\ga_1\al\ga_2,\ga_1\be\ga_2)=\big(\ga_1\al\ga_2,\ga_1(c\al)\ga_2\big)=\big(\de_1\ve_p\de_2,\de_1(c\ve_p)\de_2\big)\in\nu_{N_p}$.  This completes the proof.
\epf

From Lemma \ref{lem:Zqp_L}, we immediately deduce the following:

\begin{cor}\label{cor:Zqp_L}
If $0\leq p\leq q<\xi$, and if $H_1,H_2\leq\Fx$, then
\[
(\Z_q(H_1)^\natural)^{\downarrow_p}\leq\Z_p(H_2)^\natural \iff H_1\leq H_2.  \epfreseq
\]
\end{cor}

Here is the main result of this section.  The statement involves the relations $R_{I,\bN}^S$ defined before Theorem \ref{thm:E5}.

\begin{thm}\label{thm:L}
Let $\C$ be a non-empty set of finite-dimensional vector spaces over a field $\F$, let $\LL=\LL(\C)$ be the linear category over $\C$, and keep the above notation.
Then for any~$r\geq0$ (including $r=\om$ if $\xi=\om$), the congruences on $I_r$ are precisely the universal congruence~$\nab_{I_r}$ and those of the form
\[
R_{I_q,\bN}^{I_r} \qquad\text{for some $0\leq q<r$ and some $\bN=(N_{q+1}^\natural,\Z_{q+2}(H_{q+2})^\natural,\Z_{q+3}(H_{q+3})^\natural,\ldots)$,}
\]
where 
\[
N_{q+1}\normal\GL_{q+1}\COMMA
\Fx\geq H_{q+2}\geq H_{q+3}\geq\cdots\COMMA
\Z_{q+1}(H_{q+2})\leq N_{q+1} \text{ if $q\leq r-2$}.
\]
\end{thm}

\pf
We have already noted that all of the conditions of Theorem \ref{thm:E5} are satisfied.  The only thing left to check is that, with $\bN$ as in the statement, and writing $N_p=\Z_p(H_p)$ for $p\geq q+2$, we have $(N_t^\natural)^{\downarrow_p}\leq N_p^\natural$ for all $t>p\geq q+1$.  But this follows quickly from Corollary \ref{cor:Zqp_L}.
\epf

\subsection{Visualising the lattices}\label{subsect:visual_L}

Compared to the categories considered in Sections \ref{sect:T}--\ref{sect:J}, congruence lattices of the ideals of~$\LL$ are rather complex and harder to visualise.  
Nevertheless, we may gain some intuition by considering 
$\Cong(I_r)$ for some small $r$ and small underlying fields $\F$.

The first key is to recall the normal subgroup structure of $\GL_q$.  The next result is well known, though non-trivial to prove; see for example \cite[Theorem~4.9]{Artin1957}.

\begin{prop}\label{prop:GLq}
If $q\geq3$, or if $q=2$ and $|\F|\geq4$, then the normal subgroups of the general linear group $\GL_q=\GL_q(\F)$ are precisely
\ben
\item \label{GLq1} the subgroups of the centre $Z(\GL_q)$, and
\item \label{GLq2} the subgroups of $\GL_q$ containing the special linear group $\SL_q=\SL_q(\F)$.  \epfres
\een
\end{prop}

We have already noted that the subgroups of $Z(\GL_q)$ are of the form ${\Z_q(H)=\set{cE_q}{c\in H}}$ for $H\leq\Fx$.  It is also not hard to show that the subgroups of $\GL_q$ containing $\SL_q$ are of the form
\[
\S_q(H) = \set{A\in\GL_q}{\det(A)\in H} \qquad\text{for $H\leq\Fx$,}
\]
where as usual $\det(A)$ is the determinant of $A$.
Thus, apart from the exceptional combinations of $q$ and $|\F|$ excluded in Proposition~\ref{prop:GLq}, the lattice of normal subgroups of~$\GL_q$ is 
\[
\NN(\GL_q) = \bigset{\Z_q(H),\ \S_q(H)}{H\leq\Fx}.
\]

The sublattices $\bigset{\Z_q(H)}{H\leq\Fx}$ and $\bigset{\S_q(H)}{H\leq\Fx}$ are both clearly isomorphic to $\NN(\Fx)$, but the ways in which these sublattices are ``glued together'' to form $\NN(\GL_q)$ depends on $q$ (and in particular, which elements of $\F$ are $q$th powers).
The structure of the lattice $\NN(\Fx)$ itself depends on the field $\F$; if $|\F|=n$ is finite, $\NN(\Fx)$ is isomorphic to the lattice of divisors of $|\Fx|=n-1$, since $\Fx$ is cyclic.

As an example, consider the case in which $\F=\ZZ_7=\{0,1,\ldots,6\}$ is the field with seven elements.  The (normal) subgroups of the multiplicative group $\Fx=\{1,\ldots,6\}$ are
\[
H_1=\{1\} \COMMA H_2=\{1,6\} \COMMA H_3=\{1,2,4\} \COMMA H_4= \Fx.
\]
The Hasse diagram of $\NN(\Fx)=\{H_1,H_2,H_3,H_4\}$ is pictured in Figure \ref{fig:GLq} (in red).  Figure \ref{fig:GLq} also pictures the lattices $\NN(\GL_q(\F))$ for $1\leq q\leq 6$, illustrating the various ways the two above-mentioned sublattices may be glued together.  To avoid clutter, for each $1\leq q\leq6$ and $1\leq i\leq4$ we write $Z_{q,i}=\Z_q(H_i)$ and $S_{q,i}=\S_q(H_i)$.

\begin{figure}[ht]
\begin{center}
\scalebox{.76}{
\begin{tikzpicture}[rotate=60]
\nc\radiusfordiagram{7.4}
%
%
\node () at (0,0) {
\begin{tikzpicture}[scale=1.5]
\redcongnode000{H_1}
\redcongnode{-1}11{H_2}
\redcongnode112{H_3}
\redcongnode023{H_4}
\draw(0)--(1)--(3)--(2)--(0);
\draw(0.62,2.04) node {$=\Fx$};
\end{tikzpicture}
};
%
%
\node () at (60:\radiusfordiagram) {
\begin{tikzpicture}[scale=1.5]
\Znode1100
\Znode12{-1}1
\Znode1311
\Znode1402
\draw(Z11)--(Z12)--(Z14)--(Z13)--(Z11);
\draw(-.76,2.02) node {$\GL_1=$};
\draw(-.76,0.02) node {$\SL_1=$};
\end{tikzpicture}
};
%
%
\node () at (0:\radiusfordiagram) {
\begin{tikzpicture}[scale=1.5]
\Znode2100
\Znode22{-1}1
\Znode2311
\Znode2402
\draw(Z21)--(Z22)--(Z24)--(Z23)--(Z21);
\Snode2100{-2}2
\Snode22{-1}1{-2}2
\Snode2311{-2}2
\Snode2402{-2}2
\draw(S21)--(S22)--(S24)--(S23)--(S21);
\draw(Z22)--(S21) (Z24)--(S23);
\draw(-.76-2,2.02+2) node {$\GL_2=$};
\draw(-.76-2,2.02) node {$\SL_2=$};
\end{tikzpicture}
};
%
%
\node () at (300:\radiusfordiagram) {
\begin{tikzpicture}[scale=1.5]
\Znode3100
\Znode32{-1}1
\Znode3311
\Znode3402
\draw(Z31)--(Z32)--(Z34)--(Z33)--(Z31);
\Snode310022
\Snode32{-1}122
\Snode331122
\Snode340222
\draw(S31)--(S32)--(S34)--(S33)--(S31);
\draw(Z33)--(S31) (Z34)--(S32);
\draw(-.76+2,2.02+2) node {$\GL_3=$};
\draw(0.7+2,2.01) node {$=\SL_3$};
\end{tikzpicture}
};
%
%
\node () at (240:\radiusfordiagram) {
\begin{tikzpicture}[scale=1.5]
\Znode4100
\Znode42{-1}1
\Znode4311
\Znode4402
\draw(Z41)--(Z42)--(Z44)--(Z43)--(Z41);
\Snode4100{-2}2
\Snode42{-1}1{-2}2
\Snode4311{-2}2
\Snode4402{-2}2
\draw(S41)--(S42)--(S44)--(S43)--(S41);
\draw(Z42)--(S41) (Z44)--(S43);
\draw(-.76-2,2.02+2) node {$\GL_4=$};
\draw(-.76-2,2.02) node {$\SL_4=$};
\end{tikzpicture}
};
%
%
\node () at (180:\radiusfordiagram) {
\begin{tikzpicture}[scale=1.5]
\Znode5100
\Znode52{-1}1
\Znode5311
\Znode5402
\draw(Z51)--(Z52)--(Z54)--(Z53)--(Z51);
\fill[white] (-.5,1.5)circle(.07);
\fill[white] (.5,1.5)circle(.07);
\Snode510001
\Snode52{-1}101
\Snode531101
\Snode540201
\draw(S51)--(S52)--(S54)--(S53)--(S51);
\draw(Z51)--(S51) (Z52)--(S52) (Z53)--(S53) (Z54)--(S54);
\draw(-.76,2.02+1) node {$\GL_5=$};
\node (SL5) at (-1.3,0) {$\SL_5$};
\draw[-{latex}] (SL5) --(S51);
\end{tikzpicture}
};
%
%
\node () at (120:\radiusfordiagram) {
\begin{tikzpicture}[scale=1.5]
\Znode6100
\Znode62{-1}1
\Znode6311
\Znode6402
\draw(Z61)--(Z62)--(Z64)--(Z63)--(Z61);
\Snode610003
\Snode62{-1}103
\Snode631103
\Snode640203
\draw(S61)--(S62)--(S64)--(S63)--(S61);
\draw(Z64)--(S61) ;
\draw(-.76,2.02+3) node {$\GL_6=$};
\draw(-.76,3.02) node {$\SL_6=$};
\end{tikzpicture}
};
\end{tikzpicture}
}
\caption{Hasse diagrams of $\NN(\Fx)$, and of $\NN(\GL_q)$ for $q=1,\ldots,6$, where $\F=\ZZ_7$.  See Section~\ref{subsect:visual_L}.}
\label{fig:GLq}
\end{center}
\end{figure}

We now describe the lattices $\Cong(I_r)$, $r=1,2,3$, for $\F=\ZZ_7$ (with $\Cong(I_0)$ of course being trivial).  
To understand these, it will be useful to keep in mind Theorem \ref{thm:E6}\ref{E6ii} as well, as it will help us 
visualise the construction of $\Cong(I_{r+1})$ from $\Cong(I_r)$.

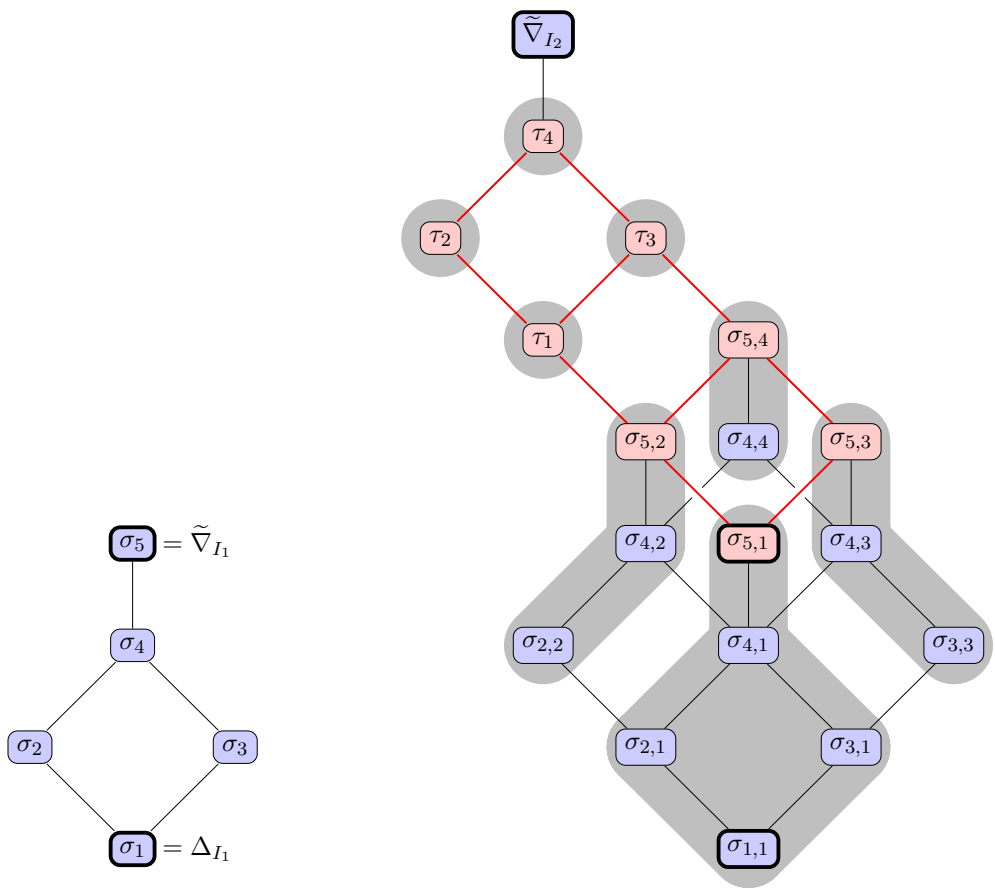
\begin{figure}[ht]
\begin{center}
\scalebox{.9}{
\begin{tikzpicture}[scale=1.5]
\begin{scope}[shift={(-6,0)}]
\thickcongnode000{\si_1}
\congnode{-1}11{\si_2}
\congnode112{\si_3}
\congnode023{\si_4}
\thickcongnode034{\si_5}
\draw (4)--(3)--(1)--(0)--(2)--(3);
\draw(0.63,0) node {$=\De_{I_1}$};
\draw(0.63,3) node {$=\nab_{I_1}$};
\end{scope}
\draw[rounded corners=2em,line width=3em,black!25,cap=round]
(0,0)--(1,1) (1,1)--(0,2) (0,2)--(0,3) (0,2)--(-1,1) (-1,1)--(0,0)
(2,2)--(1,3) (1,3)--(1,4)
(-2,2)--(-1,3) (-1,3)--(-1,4)
(0,4)--(0,5)
(-2,5)--(-2.001,5)
(-2,7)--(-2.001,7)
(-1,6)--(-1.001,6)
(-3,6)--(-3.001,6)
;
\fill[black!25](0,0)--(1,1)--(0,2)--(-1,1)--(0,0);
\thicksicongnode000011
\sicongnode10{-1}121
\sicongnode201131
\sicongnode300241
\thickredsicongnode400351
\sicongnode11{-2}222
\sicongnode31{-1}342
\redsicongnode41{-1}452
\sicongnode222233
\sicongnode321343
\redsicongnode421453
\sicongnode330444
\redsicongnode430554
\taucongnode0{-2}51
\taucongnode1{-3}62
\taucongnode2{-1}63
\taucongnode3{-2}74
\node[rounded corners,rectangle,draw,fill=blue!20, ultra thick] (nab) at (-2,8) {$\nab_{I_2}$};
\draw
(00)--(10)--(30)--(20)--(00)
(10)--(11)--(31)--(30)--(32)--(22)--(20)
(30)--(40)--(41)--(0)--(1)--(3)--(2)--(43)--(42)--(40)
(33)--(43)--(41)--(31)--(33)--(32)--(42)
(0)--(2)
(3)--(nab)
;
\fill[white] (-.5,3.5)circle(.07);
\fill[white] (.5,3.5)circle(.07);
\draw[thick,red](40)--(42)--(43)--(2)--(3)--(1)--(0)--(41)--(40) (41)--(43) (2)--(0);
\end{tikzpicture}
}
\caption{Hasse diagrams of $\Cong(I_r)$ for $r=1$ (left) and $r=2$ (right), where $I_r=I_r(\LL)$ and $\F=\ZZ_7$.  
See Section \ref{subsect:visual_L}.}
\label{fig:GLM27}
\end{center}
\end{figure}

By Theorem \ref{thm:L}, we have
\[
\Cong(I_1) = \set{R_{I_0,N^\natural}^{I_1}}{N\normal\GL_1}\cup\{\nab_{I_1}\}.
\]
Figure \ref{fig:GLM27} (left) pictures the lattice $\Cong(I_1)$, where for simplicity we denote the congruences by
\[
\si_i = R_{I_0,\Z_1(H_i)^\natural}^{I_1} \text{ for $i=1,2,3,4$} \AND
\si_5 = \nab_{I_1}.
\]
In both lattices pictured in Figure \ref{fig:GLM27}, Rees congruences are indicated by bold outline.
From this we may construct the lattice 
\[
\Cong(I_2) = \set{\si\cup\nu_N}{\si\in\Cong(I_1),\ N\normal \G_2,\ \tau_N\sub\si} \cup \{\nab_{I_2}\}.
\]
By \eqref{E4N1}, \eqref{E4N2} and Lemma \ref{lem:Zqp_L}, we see that for $i=1,2,3,4$,
\[
\text{$\tau_N = \nab_{I_1} = \si_5$ for $N=\S_2(H_i)^\natural$} \AND \text{$\tau_N = \si_i$ for $N=\Z_2(H_i)^\natural$.}
\]
From this information, we can construct the layers $\Lam(N)$, $N\normal\GL_2^\natural$, as defined in Remark \ref{rem:E3}.  We have
\[
\Lam(N)= \begin{cases}
\{\si_5\cup\nu_N\} &\text{if $N=\S_2(H_i)^\natural$ for any $i=1,2,3,4$}\\
\set{\si_i\cup\nu_N}{i=1,2,3,4,5} &\text{if $N=\Z_2(H_1)^\natural$}\\
\set{\si_i\cup\nu_N}{i=2,4,5} &\text{if $N=\Z_2(H_2)^\natural$}\\
\set{\si_i\cup\nu_N}{i=3,4,5} &\text{if $N=\Z_2(H_3)^\natural$}\\
\set{\si_i\cup\nu_N}{i=4,5} &\text{if $N=\Z_2(H_4)^\natural$.}
\end{cases}
\]
These eight layers, together with $\nab_{I_2}$, form the lattice $\Cong(I_2)$ pictured in Figure \ref{fig:GLM27} (right).  In the diagram, we write for simplicity
\[
\si_{i,j} = \si_i \cup \nu_{\Z_2(H_j)^\natural} \AND \tau_j=\nab_{I_1}\cup\nu_{\S_2(H_j)^\natural} \quad \text{ for $1\leq i\leq5$ and $1\leq j\leq4$.}
\]
Grey shaded regions in the diagram indicate the layers $\Lam(N)$, as above, and these fit together in the way described in Remark \ref{rem:E3} (cf.~Figures \ref{fig:E31} and \ref{fig:E32}).  The congruences shaded red in the diagram of $\Cong(I_2)$ are those of the form $\nab_{I_1}\cup\nu_N$ ($N\normal\GL_2^\natural$); these form a lattice isomorphic to $\NN(\GL_2)$; cf.~Figure \ref{fig:GLq}.  This calculation has been verified by GAP \cite{GAP4,GAP} in the case where~$\C$ consists of a single vector space of dimension $2$: i.e., when $\LL$ is the monoid of $2\times2$ matrices over~$\F$.

One may then construct $\Cong(I_3)$ with the same kinds of calculations, starting with the observation (again using \eqref{E4N1}, \eqref{E4N2} and Lemma \ref{lem:Zqp_L}) that for $i=1,2,3,4$,
\[
\text{$\tau_N = \nab_{I_2}$ for $N=\S_3(H_i)^\natural$} \AND \text{$\tau_N = \si_{i,i}$ for $N=\Z_3(H_i)^\natural$.}
\]
Figure \ref{fig:GLM37} shows the lattice $\Cong(I_3)$, again with the layers $\Lam(N)$, $N\normal\GL_3^\natural$, shaded, and Rees congruences denoted by filled vertices.  This time we have not labelled the congruences, but we have again indicated in red the sublattice isomorphic to $\NN(\GL_3)$; cf.~Figure \ref{fig:GLq}.  It is not so easy to verify this calculation using GAP, as even when $\C$ consists of a single vector space of dimension $3$, $\LL$ has size $7^9 = 40\ 353\ 607$.

In principle, this process could be iterated further, but we have not attempted to do this.

\begin{figure}[ht]
\begin{center}
\scalebox{.8}{
\begin{tikzpicture}[scale=1]
\nc\hori{4}
\nc\verti{1.5}
\nc\aaa{-\hori} \nc\bbb{\verti}
\nc\ccc{\hori} \nc\ddd{\verti}
\nc\eee{0} \nc\fff{2*\verti}
\draw[rounded corners=2em,line width=2em,black!20,cap=round]
(0,0)--(2,2) (2,2)--(1,3) (1,3)--(1,4) (1,4)--(-2,7) (-2,7)--(-2,8)
(0,0)--(-2,2) (-2,2)--(-1,3) (-1,3)--(-1,4) (-1,4)--(-3,6) (-3,6)--(-2,7)
;
\fill[black!20](0,0)--(2,2)--(1,3)--(1,4)--(-2,7)--(-3,6)--(-1,4)--(-1,3)--(-2,2)--(0,0);
\begin{scope}[shift={(\aaa,\bbb)}]
\draw[rounded corners=2em,line width=2em,black!20,cap=round]
(-2,2)--(0,4) (0,4)--(0,5) (0,5)--(-2,7) (-2,7)--(-2,8) (-2,7)--(-3,6) (-3,6)--(-1,4) (-1,4)--(-1,3)
;
\fill[black!20](-1,3)--(0,4)--(0,5)--(-2,7)--(-3,6)--(-1,4)--(-1,3);
\end{scope}
\begin{scope}[shift={(\ccc,\ddd)}]
\draw[rounded corners=2em,line width=2em,black!20,cap=round]
(2,2)--(0,4) (0,4)--(0,5) (0,5)--(-2,7) (-2,7)--(-2,8) (1,3)--(1,4) (1,4)--(0,5)
;
\end{scope}
\begin{scope}[shift={(\eee,\fff)}]
\draw[rounded corners=2em,line width=2em,black!20,cap=round]
(0,4)--(0,5) (0,5)--(-2,7) (-2,7)--(-2,8) 
;
\end{scope}
\begin{scope}[shift={(\eee+\ccc,\fff+\ddd)}]
\draw[rounded corners=2em,line width=2em,black!20,cap=round]
(-2,8) --(-2.001,8)
;
\end{scope}
\begin{scope}[shift={(\eee+2*\ccc,\fff+2*\ddd)}]
\draw[rounded corners=2em,line width=2em,black!20,cap=round]
(-2,8) --(-2.001,8)
;
\end{scope}
\begin{scope}[shift={(\ccc+\ccc,\ddd+\ddd)}]
\draw[rounded corners=2em,line width=2em,black!20,cap=round]
(-2,8) --(-2.001,8)
;
\end{scope}
\begin{scope}[shift={(\ccc+2*\ccc,\ddd+2*\ddd)}]
\draw[rounded corners=2em,line width=2em,black!20,cap=round]
(-2,8) --(-2.001,8)
;
\end{scope}
\foreach \a/\b in {-2/2,-1/3,0/4,-1/4,0/5,-3/6,-1/6,-2/7,-2/8} {\draw[gray!80](\a,\b)--(\a+\aaa,\b+\bbb);}
\foreach \a/\b in {0/4,0/5,-1/6,-2/7,-2/8,2/2,1/3,1/4} {\draw[gray!80](\a,\b)--(\a+\ccc,\b+\ddd);}
\foreach \a/\b in {0/4,0/5,-1/6,-2/7,-2/8} {\draw[gray!80](\a+\ccc,\b+\ddd)--(\a+\eee,\b+\fff)--(\a+\aaa,\b+\bbb);}
\foreach \a/\b in {-2/8} {\draw[thick,red!80](\a+\aaa,\b+\bbb)--(\a+\eee+2*\ccc,\b+\fff+2*\ddd);}
\foreach \a/\b in {-2/8} {\draw[thick,red!80](\a,\b)--(\a+\aaa,\b+\bbb);}
\foreach \a/\b in {-2/8} {\draw[thick,red!80](\a+\ccc,\b+\ddd)--(\a+\ccc+\aaa,\b+\ddd+\bbb);}
\foreach \a/\b in {-2/8} {\draw[thick,red!80](\a,\b)--(\a+3*\ccc,\b+3*\ddd);}
\foreach \a/\b in {-2/8} {\draw[thick,red!80](\a+3*\ccc,\b+3*\ddd)--(\a+\eee+2*\ccc,\b+\fff+2*\ddd);}
\foreach \a/\b in {-2/8} {\draw[thick,red!80](\a+\eee+\ccc,\b+\fff+\ddd)--(\a+2*\ccc,\b+2*\ddd);}
\foreach \a/\b in {-2/8} {\draw(\a+\eee+2*\ccc,\b+\fff+2*\ddd)--(\a+\eee+2*\ccc,\b+\fff+2*\ddd+1);}
\foreach \a/\b in {0/4,-1/3,-2/7} {\draw(\a+\aaa,\b+\bbb)--(\a+\aaa,\b+1+\bbb);}
\foreach \a/\b in {-2/2,-1/3} {\draw(\a+\aaa,\b+\bbb)--(\a+1+\aaa,\b+1+\bbb);}
\foreach \a/\b in {1/3,0/4,-2/7} {\draw(\a+\ccc,\b+\ddd)--(\a+\ccc,\b+1+\ddd);}
\foreach \a/\b in {2/2,1/3,1/4,0/5,-1/6} {\draw(\a+\ccc,\b+\ddd)--(\a-1+\ccc,\b+1+\ddd);}
\foreach \a/\b in {0/4,-2/7} {\draw(\a+\eee,\b+\fff)--(\a+\eee,\b+1+\fff);}
\foreach \a/\b in {0/5,-1/6} {\draw(\a+\eee,\b+\fff)--(\a-1+\eee,\b+1+\fff);}
\fill (-2,8) circle(.1);
\foreach \a/\b in {0/2,1/3,0/4,-1/3,-2/7} {\draw(\a,\b)--(\a,\b+1);}
\foreach \a/\b in {0/0,1/1,-1/1,0/2,0/3,-1/4,-2/5} {\squarelattice{\a}{\b}}
\foreach \a/\b in {-1/4,-2/5} {\squarelattice{\a+\aaa}{\b+\bbb}}
\foreach \a/\b in {0/4,0/5,-1/6,-2/7,-2/8} {\colvert{\a+\eee}{\b+\fff}{white}}
\foreach \a/\b in {-2/8} {\colvert{\a+\eee+2*\ccc}{\b+\fff+2*\ddd+1}{white}}
\foreach \a/\b in {-2/2,-1/3,0/4,-2/8} {\colvert{\a+\aaa}{\b+\bbb}{white}}
\foreach \a/\b in {2/2,1/3,0/4,1/4,0/5,-1/6,-2/7,-2/8} {\colvert{\a+\ccc}{\b+\ddd}{white}}
\foreach \a/\b in {-2/8} {\colvert{\a+\eee+\ccc}{\b+\fff+\ddd}{white}}
\foreach \a/\b in {-2/8} {\colvert{\a+\eee+2*\ccc}{\b+\fff+2*\ddd}{white}}
\foreach \a/\b in {-2/8} {\colvert{\a+2*\ccc}{\b+2*\ddd}{white}}
\foreach \a/\b in {-2/8} {\colvert{\a+3*\ccc}{\b+3*\ddd}{white}}
\foreach \x/\y in {0/0,0/3} {\colvert{\x}{\y}{black}}
\foreach \a in {3} {\colvert{-2+\aaa+\a*\ccc}{8+\bbb+\a*\ddd+1}{black}}%
\foreach \x/\y in {-2/8} {\redcolvert{\x}{\y}{red}}
\foreach \a in {1,2,3} {\redcolvert{-2+\a*\ccc}{8+\a*\ddd}{white}}
\foreach \a in {0,1,2,3} {\redcolvert{-2+\aaa+\a*\ccc}{8+\bbb+\a*\ddd}{white}}
\end{tikzpicture}
}
\caption{Hasse diagram of $\Cong(I_3)$, where $I_3=I_3(\LL)$ and $\F=\ZZ_7$.  See Section \ref{subsect:visual_L}.}
\label{fig:GLM37}
\end{center}
\end{figure}

\subsection{Congruences on ideals of $\PL$}\label{subsect:PL}

By contrast with the linear categories, the ideals $I_r$ ($2\leq r<\xi$) of the projective linear category $\PL=\PL(\C)$ satisfy 
$\Sep{I_r}$, as we observed in Remark \ref{rem:P2L}, so we may return to using Theorem~\ref{thm:E2} (and Proposition \ref{prop:small_ideal}) to describe the congruences on its ideals.

In what follows, we keep the above meaning of the parameter $\xi$, and denote $\D$-classes and ideals by $D_q=D_q(\PL)$ and $I_r=I_r(\PL)$, for appropriate~$q$ and $r$.  
For each $0\leq q<\xi$, fix some group $\H$-class $\G_q$ in $D_q$, so that $\G_q$ is isomorphic to the projective general linear group $\PGL_q=\PGL_q(\F)$.  An ideal $I_r$ ($r\geq1$) contains the IN-pairs $(I_q,N)$ for $0\leq q<r$ and $N\normal\G_{q+1}$.  As usual, these induce the congruences $R_{I_q,N}^{I_r}$, including $R_{I_q}^{I_r}$ when $N$ is trivial, and we note that $\De_{I_r}=R_{I_0}^{I_r}$, since $I_0=D_0$ is trivial.

\begin{thm}\label{thm:PL}
Let $\C$ be a non-empty set of finite-dimensional vector spaces over a field $\F$, let $\PL=\PL(\C)$ be the projective linear category over $\C$, and keep the above notation.
Then for any $r\geq0$ (including $r=\om$ if $\xi=\om$),
\[
\Cong(I_r) = \set{R_{I_q,N}^{I_r}}{0\leq q<r,\ N\normal\G_{q+1}}\cup\{\nab_{I_r}\}.
\]
\end{thm}

\pf
For $r=0$, the result says that $\Cong(I_0)=\{\nab_{I_0}\}$, which follows from $I_0$ being trivial.

For $r=1$, the result says that $\Cong(I_1)=\{\De_{I_1},\nab_{I_1}\}$.  This follows quickly from Proposition~\ref{prop:small_ideal}, keeping in mind that $I_1$ is retractable, $I_0$ is trivial, and $D_1$ is $\H$-trivial.  The technical assumption from Proposition \ref{prop:small_ideal} follows from Lemma~\ref{lem:P2L}\ref{P2Lii}.

As usual, the result then follows from Theorem \ref{thm:E2}, using Lemma \ref{lem:P2L} (cf.~Remark \ref{rem:P2L}) to establish $\Mult{I_r}$, and hence $\Sep{I_r}$, for the ideals $I_r$ ($1<r<\xi$).
\epf

\begin{rem}
Theorem \ref{thm:PL} could also be deduced from Theorem \ref{thm:L}, since $\PL=\LL/\ze$.
\end{rem}

\begin{rem}
We have not drawn any congruence lattices for ideals $I_r(\PL)$,
since for small values of $r$ and $|\F|$, these are typically chains.  Indeed, this follows from Theorem \ref{thm:PL} (cf.~Theorem \ref{thm:E2} and Figure \ref{fig:E2}), and the fact that the normal subgroups of $\PGL_q=\PGL_q(\F)$ form chains for small~$q$ and~$|\F|$, as may be readily verified using GAP \cite{GAP4}.  However, the lattices $\NN(\PGL_q)$ are not always chains, and hence neither are the lattices $\Cong(I_r(\PL))$ in general.  For example, consider the case in which $q=10$, and $\F=\ZZ_{31}$ is the field of size~$31$, and for simplicity, write $G=\GL_{10}(\ZZ_{31})$, $Z=Z(G)$ and $P=\PGL_{10}(\ZZ_{31})=G/Z$.  By the Third Isomorphism Theorem, the normal subgroups of $P$ are all of the form $N/Z$ where $N$ is a normal subgroup of~$G$ containing $Z$.  Now, $Z=\set{cE_{10}}{c\in\ZZ_{31}^\times}$, so the possible determinants of elements of $Z$ are the tenth powers of elements of $\ZZ_{31}^\times$, which are precisely~$1$,~$5$ and~$25$.  Thus, using the above notation, the only normal subgroups $N$ of $G$ containing~$Z$ are those of the form $N=\S_{10}(H)=\set{A\in\GL_{10}}{\det(A)\in H}$ where~$H$ is a subgroup of~$\ZZ_{31}^\times$ containing $\{1,5,25\}$.  The possible subgroups $H$ are of size $3$, $6$, $15$ and~$30$; let $H_1$ and~$H_2$ be the subgroups of size~$6$ and~$15$.  Since $|\S_{10}(H)|=|H|\cdot|\SL_{10}|$, it follows by Lagrange's Theorem that the normal subgroups $\S_{10}(H_1)/Z$ and $\S_{10}(H_2)/Z$ are incomparable in $\NN(P)=\NN(\PGL_{10})$.  
We thank Heiko Dietrich for suggesting this example to us.
\end{rem}

\section{Partial braid categories}\label{sect:IB}

In this final section 
we apply the theory from Section \ref{sect:H} once more,
this time to categories of partial braids.  
These categories, denoted~$\IB$, and defined in Section \ref{subsect:defnIB}, are 
categorical analogues of the inverse braid monoids of Easdown and Lavers \cite{EL2004}, and are also preimages of certain instances of the category $\I$ of partial bijections mentioned in Section~\ref{subsect:OtherT}.  We characterise Green's relations in $\IB$ and establish the requisite multiplicative properties in Section \ref{subsect:prelimIB}, where we also describe the maximum $\H$-congruence $\ze=\ze_{\IB}$ in terms of pure braid groups.  Then in Section \ref{subsect:IB} we describe the congruences on arbitrary ideals of~$\IB$; see Theorem~\ref{thm:IB}.

\subsection{Definitions and preliminaries on $\IB$}\label{subsect:defnIB}

Following \cite{Artin1947,Artin1925}, a \emph{braid of degree $n$} is a homotopy class of a collection $(s_1,\ldots,s_n)$ of strings in three-dimensional space~$\RRR^3$.  Here the $s_i$ are smooth embeddings of the unit interval $[0,1]$ into~$\RRR^3$ such that:
\bit
\item for all $i\in[n]$, $s_i(0)=(i,0,1)$ and $s_i(1)=(j_i,0,0)$ for some $j_i\in[n]$,
\item for all $i\in[n]$ and $t\in[0,1]$, the $z$-coordinate of $s_i(t)$ is $1-t$,
\item the strings do not intersect each other.
\eit
We will customarily identify $(s_1,\dots,s_n)$ and the corresponding braid (which is the homotopy class of the former).
The set $\B_n$ of all braids with $n$ strings forms the \emph{braid group of degree~$n$}; the multiplication consists of concatenation and rescaling.  For more details, see for example \cite{Birman1974,MK1999}, and Figure \ref{fig:IBn}, which gives an example product of \emph{partial} braids.  
The map $[n]\to[n]:i\mt j_i$ determined by the endpoints of the strings of $\al$ will be denoted by $\ol\al$.  This is a permutation of~$[n]$, and the map $\B_n\to\S_n:\al\mt\ol\al$ is a group homomorphism.  The kernel of this homomorphism is the \emph{pure braid group} $\P_n=\set{\al\in\B_n}{\ol\al=\id_n}$.  (We will not refer to Brauer or partition monoids/categories in this section, so we may reuse the symbols $\B$ and $\P$ without confusion.)

As in \cite{EL2004}, a \emph{partial braid of degree $n$} is a braid of degree $n$ with some number of strings removed: i.e., a homotopy class of a collection $\al=(s_{i_1},\ldots,s_{i_r})$ of strings, where ${1\leq i_1<\cdots<i_r\leq n}$, satisfying the above conditions.  
We write 
\[
\dom(\al)=\{i_1,\ldots,i_r\} \COMMA \codom(\al)=\{j_1,\ldots,j_r\} \COMMA \rank(\al)=r,
\]
where $s_{i_k}(1)=(j_k,0,0)$ for each $1\leq k\leq r$.  So $\dom(\al)$ and $\codom(\al)$ are subsets of $[n]$, and $\rank(\al)$ is an integer between $0$ and $n$, inclusive.  Note that the partial braid $\al$ induces an injective partial map $[n]\to[n]$ with domain $\dom(\al)$ and codomain $\codom(\al)$, where $i_k\mt j_k$ for each $k$; we denote this map by~$\ol\al$.

Two partial braids $\al,\be$ (of any degree) may be composed by first translating and joining the bottom of $\al$ to the top of $\be$, then removing any strings of $\al$ (respectively, $\be$) that do not become joined to a string of $\be$ (respectively, $\al$), and finally scaling the resulting object so it lies in the region $0\leq z\leq 1$ of $\RRR^3$.  An example calculation is given in Figure \ref{fig:IBn}.  For more details, see \cite{EL2004,JE2007a,Gilbert2006}.

\begin{figure}[ht]
\begin{center}
\scalebox{1}{
\begin{tikzpicture}[scale=.65]

\begin{scope}[shift={(0,0)}]
\understring21
\understring42
\overstring35
\foreach \x in {1,...,4} {\fill (\x,4)circle(.1);}
\foreach \x in {1,...,6} {\fill (\x,0)circle(.1);}
\node[left] () at (0.5,2) {$\al=$};
\draw[-{latex},thick](7,1.75)--(8,1.25);
\end{scope}

\begin{scope}[shift={(0,-6)}]
\understring51
\overstring25
\overstring63
\foreach \x in {1,...,6} {\fill (\x,4)circle(.1);}
\foreach \x in {1,...,5} {\fill (\x,0)circle(.1);}
\node[left] () at (0.5,2) {$\be=$};
\draw[-{latex},thick](7,2.25)--(8,2.75);
\end{scope}

\begin{scope}[shift={(8,-1)}]
\understring21
\understring42
\overstring35
\foreach \x in {1,...,4} {\fill (\x,4)circle(.1);}
\foreach \x in {1,...,6} {\fill (\x,0)circle(.1);}
\draw[-{latex},thick](7,1.75)--(8,1.25);
\end{scope}

\begin{scope}[shift={(8,-5)}]
\understring51
\overstring25
\overstring63
\foreach \x in {1,...,6} {\fill (\x,4)circle(.1);}
\foreach \x in {1,...,5} {\fill (\x,0)circle(.1);}
\draw[-{latex},thick](7,2.25)--(8,2.75);
\end{scope}

\begin{scope}[shift={(16,-1)}]
\halfunderstring42
\halfoverstring35
\end{scope}

\begin{scope}[shift={(16,-3)}]
\halfunderstring51
\halfoverstring25
\end{scope}

\begin{scope}[shift={(16,-3)}]
\foreach \x in {1,...,4} {\fill (\x,4)circle(.1);}
\foreach \x in {1,...,5} {\fill (\x,0)circle(.1);}
\node[right] () at (5.5,2) {$=\al\be$};
\end{scope}

\end{tikzpicture}
}
\caption{Two partial braids $\al\in\IB_{4,6}$ and $\be\in\IB_{6,5}$, and their product $\al\be\in\IB_{4,5}$.}
\label{fig:IBn}
\end{center}
\end{figure}

For natural numbers $m,n\in\N$, we write $\IB_{m,n}$ for the set of all partial braids $\al$ with $\dom(\al)\sub[m]$ and $\codom(\be)\sub[n]$.  Fix a non-empty subset $\C\sub\N$, and define
\[
\IB = \IB(\C) = \bigset{(m,\al,n)}{m,n\in\C,\ \al\in\IB_{m,n}}.
\]
For $m,n,p\in\C$, and for $\al\in\IB_{m,n}$ and $\be\in\IB_{n,p}$, we define
\[
\bd(m,\al,n)=m \COMMA \br(m,\al,n)=n \COMMA (m,\al,n)\cdot(n,\be,p)=(m,\al\be,p),
\]
where $\al\be\in\IB_{m,p}$ is as defined above.  Then
\[
(\IB,\C,\bd,\br,\cdot)
\]
is a partial semigroup, indeed an \emph{inverse category} in the sense of \cite{Kastl1979} and \cite[Section 2.3.2]{CL2002}.  For $\al\in\IB_{m,n}$, the unique inverse of $\al$ (i.e., the unique element $\be$ satisfying $\al=\al\be\al$ and $\be=\be\al\be$) will be denoted by $\al^{-1}$; it is the partial braid from $\IB_{n,m}$ obtained by reflecting $\al$ in the horizontal plane $z=\frac12$; cf.~\cite[Section 1]{EL2004}.  We call $\IB$ the \emph{partial braid category (over $\C$)}. 

As with (linear) transformations, we will avoid clutter by identifying $(m,\al,n)\in\IB$ with the partial braid $\al\in\IB_{m,n}$, but considering $m,n$ as ``encoded'' in $\al$, and writing $\bd(\al)=m$ and $\br(\al)=n$.  

An endomorphism monoid in $\IB$ is simply an \emph{inverse braid monoid} $\IB_n=\IB_{n,n}$ \cite{EL2004}.  In particular, $\IB$ contains the braid groups $\B_n$ for each $n\in\C$; these are precisely the units of $\IB$.  When $\C=\{n\}$ consists of a single natural number, $\IB=\IB_n$ is itself an inverse braid monoid.

Note that the above set-up is somewhat different from the other categories considered in Sections \ref{sect:T}--\ref{sect:J}, where $\C$ is typically a non-empty set of finite sets.  Here the subset $\C$ of $\N$ determines a family of subsets $\set{[n]}{n\in\C}$, which serve to constrain the (co)domains of partial braids in the category $\IB=\IB(\C)$.  While one could in principle replace $\C\sub\N$ with an arbitrary set of subsets of $\N$, this does not seem desirable, and would only serve to complicate notation.  On the other hand, since the bottom ideal $I_0$ of $\IB$ is still trivial for more general $\C$ (see below), the descriptions of congruences on $\IB$ and its ideals are identical, and the lattices arising are no more complicated than for the situation considered here.

It is worthwhile noting also that if $\C$ is all of $\N$ (or even a submonoid of $\N$), then $\IB$ has an additional monoidal structure; given $\al\in\IB_{m,n}$ and $\be\in\IB_{k,l}$, one forms $\al\oplus\be\in\IB_{m+k,n+l}$ by placing a translated copy of $\be$ to the right of $\al$.  This is reminiscent of the analogous operation on the partition category discussed at the start of Section \ref{subsect:prelim_J}.

\subsection{Green's relations and multiplicative properties in $\IB$}\label{subsect:prelimIB}

We keep the above notation, and let $\xi=\LSUB(\C)$.  The next result follows quickly from \cite[Theorem~5]{East2015}, which treats \emph{singular} partial braids, or  it may be proved directly using the same method of proof.

\begin{lemma}\label{lem:Green_IB}
\ben
\item \label{GIB1} For $\al,\be\in\IB$, we have
\bit
\item $(\al,\be)\in\R \iff \bd(\al)=\bd(\be)$ and $\dom(\al)=\dom(\be)$,
\item $(\al,\be)\in\L \iff \br(\al)=\br(\be)$ and $\codom(\al)=\codom(\be)$,
\item $(\al,\be)\in\J\iff(\al,\be)\in\D\iff\rank(\al)=\rank(\be)$.
\eit
\item \label{GIB2} The $\J=\D$-classes of $\IB$ are the sets
\bit
\item $D_q = D_q(\IB) = \set{\al\in\IB}{\rank(\al)=q}$ for each $0\leq q<\xi$.
\eit
These are all regular and stable, and they form a chain: $D_0<D_1<\cdots$.
\item \label{GIB3} The ideals of $\IB$ are the sets
\bit
\item $I_r = I_r(\IB) = \set{\al\in\IB}{\rank(\al)\leq r}$ for each $0\leq r<\xi$, and
\item $I_\om = I_\om(\IB) = \IB $ in the case that $\xi=\om$.
\eit
These are all regular and stable, and they form a chain: $I_0\subset I_1\subset\cdots$.
\item \label{GIB4} The $\H$-class of any idempotent from $D_q$ is isomorphic to $\B_q$, the braid group of degree~$q$.  \epfres
\een
\end{lemma}

Thus, $\IB$ is a chain of ideals, as in Definition \ref{defn:chain}.
The minimal ideal $I_0=D_0$ of~$\IB$ is trivial.

Recall that for any partial braid $\al\in\IB_{m,n}$ there is an associated injective partial map $\ol\al\in\I_{m,n}=\I_{[m],[n]}$; here $\I=\I(\C)$ denotes the category of injective partial maps $[m]\to[n]$, $m,n\in\C$, as briefly discussed in Section \ref{subsect:OtherT}.  We define the relation $\ze$ on $\IB$ by
\[
\ze = \bigset{(\al,\be)\in\IB\times\IB}{\ol\al=\ol\be}.
\]
It is easy to check that $\ze$ is a congruence; in fact it is the kernel of the functor ${\IB\to\I:\al\mt\ol\al}$.  Moreover, since $\ol\al=\ol\be$ clearly implies $\dom(\al)=\dom(\be)$ and ${\codom(\al)=\codom(\be)}$, we have $\ze\sub\H$.  We will soon see that $\ze$ is the largest $\H$-congruence on $\IB$, justifying the choice of notation.

For $A\sub\N$ with $\max(A)<\xi$, we write $\ve_A$ for any partial braid with domain and codomain both equal to $A$, and where every string of $\ve_A$ is vertical.  As noted above, when using these elements, we will have to be careful to indicate which hom-sets they belong to.  When ${m\geq\max(A)}$, clearly $\ve_A\in\IB_m$ is an idempotent.

We now verify that
ideals $I_r$ ($2\leq r<\xi$) satisfy $\Multz{I_r}$, and hence $\Sepz{I_r}$; cf.~Definitions~\ref{defn:Sep} and~\ref{defn:Mult}.

\begin{lemma}\label{lem:P2IB}
Let $m,n\in\C$, and suppose $\al,\be\in\IB_{m,n}$ with $r=\rank(\al)\geq\rank(\be)=q$.  
\ben
\item \label{P2IBi} If $r\geq2$ and $q<r$, then there exists $\ga\in I_r$ such that $\ga\al\in D_{r-1}$ and $\ga\be\in I_{r-1}\sm H_{\ga\al}$.
\item \label{P2IBii} If $q=r\geq1$ and $(\al,\be)\not\in\H$, then there exists $\ga\in I_r$ such that 
\[
\text{$[\al\ga\in D_r$ and $\be\ga\in I_{r-1}]$ \ \ \ or \ \ \ $[\ga\al\in D_r$ and $\ga\be\in I_{r-1}]$.}
\]
\item \label{P2IBiii} If $r\geq2$ and $(\al,\be)\in\H\sm\ze$, then there exists $\ga\in I_r$ such that $\ga\al\in D_{r-1}$ and $\ga\be\in I_{r-1}\sm H_{\ga\al}$.
\een
\end{lemma}

\pf
\firstpfitem{\ref{P2IBi}}  
Write $\dom(\al)=\{i_1,\ldots,i_r\}$.  Without loss of generality, we may assume that ${i_1\not\in\dom(\be)}$.  Let $\ga=\ve_A\in\IB_m$, where $A=\{i_1,\ldots,i_{r-1}\}$.  Then $\dom(\ga\al)=A$, and yet ${i_1\not\in\dom(\ga\be)}$, so $\ga\al\in D_{r-1}$ and $(\ga\al,\ga\be)\not\in\R$.

\pfitem{\ref{P2IBii}}  
By symmetry we may assume that $(\al,\be)\not\in\R$, so that $A=\dom(\al)\not=\dom(\be)$.  We then take $\ga=\ve_A\in\IB_m$, so that $\ga\al=\al\in D_r$ and $\ga\be\in I_{r-1}$.

\pfitem{\ref{P2IBiii}}  
Write $\dom(\al)=\dom(\be)=\{i_1,\ldots,i_r\}$ and $\codom(\al)=\codom(\be)=\{j_1,\ldots,j_r\}$, where $i_k\ol\al=j_k$ for all $k$.  Without loss of generality we may assume that $i_1\ol\be=j_r$.  Then with $A=\{i_1,\ldots,i_{r-1}\}$ and $\ga=\ve_A\in\IB_m$, we have $\ga\al\in D_{r-1}$, yet $j_r\in\codom(\ga\be)\sm\codom(\ga\al)$, so $(\ga\al,\ga\be)\not\in\L$.
\epf

From Lemma \ref{lem:P2IB}\ref{P2IBiii}, and the fact that $D_0$ and $D_1$ are both $\H$-trivial, we have the following:

\begin{cor}
\label{cor:zIB}
The congruence $\ze = \bigset{(\al,\be)\in\IB\times\IB}{\ol\al=\ol\be}$ is the maximum $\H$-congruence on $\IB$.  \qed
\end{cor}

\subsection{Congruences on ideals of $\IB$}\label{subsect:IB}

We now set out to apply Theorem \ref{thm:E5} and describe the congruences on the category~$\IB$ and its ideals.  
We have already noted that $I_0=D_0$ is trivial, and the remaining technical assumptions from Theorem \ref{thm:E5} follow from Lemma \ref{lem:P2IB} and Corollary \ref{cor:zIB}.  
Again, we can achieve a more explicit description of the congruences arising.

We begin as usual by fixing specific group $\H$-classes in $\IB$.  With this in mind, let $0\leq q<\xi$, and fix some $n_q\in\C$ with $n_q\geq q$.  For each braid $\al\in\B_q$, we write $\al^\natural\in\IB_{n_q}=\IB_{n_q,n_q}$ for the partial braid of degree $n_q$ with the same strings as $\al$.  For $\Om\sub\B_q$, we write $\Om^\natural=\set{\al^\natural}{\al\in\Om}$.  In particular, we take $\G_q=\B_q^\natural$, and we note that $\G_q$ is a group $\H$-class in $D_q$.

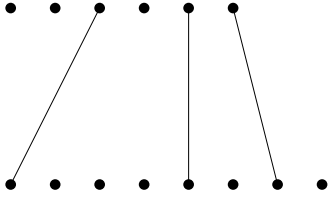
\begin{figure}[htp]
\begin{center}
\scalebox{.9}{
\begin{tikzpicture}[scale=.65]
\foreach \x/\y in {3/1,5/5,6/7} {\draw (\x,4)--(\y,0);}
\foreach \x in {1,...,6} {\fill (\x,4)circle(.12);}
\foreach \x in {1,...,8} {\fill (\x,0)circle(.12);}
\end{tikzpicture}
}
\caption{The partial braid $\ga_{A,B}\in\IB_{6,8}$, where $A=\{3,5,6\}$ and $B=\{1,5,7\}$.}
\label{fig:ga_AB}
\end{center}
\end{figure}

To describe the $\nu_N$ relations, it will first be convenient to define a family of partial braids.  
For $A,B\sub\N$ with $|A|=|B|$, say $A=\{a_1<\cdots<a_k\}$ and $B=\{b_1<\cdots<b_k\}$, we denote by $\ga_{A,B}$ any partial braid with domain $A$ and codomain $B$, and whose strings are straight lines joining $(a_i,0,1)$ to $(b_i,0,0)$ for each $i$; as usual, we need to be careful to specify which hom-set these belong to.  See Figure \ref{fig:ga_AB} for an example.

Now consider two $\H$-related elements $\al,\be\in D_q(\IB_{m,n})$, where $m,n\in\C$, and write ${A=\dom(\al)=\dom(\be)}$.  Writing $\ga = \ga_{[q],A}\in\IB_{q,m}$ (whether $q$ belongs to $\C$ or not), we then define $\phi(\al,\be) = \ga(\al\be^{-1})\ga^{-1}\in\B_q$.  Note that $\phi(\al,\be)$ can be thought of as the braid on $q$ strings obtained from $\al\be^{-1}$ by relabelling the points from $A$ by $1,\ldots,q$.  Also note that $\phi(\al,\be)$ is fully well defined, not just up to conjugacy.

\begin{lemma}\label{lem:nuN_IB}
If $0\leq q<\xi$, and if $N\normal\B_q$, then $\nu_{N^\natural} = \bigset{(\al,\be)\in\H\restr_{D_q}}{\phi(\al,\be)\in N}$.  
\end{lemma}

\pf
Let $(\al,\be)\in\H\restr_{D_q}$, say with $\al,\be\in\IB_{m,n}$, and write $A=\dom(\al)=\dom(\be)$.  Also let $\ga=\ga_{[q],A}\in\IB_{n_q,m}$.  Then since
\[
(\phi(\al,\be)^\natural,\id_{\B_q}^\natural) = (\ga\cdot\al\cdot\be^{-1}\ga^{-1},\ga\cdot\be\cdot\be^{-1}\ga^{-1})
\ANd
(\al,\be) = (\ga^{-1}\cdot\phi(\al,\be)^\natural\cdot\ga\be,\ga^{-1}\cdot\id_{\B_q}^\natural\cdot\ga\be),
\]
we have $(\al,\be)\in\nu_{N^\natural} \iff (\phi(\al,\be)^\natural,\id_{\B_q}^\natural)\in\nu_{N^\natural} \iff \phi(\al,\be)^\natural\in N^\natural \iff \phi(\al,\be)\in N$.
\epf

Next we describe the groups $Z_q$ ($0\leq q<\xi$).  Since $\nu_{Z_q}=\ze\restr_{D_q}$, Lemma~\ref{lem:nu_trick} tells us that
\[
Z_q = \set{\al^\natural}{\al\in\B_q,\ (\id_{\B_q}^\natural,\al^\natural)\in\ze} = \set{\al^\natural}{\al\in\B_q,\ \ol\al=\id_q} = \P_q^\natural
\]
is a copy of the pure braid group of degree $q$.  Thus,
\[
\ze=\nu_{\P_0^\natural}\cup\nu_{\P_1^\natural}\cup\nu_{\P_2^\natural}\cup\cdots.
\]

Lemma \ref{lem:Zqp_IB} below gives a concrete description of the $\downarrow_p$ operations.  Its statement involves the following simple construction.
For a pure braid $\al\in\P_n$, and for $m\leq n$, we denote by $\al\restr_m\in\P_m$ the braid consisting of the first $m$ strings of $\al$.  For $\Om\sub\P_n$, we write $\Om\restr_m=\set{\al\restr_m}{\al\in\Om}$.  If~$G$ is a subgroup of $\P_n$, then $G\restr_m$ is a subgroup of~$\P_m$; if $G$ is additionally a normal subgroup of~$\B_n$, then $G\restr_m$ is also a normal subgroup of~$\B_m$.  
However, it is worth noting that some normal subgroups of~$\P_n$ are not normal in~$\B_n$.
To prove Lemma \ref{lem:Zqp_IB}, we require the following technical result.

\begin{lemma}\label{lem:tauN_IB1}
Let $N\unlhd \B_q$ with $N\leq \P_q$, where $0\leq q<\xi$.  Let $\al\in N$, and suppose $\ga_1\in\IB_{s,n_q}$ and $\ga_2\in\IB_{n_q,t}$, where $s,t\in\C$.  Put $p=\rank(\ga_1\al^\natural\ga_2)$.  Then $\ga_1\al^\natural\ga_2 = \de_1\al'^\natural\de_2$ for some $\al'\in N\restr_p$, and some $\de_1\in\IB_{s,n_p}$ and $\de_2\in\IB_{n_p,t}$.
\end{lemma}

\pf
This proof is somewhat technical, though the steps have a clear geometric meaning, as illustrated in Figure \ref{fig:tauN_IB1}.  Write $\al=(s_1,\ldots,s_q)$, and define the set ${A=\codom(\ga_1)\cap[q]\cap\dom(\ga_2)}$.  The strings of $\al^\natural$ that survive in the product $\ga_1\al^\natural\ga_2$ are precisely the $s_a$ with $a\in A$; these strings (and those of $\ga_1$ and $\ga_2$ connected to them) are coloured red in Figure \ref{fig:tauN_IB1}.  Next, let $\de\in\B_q$ be any braid for which $[p]\ol\de=A$.  Then
\[
\ga_1\al^\natural\ga_2 = \ga_1(\de^{-1}\de\al\de^{-1}\de)^\natural\ga_2 = (\ga_1\de^{\natural-1}) \cdot (\de\al\de^{-1})^\natural \cdot (\de^\natural\ga_2).
\]
Note that $\de\al\de^{-1}\in N$, as $N$ is normal in $\B_q$.  Moreover, it is precisely the first $p$ strings of $\de\al\de^{-1}$ that survive in the above product.  Let
\bit
\item $\de_1\in\IB_{s,n_p}$ be the sub-braid of $\ga_1\de^{\natural-1}$ consisting of all strings whose teminal point is from~$[p]$,
\item $\de_2\in\IB_{n_p,t}$ be the sub-braid of $\de^\natural\ga_2$ consisting of all strings whose initial point is from $[p]$,
\item $\al'=(\de\al\de^{-1})\restr_p\in N\restr_p$ be the sub-braid of $\de\al\de^{-1}$ consisting of all strings whose initial (and terminal) point is from $[p]$.
\eit
Then $\ga_1\al^\natural\ga_2 = \de_1\al'^\natural\de_2$, as required.
\epf

\begin{figure}[htp]
\begin{center}
\scalebox{.9}{
\begin{tikzpicture}[scale=.65]

\begin{scope}[shift={(0,4)}]
\overstring67
\redoverstring73
\redoverstring45
\redoverstring51
\overstring24
\foreach \x in {1,...,7} {\fill (\x,4)circle(.12);}
\node[above] () at (1,4.1) {\footnotesize $1$};
\node[above] () at (7,4.1) { $s$};
\draw[|-|] (0,4)--(0,0);
\node[left] () at (0,2) {$\ga_1$};
\foreach \x in {4,5,7} {\fill[red] (\x,4)circle(.13);}
\end{scope}

\begin{scope}[shift={(0,2)}]
\redhalfoverstring14
\halfoverstring21
\halfoverstring42
\redhalfoverstring35
\redhalfoverstring53
\end{scope}
\begin{scope}[shift={(0,0)}]
\halfoverstring12
\halfoverstring24
\redhalfoverstring35
\redhalfoverstring41
\redhalfoverstring53
\foreach \x in {1,...,8} {\fill (\x,4)circle(.12);}
\foreach \x in {1,3,5} {\fill[red] (\x,4)circle(.13);}
\node[above] () at (8,4) { $\phantom{{}_q}n_q$};
\draw[-|] (0,4)--(0,0);
\node[left] () at (0,2) {$\al^\natural$};
\draw[thick,-{latex}] (10.5,2)--(12.5,2);
\end{scope}

\begin{scope}[shift={(0,-4)}]
\overstring25
\redoverstring31
\redoverstring13
\overstring64
\redoverstring58
\foreach \x in {1,...,9} {\fill (\x,0)circle(.12);}
\foreach \x in {1,...,8} {\fill (\x,4)circle(.12);}
\foreach \x in {1,3,5} {\fill[red] (\x,4)circle(.13);}
\node[above] () at (9,0.1) { $t$};
\draw[-|] (0,4)--(0,0);
\node[left] () at (0,2) {$\ga_2$};
\foreach \x in {1,3,8} {\fill[red] (\x,0)circle(.13);}
\end{scope}

\begin{scope}[shift={(16,8)}]
\overstring67
\redoverstring73
\redoverstring45
\redoverstring51
\overstring24
\foreach \x in {1,...,7} {\fill (\x,4)circle(.12);}
\node[above] () at (1,4.1) {\footnotesize $1$};
\node[above] () at (7,4.1) { $s$};
\draw[|-|] (0,4)--(0,0);
\node[left] () at (0,2) {$\ga_1$};
\foreach \x in {4,5,7} {\fill[red] (\x,4)circle(.13);}
\end{scope}

\begin{scope}[shift={(16,6)}]
\redoverstringx1210
\overstringx2240
\overstringx4250
\redoverstringx3220
\redoverstringx5230
\foreach \x in {1,...,8} {\fill (\x,2)circle(.12);}
\foreach \x in {1,3,5} {\fill[red] (\x,2)circle(.13);}
\node[above] () at (8,2) { $\phantom{{}_q}n_q$};
\draw[-|] (0,2)--(0,0);
\node[left] () at (0,1) {$\de^{\natural-1}$};
\end{scope}

\begin{scope}[shift={(16,4)}]
\redoverstringx1210
\overstringx4220
\overstringx5240
\redoverstringx2230
\redoverstringx3250
\foreach \x in {1,...,8} {\fill (\x,2)circle(.12);}
\foreach \x in {1,2,3} {\fill[red] (\x,2)circle(.13);}
\draw[-|] (0,2)--(0,0);
\node[left] () at (0,1) {$\de^\natural$};
\end{scope}

\begin{scope}[shift={(16,2)}]
\redhalfoverstring14
\halfoverstring21
\halfoverstring42
\redhalfoverstring35
\redhalfoverstring53
\end{scope}
\begin{scope}[shift={(16,0)}]
\halfoverstring12
\halfoverstring24
\redhalfoverstring35
\redhalfoverstring41
\redhalfoverstring53
\foreach \x in {1,...,8} {\fill (\x,4)circle(.12);}
\foreach \x in {1,3,5} {\fill[red] (\x,4)circle(.13);}
\draw[-|] (0,4)--(0,0);
\node[left] () at (0,2) {$\al^{\natural}$};
\end{scope}

\begin{scope}[shift={(16,-2)}]
\redoverstringx1210
\overstringx2240
\overstringx4250
\redoverstringx3220
\redoverstringx5230
\foreach \x in {1,...,8} {\fill (\x,2)circle(.12);}
\foreach \x in {1,3,5} {\fill[red] (\x,2)circle(.13);}
\draw[-|] (0,2)--(0,0);
\node[left] () at (0,1) {$\de^{\natural-1}$};
\end{scope}

\begin{scope}[shift={(16,-4)}]
\redoverstringx1210
\overstringx4220
\overstringx5240
\redoverstringx2230
\redoverstringx3250
\foreach \x in {1,...,8} {\fill (\x,2)circle(.12);}
\foreach \x in {1,2,3} {\fill[red] (\x,2)circle(.13);}
\draw[-|] (0,2)--(0,0);
\node[left] () at (0,1) {$\de^\natural$};
\end{scope}

\begin{scope}[shift={(16,-8)}]
\overstring25
\redoverstring31
\redoverstring13
\overstring64
\redoverstring58
\foreach \x in {1,...,9} {\fill (\x,0)circle(.12);}
\foreach \x in {1,...,8} {\fill (\x,4)circle(.12);}
\foreach \x in {1,3,5} {\fill[red] (\x,4)circle(.13);}
\node[above] () at (9,0.1) { $t$};
\draw[-|] (0,4)--(0,0);
\node[left] () at (0,2) {$\ga_2$};
\draw[thick,-{latex}] (5,-2)--(5,-4);
\foreach \x in {1,3,8} {\fill[red] (\x,0)circle(.13);}
\end{scope}

\begin{scope}[shift={(16,4-22)}]
\redoverstringx7431
\redoverstringx4451
\redoverstringx5411
\overstringx2441
\redoverstringx1110
\overstringx4150
\redoverstringx3120
\redoverstringx5130
\foreach \x in {1,...,7} {\fill (\x,4)circle(.12);}
\node[above] () at (1,4.1) {\footnotesize $1$};
\node[above] () at (7,4.1) { $s$};
\draw[|-|] (0,4)--(0,0);
\node[left] () at (0,2) {$\ga_1\de^{\natural-1}$};
\foreach \x in {4,5,7} {\fill[red] (\x,4)circle(.13);}
\end{scope}

\begin{scope}[shift={(16,0-22)}]
\redoverstringx1413
\overstringx4423
\overstringx5443
\redoverstringx2433
\redoverstringx3453
\redoverstringx1342
\overstringx2312
\overstringx4322
\redoverstringx3352
\redoverstringx5332
\overstringx1221
\overstringx2241
\redoverstringx3251
\redoverstringx4211
\redoverstringx5231
\redoverstringx1110
\overstringx2140
\overstringx4150
\redoverstringx3120
\redoverstringx5130
\foreach \x in {1,...,8} {\fill (\x,4)circle(.12);}
\foreach \x in {1,2,3} {\fill[red] (\x,4)circle(.13);}
\node[above] () at (8,4) { $\phantom{{}_q}n_q$};
\draw[-|] (0,4)--(0,0);
\node[left] () at (0,2) {$(\de\al\de^{-1})^\natural$};
\end{scope}

\begin{scope}[shift={(16,-4-22)}]
\redoverstringx1413
\overstringx4423
\redoverstringx2433
\redoverstringx3453
\overstringx2350
\redoverstringx3310
\redoverstringx1330
\redoverstringx5380
\foreach \x in {1,...,9} {\fill (\x,0)circle(.12);}
\foreach \x in {1,...,8} {\fill (\x,4)circle(.12);}
\foreach \x in {1,2,3} {\fill[red] (\x,4)circle(.13);}
\node[above] () at (9,0.1) { $t$};
\draw[-|] (0,4)--(0,0);
\node[left] () at (0,2) {$\de^\natural\ga_2$};
\foreach \x in {1,3,8} {\fill[red] (\x,0)circle(.13);}
\end{scope}

\begin{scope}[shift={(0,4-22)}]
\redoverstringx7431
\redoverstringx4451
\redoverstringx5411
\redoverstringx1110
\redoverstringx3120
\redoverstringx5130
\foreach \x in {1,...,7} {\fill (\x,4)circle(.12);}
\node[above] () at (1,4.1) {\footnotesize $1$};
\node[above] () at (7,4.1) { $s$};
\draw[|-|] (0,4)--(0,0);
\node[left] () at (0,2) {$\de_1$};
\foreach \x in {4,5,7} {\fill[red] (\x,4)circle(.13);}
\end{scope}

\begin{scope}[shift={(0,0-22)}]
\redoverstringx1413
\redoverstringx2433
\redoverstringx3453
\redoverstringx1342
\redoverstringx3352
\redoverstringx5332
\redoverstringx3251
\redoverstringx4211
\redoverstringx5231
\redoverstringx1110
\redoverstringx3120
\redoverstringx5130
\foreach \x in {1,...,6} {\fill (\x,4)circle(.12);}
\foreach \x in {1,2,3} {\fill[red] (\x,4)circle(.13);}
\node[above] () at (6,4) { $\phantom{{}_p}n_p$};
\draw[-|] (0,4)--(0,0);
\node[left] () at (0,2) {$\al'^\natural$};
\draw[thick,-{latex}] (12.5,2)--(10.5,2);
\end{scope}

\begin{scope}[shift={(0,-4-22)}]
\redoverstringx1413
\redoverstringx2433
\redoverstringx3453
\redoverstringx3310
\redoverstringx1330
\redoverstringx5380
\foreach \x in {1,...,9} {\fill (\x,0)circle(.12);}
\foreach \x in {1,...,6} {\fill (\x,4)circle(.12);}
\foreach \x in {1,2,3} {\fill[red] (\x,4)circle(.13);}
\node[above] () at (9,0.1) { $t$};
\draw[-|] (0,4)--(0,0);
\node[left] () at (0,2) {$\de_2$};
\foreach \x in {1,3,8} {\fill[red] (\x,0)circle(.13);}
\end{scope}

\end{tikzpicture}
}
\caption{The equality $\ga_1\al^\natural\ga_2 = \de_1\al'^\natural\de_2$ from Lemma \ref{lem:tauN_IB1}; see the proof for more details.}
\label{fig:tauN_IB1}
\end{center}
\end{figure}

The next result gives the promised description of the $\downarrow_p$ operations, and says that these are essentially equivalent to the $\restr_p$ operations.

\begin{lemma}\label{lem:Zqp_IB}
For any $0\leq p\leq q<\xi$, and for any $N\normal\B_q$ with $N\leq\P_q$, we have ${(N^\natural)^{\downarrow_p}=(N\restr_p)^\natural}$.
\end{lemma}

\pf
First consider some $\al\in (N\restr_p)^\natural$.  So $\al\in\IB_{n_p}$, and $\al=(\be\restr_p)^\natural$ for some $\be\in N$, and we note that $(\id_{\B_q}^\natural,\be^\natural)\in\nu_{N^\natural}$.  Now, $\al=\ga_1\be^\natural\ga_2$, where $\ga_1=\ve_{[p]}\in\IB_{n_p,n_q}$ and $\ga_2=\ve_{[p]}\in\IB_{n_q,n_p}$.  Since also $\ga_1\id_{\B_q}^\natural\ga_2=\id_{\B_p}^\natural$, we have $(\id_{\B_p}^\natural,\al)=(\ga_1\id_{\B_q}^\natural\ga_2,\ga_1\be^\natural\ga_2)\in\nu_{N^\natural}^\sharp\restr_{D_p}=\nu_{(N^\natural)^{\downarrow_p}}$, so it follows that $\al\in (N^\natural)^{\downarrow_p}$.  This all shows that $(N\restr_p)^\natural\leq (N^\natural)^{\downarrow_p}$.

We will show the reverse containment by showing that $\nu_{(N^\natural)^{\downarrow_p}} \sub \nu_{(N\restr_p)^\natural}$.  As in the proof of Lemma \ref{lem:Zqp_L}, $\nu_{(N^\natural)^{\downarrow_p}}$ is the transitive closure of $\Om\restr_{D_p}$, where
\[
\Om = \bigset{(\ga_1\al^\natural\ga_2,\ga_1\be^\natural\ga_2)}{\al,\be\in N,\ \ga_1,\ga_2\in \IB,\ \br(\ga_1)=\bd(\ga_2)=n_q}.
\]
And again, since $\nu_{(N\restr_p)^\natural}$ is an equivalence, the proof will be complete if we can show that $\Om\restr_{D_p}\sub\nu_{(N\restr_p)^\natural}$.  So consider some pair $(\ga_1\al^\natural\ga_2,\ga_1\be^\natural\ga_2)$ from $\Om\restr_{D_p}$.  Let $\de$, $\de_1$, $\de_2$ and $\al'$ be as in the proof of Lemma \ref{lem:tauN_IB1}.  Then with $\be'=(\de\be\de^{-1})\restr_p\in N\restr_p$, we have $\ga_1\be^\natural\ga_2=\de_1\be'^\natural\de_2$.  Thus, $(\ga_1\al^\natural\ga_2,\ga_1\be^\natural\ga_2) = (\de_1\al'^\natural\de_2,\de_1\be'^\natural\de_2) \in \nu_{(N\restr_p)^\natural}$.  As noted above, this completes the proof.
\epf

We may now prove the main result of this section.  The statement involves the relations $R_{I,\bN}^S$ defined before Theorem \ref{thm:E5}.

\begin{thm}\label{thm:IB}
Let $\C\sub\N$ be a non-empty set of natural numbers, let $\IB=\IB(\C)$ be the partial braid category over $\C$, and keep the above notation.
Then for any $r\geq0$ (including $r=\om$ if $\xi=\om$), the congruences on $I_r$ are precisely the universal congruence $\nab_{I_r}$ and those of the form
\[
R_{I_q,\bN}^{I_r} \qquad\text{for some $0\leq q<r$ and some $\bN=(N_{q+1}^\natural,N_{q+2}^\natural,\ldots)$,}
\]
where
\[
N_p\normal\B_p\ \text{for all $p\geq q+1$}\COMMA
N_p\leq\P_p\ \text{for all $p\geq q+2$}\COMMA
N_p\restr_{p-1}\leq N_{p-1} \ \text{for all $p\geq q+2$}.
\]
\end{thm}

\pf
The proof is almost identical to that of Theorem \ref{thm:L}: we apply Theorem \ref{thm:E5}, using Lemma \ref{lem:P2IB} and the above-mentioned fact that $I_0=D_0$ is trivial.  The only thing left to check is that with $\bN$ as in the statement, $(N_t^\natural)^{\downarrow_p}\leq N_p^\natural$ for all $t>p\geq q+1$.  But this follows from Lemma \ref{lem:Zqp_IB} and a simple induction.  For the $t=p+1$ case we have $(N_{p+1}^\natural)^{\downarrow_p} = (N_{p+1}\restr_p)^\natural \leq N_p^\natural$.  For $t\geq p+2$,
\[
(N_t^\natural)^{\downarrow_p} 
= (N_t\restr_p)^\natural
= ((N_t\restr_{t-1})\restr_p)^\natural
\leq (N_{t-1}\restr_p)^\natural
= (N_{t-1}^\natural)^{\downarrow_p},
\]
at which point we apply an inductive hypothesis.
\epf

\begin{rem}\label{rem:N_B}
The lattices $\Cong(I_0)=\{\nab_{I_0}\}$ and $\Cong(I_1)=\{\De_{I_1},\nab_{I_1}\}$ are of course easily understood.  
We may also visualise the lattice $\Cong(I_2)$ to some extent.  First note that $\B_1$ is trivial, and that~$\B_2$ is an infinite cyclic group generated by the standard Artin generator $\vs$ (in which the first string passes over the second, as in \cite{Artin1925,Artin1947}); thus, the subgroups of $\B_2$ are those of the form $\la\vs^k\ra$, $k\in\N$; such a subgroup is contained in $\P_2$ if and only if $k$ is even.  Thus
\[
\Cong(I_2) = \bigset{R_{I_0,(\B_1,\la\vs^{2k}\ra)}^{I_2}}{k\in\N} \cup \bigset{R_{I_1,\la\vs^{k}\ra}^{I_2}}{k\in\N} \cup \{\nab_{I_2}\}.
\]
Writing $\mid$ for the division relation in $\N$, and writing $\si_{2k}=R_{I_0,(\B_1,\la\vs^{2k}\ra)}^{I_2}=\De_{I_1}\cup\nu_{\la\vs^{2k}\ra}$ and $\tau_k=R_{I_1,\la\vs^{k}\ra}^{I_2}=\nab_{I_1}\cup\nu_{\la\vs^{k}\ra}$ for each $k\in\N$, the order in $\Cong(I_2)$ is given by
\[
\si_{2k}\sub\si_{2k'} \iff k'\mid k\COMMA
\tau_k\sub\tau_{k'} \iff k'\mid k\COMMA
\si_{2k}\sub\tau_{k'} \iff k'\mid 2k\COMMA
\si_{2k}\not\supseteq \tau_{k'},
\]
and of course every congruence is contained in $\nab_{I_2}$.  Thus, the lattice $\Cong(I_2)$ is essentially two copies of $\N$ ordered by reverse division, with one below the other, and with a top element adjoined.  

There does not seem to be much hope of visualising $\Cong(I_r)$ for $r\geq3$.
In fact, to the authors' knowledge, the normal subgroups of the braid groups $\B_n$ ($n\geq3$) have not been classified.
\end{rem}

\begin{rem}
As with the categories of transformations and partitions studied in Sections~\mbox{\ref{sect:T}--\ref{sect:J}}, we could define planar and annular reducts of $\IB$.  That is, we could consider categories such as $\OrdP{\IB}$, which consists of partial braids $\al$ with $\ol\al$ order-preserving; cf.~\cite{East2010}, where the corresponding submonoid of $\IB_n$ was denoted $\mathcal{POIB}_n$.  Note that the category $\OrdP{\I}$ of order-preserving bijections is a homomorphic image of $\OrdP{\IB}$, and also (isomorphic to) a subcategory of $\OrdP{\IB}$.  Indeed, the partial braids of the form $\ga_{A,B}$, pictured in Figure \ref{fig:ga_AB}, form a subcategory of $\IB$ isomorphic to $\OrdP{\I}$; cf.~\cite[Section 3]{East2010}.  The subcategories $\OrdPR{\IB}$, $\OriP{\IB}$ and $\OriPR{\IB}$, with analogous meanings, could also be defined.  

Congruences on these subcategories and their ideals may be described in essentially the same way as for $\IB$ itself; cf.~Theorem~\ref{thm:IB}.  Part \ref{GIB4} of Lemma \ref{lem:Green_IB} needs to be updated, to reflect the fact that group $\H$-classes are subgroups of $\B_q$ mapping onto trivial, cyclic or dihedral permutation groups as appropriate, and each occurrence of $\B_q$ in Theorem \ref{thm:IB} needs to be replaced by the relevant subgroup.  In each of these subcategories, the maximum $\H$-congruence is still of the form $\ze=\nu_{\P_0^\natural}\cup\nu_{\P_1^\natural}\cup\nu_{\P_2^\natural}\cup\cdots$.  
\end{rem}

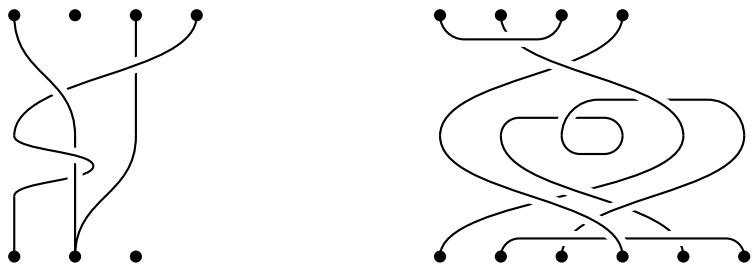
\begin{figure}[ht]
\begin{center}
\scalebox{1}{
\begin{tikzpicture}[scale=.8]

\begin{scope}[shift={(0,2)},yscale=0.5]
\overstring33
\overstring41
\overstring12
\end{scope}

\begin{scope}[shift={(0,0)},yscale=0.5]
\understring22
\overstringx14{2.3}3
\overstringx{2.3}312
\overstringx2320
\overstringx1210
\understring32
\end{scope}

\begin{scope}[shift={(0,0)}]
\foreach \x in {1,2,3,4} {\fill (\x,4)circle(.1);}
\foreach \x in {1,2,3} {\fill (\x,0)circle(.1);}
\end{scope}

\begin{scope}[shift={(7,2)}]
\overdarcx24{.3}
\overdarcx36{.6}
\halfunderstring41
\halfoverstring25
\end{scope}

\begin{scope}[shift={(7,0)}]
\overuarcx13{.4}
\halfoverstring51
\halfoverstring25
\halfoverstring63
\overdarcx26{.3}
\halfoverstring14
\end{scope}

\begin{scope}[shift={(7,-2)}]
\overuarcx34{.3}
\end{scope}

\begin{scope}[shift={(7,0)}]
\foreach \x in {1,...,4} {\fill (\x,4)circle(.1);}
\foreach \x in {1,...,6} {\fill (\x,0)circle(.1);}
\end{scope}

\end{tikzpicture}
}
\caption{A partial vine (left) and a tangle (right).}
\label{fig:vine_tangle}
\end{center}
\end{figure}

\begin{rem}
There are other categories of braid-like objects, including (partial) vines \cite{Lavers1997,JE2007b,JEpresDC} and tangles \cite{Turaev1989}; see Figure \ref{fig:vine_tangle} for examples.  It would be interesting to attempt to use the methods of this paper to describe the congruences of such categories.
\end{rem}

\subsection*{Acknowledgements}

The first author is supported by ARC Future Fellowship FT190100632.  The second author is supported by EPSRC grant EP/S020616/1.
We were assisted by GAP \cite{GAP4,GAP} in the initial stages of research, and in producing some of the lattice diagrams throughout the paper.
We thank Mikhail Volkov and Mark Sapir for useful discussions, and for drawing our attention to Klimov's paper \cite{Klimov1977}.
We also thank the referee for their detailed reading of the paper, and for their helpful suggestions.

\footnotesize
\def\bibspacing{-1.1pt}
\bibliography{biblio}

\def\cprime{$'$}
\begin{thebibliography}{100}

\bibitem{Abramsky2008}
S.~Abramsky.
\newblock Temperley-{L}ieb algebra: from knot theory to logic and computation
  via quantum mechanics.
\newblock In {\em Mathematics of quantum computation and quantum technology},
  Chapman \& Hall/CRC Appl. Math. Nonlinear Sci. Ser., pages 515--558. Chapman
  \& Hall/CRC, Boca Raton, FL, 2008.

\bibitem{Aizenstat1962}
A.~J. A{\u\i}zen{\v{s}}tat.
\newblock On homomorphisms of semigroups of endomorphisms of ordered sets.
\newblock {\em Leningrad. Gos. Ped. Inst. U\v cen. Zap.}, 238:38--48, 1962.

\bibitem{ABG2016}
J.~Ara\'{u}jo, W.~Bentz, and G.~M.~S. Gomes.
\newblock Congruences on direct products of transformation and matrix monoids.
\newblock {\em Semigroup Forum}, 97(3):384--416, 2018.

\bibitem{Artin1925}
E.~Artin.
\newblock Theorie der {Z}\"{o}pfe.
\newblock {\em Abh. Math. Sem. Univ. Hamburg}, 4(1):47--72, 1925.

\bibitem{Artin1947}
E.~Artin.
\newblock Theory of braids.
\newblock {\em Ann. of Math. (2)}, 48:101--126, 1947.

\bibitem{Artin1957}
E.~Artin.
\newblock {\em Geometric algebra}.
\newblock Interscience Publishers, Inc., New York-London, 1957.

\bibitem{ADV2012_2}
K.~Auinger, I.~Dolinka, and M.~V. Volkov.
\newblock Equational theories of semigroups with involution.
\newblock {\em J. Algebra}, 369:203--225, 2012.

\bibitem{BS2011}
J.~Baez and M.~Stay.
\newblock Physics, topology, logic and computation: a {R}osetta {S}tone.
\newblock In {\em New structures for physics}, volume 813 of {\em Lecture Notes
  in Phys.}, pages 95--172. Springer, Heidelberg, 2011.

\bibitem{Baez1994}
J.~C. Baez.
\newblock Hochschild homology in a braided tensor category.
\newblock {\em Trans. Amer. Math. Soc.}, 344(2):885--906, 1994.

\bibitem{BD1995}
J.~C. Baez and J.~Dolan.
\newblock Higher-dimensional algebra and topological quantum field theory.
\newblock {\em J. Math. Phys.}, 36(11):6073--6105, 1995.

\bibitem{BN1996}
J.~C. Baez and M.~Neuchl.
\newblock Higher-dimensional algebra. {I}. {B}raided monoidal {$2$}-categories.
\newblock {\em Adv. Math.}, 121(2):196--244, 1996.

\bibitem{BH2014}
G.~Benkart and T.~Halverson.
\newblock Motzkin algebras.
\newblock {\em European J. Combin.}, 36:473--502, 2014.

\bibitem{BH2019}
G.~Benkart and T.~Halverson.
\newblock Partition algebras {$\mathsf{P}_k(n)$} with {$2k>n$} and the
  fundamental theorems of invariant theory for the symmetric group
  {$\mathsf{S}_n$}.
\newblock {\em J. Lond. Math. Soc. (2)}, 99(1):194--224, 2019.

\bibitem{BRY2012}
S.~Bigelow, E.~Ramos, and R.~Yi.
\newblock The {A}lexander and {J}ones polynomials through representations of
  rook algebras.
\newblock {\em J. Knot Theory Ramifications}, 21(12):1250114, 18, 2012.

\bibitem{Birman1974}
J.~S. Birman.
\newblock {\em Braids, links, and mapping class groups}.
\newblock Princeton University Press, Princeton, N.J.; University of Tokyo
  Press, Tokyo, 1974.
\newblock Annals of Mathematics Studies, No. 82.

\bibitem{BDP2002}
M.~Borisavljevi{\'c}, K.~Do{\v{s}}en, and Z.~Petri{\'c}.
\newblock Kauffman monoids.
\newblock {\em J. Knot Theory Ramifications}, 11(2):127--143, 2002.

\bibitem{Brauer1937}
R.~Brauer.
\newblock On algebras which are connected with the semisimple continuous
  groups.
\newblock {\em Ann. of Math. (2)}, 38(4):857--872, 1937.

\bibitem{CH1999}
P.~M. Catarino and P.~M. Higgins.
\newblock The monoid of orientation-preserving mappings on a chain.
\newblock {\em Semigroup Forum}, 58(2):190--206, 1999.

\bibitem{CPbook}
A.~H. Clifford and G.~B. Preston.
\newblock {\em The algebraic theory of semigroups. {V}ol. {I}}.
\newblock Mathematical Surveys, No. 7. American Mathematical Society,
  Providence, R.I., 1961.

\bibitem{CPbook2}
A.~H. Clifford and G.~B. Preston.
\newblock {\em The algebraic theory of semigroups. {V}ol. {II}}.
\newblock Mathematical Surveys, No. 7. American Mathematical Society,
  Providence, R.I., 1967.

\bibitem{CL2002}
J.~R.~B. Cockett and S.~Lack.
\newblock Restriction categories. {I}. {C}ategories of partial maps.
\newblock {\em Theoret. Comput. Sci.}, 270(1-2):223--259, 2002.

\bibitem{DDE2019}
I.~{{\DJ}}ur{\dj}ev, I.~Dolinka, and J.~East.
\newblock Sandwich semigroups in diagram categories.
\newblock {\em Preprint}, 2019, {\tt arXiv:1910.10286}.

\bibitem{DE2018}
I.~Dolinka and J.~East.
\newblock Semigroups of rectangular matrices under a sandwich operation.
\newblock {\em Semigroup Forum}, 96(2):253--300, 2018.

\bibitem{DE4}
I.~Dolinka and J.~East.
\newblock Twisted {B}rauer monoids.
\newblock {\em Proc. Roy. Soc. Edinburgh Sect. A}, 148(4):731--750, 2018.

\bibitem{DEG2017}
I.~Dolinka, J.~East, and R.~D. Gray.
\newblock Motzkin monoids and partial {B}rauer monoids.
\newblock {\em J. Algebra}, 471:251--298, 2017.

\bibitem{Sandwich1}
I.~Dolinka, I.~{\DJ}ur{\dj}ev, J.~East, P.~Honyam, K.~Sangkhanan, J.~Sanwong,
  and W.~Sommanee.
\newblock Sandwich semigroups in locally small categories {I}: foundations.
\newblock {\em Algebra Universalis}, 79(3):Art. 75, 35 pp, 2018.

\bibitem{Sandwich2}
I.~Dolinka, I.~{\DJ}ur{\dj}ev, J.~East, P.~Honyam, K.~Sangkhanan, J.~Sanwong,
  and W.~Sommanee.
\newblock Sandwich semigroups in locally small categories {II}:
  transformations.
\newblock {\em Algebra Universalis}, 79(3):Art. 76, 53 pp, 2018.

\bibitem{DW2000}
W.~F. Doran, IV and D.~B. Wales.
\newblock The partition algebra revisited.
\newblock {\em J. Algebra}, 231(1):265--330, 2000.

\bibitem{DP2013}
K.~Do\v{s}en and Z.~Petri\'{c}.
\newblock Syntax for split preorders.
\newblock {\em Ann. Pure Appl. Logic}, 164(4):443--481, 2013.

\bibitem{EL2004}
D.~Easdown and T.~G. Lavers.
\newblock The inverse braid monoid.
\newblock {\em Adv. Math.}, 186(2):438--455, 2004.

\bibitem{JE2007a}
J.~East.
\newblock Braids and partial permutations.
\newblock {\em Adv. Math.}, 213(1):440--461, 2007.

\bibitem{JE2007b}
J.~East.
\newblock Vines and partial transformations.
\newblock {\em Adv. Math.}, 216(2):787--810, 2007.

\bibitem{East2010}
J.~East.
\newblock Braids and order-preserving partial permutations.
\newblock {\em J. Knot Theory Ramifications}, 19(8):1025--1049, 2010.

\bibitem{East2015}
J.~East.
\newblock Singular braids and partial permutations.
\newblock {\em Acta Sci. Math. (Szeged)}, 81(1-2):55--77, 2015.

\bibitem{JErook}
J.~East.
\newblock Presentations for rook partition monoids and algebras and their
  singular ideals.
\newblock {\em J. Pure Appl. Algebra}, 223(3):1097--1122, 2019.

\bibitem{JEpresDC}
J.~East.
\newblock Presentations for tensor categories.
\newblock {\em Preprint}, 2020, {\tt arXiv:2005.01953}.

\bibitem{EG2017}
J.~East and R.~D. Gray.
\newblock Diagram monoids and {G}raham--{H}oughton graphs: {I}dempotents and
  generating sets of ideals.
\newblock {\em J. Combin. Theory Ser. A}, 146:63--128, 2017.

\bibitem{EMRT2018}
J.~East, J.~D. Mitchell, N.~Ru\v{s}kuc, and M.~Torpey.
\newblock Congruence lattices of finite diagram monoids.
\newblock {\em Adv. Math.}, 333:931--1003, 2018.

\bibitem{ER2018}
J.~East and N.~Ru\v{s}kuc.
\newblock Congruences on infinite partition and partial {B}rauer monoids.
\newblock {\em Preprint}, 2018, {\tt arXiv:1809.07427}.

\bibitem{Fernandes2001}
V.~H. Fernandes.
\newblock The monoid of all injective order preserving partial transformations
  on a finite chain.
\newblock {\em Semigroup Forum}, 62(2):178--204, 2001.

\bibitem{FGJ2005}
V.~H. Fernandes, G.~M.~S. Gomes, and M.~M. Jesus.
\newblock Congruences on monoids of order-preserving or order-reversing
  transformations on a finite chain.
\newblock {\em Glasg. Math. J.}, 47(2):413--424, 2005.

\bibitem{FGJ2009}
V.~H. Fernandes, G.~M.~S. Gomes, and M.~M. Jesus.
\newblock Congruences on monoids of transformations preserving the orientation
  of a finite chain.
\newblock {\em J. Algebra}, 321(3):743--757, 2009.

\bibitem{FL2011}
D.~G. FitzGerald and K.~W. Lau.
\newblock On the partition monoid and some related semigroups.
\newblock {\em Bull. Aust. Math. Soc.}, 83(2):273--288, 2011.

\bibitem{FHH2009}
D.~Flath, T.~Halverson, and K.~Herbig.
\newblock The planar rook algebra and {P}ascal's triangle.
\newblock {\em Enseign. Math. (2)}, 55(1-2):77--92, 2009.

\bibitem{FY1989}
P.~J. Freyd and D.~N. Yetter.
\newblock Braided compact closed categories with applications to
  low-dimensional topology.
\newblock {\em Adv. Math.}, 77(2):156--182, 1989.

\bibitem{GMbook}
O.~Ganyushkin and V.~Mazorchuk.
\newblock {\em Classical finite transformation semigroups, an introduction},
  volume~9 of {\em Algebra and Applications}.
\newblock Springer-Verlag London, Ltd., London, 2009.

\bibitem{GAP4}
The GAP~Group.
\newblock {\em {GAP -- Groups, Algorithms, and Programming, Version 4.10.1}},
  2019.

\bibitem{Gilbert2006}
N.~D. Gilbert.
\newblock Presentations of the inverse braid monoid.
\newblock {\em J. Knot Theory Ramifications}, 15(5):571--588, 2006.

\bibitem{GH1992}
G.~M.~S. Gomes and J.~M. Howie.
\newblock On the ranks of certain semigroups of order-preserving
  transformations.
\newblock {\em Semigroup Forum}, 45(3):272--282, 1992.

\bibitem{GL1996}
J.~J. Graham and G.~I. Lehrer.
\newblock Cellular algebras.
\newblock {\em Invent. Math.}, 123(1):1--34, 1996.

\bibitem{GL1998}
J.~J. Graham and G.~I. Lehrer.
\newblock The representation theory of affine {T}emperley-{L}ieb algebras.
\newblock {\em Enseign. Math. (2)}, 44(3-4):173--218, 1998.

\bibitem{Gray2007}
R.~Gray.
\newblock Idempotent rank in endomorphism monoids of finite independence
  algebras.
\newblock {\em Proc. Roy. Soc. Edinburgh Sect. A}, 137(2):303--331, 2007.

\bibitem{Gray2008}
R.~Gray.
\newblock Hall's condition and idempotent rank of ideals of endomorphism
  monoids.
\newblock {\em Proc. Edinb. Math. Soc. (2)}, 51(1):57--72, 2008.

\bibitem{Grood2006}
C.~Grood.
\newblock The rook partition algebra.
\newblock {\em J. Combin. Theory Ser. A}, 113(2):325--351, 2006.

\bibitem{Hall1973}
T.~E. Hall.
\newblock On regular semigroups.
\newblock {\em J. Algebra}, 24:1--24, 1973.

\bibitem{HR2005}
T.~Halverson and A.~Ram.
\newblock Partition algebras.
\newblock {\em European J. Combin.}, 26(6):869--921, 2005.

\bibitem{Howie1978}
J.~M. Howie.
\newblock Idempotent generators in finite full transformation semigroups.
\newblock {\em Proc. Roy. Soc. Edinburgh Sect. A}, 81(3-4):317--323, 1978.

\bibitem{Howie}
J.~M. Howie.
\newblock {\em Fundamentals of semigroup theory}, volume~12 of {\em London
  Mathematical Society Monographs. New Series}.
\newblock The Clarendon Press, Oxford University Press, New York, 1995.
\newblock Oxford Science Publications.

\bibitem{Howie1990}
J.~M. Howie and R.~B. McFadden.
\newblock Idempotent rank in finite full transformation semigroups.
\newblock {\em Proc. Roy. Soc. Edinburgh Sect. A}, 114(3-4):161--167, 1990.

\bibitem{Jech2003}
T.~Jech.
\newblock {\em Set theory}.
\newblock Springer Monographs in Mathematics. Springer-Verlag, Berlin, 2003.
\newblock The third millennium edition, revised and expanded.

\bibitem{Jones1983_2}
V.~F.~R. Jones.
\newblock Index for subfactors.
\newblock {\em Invent. Math.}, 72(1):1--25, 1983.

\bibitem{Jones1987}
V.~F.~R. Jones.
\newblock Hecke algebra representations of braid groups and link polynomials.
\newblock {\em Ann. of Math. (2)}, 126(2):335--388, 1987.

\bibitem{Jones1994}
V.~F.~R. Jones.
\newblock The {P}otts model and the symmetric group.
\newblock In {\em Subfactors ({K}yuzeso, 1993)}, pages 259--267. World Sci.
  Publ., River Edge, NJ, 1994.

\bibitem{Jones1994_a}
V.~F.~R. Jones.
\newblock A quotient of the affine {H}ecke algebra in the {B}rauer algebra.
\newblock {\em Enseign. Math. (2)}, 40(3-4):313--344, 1994.

\bibitem{Jones2001}
V.~F.~R. Jones.
\newblock The annular structure of subfactors.
\newblock In {\em Essays on geometry and related topics, {V}ol. 1, 2},
  volume~38 of {\em Monogr. Enseign. Math.}, pages 401--463. Enseignement
  Math., Geneva, 2001.

\bibitem{JP2015}
V.~F.~R. Jones and D.~Penneys.
\newblock Infinite index subfactors and the {GICAR} categories.
\newblock {\em Comm. Math. Phys.}, 339(2):729--768, 2015.

\bibitem{JS1993}
A.~Joyal and R.~Street.
\newblock Braided tensor categories.
\newblock {\em Adv. Math.}, 102(1):20--78, 1993.

\bibitem{Kastl1979}
J.~Kastl.
\newblock Inverse categories.
\newblock In {\em Algebraische {M}odelle, {K}ategorien und {G}ruppoide},
  volume~7 of {\em Stud. Algebra Anwendungen}, pages 51--60. Akademie-Verlag,
  Berlin, 1979.

\bibitem{Kauffman1987}
L.~H. Kauffman.
\newblock State models and the {J}ones polynomial.
\newblock {\em Topology}, 26(3):395--407, 1987.

\bibitem{Kauffman1990}
L.~H. Kauffman.
\newblock An invariant of regular isotopy.
\newblock {\em Trans. Amer. Math. Soc.}, 318(2):417--471, 1990.

\bibitem{Kelly1964}
G.~M. Kelly.
\newblock On the radical of a category.
\newblock {\em J. Austral. Math. Soc.}, 4:299--307, 1964.

\bibitem{Klimov1977}
V.~N. Klimov.
\newblock Congruences of globally idempotent semigroups.
\newblock {\em Ural. Gos. Univ. Mat. Zap.}, 10(3):73--105, 217, 1977.

\bibitem{KX2001}
S.~K{\"o}nig and C.~Xi.
\newblock A characteristic free approach to {B}rauer algebras.
\newblock {\em Trans. Amer. Math. Soc.}, 353(4):1489--1505, 2001.

\bibitem{Lallement1967}
G.~Lallement.
\newblock Demi-groupes r\'{e}guliers.
\newblock {\em Ann. Mat. Pura Appl. (4)}, 77:47--129, 1967.

\bibitem{LF2006}
K.~W. Lau and D.~G. FitzGerald.
\newblock Ideal structure of the {K}auffman and related monoids.
\newblock {\em Comm. Algebra}, 34(7):2617--2629, 2006.

\bibitem{Lavers1997}
T.~G. Lavers.
\newblock The theory of vines.
\newblock {\em Comm. Algebra}, 25(4):1257--1284, 1997.

\bibitem{LZ2012}
G.~Lehrer and R.~Zhang.
\newblock The second fundamental theorem of invariant theory for the orthogonal
  group.
\newblock {\em Ann. of Math. (2)}, 176(3):2031--2054, 2012.

\bibitem{LZ2015}
G.~Lehrer and R.~B. Zhang.
\newblock The {B}rauer category and invariant theory.
\newblock {\em J. Eur. Math. Soc. (JEMS)}, 17(9):2311--2351, 2015.

\bibitem{LZ2017}
G.~I. Lehrer and R.~Zhang.
\newblock Invariants of the special orthogonal group and an enhanced {B}rauer
  category.
\newblock {\em Enseign. Math.}, 63(1-2):181--200, 2017.

\bibitem{Liber1953}
A.~E. Liber.
\newblock On symmetric generalized groups.
\newblock {\em Mat. Sbornik N.S.}, 33(75):531--544, 1953.

\bibitem{MacLane1998}
S.~Mac~Lane.
\newblock {\em Categories for the working mathematician}, volume~5 of {\em
  Graduate Texts in Mathematics}.
\newblock Springer-Verlag, New York, second edition, 1998.

\bibitem{Malcev1952}
A.~I. Mal{\cprime}cev.
\newblock Symmetric groupoids ({R}ussian).
\newblock {\em Mat. Sbornik N.S.}, 31(73):136--151, 1952.
\newblock English translation in Twelve papers in logic and algebra, Amer.
  Math. Soc. Translations Ser 2 \textbf{113}, AMS, 1979, pp. 235--250.

\bibitem{Malcev1953}
A.~I. Mal{\cprime}cev.
\newblock Multiplicative congruences of matrices.
\newblock {\em Doklady Akad. Nauk SSSR (N.S.)}, 90:333--335, 1953.

\bibitem{MSS2000}
M.~P.~O. Marques-Smith and R.~P. Sullivan.
\newblock The congruences on the semigroup of balanced transformations of an
  infinite set.
\newblock {\em J. Algebra}, 234(1):1--30, 2000.

\bibitem{Martin1991}
P.~Martin.
\newblock {\em Potts models and related problems in statistical mechanics},
  volume~5 of {\em Series on Advances in Statistical Mechanics}.
\newblock World Scientific Publishing Co., Inc., Teaneck, NJ, 1991.

\bibitem{Martin1994}
P.~Martin.
\newblock Temperley-{L}ieb algebras for nonplanar statistical mechanics---the
  partition algebra construction.
\newblock {\em J. Knot Theory Ramifications}, 3(1):51--82, 1994.

\bibitem{Martin2008}
P.~Martin.
\newblock Diagram categories, representation theory, statistical mechanics.
\newblock In {\em Noncommutative rings, group rings, diagram algebras and their
  applications}, volume 456 of {\em Contemp. Math.}, pages 99--136. Amer. Math.
  Soc., Providence, RI, 2008.

\bibitem{MM2014}
P.~Martin and V.~Mazorchuk.
\newblock On the representation theory of partial {B}rauer algebras.
\newblock {\em Q. J. Math.}, 65(1):225--247, 2014.

\bibitem{Martin2015}
P.~P. Martin.
\newblock The decomposition matrices of the {B}rauer algebra over the complex
  field.
\newblock {\em Trans. Amer. Math. Soc.}, 367(3):1797--1825, 2015.

\bibitem{Maz1998}
V.~Mazorchuk.
\newblock On the structure of {B}rauer semigroup and its partial analogue.
\newblock {\em Problems in Algebra}, 13:29--45, 1998.

\bibitem{GAP}
J.~D. Mitchell et~al.
\newblock {\em Semigroups - GAP package, Version 3.2.0}, Oct 2019.

\bibitem{MK1999}
K.~Murasugi and B.~I. Kurpita.
\newblock {\em A study of braids}, volume 484 of {\em Mathematics and its
  Applications}.
\newblock Kluwer Academic Publishers, Dordrecht, 1999.

\bibitem{NS1978}
T.~E. Nordahl and H.~E. Scheiblich.
\newblock Regular {$\ast $}-semigroups.
\newblock {\em Semigroup Forum}, 16(3):369--377, 1978.

\bibitem{Okninski1998}
J.~Okni{\'n}ski.
\newblock {\em Semigroups of matrices}, volume~6 of {\em Series in Algebra}.
\newblock World Scientific Publishing Co., Inc., River Edge, NJ, 1998.

\bibitem{Scheiblich1973}
H.~E. Scheiblich.
\newblock Concerning congruences on symmetric inverse semigroups.
\newblock {\em Czechoslovak Math. J.}, 23(98):1--10, 1973.

\bibitem{Solomon2002}
L.~Solomon.
\newblock Representations of the rook monoid.
\newblock {\em J. Algebra}, 256(2):309--342, 2002.

\bibitem{Stroppel2005}
C.~Stroppel.
\newblock Categorification of the {T}emperley-{L}ieb category, tangles, and
  cobordisms via projective functors.
\newblock {\em Duke Math. J.}, 126(3):547--596, 2005.

\bibitem{Sutov1961}
{\`E}.~G. {\v{S}}utov.
\newblock Homomorphisms of the semigroup of all partial transformations.
\newblock {\em Izv. Vys\v s. U\v cebn. Zaved. Matematika}, 1961(3
  (22)):177--184, 1961.

\bibitem{Sutov1961_2}
{\`E}.~G. {\v{S}}utov.
\newblock Semigroups of one-to-one transformations.
\newblock {\em Dokl. Akad. Nauk SSSR}, 140:1026--1028, 1961.

\bibitem{TL1971}
H.~N.~V. Temperley and E.~H. Lieb.
\newblock Relations between the ``percolation'' and ``colouring'' problem and
  other graph-theoretical problems associated with regular planar lattices:
  some exact results for the ``percolation'' problem.
\newblock {\em Proc. Roy. Soc. London Ser. A}, 322(1549):251--280, 1971.

\bibitem{Turaev1989}
V.~G. Turaev.
\newblock Operator invariants of tangles, and {$R$}-matrices.
\newblock {\em Izv. Akad. Nauk SSSR Ser. Mat.}, 53(5):1073--1107, 1135, 1989.

\bibitem{VZ1998}
F.~Van~Oystaeyen and Y.~Zhang.
\newblock The {B}rauer group of a braided monoidal category.
\newblock {\em J. Algebra}, 202(1):96--128, 1998.

\bibitem{Wenzl1988}
H.~Wenzl.
\newblock On the structure of {B}rauer's centralizer algebras.
\newblock {\em Ann. of Math. (2)}, 128(1):173--193, 1988.

\bibitem{Wilcox2007}
S.~Wilcox.
\newblock Cellularity of diagram algebras as twisted semigroup algebras.
\newblock {\em J. Algebra}, 309(1):10--31, 2007.

\bibitem{Xi1999}
C.~Xi.
\newblock Partition algebras are cellular.
\newblock {\em Compositio Math.}, 119(1):99--109, 1999.

\bibitem{YY2004}
H.~B. Yang and X.~L. Yang.
\newblock Maximal subsemigroups of finite transformation semigroups {$K(n,r)$}.
\newblock {\em Acta Math. Sin. (Engl. Ser.)}, 20(3):475--482, 2004.

\end{thebibliography}
\bibliographystyle{abbrv}

\end{document}